\documentclass[reqno]{amsart}

\usepackage{color}
\usepackage{amscd}
\usepackage{amsmath}
\usepackage{cases}
\usepackage{amssymb}
\usepackage{amsfonts}
\usepackage{amsthm}
\usepackage[all]{xy}
\usepackage{enumerate}
\usepackage{mathrsfs}

\newtheorem{Thm}{{\bf Theorem}}[section]

\newtheorem{Lem}[Thm]{{\bf Lemma}}
\newtheorem{Prop}[Thm]{{\bf Proposition}}
\newtheorem{Cor}[Thm]{{\bf Corollary}}

\theoremstyle{definition}

\newtheorem{Def}[Thm]{{\bf Definition}}
\newtheorem{Rem}[Thm]{{\bf Remark}}

\newcounter{Exami}

\numberwithin{equation}{section}

\newcommand\rank{\mathop{\rm rank}\nolimits}
\newcommand\im{\mathop{\rm Im}\nolimits}

\newcommand{\Aut}{\mathop{\rm Aut}\nolimits}
\newcommand{\Tr}{\mathop{\rm Tr}\nolimits}
\newcommand\Hom{\mathop{\rm Hom}\nolimits}
\newcommand\Ext{\mathop{\rm Ext}\nolimits}
\newcommand\End{\mathop{\rm End}\nolimits}
\newcommand\Pic{\mathop{\rm Pic}\nolimits}
\newcommand\Spec{\mathop{\rm Spec}\nolimits}

\newcommand{\res}{\mathop{\sf res}\nolimits}

\newcommand\Proj{\mathop{\rm Proj}\nolimits}

\begin{document}

\title[Moduli of framed logarithmic and parabolic connections]{Moduli spaces of framed logarithmic
and parabolic connections on a Riemann surface}

\thanks{This work was partly supported by JSPS KAKENHI: Grant Numbers JP17H06127, JP19K03422, JP19K14506,
JP22H00094. The first author is partially supported by a J. C. Bose Fellowship (JBR/2023/000003).}

\author[I. Biswas]{Indranil Biswas}

\address{Department of Mathematics, Shiv Nadar University, NH91, Tehsil
Dadri, Greater Noida, Uttar Pradesh 201314, India}

\email{indranil.biswas@snu.edu.in, indranil29@gmail.com}

\author[M. Inaba]{Michi-aki Inaba}

\address{Department of Mathematical and Physical Sciences, 
Nara Women's University, Kitauoya-nishimachi, Nara 630-8506, Japan}

\email{inaba@cc.nara-wu.ac.jp}

\author[A. Komyo]{Arata Komyo}

\address{Department of Material Science, Graduate School of Science, University of Hyogo, 
2167 Shosha, Himeji, Hyogo 671-2280, Japan}

\email{akomyo@sci.u-hyogo.ac.jp}

\author[M.-H. Saito]{Masa-Hiko Saito}

\address{Department of Data Science, Faculty of Business Administration, 
KobeGakuin University, Minatojima, Chuou-ku, Kobe, 650-8586, Japan}

\email{mhsaito@ba.kobegakuin.ac.jp}



\date{}

\begin{abstract}
We construct moduli spaces of framed logarithmic connections and also moduli spaces of framed parabolic
connections. It is shown that these moduli spaces possess a natural algebraic symplectic structure. We
also give an upper bound of the transcendence degree of the algebra of regular functions on the moduli
space of parabolic connections.
\end{abstract}

\maketitle

\tableofcontents

\section{Introduction}

Let $X$ be a compact connected Riemann surface. Since the fundamental group $\pi_1(X)$ of $X$ is finitely presented,
and $\text{GL}(r,{\mathbb C})$ is an affine algebraic group defined over $\mathbb C$,
the space of homomorphisms $\text{Hom}(\pi_1(X),\, \text{GL}(r,{\mathbb C}))$ is an affine complex algebraic variety.
The adjoint action of $\text{GL}(r,{\mathbb C})$ on itself produces an action of $\text{GL}(r,{\mathbb C})$ on
$\text{Hom}(\pi_1(X),\, \text{GL}(r,{\mathbb C}))$.
The moduli space
$${\mathcal M}_R(r)\,:=\,\text{Hom}(\pi_1(X),\, \text{GL}(r,{\mathbb C}))/\!\!/\text{GL}(r,{\mathbb C})
\,=\, \text{Spec}\, {\mathbb C}[\text{Hom}(\pi_1(X),\, \text{GL}(r,{\mathbb C}))]^{\text{GL}(r,{\mathbb C})}
$$
of equivalence classes of representations has an algebraic symplectic structure which was constructed by Goldman \cite{Go} and
Atiyah--Bott \cite{AB}. 
Let ${\mathcal M}_C(r)$ be the moduli space of holomorphic connections on $X$ of rank $r$.
This moduli space also has an algebraic symplectic structure.
The Riemann--Hilbert correspondence identifies ${\mathcal M}_R(r)$ with ${\mathcal M}_C(r)$.
The Riemann--Hilbert correspondence is only complex
analytic and not algebraic, and consequently the identification between ${\mathcal M}_R(r)$ and ${\mathcal M}_C(r)$ is
complex analytic but not algebraic. However, the transport of the symplectic form on ${\mathcal M}_R(r)$ to ${\mathcal M}_C(r)$
by this complex analytic identification actually remains algebraic.
This paper is divided into two parts.
The first part is related to the fact that ${\mathcal M}_C(r)$ has an algebraic symplectic structure.
The second part is related to the fact that the Riemann--Hilbert correspondence is not algebraic.

Now we will discuss on the first part. Fix finitely many distinct points $x_1,\, \cdots,\, x_n$ of
$X$ and denote the divisor $x_1 + \cdots +x_n$ on $X$ by $D$. Consider logarithmic connections on
$X$ of rank $r$ whose polar part is supported on $D$. The corresponding moduli space is known to have a Poisson structure. This
Poisson structure is \textit{not} symplectic if $n\, >\, 0$.

It is shown in Corollary \ref{Cor: Poisson structure on logarithmic moduli} that the Poisson structure on the 
moduli space of logarithmic connections can be elevated to a symplectic structure by introducing frames, over 
the points of $D$, of the holomorphic vector bundle underlying the logarithmic connections. This entails 
construction of the moduli space of framed logarithmic connections that occupy a large fraction of the 
article. The key theorem in the first part of this paper is Theorem \ref{Theorem: d-closedness on the moduli 
of framed connections}, which establishes the $d$-closedness of the canonical nondegenerate $2$-form on the 
moduli space of framed connections. This produces a Poisson structure on the moduli space of logarithmic 
connections; a geometric invariant theoretic construction of this moduli space was given by Nitsure 
\cite{Nit}.

In \cite{BLP0} and \cite{BLP}, generalized Higgs bundles on $X$ were considered where the Higgs fields are 
allowed to have poles along a fixed divisor $D$ on $X$. The corresponding moduli spaces have a Poisson 
structure which was constructed independently by Bottacin \cite{Bo} and Markman \cite{Mark}. It was shown in 
\cite{BLP0} and \cite{BLP} that by imposing frames of the vector bundles
underlying the Higgs bundles, over $D$, these Poisson structures can 
be enhanced to symplectic structure. The present work is an analogue of \cite{BLP} for connections, in place 
of Higgs fields.

The moduli space of logarithmic parabolic connections was constructed 
in \cite{IIS1} and \cite{Inaba-1}.
If we fix eigenvalues of residues of logarithmic parabolic connections,
then the moduli space of logarithmic parabolic connections with the fixed eigenvalues of residues has a canonical symplectic structure.
In Section \ref{Subsect:PoissonStr}, we discuss a relationship between the framed logarithmic connections 
and the logarithmic parabolic connections.
As an outcome, it is proved
that the moduli space of logarithmic parabolic connections has a canonical Poisson structure,
whose restriction to the locus of fixed eigenvalues of residues induces the symplectic structure
due to \cite{IIS1} and \cite{Inaba-1}
(Corollary \ref{Cor: Poisson structure on parabolic connection moduli}).
Moreover, this Poisson structure satisfies the condition that the forgetful map to the moduli space of
logarithmic connections --- that forgets the parabolic structure --- is a Poisson map.
The restriction of this Poisson map to the loci of fixed eigenvalues of residues 
is an isomorphism if the eigenvalues are generic, and
it produces a resolution of singularities if the eigenvalues are special.

Now we will discuss on the second part.
In this part, we focus on the algebraic moduli space of logarithmic parabolic connections such that
eigenvalues of residues are fixed. We call this moduli space the {\it de Rham moduli space}.
This moduli space is related to other moduli spaces having rich geometric structures.
First, there is the moduli space of equivalence classes of representations of the fundamental
group $\pi_1(X\setminus D)$ with fixed local monodromy data around the points of $D$,
which is known as the {\it character variety}. The relationship
between the moduli space of logarithmic parabolic connections and the character variety is given by the Riemann--Hilbert correspondence.
In the framework of \cite{IIS1} and \cite{Inaba-1}, 
the Riemann--Hilbert correspondence gives a simultaneous family of holomorphic maps
from the de Rham moduli spaces to the character varieties
over all the eigenvalues of residues. This Riemann--Hilbert morphism
is biholomorphic when the eigenvalues of residues are generic,
and it is an analytic resolution of singularities when the eigenvalues of residues are special.
Note that the characteristic variety in \cite{IIS1} and \cite{Inaba-1} is not
smooth for special eigenvalues of residues, but its singularities actually well
explain the geometry of special solutions of the isomonodromy equation
(see \cite{ST}).
Simpson introduced in \cite{Sim0} the notion of a filtered local system
which bijectively correspond to the parabolic connections
under the assumption that the eigenvalues of residues are fixed.
In \cite{Yamakawa}, Yamakawa constructed the algebraic moduli space of filtered local systems, which
is actually nonsingular. We call it the {\it Betti moduli space}.
Yamakawa proved in \cite{Yamakawa} that the Riemann--Hilbert morphism
is a biholomorphism between the de Rham moduli space
and the Betti moduli space.
Secondly, there is the moduli space of logarithmic parabolic Higgs bundles 
with fixed eigenvalues of residues together with stability data.
We call this moduli space the {\it Dolbeault moduli space}.
The relation between these moduli spaces is given by
the logarithmic version of the non-abelian Hodge theory constructed by Simpson in \cite{Sim0}.

In the case where the polar divisor $D$ is empty,
Simpson introduced in \cite{Sim1} and \cite{Sim2},
the three moduli spaces in his framework:
the {\it de Rham moduli space}, the {\it Dolbeault moduli space}, and the {\it Betti moduli space}.
These are {\it algebraic} moduli spaces
and are related to each other by 
the non-abelian Hodge theory and 
the Riemann--Hilbert correspondence.
However, the algebraic structures of these moduli spaces are very different.
In this paper, we consider the logarithmic version of these
three moduli spaces.
First our Betti moduli space is affine when the eigenvalues
of the residues are generic.
So the transcendence degree of its
affine coordinate ring is equal to the dimension of the moduli space.
On the other hand, the transcendence degree of the ring of global
algebraic functions on the Dolbeault moduli space is exactly the half
of the dimension of the moduli space,
a fact which is deduced from the properness of the Hitchin map.
In some cases, the global algebraic functions on the de Rham moduli spaces are simply the
constant scalars \cite{BR}. For general logarithmic connections, the coefficients of the
characteristic polynomial of residue at each singular point give algebraic functions on the moduli space.
The main theorem of the second part of this paper is
Theorem \ref{Thm: global algebraic functions on connection moduli}
which states that the transcendence degree
of the ring of global algebraic functions on our de Rham moduli space
is less than or equal to that for our Dolbeault moduli space.
In particular, our de Rham moduli space is not affine (this was announced in \cite{BIKS}).
To be precise, there was in fact an inadequate argument on finite generation of a graded ring
in the outline of the proof of \cite[Theorem 10]{BIKS}.
In this paper, we reconstruct a proof of it through a refinement of the statement 
(see Theorem \ref{Thm: global algebraic functions on connection moduli}).
As a consequence of Theorem \ref{Thm: global algebraic functions on connection moduli},
the Riemann--Hilbert morphism, which appears in \cite{IIS1}, \cite{Inaba-1},
is not algebraic in the logarithmic case (see Corollary \ref{Cor: logarithmic Riemann-Hilbert is not algebraic}).

Regarding the above three moduli spaces, we are mostly interested in the case where $X$ is defined over the 
field of complex numbers. However, it is also worth considering the case where the base field is of positive 
characteristic. When the base field is of positive characteristic, N. Katz introduced the notion of 
$p$-curvature in \cite{Katz}, from which Laszlo and Pauly derived the proper Hitchin map on a de Rham moduli 
space (see \cite{Las-Pa}). By the investigation of the Hitchin map on a de Rham moduli space by Groechenig in 
\cite{Groech}, the ring of global algebraic functions on the de Rham moduli space of connections without pole 
has the same transcendence degree as that of the ring of global algebraic functions on the Dolbeault moduli
space, when the characteristic of the base field is positive. So the
similar inequality as in Theorem \ref{Thm: global algebraic functions on connection 
moduli} for connections without pole becomes the equality for curves when the base field is of positive 
characteristic, while the inequality is strict for curves of higher genus defined over the field of complex numbers 
(see \cite{BR}).

Analogous to the regular case in \cite{Biswas}, we can also show, in the logarithmic case, that
the pullback, via the Riemann--Hilbert morphism, of the canonical algebraic symplectic form
on the Betti moduli space coincides with that on the de Rham moduli space.
Although not stated explicitly, it can also be found in the proof of
\cite[Proposition 7.3]{Inaba-1}. This was also proved in the earlier work in the rank two case
by Iwasaki \cite{Iw1}. In fact, the main point of \cite{Iw1} is the construction of the isomonodromic lift
of the family of symplectic forms.
A more conceptual construction of the isomonodromic lift of the family of symplectic forms
was constructed by Komyo in \cite{Komyo} --- from the moduli theoretic point of view --- by
using the cohomological description of the isomonodromic deformation given in \cite{BHH}.

P. Boalch proved the following: The monodromy map between any moduli space of unramified irregular singular connections of any rank
on a curve of genus zero and its corresponding wild character variety
is symplectic structure preserving \cite[p.~182, Theorem 6.1]{Bo1} (see also \cite{Bo2}).
The algebraic moduli space of unramified irregular singular
connections and its algebraic symplectic structure are constructed
in \cite{Inaba-Saito}.

We give a brief outline of the contents of this paper.

Section \ref{section: framed G-connections} provides general notions
of framed principal $G$-bundles on a compact Riemann surface $X$
and also of framed $G$-connections.

From Section \ref{section: construction of moduli}, we restrict to the case
of $G\,=\,\mathrm{GL}(r,\mathbb{C})$.
Subsection \ref{subsection: definition of framed moduli}
provides the formulation of moduli problem for framed connections.
Subsection \ref{subsection: construction of moduli stack}
provides the construction of the moduli space of framed 
$\mathrm{GL}(r,\mathbb{C})$-connections as a Deligne--Mumford stack
and also the irreducibility of its open substack where the underlying
framed bundles are simple.

Section \ref{section: symplectic structure} is devoted to the
construction of a canonical $2$-form on the moduli space of framed connections
and also to prove its $d$-closedness.
The main technical part is Subsection \ref{2019.11.29.10.38}.
Over the open subset where the underlying framed bundles are simple,
the canonical $2$-form on the moduli space of framed connections
becomes $d$-closed (Proposition \ref{2019.12.22.18.19}, Proposition
\ref{2019.12.22.18.19_a}).
Its proof is essentially reduced to the $d$-closedness 
of the canonical $2$-form on the character variety constructed
by Goldman in \cite{Go} when the genus of $X$ is greater than $1$.
When the genus of $X$ is zero or one, the proof of $d$-closedness is reduced to that for the
form on the moduli space of parabolic connections given in \cite{Inaba-1}.
In Subsection \ref{subsection: main d-closedness}, we prove the 
$d$-closedness of the canonical $2$-form on the entire moduli space
of simple framed connections
(see Theorem \ref{Theorem: d-closedness on the moduli of framed connections}), which is the main theorem of
the first half.
Its proof is reduced to Proposition \ref{2019.12.22.18.19} and
Proposition \ref{2019.12.22.18.19_a} through an argument for extending the polar divisor.
Subsection \ref{subsection: d-closedness for H-framing} and Subsection
\ref{Subsect:PoissonStr} are immediate consequences of
Theorem \ref{Theorem: d-closedness on the moduli of framed connections}.
We can see that the Poisson structure on several known moduli spaces
of connections can be reconstructed from the symplectic structure
on the moduli space of framed connections. 

Section \ref{section: global algebraic functions implies not affine} is devoted to
establishing an upper bound for
the transcendence degree of the ring of global algebraic functions on the moduli space of parabolic 
connections. In Subsection \ref{subsection: moduli of parabolic connections and Higgs bundles}, we recall the 
notions of parabolic connections and parabolic Higgs bundles, which work over the base field of arbitrary 
characteristic. In Subsection \ref{subsection: codimension estimation}, we prove in Proposition
\ref{proposition: codimension in de Rham moduli} that the locus of non-simple underlying quasi-parabolic bundles 
has codimension at least $2$ in the moduli space of parabolic connections. The proof is carried out by 
constructing a parameter space of non-simple quasi-parabolic bundles and a compatible connections on them. 
The essential part is to bound the dimension of the parameter space of non-simple quasi-parabolic bundles
(see Propositions \ref{dimension of non-simple parabolic bundles in higher genus case}, \ref{dimension of 
non-simple parabolic bundles when genus=1}, \ref{proposition: codimension genus=0 and n>3} and
\ref{proposition: codimension genus zero n=3 case}). Since we need to verify many cases, the proofs of these 
propositions contain a considerable amount of calculation, but each step is checked by relatively elementary 
arguments. By virtue of Proposition \ref{proposition: codimension in de Rham moduli}, the ring of global 
algebraic functions on the moduli space of parabolic connections can be replaced with that on the open loci 
where the underlying quasi-parabolic bundles are simple. Subsection \ref{subsection: transcendence degree} 
provides the main estimate for the transcendence degree of the global algebraic functions on the moduli 
space of parabolic connections. Over the moduli space of simple quasi-parabolic bundles, we construct in 
Proposition \ref{Prop: Deligne--Hitchin family} something like a relative compactification of a 
Deligne--Hitchin family, whose generic fiber is a relative compactification of the moduli space of compatible 
parabolic connections and whose special fiber is that of parabolic Higgs bundles. This family gives a family 
of sheaves of graded rings over the moduli space of simple quasi-parabolic bundles. A rough idea of the proof 
of Theorem \ref{Thm: global algebraic functions on connection moduli} is to estimate the transcendence degree 
of the ring of global sections of this sheaf of graded rings. In order to correct the flaw in the proof of 
\cite[Theorem 10]{BIKS}, we actually consider the subring generated by a suitable transcendence basis of the 
graded ring over a generic fiber and compare it with that on the special fiber.

\section{Framed $G$-connections}\label{section: framed G-connections}

Let $X$ be a compact connected Riemann surface, and let 
$x_1,\, \cdots,\, x_n$ be finitely many distinct points on $X$. Let 
\begin{equation*}
D\,=\, x_1 + \cdots +x_n 
\end{equation*}
be the reduced effective divisor on $X$. For notational convenience, the subset
$\{x_1,\,\cdots,\, x_n\}\, \subset\, X$ will also be denoted by $D$.
Denote by $K_X$ the holomorphic cotangent bundle of $X$.

\subsection{Framed principal $G$-bundles}

Let $G$ be a connected complex reductive affine algebraic group. The Lie algebra of $G$
will be denoted by $\mathfrak g$. Let
\begin{equation}\label{dp0}
p \,\colon\, E_G \,\longrightarrow\, X
\end{equation}
be a holomorphic principal $G$-bundle over $X$. For any point $x\,\in\, X$,
the fiber $p^{-1} (x) \subset E_G$ will be denoted by $(E_G)_x$. 

\begin{Def}[{See \cite[p.~5]{BLP}}]
For each point $x$ of the above subset $D$, fix a closed complex Lie proper subgroup
\begin{equation*}
H_x \,\subsetneq \,G\, .
\end{equation*}
A \textit{framing} of $E_G$ over the divisor $D$ is a map
\begin{equation*}
\phi \,\colon\, D \,\longrightarrow\, \bigcup_{x \in D} (E_G)_x/ H_x
\end{equation*}
such that $\phi (x) \,\in \,(E_G)_x/ H_x$ for every $x \,\in\, D$.
A \textit{framed principal $G$-bundle} on X is a holomorphic principal $G$-bundle $E_G$ on X 
equipped with a framing over $D$.
\end{Def}

A framing $\phi$ of $E_G$ produces a reduction of structure group
\begin{equation}\label{e2}
\textbf{H}_x\, :=\, q^{-1}_x(\phi(x))\, \subset\, (E_G)_x
\end{equation}
to $H_x$ at each point $x\, \in\, D$, where $q_x\, :\, (E_G)_x\, \longrightarrow\, (E_G)_x/ H_x$ is the
quotient map.

\subsection{Adjoint bundle for framed principal $G$-bundles}

Let $T_{E_G/X}\, \longrightarrow\, E_G$ be the relative tangent bundle for the projection $p$ in \eqref{dp0}.
Using the action of the group $G$ on $E_G$, this relative tangent bundle $T_{E_G/X}$ is identified 
with the trivial vector bundle $E_G\times \mathfrak{g} \,\longrightarrow\, E_G$ 
with fiber $\mathfrak{g} \,=\, \mathrm{Lie}(G)$. The quotient $(T_{E_G/X})/G$ is a vector bundle over $X$. 
The above identification of $T_{E_G/X}$ with $E_G\times \mathfrak{g}$ produces an identification of $(T_{E_G/X})/G$
with the vector bundle on $X$ associated to the principal $G$-bundle $E_G$ 
for the adjoint action of $G$ on $\mathfrak{g}$. 
This associated vector bundle, which is denoted by $\mathrm{ad}(E_G)$, 
is called the \textit{adjoint bundle} for $E_G$.
The fiber over any $x \,\in\, X$ for the natural projection $\mathrm{ad}(E_G) \,\longrightarrow \,X$ 
will be denoted by $\mathrm{ad}(E_G)_x$; it is a Lie algebra isomorphic to $\mathfrak g$.

Since the group $G$ is reductive, its Lie algebra $\mathfrak g$ admits $G$-invariant nondegenerate symmetric
bilinear forms. Fix a $G$-invariant nondegenerate symmetric bilinear form
\begin{equation}\label{eq1}
\sigma \,\,\colon\,\, \mathrm{Sym}^2( \mathfrak{g}) \,\longrightarrow\, \mathbb{C}
\end{equation}
on $\mathfrak{g}$. From the above construction of $\mathrm{ad}(E_G)$ it follows that
given any point $z\,\in\, (E_G)_y$ there is a corresponding isomorphism of Lie algebras
$I_z\, :\, \mathfrak{g} \longrightarrow \mathrm{ad}(E_G)_y$. Using $I_z$,
the form $\sigma$ in \eqref{eq1} produces a symmetric nondegenerate bilinear form 
on the fiber $\mathrm{ad}(E_G)_y$; this bilinear form on $\mathrm{ad}(E_G)_y$
constructed using $\sigma$ is actually
independent of the choice of the point $z$ because $\sigma$ is $G$-invariant. Let
\begin{equation}\label{ews}
\widehat{\sigma}\, \,\colon\, \mathrm{Sym}^2( \mathrm{ad}(E_G))\,\longrightarrow\, {\mathcal O}_X
\end{equation}
be the bilinear form constructed as above using $\sigma$.

Let $\phi$ be a framing of $E_G$ over $D$. For every $x\, \in\, D$, define the Lie subalgebra
\begin{equation}\label{ehx}
\mathcal{H}_x\,:=\, \text{ad}(\textbf{H}_x)\, \subset\, \text{ad}(E_G)_x
\end{equation}
(see \eqref{e2}).

\subsection{Framing of $G$-connections}

Take a holomorphic principal $G$-bundle $E_G$ over $X$.
Let $TE_G$ be the holomorphic tangent bundle of $E_G$.
Consider the action of $G$ on $TE_G$ given by the tautological action of $G$ on $E_G$.
The quotient 
\begin{equation*}
\mathrm{At} (E_G) \,:=\, (TE_G) / G
\end{equation*}
is a holomorphic vector bundle over $X$; it is called the \textit{Atiyah algebra} for $E_G$.
The Lie bracket operation of the vector fields on $E_G$ produces a Lie algebra structure on the coherent sheaf
associated to $\mathrm{At} (E_G)$. There is a short exact sequence of holomorphic vector bundles on $X$
\begin{equation}\label{atb}
0 \,\longrightarrow\, \mathrm{ad}(E_G)\,\longrightarrow\, 
\mathrm{At} (E_G)\,\xrightarrow{\ p_{\mathrm{At}} \ }\,
TX \,\longrightarrow\, 0\, ,
\end{equation}
where the projection $p_{\mathrm{At}}$ is given by the differential $dp$ of the map $p$
in \eqref{dp0} \cite{At}. All the homomorphisms in \eqref{atb} are compatible with the
Lie algebra structures. Define a holomorphic vector bundle $\mathrm{At}_D (E_G)$
over $X$ as $$\mathrm{At}_D (E_G)\, :=\, p_{\mathrm{At}}^{-1} (TX\otimes {\mathcal O}_X(-D))
\, \subset\, \mathrm{At} (E_G)\, .$$
Then \eqref{atb} gives the following short exact sequence of holomorphic vector bundles on $X$:
\begin{equation}\label{atb2}
0 \,\longrightarrow\, 
\mathrm{ad}(E_G) \longrightarrow 
\mathrm{At}_D (E_G) \,\xrightarrow{\ \,p_{\mathrm{At}_D}\ }\,
TX(-D)\,:=\, TX\otimes{\mathcal O}_X(-D)\, \longrightarrow 0\, ,
\end{equation}
where $p_{\mathrm{At}_D}$ is the restriction, to $\mathrm{At}_D (E_G)\, \subset\,
\mathrm{At} (E_G)$, of the homomorphism $p_{\mathrm{At}}$ in \eqref{atb}.

\begin{Def}[{\cite{At}}] 
A \textit{holomorphic connection} on $E_G$ is a holomorphic homomorphism of vector bundles
\begin{equation*}
\nabla \ \colon\ TX \ \longrightarrow\ \mathrm{At} (E_G)
\end{equation*}
such that $ p_{\mathrm{At}} \circ \nabla\,=\,\mathrm{Id}_{TX}$, where
$p_{\mathrm{At}}$ is the projection in \eqref{atb}.
A \textit{$D$-twisted holomorphic connection} on $E_G$ (also called a
logarithmic connection on $E_G$ with polar part on $D$) is 
a holomorphic homomorphism of vector bundles
\begin{equation*}
\nabla \ \colon\ TX(-D) \ \longrightarrow\ \mathrm{At}_D (E_G)
\end{equation*}
such that $ p_{\mathrm{At}_D} \circ \nabla\,=\,\mathrm{Id}_{TX(-D)}$, where $p_{\mathrm{At}_D}$ is the
homomorphism in \eqref{atb2}.
\end{Def}

For a $D$-twisted holomorphic connection $\nabla$ on $E_G$, consider the following commutative diagram
$$
\begin{matrix}
0 & \longrightarrow & \mathrm{ad}(E_G) & \longrightarrow & 
\mathrm{At}_D (E_G) & \stackrel{\nabla}{\longleftarrow} &TX(-D) & \longrightarrow &0\\
&& \Vert && ~\,~\, \Big\downarrow\iota'' && ~\,~\, \Big\downarrow\iota'\\
0 & \longrightarrow & \mathrm{ad}(E_G) & \longrightarrow & 
\mathrm{At} (E_G) & \stackrel{p_{\mathrm{At}}}{\longrightarrow} &TX & \longrightarrow &0
\end{matrix}
$$
where $\iota'$ and $\iota''$ are the natural inclusion homomorphisms. For any point
$x\, \in\, D$, the homomorphism of fibers $$\iota'(x) \,:\, TX(-D)_x\, \longrightarrow\,
T_xX$$ vanishes, and hence
$(p_{\mathrm{At}}\circ\iota''\circ\nabla)(TX(-D)_x)\,=\, 0$ by the commutativity of the
above diagram. Consequently, we have
$$(\iota''\circ\nabla)(TX(-D)_x)\,\subset\, \mathrm{ad}(E_G)_x\, .$$
Note that for any point $x\, \in\, D$, using the Poincar\'e adjunction formula it follows that
\begin{equation}\label{pi}
a_x\, :\, TX(-D)_x\, \stackrel{\sim}{\longrightarrow}\, \mathbb C.
\end{equation}
The element
$$\mathrm{res}_{x} (\nabla)\, :=\, (\iota''\circ\nabla)(1)\,\in\, \mathrm{ad}(E_G)_x$$
is called the residue of the logarithmic connection $\nabla$ at $x$. To describe this residue
explicitly, first recall that a holomorphic
connection on $E_G$ furnishes lift of holomorphic vector fields on any open subset $U$ of $X$ to 
$G$--invariant holomorphic vector fields on $E_G\big\vert_{p^{-1}(U)}$. Similarly, 
a $D$-twisted holomorphic connection $\nabla$ furnishes lift of holomorphic vector fields on any
open subset $U\, \subset\, X$, vanishing on $D\cap U$, to the $G$--invariant holomorphic vector
fields on $E_G\big\vert_{p^{-1}(U)}$. In other words, these lifts are locally
defined $G$--invariant holomorphic sections of $TE_G(-\log p^{-1}(D))$. Therefore, given a vector field $v$ defined on a
neighborhood of $x_i\, \in\, D$ of $X$, such that $v(x_i)\,=\,0$
and $a_{x_i}(v(x_i))\, \not=\, 0$ (see \eqref{pi}), its lift $\widetilde{v}$ to $E_G$ for $\nabla$ may be nonzero
on $p^{-1}(x_i)$ because $\widetilde{v}$ may be a nonzero vertical vector field on $p^{-1}(x_i)$. The residue of
$\nabla$ at $x_i$ is $\widetilde{v}(p^{-1}(x_i))/a_{x_i}(v(x_i))\, \in\,
\mathrm{ad}(E_G)_{x_i}$ (see \eqref{pi}).

For any $x \in D$, let 
\begin{equation}\label{2020.6.24.10.38}
\mathcal{H}^{\perp}_x \,\subset\, \mathrm{ad}(E_G)_x
\end{equation}
be the annihilator of $\mathcal{H}_x \,\subset \,\mathrm{ad}(E_G)_x$, defined in
\eqref{ehx}, with respect to the bilinear form $\widehat{\sigma}(x)$ in \eqref{ews}.

\begin{Def}\label{definition of framed connection}
A \textit{framed $G$-connection} is a triple of the form $(E_G,\,\nabla,\,\phi)$,
where $(E_G,\, \phi)$ is a framed principal $G$-bundle and $\nabla\,\colon\, TX(-D)
\,\longrightarrow \,\mathrm{At}_D(E_G)$ is a $D$-twisted connection such
that $\mathrm{res}_{x} (\nabla)\, \in\, \mathcal{H}^{\perp}_{x} \,\subset\, \mathrm{ad}(E_G)_{x}$
for every $x\,\in\, D$, where $\mathcal{H}^{\perp}_{x}$ is constructed in \eqref{2020.6.24.10.38}.
\end{Def}

\subsection{Infinitesimal deformations}

Consider the following $2$-term complex of sheaves on $X$:
\begin{equation}\label{2019.11.19.14.51}
\mathcal{C}_{\bullet} \,\colon\, \mathrm{ad} (E_G) (-D)\,:=\,\mathrm{ad} (E_G)\otimes {\mathcal O}_X(-D) 
\,\xrightarrow{\ \nabla \ } \, \mathrm{ad} (E_G) \otimes K_X (D)\,:=\, \mathrm{ad} (E_G) \otimes K_X
\otimes {\mathcal O}_X(D).
\end{equation}

\begin{Lem}[{See \cite[Lemma 3.5]{BLP} and \cite[Proposition 4.4]{Ch}}]\label{lemid}
Assume that $H_x\,=\, \{e \}$ for every $x\,\in\, D$. 
The infinitesimal deformations of the framed $G$-connection $(E_G,\, \nabla,\, \phi)$ 
are parametrized by the elements of the first hypercohomology $\mathbb{H}^1(\mathcal{C}_{\bullet})$
of the complex in \eqref{2019.11.19.14.51}.
\end{Lem}

Let
\begin{equation}\label{dn1}
(E_G,\, \widehat{\nabla},\, \phi)
\end{equation}
be a framed $G$-connection (see Definition \ref{definition of framed connection}). 
Consider the subspace $\mathcal{H}_x \,\subset\, \mathrm{ad}(E_G)_x$ in \eqref{ehx}. 
Let $\mathrm{ad}_{\phi}(E_G)$ and $\mathrm{ad}^n_{\phi}(E_G)$ be the holomorphic vector bundles on $X$ 
defined by the following short exact sequences of coherent analytic sheaves on $X$:
\begin{equation}\label{epq1}
0 \longrightarrow \mathrm{ad}_{\phi}(E_G) 
\longrightarrow \mathrm{ad}(E_G)
\longrightarrow \bigoplus_{x \in D} \mathrm{ad}(E_G)_x/\mathcal{H}_x
\longrightarrow 0
\end{equation}
and
\begin{equation}\label{epq2}
0 \longrightarrow \mathrm{ad}^n_{\phi}(E_G) 
\longrightarrow \mathrm{ad}(E_G)
\longrightarrow \bigoplus_{x \in D} \mathrm{ad}(E_G)_x/\mathcal{H}^{\perp}_x
\longrightarrow 0,
\end{equation}
respectively.

\begin{Lem}\label{lem2}
The $D$-twisted connection $\widehat{\nabla}$ in \eqref{dn1} gives a holomorphic differential operator
$$\nabla \,\colon\, \mathrm{ad} (E_G)\,\longrightarrow \,\mathrm{ad} (E_G) \otimes K_X (D)\,=\,
\mathrm{ad} (E_G) \otimes K_X\otimes{\mathcal O}_X (D)\, .
$$
If $\widehat{\nabla}$ is a framed $G$-connection, then
$\nabla$ sends the subsheaf $\mathrm{ad}_{\phi} (E_G)$ in \eqref{epq1} to
$\mathrm{ad}^n_{\phi} (E_G) \otimes K_X (D)$, where $\mathrm{ad}^n_{\phi}(E_G)$ is constructed in
\eqref{epq2}.
\end{Lem}

\begin{proof}
Let $s$ be a holomorphic section of $\mathrm{ad} (E_G)$ defined over an open subset $U\, \subset\, X$.
Then $s$ defines a $G$--invariant holomorphic vector field on $p^{-1}(U)\, \subset\, E_G$ which is
vertical for the projection $p$ in \eqref{dp0}; this vertical vector field on $p^{-1}(U)$ will be denoted by
$\widetilde s$. Take any $t\, \in\, H^0(U,\, TX(-D))$. Let
$$\widetilde{t}\, :=\, \widehat{\nabla}(t)\, \in\, H^0(p^{-1}(U),\, TE_G(-\log p^{-1}(D)))^G$$
be the horizontal lift of $t$ for the $D$-twisted connection $\widehat\nabla$ in \eqref{dn1}. Now
consider the Lie bracket of vector fields
$$
[\widetilde{t},\, \widetilde{s}]\, \in\, H^0(p^{-1}(U),\, TE_G)\, .
$$
Note that $[\widetilde{t},\, \widetilde{s}]$ is $G$--invariant because both $\widetilde{s}$ and $\widetilde{t}$ are so. Furthermore,
$[\widetilde{t}, \, \widetilde{s}]$ is vertical for the projection $p$, because $\widetilde{s}$ is vertical and $\widetilde{t}$
is $G$--invariant. Indeed, for any holomorphic function $f$ on $U$, evidently $\widetilde{s}(f\circ p)\,=\, 0$ (recall that
$\widetilde s$ is vertical), and also we have $\widetilde{t}(f\circ p)$ to be $G$-invariant, so $\widetilde{s}(\widetilde{t}
(f\circ p))\, =\, 0$. Consequently, $[\widetilde{t},\, \widetilde{s}]$ produces a holomorphic section
of $\mathrm{ad}(E_G)$ over $U$; this section of $\mathrm{ad} (E_G)$ over $U$ will be
denoted by $[t,\, s]'$. Next note that a holomorphic function $f$ on $U$ satisfies
$$
[\widetilde{f\cdot t},\, \widetilde{s}]\,=\, (f\circ p)\cdot [\widetilde{t},\, \widetilde{s}]
- \widetilde{s} (f\circ p)\cdot \widetilde{t}\,=\, (f\circ p)\cdot [\widetilde{t},\, \widetilde{s}]
$$
because $\widetilde{s} (f\circ p)\,=\, 0$ (recall that $\widetilde s$ is a vertical
vector field). Consequently, there is a homomorphism
$$
\nabla\,\, :\,\, \mathrm{ad} (E_G) \longrightarrow \mathrm{ad} (E_G) \otimes K_X (D)
$$
uniquely defined by the equation
$$
\langle \nabla(s),\, t\rangle \,=\, [t,\, s]'\, ,
$$
where $s$ and $t$ are locally defined holomorphic sections of $\mathrm{ad} (E_G)$ and
$TX(-D)$ respectively, while $\langle -,\, -\rangle$ is the natural pairing
$TX(-D)\otimes K_X(D)\, \longrightarrow\, {\mathcal O}_X$.

Recall from Definition \ref{definition of framed connection} that
$\mathrm{res}_{x} ( \nabla)\, \in\, \mathcal{H}^{\perp}_{x}$. Therefore, from
the property of residues mentioned earlier it follows immediately that
$\widetilde{t}(x)\, \in\, {\mathcal H}^\perp_x$ for every $x\, \in\, D$. Now if
$s$ is a locally defined holomorphic section of $\mathrm{ad}_{\phi} (E_G)$, then
$\widetilde{s}(x)\, \in\, {\mathcal H}_x$. Next note that
\begin{equation}\label{lb}
[{\mathcal H}^\perp_x,\, {\mathcal H}_x]\, \subset\, {\mathcal H}^\perp_x\, ,
\end{equation}
because
$$
\widehat{\sigma}(x)([a,\, b]\otimes c)\,=\, \widehat{\sigma}(x)(a\otimes [b,\, c])
$$
for all $a,\, b,\, c\,\in\, \text{ad}(E_G)_x$ (this is derived using the given
condition on $\sigma$ that it is $G$-invariant). As a consequence of \eqref{lb}, the homomorphism
$\nabla$ maps the subsheaf $\mathrm{ad}_{\phi} (E_G)$ to $\mathrm{ad}^n_{\phi} (E_G) \otimes K_X (D)$.
\end{proof}

In view of Lemma \ref{lem2}, the following $2$-term complex of sheaves on $X$ is obtained
\begin{equation}\label{db}
\mathcal{D}_{\bullet} \ \colon\ \mathrm{ad}_{\phi} (E_G) \
\xrightarrow{\,\,\ \nabla\,\, \ } \ \mathrm{ad}^n_{\phi} (E_G) \otimes K_X (D)\, .
\end{equation}

The next lemma is straight-forward to prove.

\begin{Lem}\label{Lemma: infinitesimal deformation of framed connection} 
The infinitesimal deformations of the framed $G$-connection $(E_G, \,\widehat{\nabla},\, \phi)$
in \eqref{dn1} are parametrized by the elements of the first hypercohomology $\mathbb{H}^1(
\mathcal{D}_{\bullet})$ of the complex in \eqref{db}.
\end{Lem}

\section{Construction of the moduli space}\label{section: construction of moduli}

We now assume that 
$G\,=\,\mathrm{GL}(r,\mathbb{C})$. Fix a closed complex algebraic proper subgroup
$H_x \,\subsetneq\, G$ for each $x\in D$, and set $H\,=\, \{ H_x \}_{x \in D}$
to be the collection of subgroups indexed by the points of $D$.
For a framed vector bundle $(E,\, \phi)$, if $E_G$ is the principal $\mathrm{GL}(r,\mathbb{C})$-bundle associated 
to the vector bundle $E$, then $\text{ad}(E_G)\,=\, \mathrm{End}(E)$. Define
\begin{equation*}
\mathrm{End}_{\phi}(E)\,:=\, \mathrm{ad}_{\phi} (E_G) 
\quad 
\text{and}
\quad
\mathrm{End}^n_{\phi}(E)\,:=\,\mathrm{ad}^n_{\phi} (E_G)
\end{equation*}
(see Lemma \ref{lem2}).

\subsection{Definition of the moduli functors}\label{subsection: definition of framed moduli}

A framed $\mathrm{GL}(r,\mathbb{C})$-connection $(E,\,\phi, \, \nabla)$ on $X$ will be called {\it simple} if
\begin{equation*}
\mathrm{ker}\left(H^0(X,\,\mathcal{E}nd_{\phi}(E) )\,
\xrightarrow{\nabla}\, H^0(X,\, \mathcal{E}nd_{\phi}^n(E)\otimes K_X(D))\right) \,=\,0\, .
\end{equation*}

\begin{Def}\label{2019.10.28.16.49}
Define a stack $\mathcal{M}^{H}_{\rm{FC}}(d)$
of simple framed $\mathrm{GL}(r,\mathbb{C})$-connections, for $H$,
by breaking into the following two cases.
\begin{itemize}
\item If $\mathbb{C}^* \cdot{\rm Id}\, \not\subset\, H_x$ for some $x\,\in\, D$,
then define a stack $\mathcal{M}^{H}_{\rm{FC}}(d)$ 
over the category of locally Noetherian schemes over $\Spec\mathbb{C}$
whose objects are quadruples $(S,\, E,\, \phi\,= \,\{ \phi_{x\times S} \}_{x \in D} ,\, \nabla)$
of the following type:
\begin{itemize}
\item[(1)] $S$ is a locally Noetherian scheme over $\Spec\mathbb{C}$,
and $E\, \longrightarrow\, X\times S$ is a vector bundle of rank $r$
and $\deg(E|_{X\times s})\,=\,d$ for any geometric point $s$ of $S$.

\item[(2)] $\phi_{x\times S}$ is a section of the structure map
\begin{equation*}
\mathrm{Isom}_{S} (\mathcal{O}_{x\times S}^{\oplus r} ,\, \, E|_{x\times S} )/ (H_x \times S)
 \longrightarrow x\times S.
\end{equation*} 
Here the action of the group scheme $H_x \times S$ over $S$
on $\mathrm{Isom}_{S} (\mathcal{O}_{x\times S}^{\oplus r},\, \, E|_{x\times S})$
is the restriction of the natural transitive action of 
the group scheme $\mathrm{GL}(r,\mathbb{C}) \times S$ over $S$ on 
$\mathrm{Isom}_{S} (\mathcal{O}_{x\times S}^{\oplus r},\, \, E|_{x\times S})$
given by the standard action of $\text{GL}(r,{\mathbb C})$ on $\mathcal{O}_{x\times S}^{\oplus r}$.
Define a $S$-scheme $\widetilde{S}$ 
and a map $\widetilde{S} \,\longrightarrow \,\mathrm{Isom}_{S} 
(\mathcal{O}_{x\times S}^{\oplus r},\,\, E|_{x\times S} )$
such that the diagram
\begin{equation*}
\xymatrix{
x\times\widetilde{S} \ar[d] \ar[r]
& \mathrm{Isom}_{S} (\mathcal{O}_{x\times S}^{\oplus r}, \,\, E|_{x\times S}) \ar[d] \\
x\times S \ar[r]^-{\phi_{x\times S}}
& \mathrm{Isom}_{S} (\mathcal{O}_{x\times S}^{\oplus r}, \,\, E|_{x\times S})/(H_x \times S)
}
\end{equation*}
is Cartesian.
Let $$\widetilde{\phi}_{x \times \widetilde{S}}\,\colon\,
 \mathcal{O}_{x\times \widetilde{S}}^{\oplus r}
\,\xrightarrow{\ \sim\ }\,
E_{\widetilde{S}}|_{x\times \widetilde{S}}
$$
be the isomorphism given by the map $\widetilde{S}\,\longrightarrow\,
\mathrm{Isom}_{S}(\mathcal{O}_{x\times S}^{\oplus r},\, \,E|_{x\times S})$.

\item[(3)] $\nabla \,\colon\, E \,\longrightarrow\, E \otimes K_{X}(D)$
is a relative connection, relative to $S$.

\item[(4)]
Let $\mathrm{res}_{x\times \widetilde{S}}(\nabla_{\widetilde{S}}) 
\,\in\, \mathrm{End}(E_{\widetilde{S}})|_{x\times \widetilde{S}}$
be the residue 
of the induced connection
$\nabla_{\widetilde{S}}\,\colon \,E_{\widetilde{S}}\,\longrightarrow\, E_{\widetilde{S}} 
\otimes K_{X}(D)$. Then
$\widetilde{\phi}^{-1}_{x \times \widetilde{S}} \circ
\mathrm{res}_{x\times\widetilde{S}}(\nabla_{\widetilde{S}})
\circ\widetilde{\phi}_{x \times\widetilde{S}}
\,\,\in\,\, \mathfrak{h}^{\perp} \otimes \mathcal{O}_{\widetilde{S}}$.

\item[(5)] For each point $s\,\in\, S$, 
the framed $\mathrm{GL}(r,\mathbb{C})$-connection $(E_s,\,\phi_s, \, \nabla_s)$ is simple.
Recall that $(E_s,\,\phi_s, \, \nabla_s)$ is simple if
\begin{equation*}
\mathrm{ker}\left(H^0(X,\,\mathcal{E}nd_{\phi_s}(E_s) )\,
\xrightarrow{\nabla_s}\, H^0( X,\, \mathcal{E}nd_{\phi_s}^n(E_s)\otimes K_X(D))\right) \,=\,0\, .
\end{equation*}

\end{itemize}
A morphism 
\begin{equation*}
(S,\, E,\, \phi,\,\nabla )
\longrightarrow (S',\, E',\, \phi',\, \nabla' )
\end{equation*}
in $\mathcal{M}^{H}_{\rm{FC}}$ is a Cartesian square 
\begin{equation*}
\xymatrix{
 E \ar[r]^-{\sigma} \ar[d] & E' \ar[d] \\
S \ar[r]^{\widetilde{\sigma}} & S'
}
\end{equation*}
such that the diagram
\begin{equation*}
\xymatrix{
 E \ar[r]^-{\nabla} \ar[d]_{\cong}^{\sigma} 
 & E \otimes K_{X}(D) \ar[d]_{\cong}^-{\sigma \otimes \mathrm{Id}} \\
 E' \times_{S'} S \ar[r]^-{\nabla' } 
 & (E' \times_{S'} S) \otimes K_{X}(D)
 }
\end{equation*}
is commutative and 
$(\widetilde{\phi}'_{x \times \widetilde{S}})^{-1}\circ \sigma_{\widetilde{S}} 
\circ \widetilde{\phi}_{x \times \widetilde{S}}\,
\in\, H_{x} \times\widetilde{S}$ for each $x \,\in\, D$.

\item If $\mathbb{C}^* e \,\subset\, H_x$ for all $x\,\in\, D$,
then define $\mathcal{M}^{H}_{\rm{FC}}(d)$ 
to be the stackification of $\text{pre-}\mathcal{M}^{H}_{\rm{FC}}(d)$ (see \cite[Theorem 4.6.5]{Olsson}).
Here $\text{pre-}\mathcal{M}^{H}_{\rm{FC}}(d)$ is the fibered category
over the category of locally Noetherian schemes over $\Spec\,\mathbb{C}$
whose objects are quadruples $(S,\, E,\, \phi\,=\, \{ \phi_{x\times S} \}_{x \in D} ,\,\nabla) $
that satisfy (1), (3) and (4) as above as well as the following $(2)'$ and $(5)'$:
\begin{itemize}
\item[$(2)'$:] 
$\phi_{x\times S}$ is a section of the structure map
\begin{equation*}
\mathrm{Isom}_{S} 
(\mathcal{O}_{x\times S}^{\oplus r},\, \, E|_{x\times S}
)/ (H_x \times S)
\, \longrightarrow\, x\times S.
\end{equation*}
Here the action of the group scheme $H_x \times S$, over $S$,
on $\mathrm{Isom}_{S} (\mathcal{O}_{x\times S}^{\oplus r} ,\, \, E|_{x\times S} )$
is the restriction of the natural transitive group action of 
the group scheme $\mathrm{GL}(r,\mathbb{C}) \times S$ over $S$ on 
$\mathrm{Isom}_{S} 
(\mathcal{O}_{x\times S}^{\oplus r},\,\, E|_{x\times S})$ given by the
standard action of $\text{GL}(r,{\mathbb C})$ on $\mathcal{O}_{x\times S}^{\oplus r}$.
Define a $S$-scheme $\widehat{S}$ 
and a map $\widehat{S} \,\longrightarrow\, \mathrm{Isom}_{S} 
(\mathcal{O}_{x\times S}^{\oplus r},\, \, E|_{x\times S})
/ (\mathbb{C}^* e \times S)$
such that the diagram
\begin{equation*}
\xymatrix{
x\times\widehat{S} \ar[d] \ar[r]
& \mathrm{Isom}_{S} (\mathcal{O}_{x\times S}^{\oplus r}, \,\, E|_{x\times S})
/ (\mathbb{C}^* e \times S) \ar[d] \\
x\times S \ar[r]^-{\phi_{x\times S}}
& \mathrm{Isom}_{S} (\mathcal{O}_{x\times S}^{\oplus r}, \,\,E|_{x\times S} )/(H_x \times S)
}
\end{equation*}
is Cartesian.
Denote by $\widehat{\phi}_{x \times \widehat{S}}\,\colon\, 
\mathbb{P}(\mathcal{O}_{x\times\widehat{S}}^{\oplus r})
\,\stackrel{\sim}{\longrightarrow}\, 
\mathbb{P}( E_{\widehat{S}}|_{x\times\widehat{S}})
$
the isomorphism given by 
the map $\widehat{S} \,\longrightarrow \,\mathrm{Isom}_{S} 
(\mathcal{O}_{x\times S}^{\oplus r}, \,\,E|_{x\times S})
/ (\mathbb{C}^* e \times S)$.

\item[$(5)'$:] $(E_s,\,\phi_s, \, \nabla_s)$ is simple for each point $s \,\in\, S$, that is, 
\begin{equation*}
\mathrm{ker}(H^0(X,\,\mathcal{E}nd_{\phi_s}(E_s) )\,
\xrightarrow{\nabla_s}\, H^0( X,\, \mathcal{E}nd_{\phi_s}^n(E_s)\otimes K_X(D))) \,=\,
\mathbb{C}.
\end{equation*}
\end{itemize}
A morphism
\begin{equation*}
(S,\, E,\, \phi,\, \nabla )
\longrightarrow (S',\, E',\, \phi',\, \nabla' )
\end{equation*}
in $\text{pre-}\mathcal{M}^{H}_{\rm{FC}}(d)$
is a triple $(\mathcal{L},\, \sigma ,\, \widetilde{\sigma})$, where
$\mathcal{L}$ is a line bundle on $S'$ and $\sigma ,\, \widetilde{\sigma}$ are
maps that fit in a Cartesian square 
\begin{equation*}
\xymatrix{
 E \ar[r]^-{\sigma} \ar[d] & E' \otimes \mathcal{L} \ar[d] \\
S \ar[r]^{\widetilde{\sigma}} & S'
}
\end{equation*}
such that the diagram
\begin{equation*}
\xymatrix{
 E \ar[r]^-{\nabla} \ar[d]_{\cong}^{\sigma} 
 & E \otimes K_{X}(D) \ar[d]_{\cong}^-{\sigma \otimes \mathrm{Id}}\\
 (E' \otimes\mathcal{L})\times_{S'} S \ar[r]^-{\nabla' \otimes\mathcal{L}} 
 & ((E' \otimes\mathcal{L}) \times_{S'}S) \otimes K_{X}(D)
}
\end{equation*}
is commutative, and 
$$(\widehat{\phi}'_{x \times \widehat{S}} )^{-1}\circ 
\overline{\sigma}_{x\times \widehat{S}} 
\circ \widehat{\phi}_{x \times \widehat{S}}
\,\in\, (H_{x}/\mathbb{C}^*e ) \times \widehat{S}$$ for each $x \,\in\, D$,
where $\overline{\sigma}_{x\times\widehat{S}}\,
\colon\,\mathbb{P}(E_{\widehat{S}}|_{x \times\widehat{S}})\,\longrightarrow\, 
\mathbb{P}(E_{\widehat{S}}'|_{x \times\widehat{S}})$ 
is induced by $\sigma$.
\end{itemize}
\end{Def}

We say that $\sigma$ is an {\it automorphism} of a framed $G$-connection $(E,\,\phi,\,\nabla)$
if $\sigma$ is a holomorphic automorphism of the vector bundle $E$ on $X$ 
such that 
the diagram
\begin{equation*}
\xymatrix{
 E \ar[r]^-{\nabla} \ar[d]^{\sigma} 
 & E \otimes K_{X}(D) \ar[d]^-{\sigma \otimes \mathrm{Id}}\\
E \ar[r]^-{\nabla} 
& E \otimes K_{X}(D) 
}
\end{equation*}
is commutative and 
$ \sigma|_{x}\circ \phi_{x}$ coincides with $\phi_{x}$ in the
quotient $\mathrm{Isom}(\mathcal{O}_{x}^{\oplus r} , \,\, E|_x)/H_x$
for each $x\,\in\, D$. Denote by $\mathrm{Aut}(E,\,\phi,\,\nabla)$
the space of all automorphisms of a framed $G$-connection $(E,\,\phi,\,\nabla)$.

\begin{Prop}\label{2020.1.22.10.51}
Assume that $\mathbb{C}^* \cdot{\rm Id}\, \not\subset\, H_x$ for some $x\,\in\, D$.
Let $(E,\,\phi ,\, \nabla)$ be a simple framed $G$-connection over $X$ (see Definition
\ref{2019.10.28.16.49}(5)). Then $\mathrm{Aut}(E,\phi,\nabla)$ is a finite group.
\end{Prop}

\begin{proof}
The space $\mathrm{Aut}(E,\,\phi,\,\nabla)$ has the structure of a group scheme 
of finite type over $\mathbb{C}$.
We can see that 
the tangent space of $\mathrm{Aut}(E,\,\phi,\,\nabla)$ at the identity element is
isomorphic to 
\begin{equation*}
\mathrm{ker}(H^0(X,\,\mathcal{E}nd_{\phi}(E) )\,
\xrightarrow{\,\,\nabla\,\,}\, H^0( X,\, \mathcal{E}nd_{\phi}^n(E)\otimes K_X (D)) ),
\end{equation*}
which is zero because $(E,\,\phi,\,\nabla)$ is simple.
Consequently, $\mathrm{Aut}(E,\,\phi,\,\nabla)$ is a finite group.
\end{proof}

\begin{Prop}\label{2019.12.24.23.58}
Assume that $H_x\,=\,\{ e\}$ for all $x\,\in\, D$.
Let $(E,\,\phi ,\,\nabla)$ be a simple framed $G$-connection over $X$
associated to $\{H_x\}_{x\in D}$.
Then $\mathrm{Aut}(E,\,\phi,\,\nabla) \,=\, \{ \mathrm{Id}_{E}\}$.
\end{Prop}

\begin{proof}
An automorphism $\sigma\,\in \,\mathrm{Aut}(E,\,\phi,\,\nabla)$
is an automorphism of the vector bundle $E$
such that $\nabla \circ \sigma \,=\, \sigma \circ \nabla$ and
$\phi_x\circ\sigma|_x\circ\phi_x^{-1}\,\in\, H_x$
for all $x \,\in\, D$.
Since $H_x\,=\,\{ e\}$ by the assumption,
it follows that $\sigma|_x\,=\,\mathrm{Id}|_{E|_x}$ for all $x\,\in\, D$.
Now set
\begin{equation*}
\widetilde{\sigma} \,:=\, \mathrm{Id}_{E} - \sigma \,\colon\, E \,\longrightarrow\, E\, .
\end{equation*}
Then $\widetilde{\sigma}|_x\,=\,0$
for all $x\,\in\, D$, and it is straightforward to check that $\nabla \circ \widetilde{\sigma} -
\widetilde{\sigma} \circ \nabla\,=\,0$, that is,
\begin{equation*}
\widetilde{\sigma} \,\in\,
\mathrm{ker}(H^0(X,\,\mathcal{E}nd_{\phi}(E))\,
\xrightarrow{\,\nabla\,}\, H^0(X,\, \mathcal{E}nd^n_{\phi}(E)\otimes K_X(D)))\, .
\end{equation*}
Since $(E,\,\phi,\,\nabla)$ is simple, it follows
that $\widetilde{\sigma}\,=\,0$, and hence $\sigma \,=\, \mathrm{Id}_{E}$.
\end{proof}

\begin{Prop}\label{2020.1.22.11.10}
Assume that $\mathbb{C}^* \cdot{\rm Id}\, \subset\, H_{x}$ for all $x\,\in\, D$.
If $(E,\,\phi ,\, \nabla)$ is a simple framed $G$-connection over $X$, 
then the quotient
$\mathrm{Aut}(E,\phi,\nabla)/(\mathbb{C}^* \cdot{\rm Id})$ is a finite group.
\end{Prop}

\begin{proof}
The tangent space of $\mathrm{Aut}(E,\,\phi,\,\nabla)/(\mathbb{C}^* \cdot{\rm Id})$ is
zero, because we have $(E,\,\phi,\,\nabla)$ to be simple. Consequently, 
$\mathrm{Aut}(E,\,\phi,\,\nabla)/(\mathbb{C}^* \cdot{\rm Id})$ is a finite group.
\end{proof}

\subsection{Representation of moduli functors as Deligne--Mumford stacks}\label{subsection: construction of moduli stack}

\begin{Prop}\label{2020.2.24.14.52}
The stack $\mathcal{M}^{H}_{\rm{FC}}(d)$ in Definition \ref{2019.10.28.16.49}
is a Deligne--Mumford stack.
\end{Prop}

\begin{proof}
Fix a very ample line bundle $\mathcal{O}_{X}(1)$ on the curve $X$.
Define a polynomial $\theta_d(m)$ in $m$
to be $$\theta_d(m)\,=\,r d_X m +d +r(1-g)\, ,$$
where $d_X\,:=\, \deg \mathcal{O}_X(1)$
and $g$ is the genus of $X$.
Let
\begin{equation}\label{sd}
\Sigma_{m_0}^{d}
\end{equation}
denote the fibered category whose objects are
simple framed $\mathrm{GL}(r,\mathbb{C})$-connections $(E,\,\phi,\,\nabla)$ on $X\times S$
such that 
\begin{itemize} 
\item $H^1(X,\, E_s (m_0-1))\,=\,0$ for each $s\,\in\, S$, and

\item $\chi(E_s(m))\,= \,\theta_d(m)$ for each $s \,\in\, S$ and all $m \,\in\, \mathbb{Z}$.
\end{itemize}
The fibered categories $\Sigma_{m_0}^{d}$ in \eqref{sd} form an open covering of $\mathcal{M}^{H}_{\rm{FC}}
(d)$. So we only have to prove that each $\Sigma_{m_0}^{d}$ is a Deligne--Mumford stack.

Let $$\mathcal{O}_{X\times Q^d_{m_0}}(-m_0)^{\oplus \theta_d(m_0)}\,\longrightarrow\, 
\mathcal{E}$$
be the universal quotient sheaf of the Quot-scheme 
$\mathrm{Quot}^{\theta_d}_{(\mathcal{O}_X(-m_0)^{\oplus \theta_d(m_0)}/X)}$.
Define the open subset $Q^d_{m_0}$ of $\mathrm{Quot}^{\theta_d}_{(\mathcal{O}_X(-m_0)^{\oplus
\theta_d(m_0)}/X)}$ by
$$
Q^d_{m_0} \, := \,\left\{
s \,\in\, \mathrm{Quot}^{\theta_d}_{(\mathcal{O}_X(-m_0)^{\oplus \theta_d(m_0)}/X)}
\ \middle| \
\begin{array}{l}
\text{(i) $h^0(X,\, \mathcal{E}_s(m_0) )\,=\,\theta_d(m_0)$;}\\
\text{(ii) $h^i(X,\, \mathcal{E}_s(m_0-i) )\,=\,0$\,\, for all $i\,>\,0$; and}\\
\text{(iii) $\mathcal{E}_s$ is locally free.}
\end{array}
\right\}.
$$
There is a locally free $\mathcal{O}_{Q^d_{m_0}}$--module $\mathscr{H}_D$
such that $\boldsymbol{V}^*(\mathscr{H}_D)\,:=\,\mathrm{Spec}(\mathrm{Sym}^*\mathscr{H}_D)$ 
represents the functor
\begin{equation*}
 S \,\,\longmapsto \,\,\bigoplus_{x \in D}
 \mathrm{Hom}_{x \times S} 
(\mathcal{O}^{\oplus r}_{x\times S},\, \,
 \mathcal{E}_{X\times S}|_{x\times S})
\,\, \in\,\, (\mathit{Sets})
\end{equation*}
for any Noetherian schemes $S$ over $Q^d_{m_0}$. There is a universal family 
\begin{equation}\label{2019.11.27.14.15}
\varphi^{x}\, \colon\,
\mathcal{O}^{\oplus r}_{x \times \boldsymbol{V}^*(\mathscr{H}_D)}
\,\longrightarrow\, 
\mathcal{E}_{X \times \boldsymbol{V}^*(\mathscr{H}_D)} |_{x \times \boldsymbol{V}^*(\mathscr{H}_D)}
\end{equation}
for every $x \,\in \,D$. Define $\widetilde{Q}^{d}_{m_0}$ as follows:
\begin{equation*}
\widetilde{Q}^{d}_{m_0}\,:=\, \left\{
s \,\in\, \boldsymbol{V}^*(\mathscr{H}_D)
\ \middle| \
\text{$\mathrm{coker} (\varphi^{x}_s)\,=\,0$} \,
\right\}.
\end{equation*}
Consider the map $\widetilde{Q}^{d}_{m_0} \, \longrightarrow \, Q^d_{m_0}$.
For each Noetherian scheme $S$ over $Q^d_{m_0}$,
the natural transitive group action of $G \times S$
on 
$$ 
\mathrm{Isom}_{x\times S}
(\mathcal{O}^{\oplus r}_{x\times S},\,\,
 \mathcal{E}_{X\times S}|_{x\times S})
$$ 
induces an action on $\widetilde{Q}^{d}_{m_0}$ of the group scheme
$\left(\prod_{x \in D}G\right)\times Q^d_{m_0}$ over $Q^d_{m_0}$.
The group scheme $\left(\prod_{x \in D}H_x\right)\times Q^d_{m_0}$ acts 
on $\widetilde{Q}^{d}_{m_0}$
by restricting this action of $(\prod_{x \in D}G)\times Q^d_{m_0}$ 
on $\widetilde{Q}^{d}_{m_0}$. 
Set 
$$
\widetilde{Q}^{d,H}_{m_0}\, :=\, 
\widetilde{Q}^{d}_{m_0}\Big/ \left((\prod_{x \in D}H_x)\times Q^d_{m_0}\right).
$$
Let $\widetilde{\mathcal{E}}$ be the pull-back of the family $\mathcal{E}$
under the map $X \times \widetilde{Q}^{d,H}_{m_0}
\, \longrightarrow \, X \times Q^d_{m_0}$.
We have a family $\widetilde{\phi}_x$ of sections of 
\begin{equation*}
\mathrm{Isom}_{\widetilde{Q}^{d,H}_{m_0}} 
(\mathcal{O}_{\widetilde{Q}^{d,H}_{m_0}}^{\oplus r}, \,\,
\widetilde{\mathcal{E}}|_{x\times \widetilde{Q}^{d,H}_{m_0}} )
/ (H_x \times \widetilde{Q}^{d,H}_{m_0})
\, \longrightarrow \, \widetilde{Q}^{d,H}_{m_0}
\end{equation*} 
induced by $\varphi^x$.
Put $\widetilde{\phi}\,:=\,\{\widetilde{\phi}_x \}_{x \in D}$.

Let
\begin{equation}\label{epi}
\pi\,\colon\, X \times \widetilde{Q}^{d,H}_{m_0} \,\longrightarrow\, 
\widetilde{Q}^{d,H}_{m_0}
\end{equation}
be the natural projection map. Let $\mathrm{At}(\widetilde{\mathcal{E}})$ be the
Atiyah bundle for $\widetilde{\mathcal{E}}$. Then there is a short exact sequence
\begin{equation*}
0 \longrightarrow 
\mathcal{E}nd(\widetilde{\mathcal{E}})\longrightarrow 
\mathrm{At}(\widetilde{\mathcal{E}})\xrightarrow{\mathrm{symb}_1}
T_{X\times 
\widetilde{Q}^{d,H}_{m_0}/\widetilde{Q}^{d,H}_{m_0}}
\longrightarrow 0.
\end{equation*}
Set $\mathrm{At}_{D} (\widetilde{\mathcal{E}}) \,\,:=\,\,
\mathrm{symb}_1^{-1}\left(T_{X\times \widetilde{Q}^{d,H}_{m_0}
/\widetilde{Q}^{d,H}_{m_0}}(-D\times \widetilde{Q}^{d,H}_{m_0})
\right)$.
The natural short exact sequences of Atiyah bundles induces an exact sequence
\begin{equation*}
\xymatrix{
0 \ar[r] 
& \mathcal{E}nd(\widetilde{\mathcal{E}}) \otimes K_{X} \ar[r]\ar[d]
&\mathrm{At}(\widetilde{\mathcal{E}}) \otimes K_{X} \ar[r]^-{\mathrm{symb}_1} \ar[d]
&\mathcal{O}_{X \times \widetilde{Q}^{d,H}_{m_0}}\ar[r] \ar[d]
&0 \\
0 \ar[r] 
& \mathcal{E}nd(\widetilde{\mathcal{E}}) \otimes K_{X}(D) \ar[r]\ar[d]^-{p'}
&\mathrm{At}_{D} (\widetilde{\mathcal{E}}) \otimes K_X(D)
\ar[r]^-{\mathrm{symb}_1^{D}} \ar[d]^{p}
& \mathcal{O}_{X \times \widetilde{Q}^{d,H}_{m_0}} \ar[r] 
&0 \\
&\mathcal{E}nd(\widetilde{\mathcal{E}}) 
\otimes K_{X}(D)|_{D\times \widetilde{Q}^{d,H}_{m_0}} \ar[r]^-q_-{\cong}
&\left(
\mathrm{At}_{D} (\widetilde{\mathcal{E}}) \otimes K_X(D)
\right)
\Big/\left( \mathrm{At}(\widetilde{\mathcal{E}}) \otimes K_{X} \right)\rlap{.}
&&
}
\end{equation*}
In particular, there are two compositions of maps
\begin{equation}\label{eq:2021_12_17_12_9}
\begin{aligned}
&\mathrm{At}_{D} (\widetilde{\mathcal{E}}) \otimes K_X(D) 
\xrightarrow{\,\,\, p\,\,} 
\left(
\mathrm{At}_{D} (\widetilde{\mathcal{E}}) \otimes K_X(D)
\right)
\Big/\left( \mathrm{At}(\widetilde{\mathcal{E}}) \otimes K_{X} \right) \\
&\qquad \xrightarrow{\, q^{-1}\, } 
\mathcal{E}nd(\widetilde{\mathcal{E}}) 
\otimes K_{X}(D)|_{D\times \widetilde{Q}^{d,H}_{m_0}}
\xrightarrow{\, \mathrm{res}_{D \times \widetilde{Q}^{d,H}_{m_0}} \,}
\mathrm{End} (\widetilde{\mathcal{E}})|_{D\times \widetilde{Q}^{d,H}_{m_0}}
\end{aligned}
\end{equation}
and
\begin{equation}\label{eq:2021_12_17_12_10}
\mathcal{E}nd(\widetilde{\mathcal{E}}) \otimes K_{X}(D)
 \xrightarrow{\, p'\, } 
\mathcal{E}nd(\widetilde{\mathcal{E}}) 
\otimes K_{X}(D)|_{D\times \widetilde{Q}^{d,H}_{m_0}}
\xrightarrow{\, \mathrm{res}_{D \times \widetilde{Q}^{d,H}_{m_0}} \,}
\mathrm{End} (\widetilde{\mathcal{E}})|_{D\times \widetilde{Q}^{d,H}_{m_0}};
\end{equation}
here $\mathrm{res}_{D \times \widetilde{Q}^{d,H}_{m_0}}$
is the residue map
\begin{equation*}
\mathrm{res}_{D \times \widetilde{Q}^{d,H}_{m_0}}
\,\,\colon\,\, \mathrm{End}(\widetilde{\mathcal{E}}) 
\otimes K_{X}(D)|_{D\times \widetilde{Q}^{d,H}_{m_0} }
\longrightarrow
\mathrm{End} (\widetilde{\mathcal{E}})|_{D\times \widetilde{Q}^{d,H}_{m_0}}.
\end{equation*}
Using the family $\widetilde{\phi}$ of framings, and 
the Lie subgroups $H_x$ for each $x\,\in \,D$, we may define 
a subsheaf $\mathcal{H}^{\perp}_{x \times \widetilde{Q}^{d,H}_{m_0}} 
\,\subset\, \mathrm{End} (\widetilde{\mathcal{E}})|_{D\times \widetilde{Q}^{d,H}_{m_0}}$
as done in \eqref{2020.6.24.10.38}. Define subsheaves 
$\mathrm{At}^{\widetilde{\phi}}_{D} (\widetilde{\mathcal{E}})
\,\subset\, \mathrm{At}_{D} (\widetilde{\mathcal{E}})$
and $\mathcal{E}nd_{\widetilde{\phi}}^n(\widetilde{\mathcal{E}}) 
\,\subset\, \mathcal{E}nd(\widetilde{\mathcal{E}})$ 
such that 
$\mathrm{At}^{\widetilde{\phi}}_{D} (\widetilde{\mathcal{E}}) \otimes K_X(D)$
is the inverse image of $\mathcal{H}^{\perp}_{x \times \widetilde{Q}^{d,H}_{m_0}}$
under the composition of maps in \eqref{eq:2021_12_17_12_9}
and $\mathcal{E}nd_{\widetilde{\phi}}^n(\widetilde{\mathcal{E}})\otimes K_X(D)$
is the inverse image of $\mathcal{H}^{\perp}_{x \times \widetilde{Q}^{d,H}_{m_0}}$
under the composition of maps in \eqref{eq:2021_12_17_12_10}. 

By \cite[Theorem 7.7.6]{Groth1}, there is a unique
coherent sheaf $\mathscr{H}$ on $\widetilde{Q}^{d,H}_{m_0}$ (up to isomorphism) such that
\begin{equation}\label{eq:2021.11.25.16.17}
(\pi_{Q'})_*\left( \left( \mathrm{At}^{\widetilde{\phi}}_{D} (\widetilde{\mathcal{E}})
\otimes K_X(D)\right)
\otimes_{\mathcal{O}_{Q'}} M \right)\,\,\cong\,\,
\mathcal{H}{\rm om}_{\mathcal{O}_{Q'}} (\mathscr{H}_{Q'} , M)
\end{equation}
for any $\widetilde{Q}^{d,H}_{m_0}$-scheme $Q'$ 
and any quasi-coherent sheaf $M$ on $Q'$.
Here $\pi_{Q'}$ is the natural projection $X \times Q' \,\longrightarrow\, Q'$.
Set 
$$\boldsymbol{V} (\mathscr{H}) \,:=\,
\mathrm{Spec}(\mathrm{Sym}^*(\mathscr{H}))\, .$$
There is a natural morphism $\varphi \,\in\,
\mathcal{H}{\rm om}_{\mathcal{O}_{\boldsymbol{V} (\mathscr{H})}} 
(\mathscr{H}_{\boldsymbol{V} (\mathscr{H})} ,\, \mathcal{O}_{\boldsymbol{V} (\mathscr{H})}) $
by the definition of $\boldsymbol{V} (\mathscr{H})$. In view of
the isomorphism in \eqref{eq:2021.11.25.16.17}, there is an element 
$\varphi' \,\in\, \pi_*( ( \mathrm{At}^{\widetilde{\phi}}_{D} (\widetilde{\mathcal{E}})
\otimes K_X(D) )_{\boldsymbol{V} (\mathscr{H})})$
corresponding to $\varphi$. The morphism $\mathrm{symb}_1^{D}$ induces a morphism 
$$
(\pi_{\boldsymbol{V} (\mathscr{H})})_*
\left( \left( \mathrm{At}^{\widetilde{\phi}}_{D} (\widetilde{\mathcal{E}}) \otimes K_X(D)
\right)_{\boldsymbol{V} (\mathscr{H})}\right)
\, \longrightarrow \, 
(\pi_{\boldsymbol{V} (\mathscr{H})})_*(\mathcal{O}_{X \times \boldsymbol{V} (\mathscr{H})})
\cong \mathcal{O}_{\boldsymbol{V} (\mathscr{H})}.
$$
Using this morphism, there is a function $f_{\mathrm{symb}_1^{D}}$ 
on $\boldsymbol{V} (\mathscr{H})$
corresponding to $\varphi'$.
Denote by $I_{\mathrm{symb}_1^{D}}$ the ideal sheaf of
$\mathcal{O}_{\boldsymbol{V} (\mathscr{H})}$ generated by 
$f_{\mathrm{symb}_1^{D}}-1$.
Put
$$Z^{d,H}_{m_0} \,:=\,
\mathrm{Spec}(\mathcal{O}_{\boldsymbol{V} (\mathscr{H})}\, /\, I_{\mathrm{symb}_1^{D}})\, .$$
Also denote by $\widetilde{\mathcal{E}}$ 
the pull-back of $\widetilde{\mathcal{E}}$ under the 
natural morphism $X \times Z^{d,H}_{m_0} \,\longrightarrow\, 
X\times \widetilde{Q}^{d,H}_{m_0,\mathrm{spl}}$,
and let
\begin{equation*}
\widetilde{\nabla}\,\, \colon\,\, \widetilde{\mathcal{E}} \longrightarrow 
\widetilde{\mathcal{E}} \otimes_{\mathcal{O}_{X \times \boldsymbol{V}^* (\mathscr{H})_{\epsilon}}}
K_{X}(D)
\end{equation*}
be a universal relative connection on $\widetilde{\mathcal{E}}$,
which is determined by $\varphi'$.
Define the open subset $(Z^{d,H}_{m_0})'$
of $Z^{d,H}_{m_0}$ by
\[
 (Z^{d,H}_{m_0})'\,\, =\,\,
\left\{ s\in Z^{d,H}_{m_0}\,\,\middle|\,\,\,
 \text{$(\widetilde{\mathcal E},\widetilde{\phi},\widetilde{\nabla})|_{X\times\{s\}}$ is a
 simple framed connection} \right\}
\]
and denote by $\big(\widetilde{\mathcal{E}},\, \widetilde{\phi}\,=\, 
\big\{[\widetilde{\phi}_{x\times (Z^{d,H}_{m_0})'}] \big\}_{x\in D},\, 
\widetilde{\nabla}\big)$
a universal family of $m_0$-regular simple framed $G$-connections
on $X\times (Z^{d,H}_{m_0})'$.
Here $m_0$-regular means 
$H^1(X,\, \widetilde{\mathcal{E}}_s (m_0-1))\,=\,0$ for each $s\,\in\, (Z^{d,H}_{m_0})'$.

Now consider the case where $\mathbb{C}^* \cdot{\rm Id}\, \not\subset\, H_x$ for some $x\,\in\, D$.
There exists an action of $\mathrm{GL}(\theta_d(m_0),\mathbb{C})$ 
on $(Z^{d,H}_{m_0})'$ given by
\begin{equation*}
\left(\mathcal{O}_{X\times S}^{\oplus \theta_{d}(m)}\,\stackrel{q}{\longrightarrow}\, E,\,
\phi ,\,\nabla\right)
\,\,\longmapsto\,\,
\left(\mathcal{O}_{X\times S}^{\oplus \theta_{d}(m)}\,\xrightarrow{\, q \circ g \,}
\, E,\, \phi , \, \nabla\right)
\end{equation*}
on $S$-points for $g \,\in\, \mathrm{GL}(\theta_d(m_0),\mathbb{C})_S$. Consider the map
\begin{equation*}
\begin{aligned}
(Z^{d,H}_{m_0})' &\,\longrightarrow\, \Sigma_{m_0}^{d} \\
(\mathcal{O}_{X\times S}^{\oplus \theta_{d}(m)}\,\xrightarrow{\,\,q\,\,}\, E, \,\phi ,\,
\nabla ) &\,\longmapsto\, (S,\, E,\, \phi ,\, \nabla).
\end{aligned}
\end{equation*}
This map gives an isomorphism
\[
\Sigma_{m_0}^{d}\,\,\cong\,\, 
\big[ \, (Z^{d,H}_{m_0})' \, \big/\, \mathrm{GL}(\theta_d(m_0),\mathbb{C}) \, \big] \,. 
\]
Here 
$\big[ \, (Z^{d,H}_{m_0})' \, \big/ \, \mathrm{GL}(\theta_d(m_0),\mathbb{C}) \, \big]$
is a quotient stack. 
Consequently, $\Sigma_{m_0}^{d}$ is an algebraic stack.
Using Proposition \ref{2020.1.22.10.51} it follows that
$\Sigma_{m_0}^{d}$ is in fact a Deligne--Mumford stack (see \cite[Corollary 8.4.2]{Olsson}).

Next we consider the case where $\mathbb{C}^* \cdot{\rm Id}\, \subset\, H_x$ for
every $x\,\in\, D$.
The $\mathbb{C}^*$-action on $(Z^{d,H}_{m_0})'$ is trivial,
because $\mathbb{C}^* \cdot{\rm Id}\, \subset\, H_x$ for all $x\,\in\, D$.
There exists a natural action of $\mathrm{PGL}(\theta_d(m_0),\,\mathbb{C})$ 
on $(Z^{d,H}_{m_0})'$.
Define a map
\begin{equation*}
\begin{aligned}
(Z^{d,H}_{m_0})' &\,\,\longrightarrow\,\, \Sigma_{m_0}^{d}\\
(\mathcal{O}_{X\times S}^{\oplus \theta_{d}(m)}\xrightarrow{\,\,q\,\,} E,\, \phi ,\,\nabla )
&\,\longmapsto\, (S,\, E,\, \phi ,\, \nabla).
\end{aligned}
\end{equation*}
It is straightforward to check that this map gives an isomorphism
$\Sigma_{m_0}^{d} \,\cong\, \big[ \, (Z^{d,H}_{m_0})' \, \big/ \,
\mathrm{PGL}(\theta_d(m_0),\mathbb{C}) \, \big] $.
Then $\Sigma_{m_0}^{d}$ is an algebraic stack.
By Proposition \ref{2020.1.22.11.10}, 
it follows that $\Sigma_{m_0}^{d}$ is in fact a Deligne--Mumford stack
(see \cite[Corollary 8.4.2]{Olsson}).
\end{proof}

\begin{Rem}
If $H_x \,=\, \{ e\}$ for all $x\,\in\, D$,
then $\mathcal{M}^{H}_{\rm{FC}}(d)$ is an algebraic space
by Proposition \ref{2019.12.24.23.58}.
\end{Rem}

\begin{Rem}
In the proof of Proposition \ref{2020.2.24.14.52},
we introduced the coherent sheaf $\mathscr{H}$
which is characterized by the property \eqref{eq:2021.11.25.16.17}.
Since $\mathscr{H}$ is not necessarily locally free,
we cannot see the irreducibility of the
moduli space $\mathcal{M}^H_{\mathrm{FC}}(d)$
immediately from its construction.
\end{Rem}

Define an open substack of $\mathcal{M}^{H}_{\rm{FC}}(d)$ as follows:
\begin{equation}\label{2021.11.26.11.19}
\mathcal{M}^{H}_{\rm{FC}}(d)^{\circ} \,:=\,
\{ (S,\,E,\,\phi,\,\nabla) \,\in \,\mathcal{M}^{H}_{\rm{FC}}(d) \,\mid\, 
\text{ $(E_s,\,\phi_s)$\, is simple for each $s \,\in\, S$}\}.
\end{equation}
Here we say that $(E_s,\,\phi_s)$ is simple if
$$
\begin{cases} 
H^0(X,\, \mathcal{E}nd_{\phi_s}( E_s) )\, =\,0
& \text{when $\mathbb{C}^* \cdot \mathrm{Id} \,\not\subset\, H_x$ for some $x \,\in\, D$} \\
H^0(X,\, \mathcal{E}nd_{\phi_s}( E_s) ) \,=\,\mathbb{C}
& \text{when $\mathbb{C}^* \cdot \mathrm{Id} \,\subset\, H_x$ for all $x \,\in\, D$.}
\end{cases}
$$
We adopt the above definition of simple framed bundle, in order that
the loci $\mathcal{M}^{H}_{\rm{FC}}(d)^{\circ}$ becomes open
in $\mathcal{M}^{H}_{\rm{FC}}(d)$.

\begin{Prop}\label{2023_4_2_13_21}
The open substack $\mathcal{M}^{H}_{\rm{FC}}(d)^{\circ}$ in \eqref{2021.11.26.11.19}
is irreducible.
\end{Prop}

\begin{proof}
Fix a very ample line bundle $\mathcal{O}_{X}(1)$ on the curve $X$.
Let $\theta_d(m)$ be a polynomial in $m$
defined as in Proposition \ref{2020.2.24.14.52}.
Let $(\Sigma_{m_0}^{d})^{\circ}$
denote the substack of $\mathcal{M}^{H}_{\rm{FC}}(d)^{\circ}$ whose objects are
framed $\mathrm{GL}(r,\mathbb{C})$-connections $(E,\,\phi,\,\nabla)$ on $X$
such that 
\begin{itemize} 
\item $(E,\,\phi)$ is simple,

\item $H^1(X,\, E (m_0-1))\,=\,0$, and

\item $\chi(E(m))\,= \,\theta_d(m)$ for all $m \,\in\, \mathbb{Z}$.
\end{itemize}
To prove the proposition it suffices to show that $(\Sigma_{m_0}^{d})^{\circ}$ is irreducible.

Let $V$ be a $\theta_d(m_0)$-dimensional vector space so that the underlying vector bundle $E$ 
of any object of $(\Sigma_{m_0}^{d})^{\circ}$ is described as the following quotient:
\begin{equation*}
V \otimes\mathcal{O}_X (-m_0)\, \longrightarrow E.
\end{equation*}
Take a subspace $V_r \,\subset\, V$ such that $\dim (V_r) \,=\,r$. Taking the dual of the above
quotient, and tensoring with $\mathcal{O}_X (-m_0)$, we have the following short exact sequence 
\begin{equation*}
0 \, \longrightarrow \, E^{\vee}(-m_0) 
\, \longrightarrow \, V_r^{\vee} \otimes \mathcal{O}_X
\, \longrightarrow \, F
\, \longrightarrow \, 0,
\end{equation*}
where $F$ is the quotient for the injective map 
$E^{\vee}(-m_0) \,\longrightarrow\, V_r^{\vee} \otimes \mathcal{O}_X$.
So for each object of $(\Sigma_{m_0}^{d})^{\circ}$, there is a point 
of $\mathrm{Quot}^N_{(V_r^{\vee}\otimes \mathcal{O}_X)/X}$ which
determines the underlying vector bundle of the object. Here $N$ is the length of $F$.
Note that $N$ remains constant for the underlying vector bundles.
We will show that $\mathrm{Quot}^N_{(V_r^{\vee}\otimes \mathcal{O}_X)/X}$ is 
irreducible.

The Quot-scheme $\mathrm{Quot}^N_{(V_r^{\vee}\otimes \mathcal{O}_X)/X}$ 
is smooth, because the obstructions to deformations
of 
$$[\, q \,\colon\, V_r^{\vee} \otimes \mathcal{O}_X \,\longrightarrow\, F \,]
\, \in \,\mathrm{Quot}^N_{(V_r^{\vee}\otimes \mathcal{O}_X)/X}$$
lie in $$\mathrm{Ext}^1 (\mathrm{Ker}\, q , \,F) \,\cong\, H^1 ((\mathrm{Ker}\, q)^{\vee} \otimes F)\,=\,0\, .$$
Define the map
\begin{equation*}
\begin{aligned}
f_N \,\,\colon\,\, \mathrm{Quot}^N_{(V_r^{\vee}\otimes \mathcal{O}_X)/X} 
\, &\longrightarrow \, 
\mathrm{Hilb}^N_{X}\\
[\, q \colon V_r^{\vee} \otimes \mathcal{O}_X \longrightarrow F \,] 
\, &\longmapsto\, 
\text{Divisor} ( \mathrm{det} (\mathrm{Ker}\, q \,\longrightarrow\, V_r^{\vee} \otimes \mathcal{O}_X)).
\end{aligned}
\end{equation*}
Let $H'$ be the Zariski open subset of $\mathrm{Hilb}^N_{X}$
which consists of distinct points on $X$, or in other words, 
$H'$ parametrizes the reduced subschemes. 
This Zariski open subset $H'$ is in fact irreducible.
The map $$f_N^{-1}(H') \ \longrightarrow\ H'$$ is a 
$(\mathbb{P}^{r-1}\times \cdots \times \mathbb{P}^{r-1})$-bundle; here 
$\mathbb{P}^{r-1}\times \cdots \times \mathbb{P}^{r-1}$
is the product of $N$-copies of $\mathbb{P}^{r-1}$.
Hence $f_N^{-1}(H')$ is irreducible.
Take a point $x \,=\, N_1 z_1+ \cdots +N_l z_l $ on $\mathrm{Hilb}^N_{X}$,
where $\sum_{i=1}^l N_i\, =\,N$ and $z_1,\,\cdots,\,z_l$ are distinct points on $X$.
A point on $f^{-1}_N(x)$ can be described as a collection
$\left( \,
q_i \,\colon\, V_r^{\vee} \otimes\mathcal{O}_{z_i,X}
\,\longrightarrow\, F_i\, \right)_{i=1}^l$ for which $\mathrm{length}(F_i)\,=\,N_i$. Consider the map 
$(\mathrm{Ker}\, q_i )_{z_i} \,\longrightarrow\, V_r^{\vee} \otimes \mathcal{O}_{z_i,X}$
corresponding to a point on $f^{-1}_N(x)$.
Note that $(\mathrm{Ker}\, q_i )_{z_i}\,\cong\, \mathcal{O}_{z_i,X}^{\oplus r} $.
We have a matrix representation of 
$(\mathrm{Ker}\, q_i )_{z_i} \,\longrightarrow \, V_r^{\vee} \otimes \mathcal{O}_{z_i,X}$
as follows:
\begin{equation*}
\begin{pmatrix}
1& & & & &\\
&\ddots & & & &\\
&& 1 & & & \\
& & & z^{l_{i,1}} & &\\
& & & & \ddots & \\
& & & & & z^{l_{i,s_i}}
\end{pmatrix},
\end{equation*}
where the maximal ideal $\mathfrak{m}_{z_i}$ is $\{z=0 \}$ and 
$1,\,z^{l_{i,1}},\,\cdots,\,z^{l_{i,s_i}} $ are invariant factors of
$$(\mathrm{Ker}\, q_i )_{z_i} \,\longrightarrow\,V_r^{\vee} \otimes \mathcal{O}_{z_i,X}\, .$$
For any tuple of complex numbers $a_1^{(i)},\,\cdots,\,a_{s_i}^{(i)}$,
there is a deformation of $(\mathrm{Ker}\, q_i )_{z_i} \,\longrightarrow\, 
 V_r^{\vee} \otimes \mathcal{O}_{z_i,X}$
 \begin{equation*}
\begin{pmatrix}
1& & & & &\\
&\ddots & & & &\\
&& 1 & & & \\
& & & z^{l_{i,1}}+t a_{1}^{(i)} & & \\
& & & & \ddots & \\
& & & & & z^{l_{i,s_i}}+ t a_{s_i}^{(i)}
\end{pmatrix},
\end{equation*}
over $X \times \mathrm{Spec}\, \mathbb{C}[t]$.
When the complex numbers $a_1^{(i)},\,\cdots,\,a_{s_i}^{(i)}$ are generic,
we have a deformation moving from a point on $f^{-1}_N(x)$,
where $x \,=\, N_1 z_1+ \cdots +N_l z_l $, to a point on $f_N^{-1}(H')$. 
Therefore, $\mathrm{Quot}^N_{(V_r^{\vee}\otimes \mathcal{O}_X)/X}$ 
is irreducible. 

Consider the open subset 
$$
Q' \ := \
\left\{ [\, q \,]
\, \in \,\mathrm{Quot}^N_{(V_r^{\vee}\otimes \mathcal{O}_X)/X}
\ \middle| \ 
\text{$E_q$ satisfies $H^1(E_q(m_0-1))\,=\,0$}
\right\} \, \subset\, \mathrm{Quot}^N_{(V_r^{\vee}\otimes \mathcal{O}_X)/X}.
$$
Here denote $E_q\,:=\,(\mathrm{ker} q)^{\vee} (-m_0)$ for a quotient $q$. By definition, 
$E_q$ is locally free and satisfies the condition $\chi (E_q (m)) = \theta_d(m)$ for all $m \,
\in\, \mathbb{Z}$. Let $\widetilde{Q}'$ be the scheme over $Q'$ which parametrizes 
quotients $q$ in $Q'$ and framings of $E_q$,
which is constructed as in the proof of Proposition \ref{2020.2.24.14.52}.
Now define $(\widetilde{Q}')^{\circ}$ as follows:
\begin{equation*}
(\widetilde{Q}')^{\circ}\,:=\, \left\{
s \,\in\, \widetilde{Q}'
\ \middle| \
\text{$(\widetilde{\mathcal{E}}_s ,\, \widetilde{\phi}_s )$ is simple}
\right\}.
\end{equation*}
Here $(\widetilde{\mathcal{E}} ,\, \widetilde{\phi} )$ is the family 
of vector bundles $E_q$ and framings of $E_q$ in $\widetilde{Q}'$
induced by the universal family of $\widetilde{Q}'$.
Since $\mathrm{Quot}^N_{(V_r^{\vee}\otimes \mathcal{O}_X)/X}$ 
is irreducible,
$(\widetilde{Q}')^{\circ}$ 
is also irreducible.
Let $(Z')^{\circ}$ be
the scheme over $(\widetilde{Q}')^{\circ}$ which
parametrizes quotients $q$ with framings of $E_q$ in $(\widetilde{Q}')^{\circ}$
and connections on $E_q$ that are compatible with the framings.
The scheme $(Z')^{\circ}$ is also
constructed as in the proof of Proposition \ref{2020.2.24.14.52}.
It is straightforward to check that $(Z')^{\circ}$ is smooth and 
each fiber of $(Z')^{\circ} \,\longrightarrow\,
(\widetilde{Q}')^{\circ}$ is an affine space,
which is isomorphic to $H^0(X,\, \mathrm{End}_{\phi}^n (E) \otimes K_X(D))$.
So $(Z')^{\circ}$ is irreducible.
Since a natural map from $(Z')^{\circ}$ to $(\Sigma_{m_0}^{d})^{\circ}$
is induced and this map is surjective, 
we conclude that $(\Sigma_{m_0}^{d})^{\circ}$ is irreducible. 
This completes the proof of the proposition.
\end{proof}

\section{Symplectic structures of the moduli spaces}\label{section: symplectic structure}

Throughout this section it is assumed that $G\,=\,\mathrm{GL}(r,\mathbb{C})$.

\subsection{Cotangent bundle of the moduli space of 
simple framed bundles}\label{2019.11.29.10.35}

In this subsection we assume that $H_x\,=\,\{ e \}\, \subset\, \mathrm{GL}(r,\mathbb{C})$
for all $x \,\in\, D$. Let $\mathcal{N}^e(d)$ be the following moduli space:
\begin{equation}\label{ned}
{\mathcal{N}^e(d)\,=\, \left\{\left(E,\,\phi\,=\, \{ \phi_x \}_{x \in D}\right)\ \middle| \
\begin{array}{l}
\text{$E$ is a vector bundle of degree $d$ and}\\
\text{$(E,\,\phi)$ is a simple framed principal}\\
\text{$G$-bundle, where $H_x \,=\, \{ e\}$ for all $x\,\in\, D$.} 
\end{array}
\right\}{\Big/}\Large{\sim}_{e}}.
\end{equation}
Here $(E,\, \phi ) \,\sim_e\, (E',\, \phi' )$ if there 
exists an isomorphism $\sigma \,\colon\, E\,\longrightarrow\, E'$ of vector bundles
such that the composition of homomorphisms 
$(\phi'_{x})^{-1} \circ \sigma|_{x}\circ \phi_{x}$
is the identity map of ${\mathbb C}^r$
for each $x\,\in\, D$. Since the tangent space of $\mathcal{N}^e(d)$ at
$(E,\,\phi)$ is $H^1(X,\,\mathcal{E}nd(E)(-D))$ \cite[Lemma 2.5]{BLP0}, using the Serre duality
it follows that the cotangent space of $\mathcal{N}^e(d)$ at $(E,\,\phi)$ is
$H^0 (X,\, \mathcal{E}nd(E) \otimes K_X(D))$. Let
$T^*\mathcal{N}^e(d)$ be the cotangent bundle of $\mathcal{N}^e(d)$.
For $\theta\, \in\, H^0 (X,\, \mathcal{E}nd(E) \otimes K_X(D))$,
define the following $2$-term complex:
\begin{equation*}
\mathcal{C}_{\bullet}^{\mathrm{Higgs}}
\,\colon\, \mathcal{C}_{0}^{\mathrm{Higgs}}\,:=\, \mathcal{E}nd(E)(-D) \,\xrightarrow{\,[\theta,\,\cdot]\,} \,
\mathcal{C}_{1}^{\mathrm{Higgs}} \,:=\, \mathcal{E}nd(E) \otimes K_X(D)\, .
\end{equation*}
The tangent space $T_{(E,\phi,\theta)}T^*\mathcal{N}^e(d)$ at $(E,\,\phi,\,\theta)$ is
$\mathbb{H}^1(\mathcal{C}_{\bullet}^{\mathrm{Higgs}})$ \cite[Lemma 2.7]{BLP0}.
Given an affine open covering $\{ U_{\alpha} \}$ of $X$, the hypercohomology
$\mathbb{H}^1(\mathcal{C}_{\bullet}^{\mathrm{Higgs}})$ 
admits a description in terms of \v{C}ech cohomology. In this description,
the $1$-cocycles are pairs $(\{ u_{\alpha\beta} \} ,\, \{ v_{\alpha} \})$, 
where $$u_{\alpha\beta} \,\in\, \mathcal{E}nd(E)(-D)(U_{\alpha}\cap U_{\beta})\ \
\text{ and }\ \ v_{\alpha } \,\in \,\mathcal{E}nd(E) \otimes K_X(D)(U_{\alpha})$$
such that $u_{\beta\gamma}-u_{\alpha\gamma}+u_{\alpha\beta}\,=\,0$ 
and $v_{\beta}-v_{\alpha}\,=\, [\theta,\, u_{\alpha\beta}]$.
The $1$-coboundaries are of the form $(\{w_{\alpha}-w_{\beta} \} ,\, \{[w_\alpha,\, \theta]\})$,
where $w_{\alpha} \,\in\, \mathcal{E}nd(E)(-D)(U_{\alpha})$.

We define a canonical $1$-form $\phi_{\mathcal{N}^e(d) }$ on the cotangent 
bundle $T^*\mathcal{N}^e(d)$ by
\begin{equation} \label{z1}
\begin{aligned}
\begin{split} 
\phi_{\mathcal{N}^e(d) } \,\colon\,
\mathbb{H}^1 (\mathcal{C}_{\bullet}^{\mathrm{Higgs}}) &\longrightarrow H^1(K_X) \\
[(\{ u_{\alpha\beta} \},\, \{ v_{\alpha} \} )] &\longmapsto 
[ \{ \mathrm{Tr}(\theta|_{U_{\alpha}} u_{\alpha\beta} ) \}].
\end{split}
\end{aligned}
\end{equation}

\begin{Lem}
Let $\Phi_{T^*\mathcal{N}^e(d)}$ be the Liouville $2$-form 
on the cotangent bundle $T^*\mathcal{N}^e(d)$, that is, 
$\Phi_{T^*\mathcal{N}^e(d)}$ is the exterior derivative of
the canonical $1$-form $\phi_{\mathcal{N}^e(d) }$ in \eqref{z1}.
The Liouville $2$-form $\Phi_{T^*\mathcal{N}^e(d)}$ coincides with the
bilinear form
\begin{equation*}
\begin{aligned}
\mathbb{H}^1 (\mathcal{C}_{\bullet}^{\mathrm{Higgs}}) \otimes 
\mathbb{H}^1 (\mathcal{C}_{\bullet}^{\mathrm{Higgs}}) &\longrightarrow H^1(K_X) \\
[(\{ u_{\alpha\beta} \},\, \{ v_{\alpha} \} )]
\otimes [(\{ u'_{\alpha\beta} \},\, \{ v'_{\alpha} \} )] &\longmapsto 
[\{\mathrm{Tr} ( v_{\alpha} u'_{\alpha\beta} )- \mathrm{Tr} (u_{\alpha\beta} v'_{\beta} ) \}]
\end{aligned}
\end{equation*}
on \v{C}ech cohomology.
\end{Lem}

\begin{proof}
Let $v$ and $v'$ be tangent vectors of $T^*\mathcal{N}^e(d)$ at $(E,\,\phi,\,\theta)
\,\in\, T^*\mathcal{N}^e(d)$.
Let $$D_v \,\,\colon\,\, \mathcal{O}_{T^*\mathcal{N}^e(d)} \,\longrightarrow\, 
\mathcal{O}_{T^*\mathcal{N}^e(d)}$$ be the derivative 
corresponding to $v$.
Take an affine open subset $U\,\subset\, T^*\mathcal{N}^e(d)$ such that $(E,\,\phi,\,\theta)
\,\in\, U$, and also take an affine open covering $\{ U_{\alpha} \} $ of $X\times U$
such that there is a trivialization
$$g_{\alpha} \,\colon \,E|_{U_{\alpha}} \,\stackrel{\sim}{\longrightarrow}
\, \mathcal{O}_{U_{\alpha}}^{\oplus r}$$
for each $U_{\alpha}$. 
Set $g_{\alpha\beta}\,:=\,g_{\alpha} \circ g^{-1}_{\beta}$
and $\theta_{\alpha}\,:=\, g_{\alpha} \circ \theta|_{U_{\alpha}} \circ g_{\alpha}^{-1}$.
We may describe the tangent vector $v$ as 
$$v\,=\,[(\{ u_{\alpha\beta}\} ,\, \{ v_{\alpha} \} )]\, ,$$
where $u_{\alpha\beta}\,:=\, g_{\alpha}^{-1} \circ ( D_{v}(g_{\alpha\beta} )g_{\alpha\beta}^{-1} ) \circ g_{\alpha} $
and $v_{\alpha}\,:=\, g_{\alpha}^{-1} \circ ( D_v \theta_{\alpha}) \circ g_{\alpha}$.
The exterior derivative of $\phi_{\mathcal{N}^e(d) }$ is computed as follows:
\begin{equation*}
\begin{aligned}
&D_{v}\phi_{\mathcal{N}^e(d)} ( v')-
D_{v'}\phi_{\mathcal{N}^e(d)} (v) + \phi_{\mathcal{N}^e(d) }([v,\, v']) \\
&=\, D_{v} \left( \mathrm{Tr}( \theta_{\alpha} D_{v'}(g_{\alpha\beta} )g_{\alpha\beta}^{-1} ) \right)
- D_{v'} \left( \mathrm{Tr} ( \theta_{\alpha} D_{v}(g_{\alpha\beta} )g_{\alpha\beta}^{-1} ) \right)
+ \left( \mathrm{Tr} ( \theta_{\alpha}
(D_{v'} \circ D_{v} - D_{v} \circ D_{v'} )(g_{\alpha\beta} ) )g_{\alpha\beta}^{-1} \right) \\
&=\,\mathrm{Tr} \left( D_{v}(g_{\alpha\beta}^{-1} \theta_{\alpha} ) D_{v'}(g_{\alpha\beta} ) \right)
-\mathrm{Tr} \left( D_{v'}( g_{\alpha\beta}^{-1}\theta_{\alpha})D_{v}(g_{\alpha\beta} ) \right) \\
&= \,\left( \mathrm{Tr}( D_{v}(\theta_{\alpha} ) D_{v'}(g_{\alpha\beta} )g_{\alpha\beta}^{-1} ) 
- \mathrm{Tr} ( D_{v'}( \theta_{\alpha})D_{v}(g_{\alpha\beta} )g_{\alpha\beta}^{-1} ) \right) \\
&\qquad-\left( \mathrm{Tr}(D_{v}(g_{\alpha\beta})g_{\alpha\beta}^{-1} 
\theta_{\alpha} D_{v'}(g_{\alpha\beta} )g_{\alpha\beta}^{-1} ) \right)
+\left( \mathrm{Tr} ( D_{v'}( g_{\alpha\beta})g_{\alpha\beta}^{-1}
\theta_{\alpha}D_{v}(g_{\alpha\beta} ) g_{\alpha\beta}^{-1} ) \right) \\
&= \,\left( \mathrm{Tr} ( v_{\alpha} u'_{\alpha\beta} )- \mathrm{Tr} (v'_{\alpha}u_{\alpha\beta} ) \right) 
+ \mathrm{Tr} ( [u'_{\alpha\beta} ,\, \theta_{\alpha}]u_{\alpha\beta} ) \\
&=\, \mathrm{Tr} ( v_{\alpha} u'_{\alpha\beta} )- \mathrm{Tr} (u_{\alpha\beta} v'_{\beta}).
\end{aligned}
\end{equation*}
This completes the proof of the lemma.
\end{proof}

\subsection{2-form on $\mathcal{M}^{e}_{\mathrm{FC}} (d)$}\label{2023.1.10.12.15}

As before, assume that $H_x\,=\,\{ e \}\, \subset\, \mathrm{GL}(r,\mathbb{C})$ for
every $x \,\in\, D$. Let $\mathbb{K}$ denote the following complex of coherent sheaves on $X$:
\begin{equation}\label{bK}
\mathbb{K}\ : \ \mathcal{O}_{X }\, \xrightarrow{\ d \ }\, K_X \,\longrightarrow\, 0,
\end{equation}
where $\mathcal{O}_{X }$ 
and $K_X$ are at the $0$-th position and $1$-position respectively, 
and $d$ is the de Rham differential.

Consider the complex $\mathcal{C}_{\bullet}$ in \eqref{2019.11.19.14.51}. Define a pairing
\begin{equation}\label{ted}
\begin{aligned}
\Theta^e \,\colon\,
\mathbb{H}^1 (\mathcal{C}_{\bullet}) \otimes \mathbb{H}^1 (\mathcal{C}_{\bullet}) 
&\,\longmapsto \,\mathbb{H}^2(\mathbb{K}) \,\cong \,\mathbb{C} \\
[(\{u_{\alpha\beta} \} ,\, \{ v_{\alpha} \})]\otimes[(\{u'_{\alpha\beta} \},\,
\{ v'_{\alpha} \})] &\,\longmapsto \,
[(\{ \mathrm{Tr} (u_{\alpha\beta}u'_{\beta\gamma}) \},\,
- \{\mathrm{Tr}(u_{\alpha\beta}v'_{\beta})-\mathrm{Tr}(v_{\alpha}u'_{\alpha\beta}) \})]
\end{aligned}
\end{equation}
in terms of the \v{C}ech cohomology with respect to an affine open covering $\{ U_{\alpha}\}$ of $X$.

\begin{Lem}\label{le2}
The pairing $\Theta^e$ in \eqref{ted} satisfies the identity
\begin{equation*}
\Theta^e(v,\,v)\,=\, 0\, .
\end{equation*}
Thus $\Theta^e$ is skew-symmetric and hence it produces a $2$-form on
$\mathcal{M}^e_{\mathrm{FC}} (d)$ (see Lemma \ref{lemid}).
\end{Lem}

\begin{proof}
Let $v\,=\, [(\{u_{\alpha\beta} \} ,\, \{ v_{\alpha} \})]$ be an element of
$\mathbb{H}^1 (\mathcal{C}_{\bullet})$. We compute $\Theta^e(v,\,v)$ as follows:
\begin{equation*}
\begin{aligned}
\Theta^e(v,\,v) &=\,
[(\{ \mathrm{Tr} (u_{\alpha\beta}u_{\beta\gamma}) \},\,
- \{\mathrm{Tr}(u_{\alpha\beta}v_{\beta})-\mathrm{Tr}(v_{\alpha}u_{\alpha\beta}) \})] \\
&=\, [(\{ \mathrm{Tr} (u_{\alpha\beta}u_{\beta\gamma}) \},\,
- \{\mathrm{Tr}(u_{\alpha\beta}(v_{\beta}-v_{\alpha})) \})] \\
&=\, [(\{ \mathrm{Tr} (u_{\alpha\beta}u_{\beta\gamma}) \},\,
- \{\mathrm{Tr}(u_{\alpha\beta}(\nabla \circ u_{\alpha\beta} - u_{\alpha\beta} \circ \nabla) ) \})] \\
&=\, \left[\left(\left\{ \mathrm{Tr} (u_{\alpha\beta}u_{\beta\gamma}) \right\},\,
- \left\{ d \left( \frac{1}{2}\mathrm{Tr}(u^2_{\alpha\beta} ) \right)\right\}\right)\right].
\end{aligned}
\end{equation*}
On the other hand,
\begin{equation*}
\begin{aligned}
\frac{1}{2}\mathrm{Tr}(u^2_{\alpha\beta} )
-\frac{1}{2}\mathrm{Tr}(u^2_{\alpha\gamma})
+\frac{1}{2}\mathrm{Tr}(u^2_{\beta\gamma}) 
&=\,
\frac{1}{2}\mathrm{Tr}((u_{\alpha\beta} -u_{\alpha\gamma} )(u_{\alpha\beta} +u_{\alpha\gamma} ))
+\frac{1}{2}\mathrm{Tr}(u^2_{\beta\gamma}) \\
&=\,
\frac{1}{2}\mathrm{Tr}((u_{\beta\gamma} )(u_{\beta\gamma} -u_{\alpha\beta} -u_{\alpha\gamma} )) \\
&=\,- \mathrm{Tr}(u_{\alpha\beta}u_{\beta\gamma}).
\end{aligned}
\end{equation*}
Combining these it follows that $\Theta^e(v,\,v)\,=\,0$ in $\mathbb{H}^2(\mathbb{K})$.
\end{proof}

\begin{Rem}
We have constructed a 2-form $\Theta^e$ on $\mathcal{M}^e_{\mathrm{FC}} (d)$
by \eqref{ted}.
On the other hand, there exists another definition of this 2-form from a differential geometric perspective;
this will be explained below.
First recall a description of $\mathbb{H}^1 (\mathcal{C}_{\bullet})$
as Dolbeault cohomology. (See the proof of Theorem 3.2 of \cite{Biswas}.)
Let $\overline{\partial}'$ and $\overline{\partial}'_1$ be
the Dolbeault operators for the holomorphic vector bundles $\mathcal{E}nd(E)(-D)$ and 
$\mathcal{E}nd(E)(D)$, respectively. Consider the Dolbeault resolution of the complex $\mathcal{C}_{\bullet}$:
\begin{equation*}
\xymatrix{
 0 \ar[d] & 0\ar[d] \\
\mathcal{E}nd(E)(-D) \ar[d] \ar[r]^-{\nabla}
 &\mathcal{E}nd(E)\otimes K_X (D) \ar[d]\\
 \Omega^{0,0}_X (\mathcal{E}nd(E)(-D))\ar[d]^-{\overline{\partial}'} \ar[r]^-{\nabla}
 &\Omega^{1,0}_X(\mathcal{E}nd(E)(D)) \ar[d]^-{\overline{\partial}_1'}\\
 \Omega^{0,1}_X (\mathcal{E}nd(E)(-D) )\ar[d] \ar[r]^-{\nabla'}
 &\Omega^{1,1}_X(\mathcal{E}nd(E)(D))\ar[d] \\
 0 & 0
}
\end{equation*}
where $\nabla'$ is constructed using $\nabla$ 
and the usual differential operator $\partial$ on $(0,\,1)$-forms on $X$.
Note that
\begin{equation*}
\overline{\partial}'_1\circ \nabla + \nabla' \circ \overline{\partial}'\, =\,0\, .
\end{equation*}
This produces the following complex of vector spaces
\begin{equation*}
\begin{aligned}
0 \longrightarrow C^{\infty}(X,\,\mathcal{E}nd(E)(-D)) 
&\,\xrightarrow{\ \overline{\partial}' \oplus \nabla \ }\, C^{\infty}(X,\,\Omega^{0,1}_X
(\mathcal{E}nd(E)(-D)))\oplus C^{\infty}(X,\, \mathcal{E}nd(E)\otimes K_X(D))\\
&\xrightarrow{\ \nabla'+\overline{\partial}'\ }\,
C^{\infty}(X,\, \Omega^{1,1}_X (\mathcal{E}nd(E)(D))) \,\longrightarrow\, 0.
\end{aligned}
\end{equation*}
Since the Dolbeault complex is a fine resolution of $\mathcal{C}_{\bullet}$, it follows
immediately that
\begin{equation*}
\mathbb{H}^1 (\mathcal{C}_{\bullet})\,\,=\,\, 
\frac{\mathrm{Ker(\nabla'+\overline{\partial}')}}{\mathrm{Im}(\overline{\partial}' \oplus \nabla)}.
\end{equation*}
Then the $2$-form $\Theta^e$ in \eqref{ted} can be described using the Dolbeault cohomology in the following way:
\begin{equation*}
(\omega_1,\,\omega_2) \otimes (\omega_1',\,\omega_2')\,\longmapsto\, 
\int_X \mathrm{Tr} (\omega_1 \wedge \omega_2')+ 
\int_X\mathrm{Tr} (\omega_2\wedge \omega_1')\, .
\end{equation*}
\end{Rem}

\subsection{Symplectic structure on $\mathcal{M}^{e}_{\mathrm{FC}} 
(d)^{\circ}$}\label{2019.11.29.10.38}

Now take the restriction of $\Theta^e$ to $\mathcal{M}^{e}_{\mathrm{FC}}(d)^{\circ}$;
here $\mathcal{M}^{e}_{\mathrm{FC}}(d)^{\circ}$ is the open substack 
of $\mathcal{M}^{e}_{\mathrm{FC}}(d)$ defined in \eqref{2021.11.26.11.19}, or in other words,
the underlying framed bundle $(E, \, \phi)$ of any
$(E,\, \phi,\, \nabla)\, \in\, \mathcal{M}^{e}_{\mathrm{FC}}(d)^{\circ}$ 
satisfies the condition that it is simple. Denote this restriction of
$\Theta^e$ to $\mathcal{M}^{e}_{\mathrm{FC}}(d)^{\circ}$
by $\Theta^e|_{\mathcal{M}^{e}_{\mathrm{FC}}(d)^{\circ}}$.
It will be shown that this restriction of $\Theta^e$ is a symplectic form.

Let
\begin{equation}\label{fm}
p_1 \,\colon\, \mathcal{M}^{e}_{\rm{FC}}(d)^{\circ} \,\longrightarrow\, \mathcal{N}^e(d)
\end{equation}
be the forgetful map that simply forgets the connection.
Take an analytic open subset $U\, \subset\, \mathcal{N}^e(d)$, which is assumed to be
small enough. Then there exist sections, over $U$, of the map $p_1$ in \eqref{fm}. Let
$$s\,\colon\, U \,\longrightarrow\, p_1^{-1}(U)$$ be a holomorphic section.
Using $s$, an isomorphism
\begin{align}
P_1\,\colon\, T^*U &\,\xrightarrow{\ \cong \ }\, p_1^{-1}(U) \label{eqt}\\
(y,\,v) &\,\longmapsto \,s(y) + v\nonumber
\end{align}
is obtained. The restriction, to $p_1^{-1}(U)$, of the form $\Theta^e|_{\mathcal{M}^{e}_{\mathrm{FC}}(d)^{\circ}}$
is denoted by $\Theta^e|_{p_1^{-1}(U)}$.

\begin{Lem}\label{2019.11.19.15.22}
Let $\Phi_U$ be the Liouville $2$-form on the cotangent bundle $T^*U$.
Then,
\begin{equation*}
\Theta^e|_{p_1^{-1}(U)} - (P_1^{-1})^* \Phi_U \,=\, p_1^* (s^* \Theta^e|_{p_1^{-1}(U)})\, ,
\end{equation*}
where $\Theta^e|_{p_1^{-1}(U)}$ is
the restriction of the form in \eqref{ted}, and $p_1$ is the projection in \eqref{fm}, while $P_1$ is
the isomorphism in \eqref{eqt}.
\end{Lem}

\begin{proof}
Take a point $z\,=\,(E,\,\phi,\,\nabla)$ of $p_1^{-1}(U)$.
Let $\nabla(E,\,\phi)$ be the connection associated to the point $s\circ p_1(z)$.
The image of $z$ under the map $P^{-1}_1$ in \eqref{eqt} is as follows:
\begin{equation*}
P_1^{-1}(z)\,=\,P^{-1}_1(E,\,\phi,\,\nabla) \,= \,(E,\,\phi,\, \nabla - \nabla(E,\phi))\, .
\end{equation*}
Let $[(\{u_{\alpha\beta} \} ,\,\{ v_{\beta} \})]$ be an element of
$\mathbb{H}^1((\mathcal{C}_{\bullet})_z)$, where $(\mathcal{C}_{\bullet})_z$ is the complex in \eqref{2019.11.19.14.51} associated 
to $z\,=\,(E,\,\phi,\, \nabla)$. Recall from
Lemma \ref{lemid} that $\mathbb{H}^1((\mathcal{C}_{\bullet})_z)$ is the tangent space
of $p_1^{-1}(U)$ at $z$.
Note that $u_{\alpha\beta}$ and $v_{\alpha}$ satisfy the equality
$$v_{\beta}-v_{\alpha}\,=\, \nabla\circ u_{\alpha\beta}-u_{\alpha\beta} \circ \nabla.$$
Let $[(\{u_{\alpha\beta}\},\, \{ v^s_{\alpha} \})]$ be the element of 
$\mathbb{H}^1((\mathcal{C}_{\bullet})_{s\circ p_1(z)})$ such that 
$$(s \circ p_1)_{*} ([(\{u_{\alpha\beta} \} ,\,\{ v_{\beta} \})])
\,=\, [(\{u_{\alpha\beta} \} ,\, \{ v^s_{\beta} \})]\, . $$
Note that $u_{\alpha\beta}$ and $v^s_{\alpha}$ satisfy the equality
$$v^s_{\beta}-v^s_{\alpha}\,=\, \nabla(E,\phi)\circ u_{\alpha\beta}-u_{\alpha\beta} \circ
\nabla(E,\phi)\, .$$
Since
\begin{equation*}
(v_{\beta} - v^s_{\beta})- (v_{\alpha} - v^s_{\alpha})
\,=\, [ \nabla- \nabla(E,\phi),\, u_{\alpha\beta}]\, ,
\end{equation*}
it follows that $[(\{ u_{\alpha\beta} \} ,\,\{v_{\alpha}-v_{\alpha}^s \})]$ is an element of $\mathbb{H}^1
(\mathcal{C}^{\mathrm{Higgs}}_{\bullet})$; recall that $\mathbb{H}^1
(\mathcal{C}^{\mathrm{Higgs}}_{\bullet})$ is the tangent space of $T^*U$ at $P_1^{-1}(z)$ (see \eqref{eqt}).
There is a map
\begin{equation*}
\begin{aligned}
 \mathbb{H}^1 ((\mathcal{C}_{\bullet})_z) &\longrightarrow 
\mathbb{H}^1 (\mathcal{C}^{\mathrm{Higgs}}_{\bullet}) \\ 
[(\{u_{\alpha\beta} \},\,\{ v_{\alpha} \})] &\longmapsto 
[(\{ u_{\alpha\beta} \} ,\,\{v_{\alpha}-v_{\alpha}^s \})] .
\end{aligned}
\end{equation*}
This map coincides with $(P_1^{-1})_* \,\colon\,\mathbb{H}^1 ((\mathcal{C}_{\bullet})_z) \,
\longrightarrow \,\mathbb{H}^1 (\mathcal{C}^{\mathrm{Higgs}}_{\bullet})$.

Now we compute the map
\begin{equation*}
(\Theta^e - (P_1^{-1})^* \Phi_U) \,\colon \,
 \mathbb{H}^1 ((\mathcal{C}_{\bullet})_z) \otimes 
 \mathbb{H}^1 ((\mathcal{C}_{\bullet})_z) \,\longrightarrow \,
\mathbb{H}^2(\mathbb{K})
\end{equation*}
 as follows:
\begin{equation*}
\begin{aligned}
(\Theta^e - (P_1^{-1})^* \Phi_U)(v,v') &\,=\,
[(\{ \mathrm{Tr} (u_{\alpha\beta}u'_{\beta\gamma}) \},\,
- \{ \mathrm{Tr}(u_{\alpha\beta}v'_{\beta})-\mathrm{Tr}(v_{\alpha}u'_{\alpha\beta}) \})] \\
&\qquad- [ \{ 0 \},\,- \{ \mathrm{Tr}(u_{\alpha\beta}(v'_{\beta}-(v_{\beta}^s)'))
-\mathrm{Tr}((v_{\alpha}-v_{\alpha}^s) u'_{\alpha\beta}) \} ]\\
&=\,[(\{ \mathrm{Tr} (u_{\alpha\beta}u'_{\beta\gamma}) \},\,
- \{\mathrm{Tr}(u_{\alpha\beta}(v^s_{\beta})')-\mathrm{Tr}(v^s_{\alpha}u'_{\alpha\beta}) \})] 
\,\in \,\mathbb{H}^2(\mathbb{K}).
\end{aligned}
\end{equation*}
On the other hand, we compute 
\begin{equation*}
p_1^* (s^* \Theta^e) \,\colon\, 
 \mathbb{H}^1 ((\mathcal{C}_{\bullet})_z) \otimes 
\mathbb{H}^1 ((\mathcal{C}_{\bullet})_z) \xrightarrow{(s\circ p_1)_* \otimes (s\circ p_1)_* } 
 \mathbb{H}^1 ((\mathcal{C}_{\bullet})_{s\circ p_1(z)}) \otimes 
 \mathbb{H}^1 ((\mathcal{C}_{\bullet})_{s\circ p_1(z)}) \xrightarrow{\ \Theta^e \ } 
\mathbb{H}^1(\mathbb{K})
\end{equation*}
as follows:
\begin{equation*}
\begin{aligned}
p_1^* (s^* \Theta^e)(v,\,v')&= 
p_1^* (s^* \Theta^e)([(\{u_{\alpha\beta}\},\, \{ v^s_{\alpha} \})]
,[(\{u'_{\alpha\beta}\}, \,\{ (v^s_{\alpha})' \})]) \\
&=[(\{ \mathrm{Tr} (u_{\alpha\beta}u'_{\beta\gamma}) \},\,
- \{\mathrm{Tr}(u_{\alpha\beta}(v^s_{\beta})')-\mathrm{Tr}(v^s_{\alpha}u'_{\alpha\beta}) \})] 
\,\in\, \mathbb{H}^1(\mathbb{K}).
\end{aligned}
\end{equation*}
Therefore, we have the equality 
$\Theta^e|_{p_1^{-1}(U)} - (P_1^{-1})^* \Phi_U \,=\, p_1^* (s^* \Theta^e|_{p_1^{-1}(U)})$.
\end{proof}

It will now be shown that the restriction of
$\Theta^e$ to $\mathcal{M}^{e}_{\mathrm{FC}}(d)^{\circ}$ is nondegenerate:

\begin{Cor}\label{cor1}
The $2$-form $\Theta^e|_{\mathcal{M}^{e}_{\mathrm{FC}}(d)^{\circ}}$ is nondegenerate.
\end{Cor}

\begin{proof}
For any point $(E,\,\phi,\,\nabla)\, \in\, \mathcal{M}^{e}_{\rm{FC}}(d)^{\circ}$, and any tangent vectors
$v,\, w\, \in\, T_{(E,\phi,\nabla)} \mathcal{M}^{e}_{\rm{FC}}(d)^{\circ}$, we have
$$
p_1^* s^* \Theta^e (E,\phi,\nabla)(v,\, w)\,=\, 0
$$
when one of $v$ and $w$ is vertical for the projection $p_1$ in \eqref{fm}. So, if
$w$ is vertical, from Lemma \ref{2019.11.19.15.22} it follows that
\begin{equation}\label{ii0}
\Theta^e(E,\,\phi,\,\nabla)(v,\, w) \,=\, (P_1^{-1})^* \Phi_U (E,\,\phi,\,\nabla)(v,\, w)\, .
\end{equation}
Since $\Phi_U$ is a symplectic form, there is a tangent vector $v\, \in\, 
T_{(E,\phi,\nabla)} \mathcal{M}^{e}_{\rm{FC}}(d)^{\circ}$ such that
$$(P_1^{-1})^* \Phi_U (E,\,\phi,\,\nabla)(v,\, w)\, \not=\, 0.$$ Now from
\eqref{ii0} it follows that $\Theta^e(E,\,\phi,\,\nabla)(v,\, w)\, \not=\, 0$.

Since the vertical tangent spaces for the projection $T^*U\, \longrightarrow\, U$ are
Lagrangian for the Liouville $2$-form $\Phi_U$, given any non-vertical tangent
vector
$$
v\, \in\, T_{(E,\phi,\nabla)} \mathcal{M}^{e}_{\rm{FC}}(d)^{\circ}
$$
for the projection $T^*U\, \longrightarrow\, U$,
there is a vertical tangent vector
$$
w\, \in\, T_{(E,\phi,\nabla)} \mathcal{M}^{e}_{\rm{FC}}(d)^{\circ}
$$
for the projection $T^*U\, \longrightarrow\, U$,
such that $(P_1^{-1})^* \Phi_U (E,\,\phi,\,\nabla)(v,\, w)\,\not=\, 0$. Now from \eqref{ii0}
it follows that $\Theta^e(E,\,\phi,\,\nabla)(v,\, w)\, \not=\, 0$.
Consequently, the form $\Theta^e\vert_{p^{-1}(U)}$ is nondegenerate.
\end{proof}

\begin{Rem}
It was shown above that the restriction of
$\Theta^e$ to $\mathcal{M}^{e}_{\mathrm{FC}}(d)^{\circ}$
is nondegenerate by using Lemma \ref{2019.11.19.15.22}.
On the other hand, we will show that the 2-form $\Theta^e$ 
on $\mathcal{M}^{e}_{\mathrm{FC}}(d)$ is nondegenerate by using the Serre duality 
(Proposition \ref{2023.1.10.23.18} below).
So it can be shown, without using Lemma \ref{2019.11.19.15.22}, that the restriction of 
$\Theta^e$ to $\mathcal{M}^{e}_{\mathrm{FC}}(d)^{\circ}$ is nondegenerate.
Nevertheless, we have discussed nondegeneracy of the restriction 
of $\Theta^e$ by using this lemma, because this argument highlights another important perspective. 
On the other hand, Lemma \ref{2019.11.19.15.22} will be used below
in the proof of the $d$-closedness of the restriction of $\Theta^e$.
Moreover, the $d$-closedness of the restriction of $\Theta^e$ will be used below
in the proof of the $d$-closedness of $\Theta^e$ on $\mathcal{M}^{e}_{\rm{FC}}(d)$.
\end{Rem}

\begin{Prop}\label{2019.12.22.18.19}
Assume that $g\,\ge\, 2$. 
Let ${\mathcal M}^e_{\mathrm{FC}}(d)^{\circ}$
be the open subspace of $\mathcal{M}^{e}_{\mathrm{FC}}(d)$
defined in \eqref{2021.11.26.11.19} for $H=\{e\}$.
Then the restriction $\Theta^e|_{\mathcal{M}^{e}_{\mathrm{FC}}(d)^{\circ}}$ 
of the nondegenerate $2$-form $\Theta^e$ in
\eqref{ted} is $d$-closed.
\end{Prop}

\begin{proof}
The moduli space $\mathcal{N}^e(d)$ in \eqref{ned} has the open subset
$\mathcal{N}^e(d)^{\circ}$ defined by
\begin{equation*}
\mathcal{N}^e(d)^{\circ}\, :=\, \left\{ (E,\,\phi ) \,\in\, \mathcal{N}^e(d) \ \big\vert \ 
\text{$E$ is a stable vector bundle}\right\}.
\end{equation*}
Also, let $\mathcal{M}_{\mathrm{FC}}^{e}(d)^{\circ\circ}\, \subset\,
\mathcal{M}_{\mathrm{FC}}^{e}(d)^{\circ}$ (see Definition \ref{2019.10.28.16.49}) be the open subset
\begin{equation}\label{dm1}
\mathcal{M}_{\mathrm{FC}}^{e}(d)^{\circ\circ}\,:=\, \left\{ (E,\,\phi ,\,\nabla)
\,\in\, \mathcal{M}_{\mathrm{FC}}^{e}(d)^{\circ}
 \ \big\vert \ 
\text{\,$E$ is a stable vector bundle}\right\}.
\end{equation}
The openness of both $\mathcal{N}^e(d)^{\circ}$ and $\mathcal{M}_{\mathrm{FC}}^{e}
(d)^{\circ\circ}$ follows from \cite[p.~635, Theorem 2.8(B)]{Ma}.
The moduli spaces $\mathcal{N}^e(d)^{\circ}$ 
and $\mathcal{M}_{\mathrm{FC}}^{e}(d)^{\circ\circ}$ 
are non-empty because $g\,\ge\, 2$.

To prove that the form $\Theta^e|_{\mathcal{M}^{e}_{\mathrm{FC}}(d)^{\circ}}$ on 
$\mathcal{M}_{\mathrm{FC}}^{e}(d)^{\circ}$ is closed, it suffices to show that the restriction of
$\Theta^e|_{\mathcal{M}^{e}_{\mathrm{FC}}(d)^{\circ}}$ to $\mathcal{M}_{\mathrm{FC}}^{e}(d)^{\circ\circ}$ is closed.

Let $p_{1,0}\, :\, \mathcal{M}^{e}_{\rm{FC}}(d)^{\circ\circ} \,\longrightarrow \,\mathcal{N}^e(d)^{\circ}$
be the restriction of the forgetful map $p_1$ in \eqref{fm}.
Take a sufficiently small analytic open subset $U
\, \subset\, \mathcal{N}^e(d)^{\circ}$ such that there is a holomorphic section
$$s\,\colon\, U \,\longrightarrow \, p_{1,0}^{-1}(U)\, ,$$ over $U$, of $p_{1,0}$. Now
Lemma \ref{2019.11.19.15.22} says that
\begin{equation*}
\Theta^e - (P_1^{-1})^* \Phi_U \,= \,p^*_{1,0} (s^* \Theta^e)
\end{equation*}
on $p_{1,0}^{-1}(U)$.
This implies that
\begin{equation}\label{cf}
d \Theta^e\,=\, p^*_{1,0} d (s^* \Theta^e)\, 
\end{equation}
on $p_{1,0}^{-1}(U)$, because the Liouville $2$-form is $d$-closed.

In view of \eqref{cf}, to prove the theorem it suffices to show the existence of a local
holomorphic section $s\,\colon\, U \,\longrightarrow\, p^{-1}_{1,0}(U)$ of the map $p^{-1}_{1,0}$
such that $d(s^*\Theta^e)\,=\, 0$.

We shall construct a holomorphic section $s\,\colon\, U \,\longrightarrow\, p^{-1}_{1,0}(U)$ such that
$$d(s^*\Theta^e)\,=\, 0\, .$$ For that, first define a moduli space
\begin{equation*}
\mathcal{M}^{e}_{\rm{FC}}(d)^{\circ\circ}_{0}\,:= \,\left\{ (E,\phi,\nabla) \ \middle| \ 
\begin{array}{l}
\text{$E$ is a stable vector bundle of degree $d$, and}\\
\text{$(E, \phi,\nabla)$ is a framed $\mathrm{GL}(r,\mathbb{C})$-connection such that} \\
\text{$\mathrm{res}_{x_i}(\nabla)\,=\,0$ for $i=1,\cdots,n-1$ and $\mathrm{res}_{x_n}( \nabla)=-\frac{d}{r}e$}
\end{array}
\right\}{\Big/}\Large{\sim},
\end{equation*}
where $e$ is the identity matrix. There is the natural inclusion map
\begin{equation}\label{dio}
\iota\,\,\colon\,\, \mathcal{M}^{e}_{\rm{FC}}(d)^{\circ\circ}_{0}\,\,\hookrightarrow\,\,
\mathcal{M}^{e}_{\rm{FC}}(d)^{\circ\circ}\, ,
\end{equation}
where $\mathcal{M}^{e}_{\rm{FC}}(d)^{\circ\circ}$ is defined in \eqref{dm1}. Also, define two moduli spaces
\begin{equation*}
\mathcal{M}(d)_{0}^{\circ\circ}\,:= \,\left\{ (E,\nabla) \ \middle| \ 
\begin{array}{l}
\text{$E$ is a stable vector bundle of rank $r$ and degree $d$, and}\\
\text{$\nabla \,\colon\, E \,\longrightarrow \,E\otimes K_X(D)$ is a connection such that } \\
\text{$\mathrm{res}_{x_i}( \nabla)\,=\,0$ for $i\,=\,1,\cdots,n-1$ and
$\mathrm{res}_{x_n}(\nabla)\,=\,-\frac{d}{r}e$}
\end{array}
\right\}{\Big/}\Large{\sim}
\end{equation*}
and
\begin{equation*}
\mathcal{N}(d)^{\circ}\,=\, \{ E \mid 
\text{$E$ is a stable vector bundle of rank $r$ and degree $d$}
\}/\sim.
\end{equation*}
There are the forgetful maps 
\begin{equation}\label{fm1}
q_1 \,\colon\, \mathcal{M}(d)^{\circ\circ}_{0}\,\longrightarrow \,\mathcal{N}(d)^{\circ},\
q_2 \,\colon\, \mathcal{M}^{e}_{\rm{FC}}(d)^{\circ\circ}_{0}\,\longrightarrow\,
\mathcal{M}(d)^{\circ\circ}_{0}\ \text{and}\
p_2 \,\colon\, \mathcal{N}^e(d)^{\circ} \,\longrightarrow \,\mathcal{N}(d)^{\circ}\, ,
\end{equation}
where $q_2$ and $p_2$ forget the framing while $q_1$ forgets the connection.

Take an analytic open subset $U_0\,\subset\, \mathcal{N}(d)^{\circ}$. Assume that
$U_0$ is small enough and that
the image of $U$ under the forgetful map $p_2 \,\colon\, \mathcal{N}^e (d)^{\circ}\,\longrightarrow\,
\mathcal{N}(d)^{\circ}$ is contained in $U_0$, by shrinking sufficiently the analytic open subset $U$.
Take a holomorphic section
\begin{equation*}
\begin{aligned}
s_0 \,\colon\,\, U_0 &\,\longrightarrow\, q_1^{-1}(U_0) \\
E &\,\longmapsto\, (E,\,\nabla(E))
\end{aligned}
\end{equation*}
of the forgetful map 
$q_1 \,\colon\, \mathcal{M}(d)^{\circ\circ}_{0}\,\longrightarrow\, \mathcal{N}(d)^{\circ}$.
Since $H_x\,=\,\{ e\}$ for all $x\in D$, we may define a section $\widetilde{s}$ on $U$ 
\begin{equation*}
\begin{aligned}
\widetilde{s} \,\colon\,\, U &\,\longrightarrow\, 
\mathcal{M}^{e}_{\rm{FC}}(d)^{\circ\circ}_{0} \\
(E , \,\phi ) &\,\longmapsto\, (E, \,\phi,\, \nabla(E))
\end{aligned}
\end{equation*}
using the section $s_0$.
Define the section $s$ on $U$ of $p_{1,0} \,\colon\, p_1^{-1}(U)\,\longrightarrow\, U$ by
$$s \,=\, \iota \circ \widetilde{s}\, .$$

Now we shall prove that 
\begin{equation}\label{ts1}
d (s^* \Theta^e)\,=\,0
\end{equation}
for such a section.

To prove \eqref{ts1}, first
recall that the moduli space $\mathcal{M}(d)^{\circ\circ}_{0}$ is equipped with
a natural symplectic structure.
We briefly describe this symplectic structure on $\mathcal{M}(d)^{\circ\circ}_{0}$.
The tangent space to $\mathcal{M}(d)^{\circ\circ}_{0}$ at any point $(E,\,\nabla)$ is isomorphic to
the first hypercohomology $\mathbb{H}^1(\mathcal{C}^0_{\bullet})$, where
\begin{equation}\label{cc0}
\mathcal{C}^0_{\bullet}\, \,\colon\,\, \mathcal{C}^0_0\,=\,
\mathcal{E}nd (E) \, \xrightarrow{\ \nabla \ }\, 
\mathcal{C}^0_1\,=\,\mathcal{E}nd(E)\otimes K_X\, .
\end{equation}
Define a nondegenerate $2$-form $\Theta_0$ on $\mathcal{M}(d)^{\circ\circ}_{0}$
$$
\Theta_0(E,\,\nabla) \,\colon\, 
\mathbb{H}^1(\mathcal{C}^0_{\bullet})
\otimes \mathbb{H}^1(\mathcal{C}^0_{\bullet}) 
\,\longrightarrow \,\mathbb{C}
$$
exactly as done in \eqref{ted}. 
This $2$-form $\Theta_0$ is $d$-closed, which is proved in \cite{Go}.

As the second step to prove \eqref{ts1}, it will be shown that 
\begin{equation}\label{cl1}
\iota^* \Theta^e \,=\, q_2^* \Theta_0\, ,
\end{equation}
where $q_2$ and $\iota$ are the maps in \eqref{fm1} and \eqref{dio} respectively.
To prove \eqref{cl1}, note that the tangent space of 
$\mathcal{M}^{e}_{\rm{FC}}(d)^{\circ\circ}_{0}$ at $(E,\,\phi,\,\nabla)$ is isomorphic to
the first hypercohomology $\mathbb{H}^1(\mathcal{C}'_{\bullet})$ of the complex
\begin{equation}\label{2019.11.30.13.52}
\mathcal{C}'_{\bullet} \,\colon\, 
\mathcal{C}'_0\,=\,\mathcal{E}nd (E) (-D) \,\xrightarrow{\ \nabla \ }\, 
\mathcal{C}'_1\,=\,\mathcal{E}nd (E) \otimes K_X\, .
\end{equation}
For $[(\{u_{\alpha\beta}\} ,\{ v_{\alpha}\})]\,\in\, \mathbb{H}^1(\mathcal{C}'_{\bullet})$,
we have that
\begin{equation*}
\iota_*[(\{ u_{\alpha\beta}\} , \,\{ v_{\alpha}\} ) ]
\,=\,[ (\{ u_{\alpha\beta}\} , \, \{ v_{\alpha} \} ) ]\,\in\, \mathbb{H}^1(\mathcal{C}_{\bullet})
\ \text{ and } \ 
(q_2)_*[ (\{ u_{\alpha\beta} \} ,\, \{ v_{\alpha}\} ) ]
\,=\,[ (\{ u_{\alpha\beta}\} , \,\{ v_{\alpha}\} ) ] \,\in\, \mathbb{H}^1(\mathcal{C}^0_{\bullet}).
\end{equation*}
Therefore $\iota^* \Theta^e$ and $q_2^* \Theta_0$ have the following identical description:
\begin{equation*}
\begin{aligned}
\mathbb{H}^1(\mathcal{C}'_{\bullet}) \otimes
\mathbb{H}^1(\mathcal{C}'_{\bullet}) \, &\longrightarrow \, \mathbb{H}^2(\mathbb{K})
\,\cong\, \mathbb{C} \\
[(\{ u_{\alpha\beta} \} ,\,\{ v_{\alpha}\}) ] \otimes
[(\{ u_{\alpha\beta}' \} ,\,\{ v_{\alpha}'\}) ]\, &\longmapsto \, 
[(\{ \mathrm{Tr} (u_{\alpha\beta}u'_{\beta\gamma}) \},\,
- \{\mathrm{Tr}(u_{\alpha\beta}v'_{\beta})-\mathrm{Tr}(v_{\alpha}u'_{\alpha\beta}) \})].
\end{aligned}
\end{equation*}
This proves \eqref{cl1}.

Thirdly, by the equality $\iota^* \Theta^e \,=\, q_2^* \Theta_0$ in \eqref{cl1} we have
\begin{equation*}
s^* \Theta^e \,=\, \widetilde{s}^* (q_2^* \Theta_0)\, .
\end{equation*}
Since $\Theta_0$ is $d$-closed, it follows that $d(s^* \Theta^e)\,=\,0$, proving
\eqref{ts1}.

Finally, from the combination of \eqref{ts1} and the equality $d \Theta^e\,=\, p^*_{1,0} d (s^* \Theta^e)$
(see \eqref{cf}), it follows that $d \Theta^e\,=\,0$ on $p^{-1}_{1,0}(U)$.
This implies that the $2$-form $\Theta^e$ is $d$-closed on 
$\mathcal{M}_{\rm{FC}}^{e}(d)^{\circ}$. As noted before, this proves the theorem.
\end{proof}

Next the $d$-closedness of $\Theta^e$ in \eqref{ted} will be shown when $g\,=\,0$ and $g\,=\,1$.
For this purpose, we recall the definition of \textit{parabolic connections}.
Let $$(X,\,\boldsymbol{x})\, := \,(X,\, (x_1,\, \cdots,\,x_n))$$
be an $n$-pointed smooth projective curve of genus $g$ 
over $\mathbb{C}$, where $x_1, \,\cdots,\,x_n$ are distinct points of $X$.
Denote the reduced divisor $x_1+\cdots+x_n$ on $X$ by $D(\boldsymbol{x})$
or simply by $D$ if there is no possibility of confusion.
Take a positive integer $r$.

\begin{Def}\label{2020.3.27.15.05} 
A \textit{$\boldsymbol{x}$-quasi-parabolic bundle} of rank $r$ and degree $d$
is a pair $(E,\, \boldsymbol{l}\, =\, \{ l^{(i)}_* \}_{1\le i \le n})$, where
\begin{itemize}
\item[(1)] $E$ is an algebraic vector bundle on $X$ of rank $r$ and degree $d$, and

\item[(2)] $l^{(i)}_* $ is a filtration of subspaces
$E\big\vert_{x_i} \,=\, l_0^{(i)} \,\supset\, l_1^{(i)}\,\supset \,\cdots \,\supset\,
l_r^{(i)}\,=\,0$ for every $1\, \leq\, i\, \leq\, n$ 
such that $\dim (l_j^{(i)}/l_{j+1}^{(i)})\,=\,1$.
\end{itemize}
\end{Def}

Let $\boldsymbol{\alpha}$ be a tuple
$(\alpha^{(i)}_j)^{1\leq i\leq n}_{1\leq j\leq r}$ of real numbers which satisfy the condition
\[
 0\,<\,\alpha^{(i)}_1\,<\,\alpha^{(i)}_2\,<\, \cdots\,<\, \alpha^{(i)}_r\,<\,1
\]
for each $1\leq i\,\leq\, n$, and
$\alpha^{(i)}_j\,\neq\, \alpha^{(i')}_{j'}$ for all $(i,\,j)\,\neq\,(i',\,j')$.
We call the tuple $\boldsymbol{\alpha}$ a parabolic weight.
Take an element
$$\boldsymbol{\nu}\,=\,(\nu^{(i)}_j )^{1\le i \le n}_{0\le j \le r-1} \,\in\, \mathbb{C}^{nr}$$
such that $\sum_{i,j} \nu^{(i)}_j \,=\, -d \,\in\, \mathbb{Z}$.

\begin{Def}\label{def: parabolic connection}
A quadruple
$(E,\,\nabla,\, \boldsymbol{l}\,=\,\{ l^{(i)}_* \}_{1\le i \le n},\boldsymbol{\alpha})$
is called 
a \textit{$(\boldsymbol{x},\, \boldsymbol{\nu})$-parabolic connection}
of rank $r$ and degree $d$ if 
\begin{itemize}
\item[(1)] $(E,\, \boldsymbol{l}\, =\, \{ l^{(i)}_* \}_{1\le i \le n})$ is
a $\boldsymbol{x}$-quasi-parabolic bundle of rank $r$ and degree $d$, and

\item[(2)] $\nabla\,\colon\, E \,\longrightarrow \,E\otimes K_{X}(D)$ is a
logarithmic connection whose residue
$\mathrm{res}_{x_i}(\nabla)\,\colon\, E|_{x_i}\,\longrightarrow\, E|_{x_i}$
at each point $x_i$ for $1\,\leq\, i\,\leq\, n$ satisfies the condition
$(\mathrm{res}_{x_i}(\nabla)-\nu^{(i)}_j \mathrm{Id}_{E|_{x_i}}) (l^{(i)}_j)
\,\subset\, l_{j+1}^{(i)}$ for all $j\,=\,0,\,\cdots,\,r-1$.
\end{itemize}
\end{Def}

\begin{Def}\label{def: stability of parabolic connection}
A $(\boldsymbol{x},\ \boldsymbol{\nu})$-parabolic connection
$(E,\ \nabla,\ \boldsymbol{l},\,\boldsymbol{\alpha})$ is said to be
$\boldsymbol{\alpha}$-stable if
the inequality
\[
 \frac {\deg F+\sum_{i=1}^n\sum_{j=1}^r \alpha^{(i)}_j
 \dim ( (F|_{x_i}\cap l^{(i)}_{j-1})/(F|_{x_i}\cap l^{(i)}_j)) } {\rank F}
\, <\,
 \frac {\deg E+\sum_{i=1}^n\sum_{j=1}^r \alpha^{(i)}_j
 \dim ( l^{(i)}_{j-1}/l^{(i)}_j ) } {\rank E}
\]
holds for every subbundle $0\,\neq\, F\,\subsetneq\, E$ for which
$\nabla(F)\,\subset\, F\otimes\Omega^1_X(D)$.
We say that $(E,\ \nabla,\ \boldsymbol{l},\,\boldsymbol{\alpha})$
is $\boldsymbol{\alpha}$-semistable
if the weaker inequality ``$\leq$'' holds (instead of ``$<$'').
\end{Def}

\begin{Rem} 
In the non-abelian Hodge correspondence (see \cite{Sim0}),
the parabolic weight $\boldsymbol{\alpha}$ is an important
datum needed to connect to the parabolic Higgs bundles.
Since we focus on the algebraic moduli spaces,
we omit the parabolic weight $\boldsymbol{\alpha}$
to denote the parabolic connection.
So we denote a parabolic connection by $(E,\,\nabla,\, \boldsymbol{l})$, 
even though there is the parabolic weight $\boldsymbol{\alpha}$ in the background.
\end{Rem}

In the inequality for the stability condition
in Definition \ref{def: stability of parabolic connection}, we may replace the parabolic weight with a tuple of
rational numbers which is very close to $\boldsymbol{\alpha}$.
We have the following.

\begin{Thm}[{\cite[Theorem 2.2] {Inaba-1}}]
The moduli space
${\mathcal M}_{\mathrm{PC}}^{\boldsymbol{\alpha}}(\boldsymbol{\nu})$ of
$\boldsymbol{\alpha}$-stable $(\boldsymbol{x},\boldsymbol{\nu})$-parabolic connections exists as a quasi-projective scheme over $\Spec\mathbb{C}$.
\end{Thm}

Let $(E,\, \boldsymbol{l})$ be a $\boldsymbol{x}$-quasi-parabolic bundle.
Set
\[
 \End(E,\,\boldsymbol{l})\,\,:=\,\,
 \left\{ u\,\in\, \Hom_{{\mathcal O}_X}(E,\,E) \,\,\middle|\,\,\,
 \text{$u|_{x_i}(l^{(i)}_j)\,\subset\, l^{(i)}_j$ for any $i,\,j$}\right\}.
\]
We denote the invertible elements of $\End(E,\,\boldsymbol{l})$
by $\Aut(E,\,\boldsymbol{l})$.

\begin{Def}
A $\boldsymbol{x}$-quasi-parabolic bundle $(E,\, \boldsymbol{l})$ is said to be \textit{simple} 
if $\End(E,\, \boldsymbol{l}) \,=\,\mathbb{C}$,
which is equivalent to the condition that $\Aut(E,\, \boldsymbol{l}) \,=\,\mathbb{C}^{*}$.
\end{Def}

\begin{Rem} 
For each $x_i$, let $\mathbf{H}_{x_i}\subset\mathrm{GL}_r(\mathbb{C})$
be the Borel subgroup consisting of the upper triangular matrices.
Then a framed $\mathrm{GL}_r(\mathbb{C})$--bundle with respect to
the structure subgroups $(\mathbf{H}_{x_i}\subset \mathrm{GL}_r(\mathbb{C}))_{1\leq i\leq n}$
is equivalent to a $\boldsymbol{x}$-quasi-parabolic bundle.
The above definition of simple quasi-parabolic bundle is
equivalent to that of simple framed bundle with this structure subgroup
in the sense of Definition \ref{2019.10.28.16.49}.
A framed $\mathrm{GL}_r(\mathbb{C})$--connection with respect to
the structure subgroups $(\mathbf{H}_{x_i}\subset \mathrm{GL}_r(\mathbb{C}))_{1\leq i\leq n}$
is equivalent to a $(\boldsymbol{x},\ \boldsymbol{0})$-parabolic connection,
where $\mathbf{0}\,\in\,\mathbb{C}^{nr}$ is defined by
$\nu^{(i)}_j\,=\,0$ for any $i,\,j$.
\end{Rem}

For a $(\boldsymbol{x},\, \boldsymbol{\nu})$-parabolic connection
$(E,\,\nabla,\, \boldsymbol{l})$, set
\[
 \End(E,\,\nabla,\,\boldsymbol{l})\,:=\,
 \left\{ u\in \End(E,\,\boldsymbol{l}) \,\,\middle|\,\,\,
 \nabla\circ u\,=\,(u\otimes\mathrm{id})\circ\nabla \right\},
\]
and denote by $\Aut(E,\,\nabla,\,\boldsymbol{l})$
the invertible elements in $\End(E,\,\nabla,\,\boldsymbol{l})$.

\begin{Def}
A $(\boldsymbol{x},\, \boldsymbol{\nu})$-parabolic connection
$(E,\,\nabla,\, \boldsymbol{l})$ is said to be \textit{simple} 
if $\End (E,\,\nabla,\, \boldsymbol{l}) \,=\,\mathbb{C}$,
which is equivalent to the condition that
$\Aut (E,\,\nabla,\, \boldsymbol{l}) \,=\,\mathbb{C}^*$.
\end{Def}

An argument similar to the one in Proposition \ref{2020.2.24.14.52} proves 
the following proposition.

\begin{Prop}\label{Prop: moduli of simple parabolic connections}
The moduli space 
${\mathcal M}_{\mathrm{PC}}(\boldsymbol{\nu})$
of simple $(\boldsymbol{x},\boldsymbol{\nu})$-parabolic connections
exists as an algebraic space.
The moduli space 
${\mathcal M}_{\mathrm{PC}}^{\boldsymbol{\alpha}}(\boldsymbol{\nu})$ of
$\boldsymbol{\alpha}$-stable $(\boldsymbol{x},\boldsymbol{\nu})$-parabolic connections is a Zariski open subspace of
${\mathcal M}_{\mathrm{PC}}(\boldsymbol{\nu})$.
\end{Prop}

\begin{Prop}\label{2019.12.22.18.19_a}
Assume that either $g\,=\,0$ or $g\,=\,1$ and
\begin{itemize}
\item
$nr-2r-2>0$ if $g=0$,
\item
$n\geq 2$ if $g=1$.
\end{itemize}
Let ${\mathcal M}^e_{\mathrm{FC}}(d)^{\circ}$
be the open subspace of $\mathcal{M}^{e}_{\mathrm{FC}}(d)$
defined in \eqref{2021.11.26.11.19} for $H=\{e\}$.
Then the restriction $\Theta^e|_{\mathcal{M}^{e}_{\mathrm{FC}}(d)^{\circ}}$ 
of the nondegenerate $2$-form $\Theta^e$ in
\eqref{ted} is $d$-closed.
\end{Prop}

\begin{proof}
Consider the forgetful map $p_1\,\colon\, \mathcal{M}^{e}_{\rm{FC}}(d)^{\circ} \,\longrightarrow\,
\mathcal{N}^e(d)$ in \eqref{fm}, and take a sufficiently small analytic open subset $U\, \subset\,
\mathcal{N}^e(d)$. For a holomorphic section $s\,\colon\, U \,\longrightarrow\, p_1^{-1}(U)$
of $p_1$ such that 
\[
d(s^*\Theta^e|_{\mathcal{M}^{e}_{\mathrm{FC}}(d)^{\circ}})\,=\,0,
\]
we have $d \:\! \Theta^e|_{\mathcal{M}^{e}_{\mathrm{FC}}(d)^{\circ}}\,=\,0$
by the same argument as in the proof of
Theorem \ref{2019.12.22.18.19}. We will now construct such a section $s$.

Let $\mathcal{N}_{\text{par}} (d)$
be the moduli space of simple $\boldsymbol{x}$-quasi-parabolic bundles
of rank $r$ and degree $d$.
For each $x \,\in\, D$, set the complex Lie proper subgroup $H_x$
to be the subgroup of $\mathrm{GL}(r,\mathbb{C})$ 
consisting of the upper triangular matrices. It may be mentioned 
that a $\boldsymbol{x}$-quasi-parabolic bundle is the same as 
a framed bundle with respect to $\{ H_x\}_{x \in D}$. For a framed bundle
$(E,\,\phi)$, we can associate a quasi-parabolic bundle $(E,\,\boldsymbol{l})$
whose filtration $l^{(i)}_*$ on $E|_{x_i}$
is induced by the framing $\phi_{x_i}$ of $E|_{x_i}$ for each $1\,\leq\, i\,\leq\, n$.
Setting
\[
 {\mathcal N}^e(d)^{\circ}\,:=\,
 \left\{(E,\phi)\,\,\middle|\,\,
 \text{The quasi-parabolic bundle $(E,\,\boldsymbol{l})$ induced by the framing $\phi$
 is simple} \right\},
\]
there is a natural morphism
\begin{equation}
 q_B\ \colon\ {\mathcal N}^e(d)^{\circ}\ \longrightarrow\
 {\mathcal N}^{\mathrm{par}}(d)^{\circ}.
\end{equation}

Take an element $$\boldsymbol{\nu}\,=\,(\nu_{j}^{(i)} )^{1\le i \le n}_{0\le j \le r-1}
\,\in\, \mathbb{C}^{nr}$$
such that $\sum_{i,j} \nu_{j}^{(i)} \,=\, -d$. 
Let $\mathcal{M}_{\text{PC}} (\boldsymbol{\nu})^{\circ}$ 
be the moduli space defined by
\begin{equation*}
\mathcal{M}_{\text{PC}} (\boldsymbol{\nu})^{\circ}
\,:=\,\left\{ (E,\,\nabla,\, \boldsymbol{l})\in \mathcal{M}_{\text{PC}} (\boldsymbol{\nu}) \, \middle| \,
\begin{array}{l}
\text{$(E,\, \boldsymbol{l})$ is a simple $\boldsymbol{x}$-quasi-parabolic bundle}
\end{array}
\right\}{\Big/}\Large{\sim}.
\end{equation*}
Define the locally closed subspace
$\mathcal{M}_{\text{FC}}^{e} (\boldsymbol{\nu})^{\circ}$
of $\mathcal{M}_{\text{FC}}^{e} (d)^{\circ}$ by
\begin{equation*}
\mathcal{M}_{\text{FC}}^{e} (\boldsymbol{\nu})^{\circ} \ :=
\left\{ (E,\, \phi ,\,\nabla) \,\in\, \mathcal{M}_{\text{FC}}^{e} (d)^{\circ} \ \middle| \ 
\begin{array}{l}
\text{the framed bundle $(E,\phi)$ belongs to ${\mathcal N}^e(d)^{\circ}$} \\
\text{and $(E,\,\boldsymbol{l}):=q_B(E,\,\phi)$ satisfies} \\
\text{$(\mathrm{res}_{x}(\nabla)-\nu^{(i)}_{j} \mathrm{Id}_{E|_{x_i}}) (l^{(i)}_{j})\, \subset\, l_{j+1}^{(i)}$ 
for $0\leq j\leq r-1$}
\end{array}
\right\}{\Big/}\Large{\sim}.
\end{equation*}
In the above definition we have $\phi\,=\,\{\phi_x\}_{x\in D}$, where $\phi_x\,\colon\,
\mathcal{O}_X^{\oplus r}|_x \,\longrightarrow\, E\vert_x$ 
are isomorphisms defining a framing of $E$ over $D$. Since a
framing defines a parabolic structure, there is a natural map 
\begin{equation}
q_1^{\text{par}} \ \colon\ \mathcal{M}_{\text{FC}}^{e} (\boldsymbol{\nu})^{\circ}
\ \longrightarrow\ \mathcal{M}_{\text{PC}} (\boldsymbol{\nu})^{\circ}\, .
\end{equation}
Notice that $\mathcal{M}_{\text{PC}} (\boldsymbol{\nu})^{\circ}$
is non-empty by virtue of the assumption in the proposition,
and so is $\mathcal{M}_{\text{FC}}^{e} (\boldsymbol{\nu})^{\circ}$.
Consider the complex 
\begin{equation*}
\mathcal{D}_{\bullet}^{\text{par}} \ \colon\
\mathrm{ad}_{\phi} (E_G) \
\xrightarrow{\,\ [\, \nabla,\, \cdot \, ] \,\, \ } \
\mathrm{ad}^n_{\phi} (E_G) \otimes K_X (D)\, 
\end{equation*}
for $\{ H_x\}_{x \in D}$.
Here $\mathrm{ad}_{\phi} (E_G)$ and
$\mathrm{ad}^n_{\phi} (E_G) \otimes K_X (D)$ are defined as in \eqref{db}.
The tangent space of $\mathcal{M}_{\text{PC}} (\boldsymbol{\nu})^{\circ}$ at
$(E,\,\nabla, \,\boldsymbol{l})$ is $\mathbb{H}^1(\mathcal{D}^{\text{par}}_{\bullet})$.
There is also a natural morphism
\[
 p_0^{\mathrm{par}}\ \colon\ 
 \mathcal{M}_{\text{PC}} (\boldsymbol{\nu})^{\circ}
\ \longrightarrow\ {\mathcal N}_{\mathrm{par}}(d)^{\circ}
\]
which is \'etale locally an affine space bundle whose fiber is isomorphic to
$H^0(X,\,\mathrm{ad}^n_{\phi}(E_G)\otimes K_X(D))$. So there is a non-empty analytic open subset 
$U\,\subset\, {\mathcal N}^{\mathrm{par}}(d)^{\circ}$ with a local section
$s^{\mathrm{par}}\,\colon\, U \,\longrightarrow\, (p_0^{\mathrm{par}})^{-1}(U)$
of $p_0^{\mathrm{par}}$. Consider the following commutative diagram
\[
\begin{CD}
\mathcal{M}_{\text{PC}} (\boldsymbol{\nu})^{\circ}
@< q_1^{\mathrm{par}} <<
\mathcal{M}_{\text{FC}}^{e} (\boldsymbol{\nu})^{\circ}
& \; \xrightarrow[\;\subset\;]{\iota} \;& \mathcal{M}_{\text{FC}}^{e} (d)^{\circ}
\\
@V p_0^{\mathrm{par}} VV
@V V p_1|_{\mathcal{M}_{\text{FC}}^{e} (\boldsymbol{\nu})^{\circ}} V
\\
 {\mathcal N}^{\mathrm{par}}(d)^{\circ}
 @< q_B << {\mathcal N}^e(d)^{\circ} .
\end{CD}
\]
whose left square is Cartesian. The local section
$s^{\mathrm{par}}$ of $p_0^{\mathrm{par}}$
produces a local section
$s_1\,\colon\, q_B^{-1}(U)\,\longrightarrow\, (q_B\circ p_1)^{-1}(U)$
of $p_1|_{\mathcal{M}_{\text{FC}}^{e} (\boldsymbol{\nu})^{\circ}}$.
Let $\Theta_{\text{par}}$ be the symplectic structure 
on $\mathcal{M}_{\text{PC}}^{\boldsymbol{\alpha}} (\boldsymbol{\nu})$ 
constructed in \cite{Inaba-1}.
This symplectic form $\Theta_{\text{par}}$ is 
described as follows:
\begin{equation}\label{2020.1.26.14.17}
\begin{aligned}
\Theta_{\text{par}} \,\colon\,
\mathbb{H}^1 (\mathcal{D}^{\text{par}}_{\bullet}) \otimes \mathbb{H}^1 (\mathcal{D}^{\text{par}}_{\bullet}) 
&\,\longrightarrow\, \mathbb{H}^2(\mathbb{K}) \,\cong\, \mathbb{C} \\
[(\{u_{\alpha\beta} \} ,\, \{ v_{\alpha} \})]\otimes[(\{u'_{\alpha\beta} \} ,\, \{ v'_{\alpha} \})]
&\,\longmapsto\, 
[(\{ \mathrm{Tr} (u_{\alpha\beta}u'_{\beta\gamma}) \},\,
- \{\mathrm{Tr}(u_{\alpha\beta}v'_{\beta})-\mathrm{Tr}(v_{\alpha}u'_{\alpha\beta}) \})]
\end{aligned}
\end{equation}
in terms of the \v{C}ech cohomology constructed using an affine open covering $\{ U_{\alpha}\}$
(see \cite[Proposition 7.2]{Inaba-1}).
The symplectic form $\Theta_{\text{par}}$ is $d$-closed \cite[Proposition 7.3]{Inaba-1}.
Since the images of $\Theta^e$ and $\Theta_{\text{par}}$ in $\mathbb{H}^2(\mathbb{K})$
have the same description in terms of \v{C}ech cohomology
(see \eqref{ted} and \eqref{2020.1.26.14.17}), it follows that
\begin{equation*}
(q_1^{\text{par}})^*\Theta_{\text{par}}
|_{\mathcal{M}_{\text{PC}} (\boldsymbol{\nu})^{\circ}}
\, =\, \iota^* \Theta^e|_{\mathcal{M}_{\text{FC}}^{e} (d)^{\circ} }
\, .
\end{equation*}
Since $\Theta^{\mathrm{par}}$ is $d$-closed, so is
$\iota^* \Theta^e|_{\mathcal{M}_{\text{FC}}^{e} (d)^{\circ} }$.
Set $s\,:=\,\iota\circ s_1\,\colon\, U\,\longrightarrow\, p_1^{-1}(U)$,
which is a local section of $U_1$. Then the pullback $s^*(\Theta^e)\,=\,
s_1^*\iota^*\Theta^e|_{\mathcal{M}_{\text{FC}}^{e} (d)^{\circ} }$
is $d$-closed and so is $\Theta^e|_{\mathcal{M}_{\text{FC}}^{e} (d)^{\circ} }$
by the first remark in this proof.
\end{proof}

\subsection{Symplectic structure on $\mathcal{M}^{e}_{\mathrm{FC}} (d)$}\label{subsection: main d-closedness}

In Section \ref{2023.1.10.12.15} a 2-form $\Theta^e$ on $\mathcal{M}^{e}_{\mathrm{FC}} (d)$ was constructed. 
In the previous section, we considered the restriction of $\Theta^e$ on $\mathcal{M}^{e}_{\mathrm{FC}} 
(d)^{\circ} \,\subset\, \mathcal{M}^{e}_{\mathrm{FC}} (d)$. It was shown that this restriction is a 
symplectic form. Note that in the proof of the $d$-closedness of this restriction, we used irreducibility of 
$\mathcal{M}^{e}_{\mathrm{FC}} (d)^{\circ}$ (Proposition \ref{2023_4_2_13_21}) implicitly. In this section, 
we shall show that the 2-form $\Theta^e$ on $\mathcal{M}^{e}_{\mathrm{FC}} (d)$ is a symplectic form. In the 
proof of the $d$-closedness of $\Theta^e$, we will use the $d$-closedness of 
$\Theta^e|_{\mathcal{M}^{e}_{\mathrm{FC}} (d)^{\circ}}$ on $\mathcal{M}^{e}_{\mathrm{FC}} (d)^{\circ}$ for 
another effective divisor $\widetilde{D}$, instead of any argument on irreducibility of 
$\mathcal{M}^{e}_{\mathrm{FC}} (d)$.

\begin{Prop}\label{2023.1.10.23.18}
The $2$-form $\Theta^e$ on $\mathcal{M}^{e}_{\mathrm{FC}} (d)$ is nondegenerate.
\end{Prop}

\begin{proof}
Recall that the 2-form $\Theta^e$ is defined in \eqref{ted}. 
Let 
$\xi_{\Theta^H} \,\colon\, 
\mathbb{H}^1({\mathcal C}_{\bullet})\,\longrightarrow\, \mathbb{H}^1({\mathcal C}_{\bullet})^*$
be the homomorphism induced by $\Theta^H$.
Set $\mathcal{C}_0 \,:= \,{\mathcal E}nd(E)(-D)$ and 
$\mathcal{C}_1 \,:= \,{\mathcal E}nd(E)\otimes K_X(D)$.
For the above defined map $\xi_{\Theta^H}$, we have
the following commutative diagram whose rows are exact:
\[
\xymatrix{
H^0(\mathcal{C}_0) \ar[r] \ar[d]^-{b_1}
& H^0(\mathcal{C}_1)\ar[r]\ar[d]^-{b_2}
& \mathbb{H}^1({\mathcal C}_{\bullet}) \ar[d]^-{\xi_{\Theta^H}}\ar[r] 
& H^1(\mathcal{C}_0)\ar[r]\ar[d]^-{b_3}
& H^1(\mathcal{C}_1) \ar[d]^-{b_4} \\
H^1(\mathcal{C}_1)^* \ar[r]
& H^1(\mathcal{C}_0)^* \ar[r]&
\mathbb{H}^1({\mathcal C}_{\bullet})^* \ar[r]
& H^0(\mathcal{C}_1)^*\ar[r]
& H^0(\mathcal{C}_0)^* ,
}
\]
where $b_1,\,b_2,\,b_3,\,b_4$ are Serre duality isomorphisms.
So from the five lemma it follows that $\xi_{\Theta^H}$ is an isomorphism.
In other words, the $2$-form $\Theta^H$ is nondegenerate.
\end{proof}

Next we shall investigate the $d$-closedness of the $2$-form $\Theta^e$ on 
$\mathcal{M}^{e}_{\mathrm{FC}} (d)$.

\begin{Lem}\label{2023.1.9.18.44}
Let $(E_0,\, \phi_0 , \, \nabla_0)$ be a point on $\mathcal{M}^{e}_{\mathrm{FC}} (d)$.
For this point on $\mathcal{M}^{e}_{\mathrm{FC}} (d)$,
there exist a reduced effective divisor $\widetilde{D}$
and an isomorphism $\widetilde{\phi}_0 \,\colon\, \mathcal{O}_{\widetilde{D}}^{\oplus r}
\,\longrightarrow\, E_0|_{\widetilde{D}}$ such that 
$\widetilde{D} \,\supset\, D$,
$\widetilde{\phi}_0 |_D \,=\, \phi_0$, and 
$(E_0,\,\widetilde{\phi}_0)$ is simple.
\end{Lem}

\begin{proof}
Take a reduced effective divisor $\widetilde{D}$ 
such that $\widetilde{D} \,\supset\, D$ and 
$H^0(X, \,\mathcal{E}nd (E_0) (-\widetilde{D}) ) \,=\,0$. 
Moreover, take an isomorphism 
$\widetilde{\phi}_0 \,\colon\, \mathcal{O}_{\widetilde{D}}^{\oplus r}
\,\longrightarrow\, E_0|_{\widetilde{D}}$ such that $\widetilde{\phi}_0 |_D \,=\, \phi_0$.
We will show that $(E_0,\,\widetilde{\phi}_0)$ is simple.
For that, let $\textbf{g}$ be an automorphism of $(E_0,\,\widetilde{\phi}_0)$,
that is, $\textbf{g}$ is an automorphism of $E_0$ such that the following diagram
$$
\xymatrix{
\mathcal{O}_{\widetilde{D}}^{\oplus r} \ar[r]^-{\widetilde{\phi}_0}
\ar[rd]_-{\widetilde{\phi}_0}
 & E_0|_{\widetilde{D}}\ar[d]^-{\textbf{g}|_{\widetilde{D}}}\\
 & E_0|_{\widetilde{D}}
}
$$
is commutative. So the restriction $\textbf{g}|_{\widetilde{D}}$ is the identity map. 
Therefore, we have $$\textbf{g}- \mathrm{Id}_{E_0} \,\in\, H^0(X,\, \mathcal{E}nd (E_0) (-\widetilde{D}) ),$$
Since $H^0(X, \,\mathcal{E}nd (E_0) (-\widetilde{D}) )\, =\,0$,
it follows that $\textbf{g}\,=\, \mathrm{Id}_{E_0}$.
In other words, $(E_0,\,\widetilde{\phi}_0)$ is simple.
\end{proof}

Take an open covering
\begin{equation}\label{oc}
\mathcal{M}^{e}_{\rm{FC}}(d) \,=\, \bigcup_{m_0} \Sigma_{m_0}^{d},
\end{equation}
where each $\Sigma_{m_0}^{d}$ is the open substack of $\mathcal{M}^{e}_{\mathrm{FC}} (d)$
defined in \eqref{sd}. Recall that a very ample line bundle $\mathcal{O}_{X}(1)$ on the curve $X$ is fixed; set
$\theta_d(m)\,:=\,r d_X m +d +r(1-g)$, where $d_X\,:=\, \deg \mathcal{O}_X(1)$ and $g$ is the genus of $X$.
The above open substack $\Sigma_{m_0}^{d}$ is the fibered category whose objects are
simple framed $\mathrm{GL}(r,\mathbb{C})$--connections $(E,\,\phi,\,\nabla)$ on $X\times S$
such that 
\begin{itemize} 
\item $H^1(X,\, E_s (m_0-1))\,=\,0$ for each $s\,\in\, S$, and

\item $\chi(E_s(m))\,= \,\theta_d(m)$ for each $s \,\in\, S$ and all $m \,\in\, \mathbb{Z}$.
\end{itemize}
By the argument in the proof of Proposition \ref{2020.2.24.14.52},
the substack $\Sigma_{m_0}^{d}$ is of finite type.

\begin{Lem}\label{2023.1.10.12.40}
There exists a reduced effective divisor $\widetilde{D} \,\supset\, D$
such that for any points $(E,\, \phi , \, \nabla) \, \in \,\Sigma_{m_0}^{d}$,
there is an isomorphism $\widetilde{\phi}\, \colon\, \mathcal{O}_{\widetilde{D}}^{\oplus r}
\,\longrightarrow \,E|_{\widetilde{D}}$ satisfying the conditions that 
$\widetilde{\phi}|_D \,=\, \phi$ and 
$(E,\,\widetilde{\phi})$ is simple.
\end{Lem}

\begin{proof}
Take a point $s_0 \,=\, (E,\,\phi , \, \nabla) \, \in\, \Sigma_{m_0}^{d}$.
By Lemma \ref{2023.1.9.18.44}, there exists a reduced effective divisor 
$\widetilde{D}_{s_0}$
together with an isomorphism $\widetilde \phi \,\colon\, \mathcal{O}_{\widetilde{D}_{s_0}}^{\oplus r}
\,\longrightarrow\, E|_{\widetilde{D}_{s_0}}$ satisfying the following three conditions: 
$\widetilde{D}_{s_0} \,\supset\, D$,
$\widetilde \phi |_D\, =\, \phi$ and 
$(E,\,\widetilde \phi)$ is simple.

Take an open substack $U_s \,\subset\, \Sigma_{m_0}^{d}$, where ${s_0} \,\in\, U_{s_0}$ 
and $U_{s_0}$ is small enough,
and take a universal family $(\widetilde{E}, \, \psi,\, \widetilde \nabla)$ over $X \times U_{s_0}$.
Since $\widetilde E$ is locally trivial, 
we may take a lift $\widetilde{\psi} \,\colon\, \mathcal{O}_{\widetilde{D}_{s_0} \times U_{s_0}}^{\oplus r}
\,\longrightarrow \,\widetilde{E}|_{\widetilde{D}_{s_0}\times U_{s_0}}$
such that $\widetilde \psi |_{D\times U_{s_0}} = \psi$.
Note that $(\widetilde{E}, \, \widetilde \psi)|_{X \times {s_0} } \cong (E, \, \widetilde \phi)$,
which is simple. Since the requirement that $H^0(X \times s ,\, \mathcal{E}nd (\widetilde E|_{X \times s}) (-\widetilde D_{s_0}) )\,=\,0$
is an open condition, we may assume that $(\widetilde{E}, \, \widetilde \psi)$ is a family of simple framed bundles.
Consider an open covering $ \Sigma_{m_0}^{d} \,=\, \bigcup_{s_0} U_{s_0}$.
Since $\Sigma_{m_0}^{d}$ is of finite type, 
we may cover $ \Sigma_{m_0}^{d}$ by a finite number of the open substacks $\{U_{s_0}\}_{s_0} $:
$$\Sigma_{m_0}^{d} \,=\, \bigcup_{i=1}^{m} U_{s_{i}},$$
where $s_1,\,\cdots,\,s_m$ are points on $ \Sigma_{m_0}^{d}$. Now take
$$
\widetilde{D} \,:=\, \bigcup_{i=1}^{m} \widetilde D_{s_i}.
$$
Then, by the construction of $\widetilde{D}$,
for any points $(E,\,\phi , \, \nabla)\,\in \, \Sigma_{m_0}^{d}$,
there exists an isomorphism $\widetilde{\phi} \,\colon\, \mathcal{O}_{\widetilde{D}}^{\oplus r}
\,\longrightarrow\, E|_{\widetilde{D}}$ such that 
$\widetilde \phi |_D \,=\, \phi$ and $(E,\,\widetilde{\phi})$ is simple.
\end{proof}

\begin{Thm}\label{Theorem: d-closedness on the moduli of framed connections}
The nondegenerate $2$-form $\Theta^e$ on 
$\mathcal{M}^{e}_{\mathrm{FC}} (d)$ defined by
\eqref{ted} is $d$-closed.
\end{Thm}

\begin{proof}
Consider the open covering 
$\mathcal{M}^{e}_{\mathrm{FC}} (d)\,=\, \bigcup_{m_0} \Sigma^d_{m_0}$ in \eqref{oc}.
It is enough to prove that the restriction
$\Theta^e|_{\Sigma^d_{m_0}}$ is $d$-closed for each $m_0$.
Take a reduced effective divisor $\widetilde D$ as in Lemma \ref{2023.1.10.12.40}.
Let $\mathcal{M}^{e}_{\mathrm{FC}} (d,\,\widetilde{D})$ 
be the Deligne-Mumford stack constructed in Proposition \ref{2020.2.24.14.52} for $\widetilde D$.
Let $\mathcal{M}^{e}_{\mathrm{FC}} (d,\widetilde{D})^{\circ}$ 
be the Deligne-Mumford stack whose objects are objects of 
$\mathcal{M}^{e}_{\mathrm{FC}} (d,\widetilde{D})$ such that 
the underlying framed bundles are simple.
In other words, we have
\begin{equation*}
\mathcal{M}^{e}_{\mathrm{FC}} (d,\widetilde{D})^{\circ}\,=\, 
\left\{\left( \widetilde E,\, \widetilde{\phi} ,\, \widetilde{\nabla}\right)\ \middle| \
\begin{array}{l}
\text{$\widetilde E$ is a vector bundle of degree $d$,}\\
\text{$\widetilde \phi \,\colon\, \mathcal{O}_{\widetilde{D}}^{\oplus r}
\,\longrightarrow\, E|_{\widetilde{D}}$ is an isomorphism,}\\
\text{$\widetilde\nabla \,\colon\, \widetilde{E} \,\longrightarrow\, \widetilde{E}
\otimes K_{X}(\widetilde{D})$ is a connection, and} \\
\text{$(\widetilde E,\, \widetilde{\phi})$ is simple} 
\end{array}
\right\} {\Big/}\Large{\sim}_{e}.
\end{equation*}
Taking the degree of $\widetilde{D}$ to be sufficiently large,
the canonical $2$-form 
$\Theta^e|_{\mathcal{M}^{e}_{\mathrm{FC}} (d,\widetilde{D})^{\circ}}$
on $\mathcal{M}^{e}_{\mathrm{FC}} (d,\widetilde{D})^{\circ}$
is $d$-closed by Proposition \ref{2019.12.22.18.19} and Proposition \ref{2019.12.22.18.19_a}.
Define a moduli space $\mathcal{M}^{e}_{\mathrm{FC}} (d,\widetilde{D},D)$ as follows:
\begin{equation*}
{\mathcal{M}^{e}_{\mathrm{FC}} (d,\widetilde{D},D)\,=\, 
\left\{\left( \widetilde E,\, \widetilde{\phi} ,\, \widetilde{\nabla}\right) \,\in \,\mathcal{M}^{e}_{\mathrm{FC}} (d,\widetilde{D})^{\circ}
 \ \middle| \
\begin{array}{l}
\text{$\widetilde\nabla$ is regular on $\widetilde{D} \setminus D$, and} \\
\text{$(\widetilde E,\, \widetilde{\phi} |_D ,\, \widetilde \nabla )$ is simple} 
\end{array}
\right\}{\Big/}\Large{\sim}_{e}}.
\end{equation*}
Let $\iota\,:\, \mathcal{M}^{e}_{\mathrm{FC}} (d,\widetilde{D},D)\, \longrightarrow\,
\mathcal{M}^{e}_{\mathrm{FC}} (d,\widetilde{D})^{\circ}$ be the natural inclusion map 
and $\pi$ the natural map 
from $\mathcal{M}^{e}_{\mathrm{FC}} (d,\widetilde{D},D)$ to 
$\mathcal{M}^{e}_{\mathrm{FC}} (d)$ 
induced by the restriction of framings to $D$:
$$
\begin{aligned}
\pi \,\, \colon\,\, \mathcal{M}^{e}_{\mathrm{FC}} (d,\widetilde{D},D) 
&\,\longrightarrow\,\, \mathcal{M}^{e}_{\mathrm{FC}} (d) \\
(\widetilde E,\, \widetilde{\phi},\, \widetilde{\nabla}) &\,\longmapsto\,\, 
(\widetilde E,\, \widetilde{\phi}|_D,\, \widetilde{\nabla}).
\end{aligned}
$$ 
This map $\pi$ is smooth.
By Lemma \ref{2023.1.10.12.40}, 
the open substack $\Sigma_{m_0}^{d}$ is 
contained in the image of $\pi$.
We consider the following maps
\begin{equation*}
\xymatrix{
& \mathcal{M}^{e}_{\mathrm{FC}} (d,\widetilde{D},D)
\ar[r]^-{\iota} \ar[d]^-{\pi}
& \mathcal{M}^{e}_{\mathrm{FC}} (d,\widetilde{D})^{\circ} \\
\Sigma_{m_0}^{d} \ar[r]^-{\subset} & \mathcal{M}^{e}_{\mathrm{FC}} (d) \rlap{.} &
}
\end{equation*}
Let $\Theta^e_{\widetilde{D}}$ be the 2-form on $\mathcal{M}^{e}_{\mathrm{FC}} (d,\widetilde{D})^{\circ}$ 
defined in \eqref{ted}.
By the definition of $\Theta^e$ and $\Theta^e_{\widetilde{D}}$, 
which are described by the same formula via the \v{C}ech cohomology,
we have 
$$
\pi^* \Theta^e\ =\ \iota^* \Theta^e_{\widetilde{D}}.
$$
As $\Theta^e_{\widetilde{D}}$ is $d$-closed by Proposition \ref{2019.12.22.18.19} 
and Proposition \ref{2019.12.22.18.19_a},
we conclude that $\pi^* \Theta^e$ is $d$-closed.
Since $\pi$ is smooth, and the image of $\pi$ contains the open substack $\Sigma_{m_0}^{d}$,
it follows that $\Theta^e|_{\Sigma^d_{m_0}}$ is $d$-closed.
\end{proof}

\subsection{Symplectic structure on $\mathcal{M}^{H}_{\mathrm{FC}} (d)$}\label{subsection: d-closedness for H-framing}

Fix a complex Lie proper subgroup
$H_x\,\subsetneq\, \mathrm{GL}(r,\mathbb{C})$ for each $x\,\in\, D$.

Consider the complexes $\mathcal{D}_{\bullet}$ and $\mathbb{K}$ constructed in \eqref{db} and \eqref{bK}
respectively. Note that the pairing $\mathrm{ad}(E_G)\otimes \mathrm{ad} (E_G)\, \longrightarrow\,
{\mathcal O}_X$ in \eqref{ews} produces a pairing
$$
\mathrm{ad}_{\phi} (E_G)\otimes (\mathrm{ad}^n_{\phi} (E_G) \otimes K_X (D))
\, \longrightarrow\, K_X\, .
$$
The restriction of the pairing $\widehat{\sigma}$ (see \eqref{ews})
$$
\mathrm{ad}_{\phi} (E_G)\otimes \mathrm{ad}_{\phi} (E_G)\, \longrightarrow\,
{\mathcal O}_X\, ,
$$
and the homomorphism
$$
(\mathrm{ad}_{\phi} (E_G)\otimes (\mathrm{ad}^n_{\phi} (E_G) \otimes K_X (D)))\oplus
((\mathrm{ad}^n_{\phi} (E_G) \otimes K_X (D))\otimes \mathrm{ad}_{\phi} (E_G))\,\longrightarrow\, K_X
$$
constructed using $\widehat\sigma$, together produce a homomorphism
$$
\mathcal{D}_{\bullet}\otimes \mathcal{D}_{\bullet}\, \longrightarrow\, \mathbb{K}
$$
of complexes. Let
$$
\mathbb{H}^2 (\mathcal{D}_{\bullet} \otimes \mathcal{D}_{\bullet}) 
\,\longrightarrow\, \mathbb{H}^2(\mathbb{K})
$$
be the homomorphism of hypercohomologies induced by this homomorphism of complexes. Now the composition of
the natural homomorphism
$$
\mathbb{H}^1 (\mathcal{D}_{\bullet}) \otimes \mathbb{H}^1 (\mathcal{D}_{\bullet})
\,\longrightarrow\, \mathbb{H}^2 (\mathcal{D}_{\bullet} \otimes \mathcal{D}_{\bullet})
$$
with the above homomorphism of hypercohomologies produces a pairing
\begin{equation}\label{DP}
\Theta^H \,\colon\,
\mathbb{H}^1 (\mathcal{D}_{\bullet}) \otimes \mathbb{H}^1 (\mathcal{D}_{\bullet}) 
\,\longrightarrow\, \mathbb{H}^2(\mathbb{K}) \,=\, \mathbb{C}\, .
\end{equation}
In terms of the \v{C}ech cohomology with respect to an affine open covering $\{ U_{\alpha}\}$, the pairing
$\Theta^H$ in \eqref{DP} is of the form
$$
[(\{u_{\alpha\beta} \} ,\, \{ v_{\alpha} \})]\otimes[(\{u'_{\alpha\beta} \} ,\, \{ v'_{\alpha} \})]
\, \longmapsto\, 
[(\{ \mathrm{Tr} (u_{\alpha\beta}u'_{\beta\gamma}) \},\,
- \{\mathrm{Tr}(u_{\alpha\beta}v'_{\beta})-\mathrm{Tr}(v_{\alpha}u'_{\alpha\beta}) \})]\, .
$$

This pairing in \eqref{DP} gives a 2-form on $\mathcal{M}^{H}_{\mathrm{FC}} (d)$.
We also denote by $\Theta^H $ this 2-form on $\mathcal{M}^{H}_{\mathrm{FC}} (d)$.
Then $\Theta^H $ is nondegenerate by the argument as after
\cite[Theorem 5]{BIKS} by applying \cite[Proposition 4.1]{BLP}.
Now it will be shown that $\Theta^H $ is $d$-closed.

\begin{Def}
Let $\mathcal{M}^{e}_{\rm{FC}}(d)_{\mathfrak{h}^{\perp}}$ be the stack
over the category of locally Noetherian schemes whose 
objects are quadruples $(S,\, E, \,\phi=\{\phi_{x\times S} \}_{x\in D} ,\,\nabla)$
that satisfy (1), (3) and (5) in Definition \ref{2019.10.28.16.49} 
and the following $(2)''$ and $(4)''$:
\begin{itemize}
\item[$(2)''$] $\phi_{x\times S}$ is a section of the structure map
\begin{equation*}
\mathrm{Isom}_{S} 
(\mathcal{O}_{x\times S}^{\oplus r}, \, E\big\vert_{x\times S})
\,\longrightarrow \, x\times S\, .
\end{equation*} 
Denote by $$\varphi_{x \times S}\,\,\colon\,\, 
\mathcal{O}_{x\times S}^{\oplus r}\,\xrightarrow{\ \sim\ }\, 
E\big\vert_{x\times S}$$
the isomorphism given by 
the map $x\times S\,\longrightarrow\, \mathrm{Isom}_{S} 
(\mathcal{O}_{x\times S}^{\oplus r},
\, E\big\vert_{x\times S})$.

\item[$(4)''$]
Let $\mathrm{res}_{x\times S}(\nabla) \,\in\, \mathrm{End}(E)\big\vert_{x\times S}$
be the residue matrix of the connection $\nabla$ along $x\times S$.
Then $\phi^{-1}_{x \times S} \circ \mathrm{res}_{x\times S}(\nabla)
\circ \phi_{x \times S} \,\in\, \mathfrak{h}^{\perp} \otimes \mathcal{O}_{S}$.
\end{itemize}
A {\it morphism}
\begin{equation*}
(S,\, E,\, \phi,\,\nabla )\,\longrightarrow \, (S',\, E',\, \phi',\,\nabla')
\end{equation*}
in $\mathcal{M}^{e}_{\rm{FC}}(d)_{\mathfrak{h}^{\perp}}$ 
is a Cartesian square 
\begin{equation*}
\xymatrix{
 E \ar[r]^-{\sigma} \ar[d] & E' \ar[d] \\
S \ar[r]^{\widetilde{\sigma}} & S'
}
\end{equation*}
such that the diagram
\begin{equation*}
\xymatrix{
 E \ar[r]^-{\nabla} \ar[d]_{\cong}^{\sigma} 
 & E \otimes K_{X}(D) \ar[d]_{\cong}^-{\sigma} \\
 E' \ar[r]^-{\nabla'} 
 & E' \otimes K_{X}(D)
}
\end{equation*}
commutes and the composition 
$(\phi'_{x\times S})^{-1} \circ \sigma|_{x\times S}\circ \phi_{x\times S}$
coincides with the identity map of $\mathcal{O}_{x\times S}^{\oplus r}$ for each $x\,\in\, D$.
\end{Def}

\begin{Thm}\label{Thm: d-closedness on H-framed connection moduli}
The nondegenerate $2$-form $\Theta^H$
on $\mathcal{M}^{H}_{\mathrm{FC}} (d)$ 
defined by \eqref{DP} is $d$-closed.
\end{Thm}

\begin{proof}
Consider the diagram
\begin{equation*}
\xymatrix{
\mathcal{M}^{e}_{\mathrm{FC}}(d)_{\mathfrak{h}^{\perp}}
 \ar[r]^-{\pi_1} \ar[d]_-{\pi_2} &\mathcal{M}^{e}_{\mathrm{FC}} (d) \\
\mathcal{M}_{\mathrm{FC}}^{H}(d)&
}
\end{equation*}
where $\pi_1$ and $\pi_2$ are the natural maps. It is straightforward to check that 
\begin{equation*}
\pi_1^* \Theta^e\,\,=\,\,\pi_2^*\Theta^H.
\end{equation*}
Since $\Theta^e$ is $d$-closed, the form $\pi_2^*\Theta^H$ is also $d$-closed.
This implies that $\Theta^H$ is $d$-closed, because the map
$\pi_2$ is dominant. 
\end{proof}

\subsection{Poisson structure}\label{Subsect:PoissonStr}

In this subsection, we will see the details of the Poisson structure 
mentioned in the introduction. This is influenced by a construction done in
\cite{BBG}.

Let $\mathcal{M}_{\rm{C}}(d)$ be the moduli space of pairs $(E,\, \nabla)$,
where $E$ is a holomorphic vector bundle on $X$ of rank $r$ and degree $d$, and $\nabla$
is a logarithmic connection on $E$ whose singular part is contained in $D$,
such that $(E,\,\nabla)$ is simple in the sense that
the endomorphisms of $E$ preserving $\nabla$ are just the constant scalar multiplications. 
In \cite{Nit}. Nitsure constructed the moduli space
${\mathcal M}^{\mathrm{ss}}_{\mathrm{C}}(d)$
of semistable logarithmic connections,
which contains the moduli space of
stable logarithmic connections
${\mathcal M}^{\mathrm{s}}_{\mathrm{C}}(d)$
as a Zariski open subset. By its definition, our moduli space
$\mathcal{M}_{\rm{C}}(d)$ contains
${\mathcal M}^{\mathrm{s}}_{\mathrm{C}}(d)$
as a Zariski open subspace. Recall that a description of the tangent space of this moduli space
is given in \cite{Nit}.
For 
$(E,\, \nabla)\, \in\,\mathcal{M}_{\rm{C}}(d)$,
the tangent space of $\mathcal{M}_{\rm{C}}(d)$ at
$(E,\nabla)$ is
$$
T_{(E,\nabla)}\mathcal{M}_{\rm{C}}(d)
\, \,=\, \,
{\mathbb H}^1\big( {\mathcal End}(E) \to {\mathcal C}_1 \big),
$$
where ${\mathcal C}_0\,=\,{\mathcal End}(E)(-D)$,
${\mathcal C}_1\,=\,{\mathcal End}(E)\otimes K_X(D)$
and the map ${\mathcal End}(E) \,\longrightarrow\, {\mathcal C}_1$ is defined by
$u \,\longmapsto \,\nabla \circ u- u \circ \nabla$.
The cotangent space is 
\[
 T_{(E,\nabla)}^*\mathcal{M}_{\rm{C}}(d)
\,\, =\,\,
 {\mathbb H}^1\big( {\mathcal End}(E) \to {\mathcal C}_1 \big)^*
\,\, \cong\,\,
 {\mathbb H}^1\big( {\mathcal C}_0 \to {\mathcal End}(E)\otimes K_X \big),
\]
over which there is a canonical pairing
\begin{align}\label{Poisson bracket on connection moduli}
\begin{split}
 &T_{(E,\nabla)}^*\mathcal{M}_{\rm{C}}(d)
 \otimes
 T_{(E,\nabla)}^*\mathcal{M}_{\rm{C}}(d)
 \\
 &
 =\,
 {\mathbb H}^1\big( {\mathcal C}_0 \to {\mathcal End}(E)\otimes K_X \big)
 \otimes
 {\mathbb H}^1\big( {\mathcal C}_0 \to {\mathcal End}(E)\otimes K_X \big)
 \,\longrightarrow\,
 \mathbb{H}^2(\Omega_X^{\bullet})
 \,\cong \,\mathbb{C}.
\end{split}
\end{align}
Consider the open subspace
\begin{equation*}
 {\mathcal M}^e_{\rm FC}(d)'
 \,\,=\,\,
 \left\{(E,\,\nabla,\,\phi) \, \in \, {\mathcal M}^e_{\rm FC}(d) \ \middle| \
 \text{$(E,\,\nabla)$\, is simple} \right\} 
\end{equation*}
of the moduli space ${\mathcal M}^e_{\rm FC}(d)$ of simple framed connections.
Then there is a natural forgetful map
\begin{equation}\label{Poisson map from framed moduli to simple moduli}
\pi\,\, \colon\,\, {\mathcal M}^e_{\rm FC}(d)'
\, \longrightarrow \, \mathcal{M}_{\rm{C}}(d),
\end{equation}
and the induced map $\pi^*$ on the cotangent spaces makes the diagram
\[
\begin{CD}
T_{(E,\nabla)}^*\mathcal{M}_{\rm{C}}(d)
 \times
 T_{(E,\nabla)}^*\mathcal{M}_{\rm{C}}(d)
 @>>> \mathbb{H}^2(\Omega_X^{\bullet})
 \cong \mathbb{C}
 \\
 @V\pi^*\times\pi^* VV @VVV \\
 T_{(E,\nabla,\phi)}^*\mathcal{M}^e_{\rm{FC}}(d)
 \times
 T_{(E,\nabla,\phi)}^*\mathcal{M}^e_{\rm{FC}}(d)
 @>>> \mathbb{H}^2(\Omega_X^{\bullet})
 \cong \mathbb{C}
 \end{CD}
\]
commutative.
The bottom horizontal arrow satisfies the Jacobi identity,
because it corresponds to the symplectic form on the moduli space
${\mathcal M}^e_{\rm FC}(d)$ given in 
Theorem \ref{Theorem: d-closedness on the moduli of framed connections}.
So the pairing in \eqref{Poisson bracket on connection moduli}
is also skew-symmetric and satisfies the Jacobi identity.
Thus the following corollary is obtained.

\begin{Cor}\label{Cor: Poisson structure on logarithmic moduli}
The moduli space $\mathcal{M}_{\rm{C}}(d)$ has a Poisson structure
defined by the Poisson bracket in \eqref{Poisson bracket on connection moduli}.
Furthermore, the morphism 
$\pi$ in \eqref{Poisson map from framed moduli to simple moduli} becomes a Poisson map.
\end{Cor}

We will see a slightly different view of the Poisson structure on the moduli space
$\mathcal{M}_{\rm{C}}(d)$.
Set
\[
 A\,:=\, \left\{ \boldsymbol{a}\,=\,(a^{(i)}_j)^{1\leq i\leq n}_{0\leq j\leq r-1}\,\, \middle|\,\,
 \sum_i a^{(i)}_{r-1}\,=\,d \right\}.
\]
By associating the coefficients of the 
characteristic polynomial of $\res_{x_i}(\nabla)$
at each point $x_i\in D$, we can define a morphism
\begin{equation}\label{map associating coefficients of characteristic polynomial of residue}
 {\mathcal M}_{\rm C}(d) \, \longrightarrow \, A
\end{equation}
whose fiber ${\mathcal M}_{\rm C}(\boldsymbol{a})$ over $\boldsymbol{a}\,\in\, A$
is smooth for generic $\boldsymbol{a}$ but
it has singularities for special $\boldsymbol{a}$.
Consider the moduli space of simple parabolic connections
\[
 {\mathcal M}_{\rm PC}(d)\,=\,\left.
 \left\{ (E,\,\nabla,\,\boldsymbol{l}) \ \middle| \
 \begin{array}{l}
 \text{$(E,\,\boldsymbol{l}=(l^{(i)}_j))$ is a quasi-parabolic bundle of rank $r$ and degree $d$,}
 \\
 \text{$\nabla\,\colon\, E\,\longrightarrow \,E\otimes K_X(D)$
 is a connection satisfying}
 \\
 \text{$\res_{x_i}(\nabla)(l^{(i)}_j)\,\subset\, l^{(i)}_j$ for any $i,\,j$
 and $(E,\,\nabla,\,\boldsymbol{l})$ is simple.}
 \end{array}
 \right\} \right/ \sim.
\]
For the open subspace
\[
{\mathcal M}_{\rm PC}(d)'\,\,=\,\,
\left\{ (E,\,\nabla,\,\boldsymbol{l}) \, \in \, {\mathcal M}_{\rm PC}(d) \ \middle| \
\text{$(E,\,\nabla)$\, is simple} \right\}
\]
of ${\mathcal M}_{\rm PC}(d)$,
there is a canonical morphism
\begin{equation}\label{morphism from parabolic connection to logarithmic connection}
 {\mathcal M}_{\rm PC}(d)'\,\longrightarrow\,
 \mathcal{M}_{\rm{C}}(d)
\end{equation}
which is generically finite.
Set 
$\Lambda\,:=\,\left\{ (\nu^{(i)}_j)^{1\leq i\leq n}_{0\leq j\leq r-1} \, \in \, \mathbb{C}^{nr}
\ \middle| \ d+\sum_{i,j}\nu^{(i)}_j\,=\,0 \right\}$.
Then we have a smooth morphism
\begin{equation}\label{map associating exponents}
 {\mathcal M}_{\rm PC}(d) \, \longrightarrow \, 
 \Lambda
\end{equation}
whose fiber over any $\boldsymbol{\nu}\,\in\,\Lambda$ is the moduli space
${\mathcal M}_{\rm PC}(\boldsymbol{\nu})$ of 
$\boldsymbol{\nu}$-parabolic connections.
The morphism in \eqref{morphism from parabolic connection to logarithmic connection}
induces a map between the fibers of \eqref{map associating coefficients of characteristic polynomial of residue}
and \eqref{map associating exponents}
\[
{\mathcal M}_{\rm PC}(\boldsymbol{\nu})'\,:=\,
{\mathcal M}_{\rm PC}(\boldsymbol{\nu})\cap {\mathcal M}_{\rm PC}(d)'
\, \longrightarrow \, {\mathcal M}_{\rm C}(\boldsymbol{a})
\]
which is an isomorphism for generic $\boldsymbol{a}$
and it is a resolution of singularities of ${\mathcal M}_{\rm C}(\boldsymbol{a})$
for special $\boldsymbol{a}$, where $\boldsymbol{a}\,=\,(a^{(i)}_j)$
is determined by $\boldsymbol{\nu}\,=\,(\nu^{(i)}_j)$ as follows:
$$\prod_{j=0}^{r-1}(t-\nu^{(i)}_j)\,=\,t^r+a^{(i)}_{r-1}t^{r-1}+\cdots+a^{(i)}_1t+a^{(i)}_0.$$
Roughly speaking, the moduli space ${\mathcal M}_{\rm C}(\boldsymbol{a})$
for special $\boldsymbol{a}$ gives a partial resolution of singularities
of the corresponding character variety which we will define precisely later in
\eqref{character variety with fixed local monodromy}.
The meaning of the singularities of character varieties and their exceptional loci
in the moduli space ${\mathcal M}_{\rm PC}(\boldsymbol{\nu})$
(or precisely ${\mathcal M}_{\rm PC}^{\boldsymbol{\alpha}}(\boldsymbol{\nu})$)
is explained in \cite{Iw2} and \cite{IIS1} from the viewpoint of the isomonodromic deformation,
and their classification in the case of Painlev\'e equations is given in \cite{ST}.

Setting
\begin{align*}
 \widetilde{\mathcal D}^{\rm par}_0
 &=
 \left\{ u \, \in \, {\mathcal End}(E) \ \middle| \
 \text{$u|_{x_i}(l^{(i)}_j) \subset l^{(i)}_{j+1}$ for any\, $i,\,j$}
 \right\}
 \ \subset \ {\mathcal D}^{\rm par}_0
 \\
 \widetilde{\mathcal D}^{\rm par}_1
 &=
 \left\{ v \, \in \, {\mathcal End}(E)\otimes K_X(D) \ \middle| \
 \text{$\res_{x_i}(v)(l^{(i)}_j) \subset l^{(i)}_j$ for any\, $i,\,j$}
 \right\},
 \\
\end{align*}
we can define a complex ${\mathcal D}^{\rm par}_0\, \longrightarrow\, 
\widetilde{\mathcal D}^{\rm par}_1$,\,
$u \,\longmapsto\, \nabla \circ u-u \circ \nabla$, which induces complexes
$\widetilde{\mathcal D}^{\rm par}_0 \,\longrightarrow\, {\mathcal D}^{\rm par}_1$
and $\widetilde{\mathcal D}^{\rm par}_0\,\longrightarrow\, \widetilde{\mathcal D}^{\rm par}_1$.
The tangent space of the moduli space ${\mathcal M}_{\rm PC}(d)$ is
$T{{\mathcal M}_{\rm PC}(d)}\,=\,
\mathbb{H}^1({\mathcal D}^{\rm par}_0 \to \widetilde{\mathcal D}^{\rm par}_1)$
and the cotangent space is its dual
\[
T^*{{\mathcal M}_{\rm PC}(d)}\,=\,
\mathbb{H}^1({\mathcal D}^{\rm par}_0 \to \widetilde{\mathcal D}^{\rm par}_1)^*
\,\cong\, \mathbb{H}^1(\widetilde{\mathcal D}^{\rm par}_0\to {\mathcal D}^{\rm par}_1).
\]
So we can define a canonical pairing
\begin{equation}\label{Poisson product on parabolic connections}
 (T^*{{\mathcal M}_{\rm PC}(d)})\otimes (T^*{{\mathcal M}_{\rm PC}(d)})
 \,=\,
 \mathbb{H}^1(\widetilde{\mathcal D}^{\rm par}_0\to {\mathcal D}^{\rm par}_1) \otimes
 \mathbb{H}^1(\widetilde{\mathcal D}^{\rm par}_0\to {\mathcal D}^{\rm par}_1)
 \, \longrightarrow \,
 \mathbb{H}^2(\Omega^{\bullet}_X)\,\cong\,\mathbb{C}.
\end{equation}
Let $B$ be the Borel subgroup of $\mathrm{GL}(r,\mathbb{C})$ consisting of
upper triangular matrices,
and let $U$ be the subgroup of $B$ consisting of matrices whose diagonal entries are $1$.
Consider the open subspace
\[
 {\mathcal M}^U_{\rm FC}(d)'
\, =\,
 \left\{ (E,\,\nabla,\,[\phi]) \, \in \, {\mathcal M}^U_{\rm FC}(d) \ \middle| \
 \begin{array}{l}
 \text{the parabolic connection $(E,\,\nabla,\,\boldsymbol{l})$} \\
 \text{induced from $(E,\,\nabla,\,[\phi])$ is simple}
 \end{array}
 \right\}
\]
of ${\mathcal M}^U_{\rm FC}(d)$, which is the moduli space of framed connections
in Definition \ref{2019.10.28.16.49} with $H\,=\,U$.
Associating the corresponding parabolic connection,
we can define a morphism
\begin{equation}\label{Poisson map from U-framed to parabolic}
{\mathcal M}^U_{\rm FC}(d)' \, \longrightarrow \, {\mathcal M}_{\rm PC}(d)
\end{equation}
which becomes a $\big(\prod_{D}B/U\big)\big/\mathbb{C}^*$-bundle. 
By construction, the diagram
\[
\begin{CD}
(T^*{{\mathcal M}_{\rm PC}(d)})\otimes (T^*{{\mathcal M}_{\rm PC}(d)})
\,=\,
 \mathbb{H}^1(\widetilde{\mathcal D}^{\rm par}_0\to {\mathcal D}^{\rm par}_1) \otimes
 \mathbb{H}^1(\widetilde{\mathcal D}^{\rm par}_0\to {\mathcal D}^{\rm par}_1)
 @>>> 
 \mathbb{H}^2(\Omega^{\bullet}_X)\,\cong\,\mathbb{C}
 \\
 @VVV @VVV
 \\
 (T^*{{\mathcal M}^U_{\rm FC}(d)'})\otimes (T^*{{\mathcal M}^U_{\rm FC}(d)'}) 
 \,=\,
 \mathbb{H}^1(\widetilde{\mathcal D}^{\rm par}_0\to \widetilde{\mathcal D}^{\rm par}_1) \otimes
 \mathbb{H}^1(\widetilde{\mathcal D}^{\rm par}_0\to \widetilde{\mathcal D}^{\rm par}_1) 
 @>>>
 \mathbb{H}^2(\Omega^{\bullet}_X)\,\cong\,\mathbb{C}
 \end{CD}
\]
is commutative.
The lower horizontal arrow is the Poisson bracket
corresponding to the symplectic form on the moduli space
${\mathcal M}^U_{\rm FC}(d)$ given by 
Theorem \ref{Thm: d-closedness on H-framed connection moduli}.
So the pairing in \eqref{Poisson product on parabolic connections}
defines a Poisson structure on the moduli space
${\mathcal M}_{\rm PC}(d)$
and the morphism in \eqref{Poisson map from U-framed to parabolic} is a Poisson map.

We can also see that the pairing in \eqref{Poisson product on parabolic connections}
commutes with the Poisson bracket on ${\mathcal M}_{\rm PC}(\boldsymbol{\nu})$
corresponding to the symplectic form.
So the canonical inclusion
${\mathcal M}_{\rm PC}(\boldsymbol{\nu})
\, \hookrightarrow \,{\mathcal M}_{\rm PC}(d)$ is also a Poisson map.

The canonical map
$\mathbb{H}^1({\mathcal C}_0\to{\mathcal End}(E)\otimes K_X)
\,\longrightarrow\, 
\mathbb{H}^1(\widetilde{\mathcal D}^{\rm par}_0 \,\rightarrow\, {\mathcal D}^{\rm par}_1)$
coincides with the map
\[
T^*{\mathcal{M}_{\rm{C}}(d)}
\,\longrightarrow\,
T^*{\mathcal M}_{\rm PC}(\boldsymbol{\nu})
\]
on the cotangent spaces induced by the morphism in
\eqref{morphism from parabolic connection to logarithmic connection},
which means that 
the Poisson brackets in \eqref{Poisson bracket on connection moduli} commutes with
that in \eqref{Poisson product on parabolic connections}.
Combining the above, the following corollary is obtained.

\begin{Cor}\label{Cor: Poisson structure on parabolic connection moduli}
The moduli space ${\mathcal M}_{\rm PC}(d)$ of parabolic connections
has a Poisson structure defined by the Poisson bracket given in
\eqref{Poisson product on parabolic connections}.
Furthermore, the morphism
${\mathcal M}_{\rm PC}(d)' \, \longrightarrow \,
\mathcal{M}_{\rm{C}}(d)$ in
\eqref{morphism from parabolic connection to logarithmic connection}
becomes a Poisson map for this Poisson structure.
\end{Cor}

\section{The moduli space of parabolic connections is not affine}\label{section: global algebraic functions implies not affine}

\subsection{Moduli space of parabolic connections and parabolic Higgs bundles}\label{subsection: moduli of parabolic connections and Higgs bundles}

Throughout this section, we assume that $k$ is an algebraically closed field of
arbitrary characteristic.

Let $$(X,\,\boldsymbol{x})\, := \,(X,\, (x_1,\, \cdots,\,x_n))$$
be an $n$-pointed smooth projective curve of genus $g$ 
over $k$, where $x_1, \,\cdots,\,x_n$ are distinct $k$-valued points of $X$.
Denote the reduced divisor $x_1+\cdots+x_n$ on $X$ by $D$.
Take a positive integer $r$ which is not divisible by the characteristic of $k$ and also
take an integer $d$ and an element 
$$\boldsymbol{\nu}\,=\,(\nu^{(i)}_j )^{1\le i \le n}_{0\le j \le r-1}
\,\in\, k^{nr}$$
such that the equality $\sum_{i,j} \nu^{(i)}_j \,=\, -d$
holds in $k$.
Take a collection of rational numbers
$$\boldsymbol{\alpha}\,=\,(\alpha^{(i)}_j)^{1\leq i\leq n}_{1\leq j\leq r}\, \in\,
{\mathbb Q}^{rn}$$
satisfying the conditions
\begin{itemize}
\item $0\,<\,\alpha^{(i)}_1\,<\,\cdots\,<\,\alpha^{(i)}_r\,<\,1$, and

\item $\alpha^{(i)}_j\,\neq\,\alpha^{(i')}_{j'}$ for $(i,\,j)\,\neq\, (i',\,j')$.
\end{itemize}

A $(\boldsymbol{x},\ \boldsymbol{\nu})$-parabolic connection
on $X$ is defined exactly in the same way as
Definition \ref{def: parabolic connection}.
Although a parabolic connection includes the data of a parabolic weight,
we omit it and simply write $(E,\,\nabla ,\, \boldsymbol{l})$.
The definition of 
$\boldsymbol{\alpha}$-stability of parabolic connection
is also defined in the same way as 
Definition \ref{def: stability of parabolic connection}.

In the proof of the existence of the moduli space of
stable parabolic connections in \cite [Theorem 2.2] {Inaba-1},
we used the embedding to the moduli space of
parabolic $\Lambda^1_D$-triples (\cite[Theorem 5.1] {IIS1}; this
argument also works over a field of arbitrary characteristic.
So we have the following theorem;

\begin{Thm}
There exists a coarse moduli scheme
$\mathcal{M}_{\mathrm{PC}}^{\boldsymbol{\alpha}}(\boldsymbol{\nu})$
of $\boldsymbol{\alpha}$-stable
$(\boldsymbol{x},\, \boldsymbol{\nu})$-parabolic connections
on a smooth projective curve $X$ over $k$.
Furthermore, $\mathcal{M}_{\mathrm{PC}}^{\boldsymbol{\alpha}}(\boldsymbol{\nu})$
is quasi-projective over $k$.
\end{Thm}

\begin{Def}[{\cite[Lecture 14, page 99]{Mu}}]
Let $Y$ be a projective variety over $k$, and let 
${\mathcal O}_Y(1)$ be a very ample line bundle on $Y$.
Take an integer $n_0$.
A coherent sheaf $E$ on $Y$ is called $n_0$-regular if
\[
 H^i(Y,\, E\otimes{\mathcal O}_Y (n_0-i))\ =\ 0
\]
holds for all $i\,>\,0$.
\end{Def}

We will denote $E\otimes{\mathcal O}_Y(m)$ by $E(m)$ for an integer $m$.

\begin{Def} 
Let $Y$ be a projective variety over $k$. A set ${\mathcal T}$ of coherent sheaves on $Y$ is
called bounded if there is a scheme $S$ of finite type over $k$,
and a coherent sheaf ${\mathcal E}$ on $Y\times S$, such that for any member $E\,\in\, {\mathcal T}$,
there is a $k$-valued point $s\,\in\, S$ such that ${\mathcal E}|_{Y\times\{s\}}\,\cong\, E$.
\end{Def}

The following Lemma is a useful tool to show the boundedness of
a family of coherent sheaves.

\begin{Lem}[{\cite[Theorem 1.13]{Kl}}]\label{criterion of boundedness}
Let $Y$ be a projective variety over $k$, and let ${\mathcal O}_Y(1)$ be a very ample line bundle on $Y$.
Then a set ${\mathcal T}$ of coherent sheaves on $Y$ is bounded
if and only if there is an integer $n_0$ such that
all the members of ${\mathcal T}$ are $n_0$-regular
and the set
$$\left\{ \chi(E(m))\,=\, \sum_i (-1)^i \dim H^i(X,\, E(m))
\,\, \middle|\,\,\, E\,\in\,{\mathcal T}\right\}$$
of Hilbert polynomials $\chi(E(m))$
in $m$ of the members $E$ of ${\mathcal T}$ is finite.
\end{Lem}

In the same way as Proposition \ref{Prop: moduli of simple parabolic connections},
the moduli space of simple $(\boldsymbol{x}, \,\boldsymbol{\nu})$-parabolic connections
$\mathcal{M}_{\mathrm{PC}}(\boldsymbol{\nu})$
is an algebraic space over $k$, and the moduli space
$\mathcal{M}_{\mathrm{PC}}^{\boldsymbol{\alpha}}(\boldsymbol{\nu})$
of $\boldsymbol{\alpha}$-stable
$(\boldsymbol{x},\, \boldsymbol{\nu})$-parabolic connections
is a Zariski open subspace of 
$\mathcal{M}_{\mathrm{PC}}(\boldsymbol{\nu})$.
Since $\mathcal{M}_{\mathrm{PC}}^{\boldsymbol{\alpha}}(\boldsymbol{\nu})$
is quasi-projective over $k$,
we can take an integer $n_0$ such that for all
$(E,\,\nabla,\,\boldsymbol{l})\,\in\, \mathcal{M}_{\mathrm{PC}}^{\boldsymbol{\alpha}}(\boldsymbol{\nu})$,
the underlying vector bundle $E$ is $n_0$-regular.

Fix a line bundle $L$ on $X$ and a logarithmic connection
\[
 \nabla_L\,\colon\, L\, \longrightarrow\, L\otimes K_X(D)\, 
\]
such that
$\mathrm{res}_{x_i}(\nabla_L)\,=\,\sum_{j=0}^{r-1}\nu^{(i)}_j$
for all $1\,\leq\, i\,\leq\, n$. Set
\begin{align}
 \mathcal{M}_{\mathrm{PC}}(\boldsymbol{\nu},\nabla_L)
 \,:&=\,
 \left\{
 (E,\,\nabla,\,\boldsymbol{l})\,\in\,
 \mathcal{M}_{\mathrm{PC}}(\boldsymbol{\nu})
 \,\,\big\vert\,
 \det(E,\,\nabla)\,\cong\, (L,\,\nabla_L)
 \right\}, \label{ne1}
 \\
 \mathcal{M}_{\mathrm{PC}}^{\boldsymbol{\alpha}}(\boldsymbol{\nu},\nabla_L)
 \,:&=
 \left\{
 (E,\,\nabla,\,\boldsymbol{l})\,\in\,
 \mathcal{M}_{\mathrm{PC}}^{\boldsymbol{\alpha}}(\boldsymbol{\nu})
 \,\,\big\vert\, \det(E,\,\nabla)\,\cong\, (L,\,\nabla_L)
 \right\}.\label{ne3}
\end{align}
These are closed subspaces of
$\mathcal{M}_{\mathrm{PC}}(\boldsymbol{\nu})$ and
$ \mathcal{M}_{\mathrm{PC}}^{\boldsymbol{\alpha}}(\boldsymbol{\nu})$,
respectively.
Setting
\[
\mathcal{M}^{n_0\text{\rm -reg}}_{\mathrm{PC}}(\boldsymbol{\nu})\,:=\,
 \left\{ (E,\,\nabla,\,\boldsymbol{l})\,\in \,\mathcal{M}_{\mathrm{PC}}(\boldsymbol{\nu})
 \,\, \big\vert \, \text{$E$ is $n_0$-regular} \right\},
\]
there is a canonical open immersion
\[
 \iota\,\,\colon\,\, 
 \mathcal{M}_{\mathrm{PC}}^{\boldsymbol{\alpha}}(\boldsymbol{\nu})
\,\, \hookrightarrow\,\,
\mathcal{M}_{\mathrm{PC}}^{n_0\text{\rm -reg}}(\boldsymbol{\nu}).
\]
Set
\[
\mathcal{M}_{\mathrm{PC}}^{n_0\text{\rm -reg}}(\boldsymbol{\nu},\nabla_L)
 \,\,:=\,\,
 \left\{
 (E,\,\nabla,\,\boldsymbol{l})\,\in\,
 \mathcal{M}_{\mathrm{PC}}^{n_0\text{\rm -reg}}(\boldsymbol{\nu})
 \,\,\big\vert\,
 \det(E,\,\nabla)\,\cong\, (L,\,\nabla_L)
 \right\}.
\]
Then $\mathcal{M}_{\mathrm{PC}}^{n_0\text{\rm -reg}}(\boldsymbol{\nu},\nabla_L)$
is a closed subspace of
$\mathcal{M}^{n_0\text{\rm -reg}}_{\mathrm{PC}}(\boldsymbol{\nu})$
and it contains
$\mathcal{M}_{\mathrm{PC}}^{\boldsymbol{\alpha}}(\boldsymbol{\nu},\nabla_L)$
as a Zariski open subspace.
 
Under the assumption that the rank $r$ is not divisible by the characteristic of $k$,
the proof of the smoothness of the moduli space given in \cite [Theorem 2.1]{Inaba-1}
works because the assumption ensures that the Killing form on $\text{sl}(r,k)$
remains nondegenerate. This is elaborated in the following proposition.

\begin{Prop}\label{Prop: smoothness in arbitrary characteristic}
Assume that the $r$ is not divisible by the characteristic of $k$.
Then the moduli space 
$\mathcal{M}_{\mathrm{PC}}(\boldsymbol{\nu},\nabla_L)$
is smooth over $k$ and so is its open subspace
$\mathcal{M}_{\mathrm{PC}}^{\boldsymbol{\alpha}}(\boldsymbol{\nu},\nabla_L)$.
\end{Prop}

\begin{proof}
We use the criterion of smoothness
in \cite[Proposition 17.14.2]{Groth2}.
Let $A$ be an Artinian local ring over $k$ with the maximal ideal
$\mathfrak{m}$, and let $I$ be an ideal of $A$ such that $\mathfrak{m}I\,=\,0$.
Suppose that we are given a morphism 
$\Spec A/I\,\longrightarrow \,\mathcal{M}_{\mathrm{PC}}(\boldsymbol{\nu},\nabla_L)$
which corresponds to a flat family $(E,\,\nabla,\,\boldsymbol{l})$ of parabolic connections on
$X\times\Spec A/I$ over $A/I$. It suffices to construct a flat family
$(\widetilde{E},\,\widetilde{\nabla},\,\widetilde{\boldsymbol{l}})$
of $\boldsymbol{\nu}$-parabolic connections on $X\times\Spec A$ over
$\Spec A$ which is a lift of $(E,\,\nabla,\,\boldsymbol{l})$.

There is an isomorphism 
$\varphi\,\colon\, \det E\,\xrightarrow{\,\,\sim\,\,\,}\, L\otimes A/I$
such that $(\nabla_L\otimes A/I)\circ\varphi\,=\,(\varphi\otimes\mathrm{id})\circ\Tr(\nabla)$. 
Take an affine open covering $\{U_{\alpha}\}$ of $X$ satisfying the condition that there is an isomorphism
$$\phi_{\alpha}\,\colon\, E|_{U_{\alpha}\times\Spec A/I}\,
\xrightarrow{\,\,\sim\,\,\,}\, {\mathcal O}_{U_{\alpha}\times\Spec A/I}^{\oplus r}.$$
Set
$\overline{\phi_{\alpha}}\,:=\,\phi_{\alpha}\otimes A/\mathfrak{m}$
and $\overline{\varphi}\,:=\,\varphi\otimes A/\mathfrak{m}$.
After replacing $\phi_{\alpha}$ with
$(1+r^{-1}a)\phi_{\alpha}$ for some $a\,\in\, I{\mathcal O}_{U_{\alpha}}$,
we may assume that $$\det(\phi_{\alpha})\circ\varphi^{-1}\ =\
(\det(\overline{\phi_{\alpha}})\circ\overline{\varphi}^{-1})\otimes \mathrm{id}_{A/I}$$ as maps from
$L\otimes A/I$ to ${\mathcal O}_{U_{\alpha}\times\Spec A/I}$.
Set $E_{\alpha}\,:=\,{\mathcal O}_{U_{\alpha}\times\Spec A}^{\oplus r}$
and put
$$\varphi_{\alpha}\,:=\, (\overline{\varphi}\circ\det(\overline{\phi_{\alpha}})^{-1})\otimes A
\, \colon\,\, \det(E_{\alpha})\, \xrightarrow{\,\,\sim\,\,\,}\, L\otimes \mathrm{id}_A.$$
Choose a lift
\[
\theta_{\beta\alpha}\,\colon\,\, E_{\alpha}|_{U_{\alpha\beta}\times\Spec A}
\, \xrightarrow{\,\,\sim\,\,\,}\, E_{\beta}|_{U_{\alpha\beta}\times\Spec A}
\]
of $\phi_{\beta}\circ\phi_{\alpha}^{-1}$. Replacing $\theta_{\beta\alpha}$ with
$(1+r^{-1}b)\theta_{\beta\alpha}$ for some $b\,\in\, I{\mathcal O}_{U_{\alpha\beta}\times \Spec A}$,
we may assume that coincidence
$$\det(\theta_{\beta\alpha})\ =\
\varphi_{\beta}^{-1}\circ\varphi_{\alpha}$$
as maps from $\det(E_{\alpha})|_{U_{\alpha\beta}\times\Spec A}$ to
$\det(E_{\beta})|_{U_{\alpha\beta}\times\Spec A}$.
If $x_i\,\in\, U_{\alpha}$, then we take a quasi-parabolic structure $l^{\alpha}_*$ on $E_{\alpha}$
at $x_i\times\Spec A$ which is a lift of $l^{(i)}_*$.
Take a relative connection
\[
\nabla_{\alpha}\ \colon\ E_{\alpha} \ \longrightarrow \ E_{\alpha}\otimes
\Omega_{X\times \Spec A/\Spec A}(D\times\Spec A)
\]
such that 
$\nabla_{\alpha}\otimes A/I\,=\,
\phi_{\alpha}\circ\nabla|_{U_{\alpha}\times\Spec A/I}\circ\phi_{\alpha}^{-1}$ and
$\big(\mathrm{res}_{x_i\times\Spec A}(\nabla_{\alpha})-\nu^{(i)}_j\big)
(l^{\alpha}_j)\,\subset\, l^{\alpha}_{j+1}$ for all $0\,\leq\, j\,\leq\, r-1$.
After replacing $\nabla_{\alpha}$ with
$\nabla_{\alpha}+r^{-1}\eta\otimes\mathrm{id}_{E_{\alpha}}$ for some
$\eta\,\in\, I \, \Omega^1_{U_{\alpha}\times\Spec A/\Spec A}$,
we may assume that
$\varphi_{\alpha}\Tr(\nabla_{\alpha})\varphi_{\alpha}^{-1}\,=\,\nabla_L\otimes\mathrm{id}_A$.
Put
$(\overline{E},\,\overline{\nabla},\,\overline{\boldsymbol{l}})\,:=\,
(E,\,\nabla,\,\boldsymbol{l})\otimes A/\mathfrak{m}$ and set
\begin{align*}
 {\mathcal D}^{\mathrm{par}}_{\mathfrak{sl},0}
 &\,=\,
 \left\{ u\,\in\,{\mathcal End}(\overline{E}) \, \middle| \,\,
 \text{$\Tr(u)\,=\,0$ and $\mathrm{res}_{x_i}(u)(\overline{l}^{(i)}_j)\,\subset\, \overline{l}^{(i)}_j$
 for any $i,\,j$} \right\}
 \\
{\mathcal D}^{\mathrm{par}}_{\mathfrak{sl},1}
 &\,=\,
 \left\{ v\,\in\,{\mathcal End}(\overline{E})\otimes K_X(D) \, \middle| \,\,
 \text{$\Tr(v)\,=\,0$ and $\mathrm{res}_{x_i}(v)(\overline{l}^{(i)}_j)\,\subset\, \overline{l}^{(i)}_{j+1}$
 for any $i,\,j$} \right\}
 \\
\nabla_{{\mathcal D}^{\mathrm{par}}_{\mathfrak{sl},\bullet}} &\,\colon\,
{\mathcal D}^{\mathrm{par}}_{\mathfrak{sl},0} \, \longrightarrow\,
{\mathcal D}^{\mathrm{par}}_{\mathfrak{sl},1}, \ \ \,u \, \longmapsto\,
\overline{\nabla}\circ u-u\circ\overline{\nabla}.
\end{align*}
Then we get a cohomology class
$[\{\theta_{\gamma\alpha}^{-1}\theta_{\gamma\beta}\theta_{\beta\alpha}-\mathrm{id}\},\,
\{ \theta_{\beta\alpha}^{-1}\circ\nabla_{\beta}\circ\theta_{\beta\alpha}-\nabla_{\alpha}\}]
\,\in\, \mathbb{H}^2({\mathcal D}^{\mathrm{par}}_{\mathfrak{sl},\bullet})
\otimes I$
whose vanishing is equivalent to the existence of a lift
$(\widetilde{E},\,\widetilde{\nabla},\,\widetilde{\boldsymbol{l}})
\,\in\,{\mathcal M}_{\mathrm{PC}}(\boldsymbol{\nu},\nabla_L)(A)$
of $(E,\,\nabla,\,\boldsymbol{l})$. There is a commutative diagram with exact rows
\[
\begin{CD}
H^1( {\mathcal D}^{\mathrm{par}}_{\mathfrak{sl},0} )
@>>>
H^1( {\mathcal D}^{\mathrm{par}}_{\mathfrak{sl},1} )
@>>>
\mathbb{H}^2({\mathcal D}^{\mathrm{par}}_{\mathfrak{sl},\bullet})
@>>> 0
\\
@V\cong VV @V\cong VV @V\cong VV \\
H^0({\mathcal D}^{\mathrm{par}}_{\mathfrak{sl},1} )^{\vee}
@>>>
H^0( {\mathcal D}^{\mathrm{par}}_{\mathfrak{sl},0} )^{\vee}
@>>>
\mathbb{H}^0({\mathcal D}^{\mathrm{par}}_{\mathfrak{sl},\bullet})^{\vee}
@>>> 0
\end{CD}
\]
induced by the Serre duality.
Take any member 
$$u\,\in\, \mathbb{H}^0({\mathcal D}^{\mathrm{par}}_{\mathfrak{sl},\bullet})
\,=\,\ker\big( H^0( {\mathcal D}^{\mathrm{par}}_{\mathfrak{sl},0} )
\, \xrightarrow{\,\, \nabla_{{\mathcal D}^{\mathrm{par}}_{\mathfrak{sl},\bullet}}\,\,\,}\,
\mathbb{H}^0({\mathcal D}^{\mathrm{par}}_{\mathfrak{sl},1}) \big).$$
Since $(\overline{E},\,\overline{\nabla},\,\overline{\boldsymbol{l}})$ is simple, 
we can write $u\,=\,c\cdot\mathrm{id}_{\overline{E}}$ for some $c\,\in\, k$.
By the definition of ${\mathcal D}^{\mathrm{par}}_{\mathfrak{sl},0}$,
we have $0\,=\,\Tr(u)\,=\,\Tr(c\,\mathrm{id}_{\overline{E}})\,=\,rc$.
Since $r^{-1}\,\in \,k^{\times}$ by the assumption, we have $c\,=\,0$.
Thus $u\,=\,0$ and we have $\mathbb{H}^0({\mathcal D}^{\mathrm{par}}_{\mathfrak{sl},\bullet})\,=\,0$.
So the obstruction space
$\mathbb{H}^2({\mathcal D}^{\mathrm{par}}_{\mathfrak{sl},\bullet})
\,\cong\, \mathbb{H}^0({\mathcal D}^{\mathrm{par}}_{\mathfrak{sl},\bullet})^{\vee}$
vanishes, and there is a lift
$(\widetilde{E},\,\widetilde{\nabla},\,\widetilde{l})\,\in\,
{\mathcal M}_{\mathrm{PC}}(\boldsymbol{\nu},\nabla_L)(A)$
of $(E,\nabla,l)$. This means that ${\mathcal M}_{\mathrm{PC}}(\boldsymbol{\nu},\nabla_L)$
is smooth.
\end{proof}

Using Proposition \ref{Prop: smoothness in arbitrary characteristic}
and a similar calculation as done in 
Lemma \ref{Lemma: infinitesimal deformation of framed connection},
we have the following proposition.

\begin{Prop}[{\cite[Theorem 2.1, Proposition 5.1, Proposition 5.2 and Proposition 5.3]{Inaba-1}}]
\label{Proposition: dimension of the moduli space of parabolic connections}
The dimension of the moduli space
$\mathcal{M}_{\mathrm{PC}}(\boldsymbol{\nu},\nabla_L)$
is 
$2(r^2-1)(g-1)+nr(r-1)$
which is same as the dimension of its open subspace
$\mathcal{M}_{\mathrm{PC}}^{\boldsymbol{\alpha}}(\boldsymbol{\nu},\nabla_L)$.
\end{Prop}

We can similarly define the Higgs bundles. As before,
$\boldsymbol{\alpha}$ is a parabolic weight. Take a tuple
$\boldsymbol{\mu}\,=\,(\mu^{(i)}_j)^{1\leq i\leq n}_{0\leq j\leq r-1}\,\in\,k^{nr}$
satisfying the following condition:
\[
 \sum_{i=1}^n\sum_{j=0}^{r-1} \mu^{(i)}_j\ =\ 0.
\]
We say that a tuple
$(E,\,\Phi, \,\boldsymbol{l}\,=\,\{l^{(i)}_*\}_{1\leq i\leq n})$
(equipped with a parabolic weight $\boldsymbol{\alpha}$)
is a $(\boldsymbol{x},\,\boldsymbol{\mu})$-parabolic Higgs bundle
if
\begin{itemize}
\item[(1)] $E$ is an algebraic vector bundle on $X$ of rank $r$ and degree $d$,
 
\item[(2)] $\Phi\,\colon\, E\,\longrightarrow\, E\otimes K_{X}(D)$ is an
${\mathcal O}_X$-linear homomorphism, and

\item[(3)] $l^{(i)}_* $ is a filtration $$E\big\vert_{x_i} \,=\, l_0^{(i)} \,\supset\, l_1^{(i)}
\,\supset \,\cdots \,\supset\, l_r^{(i)}\,=\,0$$
for every $x_i$ such that $\dim (l_j^{(i)}/l_{j+1}^{(i)})\,=\,1$
and $(\mathrm{res}_{x_i}(\Phi)-\mu^{(i)}_j)(l^{(i)}_j)
\,\subset\, l_{j+1}^{(i)}$ for all $j\,=\,0,\,\cdots,\,r-1$.
\end{itemize}

A $(\boldsymbol{x},\,\boldsymbol{\mu})$-parabolic Higgs bundle
$(E,\,\Phi,\,\boldsymbol{l})$ is said to be simple if every endomorphism
$f\,\colon\, E\,\longrightarrow\, E$ that commutes with $\Phi$ and preserves $\boldsymbol{l}$
is a constant scalar multiplication $f\,=\,c\,\mathrm{Id}_E$ for some $c\,\in\,k$.
Denote by $\mathcal{M}_{\mathrm{Higgs}}(\boldsymbol{\mu})$ the moduli space of
simple $\boldsymbol{\mu}$-parabolic Higgs bundles.
Define $\boldsymbol{\alpha}$-stability for
parabolic Higgs bundles analogous to Definition \ref{def: stability of parabolic connection}.
If we replace $n_0$ by a sufficiently large integer, we may assume that
for every $\boldsymbol{\alpha}$-stable $(\boldsymbol{x},\,\boldsymbol{\mu})$-parabolic
Higgs bundle $(E,\,\Phi,\,\boldsymbol{l})$, the underlying vector bundle $E$ is $n_0$-regular.

Fix a line bundle $L$ on $X$ together with a homomorphism
$\Phi_L\, \colon\, L\, \longrightarrow\, L\otimes K_X$
of ${\mathcal O}_X$--modules such that
$\res_{x_i}(\Phi_L)\,=\,\sum_{j=0}^{r-1}\mu^{(i)}_j$ for any $i$.
Set
\begin{align*}
\mathcal{M}_{\mathrm{Higgs}}(\boldsymbol{\mu},\,\Phi_L)
 \,&:=\,
 \left\{
 (E,\,\Phi,\,\boldsymbol{l})\,\in\,
\mathcal{M}_{\mathrm{Higgs}}(\boldsymbol{\mu})
 \,\,\big\vert\,\,
 (\det(E),\,\Tr(\Phi))\,\cong \,(L,\,\Phi_L)
 \right\},
 \\
\mathcal{M}^{n_0\text{-reg}}_{\mathrm{Higgs}}(\boldsymbol{\mu},\,\Phi_L)
 \,&:=\,
 \left\{
 (E,\,\Phi,\,\boldsymbol{l})\,\in\,
 \mathcal{M}_{\mathrm{Higgs}}(\boldsymbol{\mu},\,\Phi_L)
 \,\,\big\vert\,\,
 \text{$E$ is $n_0$-regular}
 \right\},
 \\
 \mathcal{M}^{\boldsymbol{\alpha}}_{\mathrm{Higgs}}(\boldsymbol{\mu},\,\Phi_L)
 \,&:=\,
 \left\{
 (E,\,\Phi,\,\boldsymbol{l})\,\in\,
 \mathcal{M}_{\mathrm{Higgs}}(\boldsymbol{\mu},\,\Phi_L)
 \,\,\big\vert\,\,
 \text{$(E,\,\Phi,\,\boldsymbol{l})$ is $\boldsymbol{\alpha}$-stable}
 \right\}.
\end{align*}

The same calculations as done in the proof of 
Proposition \ref{Prop: smoothness in arbitrary characteristic} and
Proposition \ref{Proposition: dimension of the moduli space of parabolic connections}
yield the following proposition.

\begin{Prop}[{\cite[Section 2.1]{B-Y}}, {\cite[Theorem 2.8]{Yokogawa}}]\label{Prop: smoothness and dimension of moduli of Hiiggs bundles}
Assume that $r$ is not divisible by the characteristic of $k$.
Then the moduli space
$ \mathcal{M}_{\mathrm{Higgs}}(\boldsymbol{\mu},\,\Phi_L)$ is smooth
and $\dim \mathcal{M}_{\mathrm{Higgs}}(\boldsymbol{\mu},\,\Phi_L)
\,=\,2(r^2-1)(g-1)+nr(r-1)$.
Furthermore, the open subspace
$\mathcal{M}^{\boldsymbol{\alpha}}_{\mathrm{Higgs}}(\boldsymbol{\mu},\,\Phi_L)$
of $\mathcal{M}_{\mathrm{Higgs}}(\boldsymbol{\mu},\,\Phi_L)$
consisting of $\boldsymbol{\alpha}$-stable parabolic Higgs bundles
is quasi-projective.
\end{Prop}

It is known that there is no non-constant global algebraic function
on the moduli space of logarithmic connections 
with central residues on a compact Riemann surface of genus at least $3$ \cite{BR}.
In the logarithmic case, the same statement was proved in \cite{Ari}
in a very special case when $g\,=\,0$, $r\,=\,2$ and $n\,=\,4$.
In \cite{BR}, the Betti number of the moduli space of
stable vector bundles assumed one of the key roles.
A similar result is proved in \cite{Singh}.
We will prove, in this section, a weaker result that the moduli space
$\mathcal{M}_{\mathrm{PC}}^{\boldsymbol{\alpha}}(\boldsymbol{\nu},\,\nabla_L)$
of $(\boldsymbol{x},\,\boldsymbol{\nu})$-parabolic connections
is not affine for any genus, except for several special cases.
We use a part of the ideas in \cite{BR}
and compare the transcendence degree, of the
ring of global algebraic functions on the moduli space 
$\mathcal{M}_{\mathrm{PC}}^{\boldsymbol{\alpha}}(\boldsymbol{\nu},\,\nabla_L)$
of parabolic connections,
with that on the moduli space
$\mathcal{M}_{\mathrm{Higgs}}^{\boldsymbol{\alpha}}(\boldsymbol{0},\,0)$
of parabolic Higgs bundles.
Our argument also works over the base field of positive characteristic,
which is consistent with the existence of the Hitchin map
on the moduli space of connections
(\cite{Las-Pa}, \cite{Groech}).

\subsection{Codimension estimation for non-simple underlying bundle}\label{subsection: codimension estimation}

This subsection provides an improvement of the result 
of \cite[Section 5]{Inaba-1}.
Throughout this subsection, $k$ is assumed to be an algebraically closed field
of arbitrary characteristic.

Now let $X$ be a smooth projective irreducible curve over $\Spec k$ of genus is $g$, and let $D\,=\,x_1+
\ldots+x_n$ be a reduced divisor on $X$. Fix a line bundle $L$ of degree $d$ on $X$. Consider the set
\[
 \left| {\mathcal N}^{n_0\text{\rm -reg}}_{\mathrm{par}}(L) \right|
 \ =\ \left\{
 (E,\,\boldsymbol{l}) \right\}/\!\cong
\]
of isomorphism classes of quasi-parabolic bundles $(E,\,\boldsymbol{l})$ on $(X,\,D)$ such that
\begin{itemize}
\item[(i)] $E$ is an algebraic vector bundle on $X$ of rank $r$ with $\det E\,\cong \,L$,

\item[(ii)] $\boldsymbol{l}$ is a quasi-parabolic structure consisting of
filtrations $$E|_{x_i}\,=\,l^{(i)}_0\,\supsetneq\, l^{(i)}_1\,\supsetneq\,\cdots\,\supsetneq\,
l^{(i)}_{r-1}\,\supsetneq\, l^{(i)}_r\,=\,0$$
for every $x_i\,\in\, D$ and

\item[(iii)] $E_0$ is $n_0$-regular.
\end{itemize}
By virtue of Lemma \ref{criterion of boundedness}, there is a scheme $Z$ of finite type
over $\Spec k$ and a flat family $(\widetilde{E},\,\widetilde{\boldsymbol{l}})$
of quasi-parabolic bundles on $X\times Z$ over $Z$
such that every member $(E,\,\boldsymbol{l})\,\in\, {\mathcal N}^{n_0\text{\rm -reg}}_{\mathrm{par}}(L)$
is isomorphic to $(\widetilde{E},\,\widetilde{\boldsymbol{l}})|_{X\times \{p\}}$
for some point $p\,\in\, S$. Consider the subset
\[
\left| {\mathcal N}^{n_0\text{\rm -reg}}_{\mathrm{par}}(L)^{\circ} \right|
\,\,=\,\,
\left\{ (E,\,\boldsymbol{l})\,\in\, \left| {\mathcal N}^{n_0\text{\rm -reg}}_{\mathrm{par}}(L) \right| 
\,\, \middle|\,\,\, \dim \End (E,\,\boldsymbol{l}) \,=\,1 \right\}
\]
of $|{\mathcal N}^{n_0\text{\rm -reg}}_{\mathrm{par}}(L)|$
consisting of simple quasi-parabolic bundles,
where $\End(E,\,\boldsymbol{l})$ is defined by
\[
 \End(E,\,\boldsymbol{l})\,=\, \left\{
 u\,\in\, \Hom_{{\mathcal O}_X}(E,E) \,\,\middle|\,\,\,
 \text{$u|_x\colon E|_x\rightarrow E|_x$ satisfies
 $u|_x(l^x_j)\subset l^x_j$ for any $x\in D$ and any $j$} \right\}.
\]

\begin{Def}
Let $X$ be a smooth projective curve over $k$.
For a vector bundle $E$ on $X$,
we set $\mu(E)\,:=\,\deg(E)/\rank(E)$, and call it
the slope of $E$.
\end{Def}

We will construct a parameter space of 
$\left| {\mathcal N}^{n_0\text{\rm -reg}}_{\mathrm{par}}(L) \right|
\setminus \left| {\mathcal N}^{n_0\text{\rm -reg}}_{\mathrm{par}}(L)^{\circ} \right|$
whose dimension is at most
$(r^2-1)(g-1)+nr(r-1)/2-2$.
For its proof, we need the following lemma.

\begin{Lem}\label{key inequality coming from Riemann--Roch}
Let $X$ be a smooth projective curve over $k$ of genus $g\,\geq\, 2$,
and let $E$ and $F$ be semistable vector bundles on $X$
satisfying the condition that $\mu(E)\,>\,\mu(F)$.
Then the following inequality holds:
\[
\dim\Ext^1_X(F,\,E)\ \leq\
\max \big\{ \rank(E)\rank(F)(2g-3) ,\, \rank(E)\rank(F)g-1 \big\}.
\]
\end{Lem}

\begin{proof}
By the Serre duality,
we have
$\dim \Ext^1_X(F,\,E)\,=\,
\dim \Hom(E,\,F\otimes K_X)$.
Choose a general point $x\,\in\, X$.

First consider the case where $\deg(E^{\vee}\otimes F\otimes{\mathcal O}_X(x))\,>\,0$.
In this case, we have $\Hom(F,\, E\otimes{\mathcal O}_X(-x))\,=\,0$.
Note that we have
$\deg(F^{\vee}\otimes E\otimes{\mathcal O}_X(-x))\,>\,-\rank(E)\rank(F)$, because $\mu(E)\,>\,\mu(F)$.
By the Riemann--Roch theorem, we have
\begin{align*}
 \dim\Ext^1_X(F,\,E)\, =\,
 \dim \Hom(E,\,F\otimes K_X)
 &\,\leq\,
 \dim \Hom(E,\,F\otimes K_X(x))
 \\
 &=\,
 \dim\Ext^1(F,\,E\otimes{\mathcal O}_X(-x))
 \\
 &=\,
 -\rank(E)\rank(F)(1-g)-\deg(F^{\vee}\otimes E\otimes{\mathcal O}_X(-x))
 \\
 &\leq\, \rank(E)\rank(F)g-1.
\end{align*}

Secondly, consider the case where
$\deg(E^{\vee}\otimes F\otimes{\mathcal O}_X(x))\,<\,0$.
Take general points $x_1,\,\cdots,\,x_{2g-3}$ of $X$.
Then we get exact sequences
$$
0 \,\longrightarrow\, H^0(X,\,E^{\vee}\otimes F\otimes K_X(-x_1-\cdots-x_i))
\,\longrightarrow\, H^0(X,\,E^{\vee}\otimes F\otimes K_X(-x_1-\cdots-x_{i-1}))
$$
$$
\longrightarrow\, E^{\vee}\otimes F\otimes K_X(-x_1-\cdots-x_{i-1})\big|_{x_i}
$$
for $i\,=\,1,\,\cdots,\,2g-3$. Note that the condition $\deg(E^{\vee}\otimes F(x))\,<\,0$
implies that $\mu(E)\,>\,\mu(F\otimes K_X(-x_1-\cdots-x_{2g-3}))$,
which yields the following
$$\Hom(E,\, F\otimes K_X(-x_1-\cdots-x_{2g-3}))\ =\ 0,$$
because $E$ and $F$ are semistable.
So we have
\begin{align*}
 \dim \Ext^1_X(F,\,E)
 \,=\,
 \dim H^0(X,\,E^{\vee}\otimes F\otimes K_X)
 &\,\leq\,
 \sum_{i=1}^{2g-3} \dim_{\mathbb{C}} \left(E^{\vee}\otimes F\otimes K_X(-x_1-\cdots-x_{i-1})\big|_{x_i}\right)
 \\
 &=\,
 \rank(E)\rank(F)(2g-3).
\end{align*}

Consider the remaining case where
$\deg ( E^{\vee}\otimes F(t))\,=\,0$.
Take general points $x_1,\,\cdots,\,x_{2g-3}\,\in\, X$.
Then we have $\mu(E)\,=\,\mu(F\otimes K_X(-x_1-\cdots-x_{2g-3}))$. 
We can write $\mathrm{gr}(E)\,=\,\bigoplus_i E_i$
and $\mathrm{gr}(F)\,=\,\bigoplus_jF_j$
for stable vector bundles $E_i$ and $F_j$ such that
$\mu(E_i)\,=\,\mu(E)\,=\,\mu(F)\,=\,\mu(F_j)$ for any $i,\,j$.
If we take $x_1,\,\cdots,\,x_{2g-3}$ sufficiently generic,
then we may assume $E_i\,\not\cong\, F_j\otimes K_X(-x_1-\cdots-x_{2g-3})$
for any $i,\,j$. Then we have
$\Hom(E,\,F\otimes K_X(-x_1-\cdots-x_{2g-3}))\,=\,0$.
By the same argument as before, we have the inequality
$\dim\Ext^1_X(E,F)\,\leq\, \rank(E)\rank(F)(2g-3)$.
\end{proof}

\begin{Prop}\label{dimension of non-simple parabolic bundles in higher genus case}
Let $X$ be a smooth projective curve over $k$ of genus $g\,\geq\, 2$,
and let $L$ be a line bundle of degree $d$ on $X$.
Assume that the integers $r$ and $n$ satisfy the conditions $r\,\geq\, 2$ and $n\,\geq\, 1$.
Then there exists a scheme $Z$ of finite type over $\Spec k$
and a flat family $(\mathcal{E},\,\ell)$ of quasi-parabolic bundles on $X\times Z$ over $Z$ such that
\begin{itemize}
\item
$\dim Z \,\,\leq\,\, (r^2-1)(g-1)+r(r-1)n/2-2$, 
\item
$\dim \End\left((\mathcal{E},\ell)\big\vert_{X\times\{z\}}\right)\,\,\geq\,\, 2$
for any $z\,\in\, Z$,
\end{itemize}
and each member of the complement
$\left| {\mathcal N}^{n_0\text{\rm -reg}}_{\mathrm{par}}(L) \right|
\setminus \left| {\mathcal N}^{n_0\text{\rm -reg}}_{\mathrm{par}}(L)^{\circ} \right|$
is isomorphic to $(\mathcal{E},\,\ell)|_{X\times\{z\}}$
for some point $z\,\in\, Z$.
\end{Prop}

\begin{proof}
Take a quasi-parabolic bundle $(E,\,\boldsymbol{l})$ on $(X,\,D)$ with $\det E\, \cong\, L$.
Choose a point $x_i\,\in\, D$ and $l^{(i)}_j\,\subset\, E|_{x_i}$.
Then $E'\,:=\,\ker (E\,\rightarrow\, E|_{x_i}/l^{(i)}_j)$ has a canonical 
quasi-parabolic structure $\boldsymbol{l'}$ induced by $\boldsymbol{l}$.
The correspondence
$(E,\,\boldsymbol{l})\,\longmapsto\, (E',\,\boldsymbol{l'})$ gives a bijection between
the set of isomorphism classes of quasi-parabolic bundles;
it is called an elementary transformation or a Hecke modification.
After applying a finite number of elementary transformations,
it may be assumed that $r$ and $d$ are coprime.

Take a member $(E,\,\boldsymbol{l})\,\in\,
|{\mathcal N}^{n_0\text{\rm -reg}}_{\mathrm{par}}(L)|
\setminus |{\mathcal N}^{n_0\text{\rm -reg}}_{\mathrm{par}}(L)^{\circ}|$.
Since $\dim \End(E,\,\boldsymbol{l})\,>\,1$ by the definition,
we have $\dim \End(E)\,>\,1$ and $E$ is not a semistable vector bundle.
Let
\[
 E_1\,\subset\, E_2\,\subset\,\cdots\,\subset\, E_m\,=\,E
\]
be the Harder--Narasimhan filtration of $E$; note that $m\,\geq\, 2$ because $E$
is not semistable. Set $\overline{E}_1\,:=\,E_1$, $\overline{E}_s\,:=\,E_s/E_{s-1}$ for $s\,\geq\, 2$
and $r_s\,:=\,\rank\,\overline{E}_s$.
By the definition of a Harder--Narasimhan filtration,
each $\overline{E}_s$ is semistable for $1\,\leq\, s\,\leq\, m$
and the inequalities $\mu(\overline{E}_1)\,>\,\mu(\overline{E}_2)\,>\,\cdots\,>\,\mu(\overline{E}_m)$
hold. Each semistable vector bundle $\overline{E}_s$ has a Jordan--H\"{o}lder filtration
\[
 0\,\subset\, E^{(1)}_s\,\subset\, E^{(2)}_s\,\subset\,\cdots\,\subset\,
 E^{(\gamma_s)}_s \,=\, \overline{E}_s
\]
with $\gamma_s\,\geq \,1$. Set $\overline{E}^{(1)}_s\,:=\,E^{(1)}_s$,
$\overline{E}^{(i)}_s\,:= \,E^{(i)}_s/E^{(i-1)}_s$ for $2\,\leq\, i\,\leq \,\gamma_s$,
$r^{(i)}_s\,:=\,\rank\overline{E}^{(i)}_s$
and $d^{(i)}_s\,:=\,\deg\overline{E}^{(i)}_s$,
Then each $\overline{E}^{(i)}_s$ is a stable bundle on $X$ and
\[
 L \ \cong \ \det \left(\bigoplus_{s=1}^m \bigoplus_{i=1}^{\gamma_s} \overline{E}^{(i)}_s\right)
\]
holds.

Let us consider the converse. If stable bundles
$\{\overline{E}^{(i)}_s\}$ are given,
$\{\overline{E}_s\}$ are given by successive extensions
\begin{equation}\label{fundamental piece of extensions of stable bundles}
 0\, \longrightarrow\, E^{(i-1)}_s \,\longrightarrow\, E^{(i)}_s
\,\longrightarrow\, \overline{E}^{(i)}_s \,\longrightarrow\, 0
 \qquad
 (2\,\leq\, i \,\leq \,\gamma_s)
\end{equation}
with $E^{(\gamma_s)}_s\,=\,\overline{E}_s$.
If $\{\overline{E}_s\}$ are given, then $E$ is given by
successive extensions
\begin{equation}\label{fundamental piece of extensions of semistable sheaves}
 0\,\longrightarrow\, E_{s-1} \,\longrightarrow\, E_s \,\longrightarrow \,\overline{E}_s
\, \longrightarrow\, 0
 \qquad
 (2\,\leq\, s\,\leq\, m)
\end{equation}
with $E_m\,=\,E$. By its definition, $\boldsymbol{l}$ is given by a filtration
$E\big\vert_{x_i}\,=\,l^{(i)}_0\,\supset\, l^{(i)}_1\,\supset\,\cdots\,\supset\, l^{(i)}_{r-1}\,\supset\, l^{(i)}_r\,=\,0$
for each $1\,\leq\, i\,\leq\, n$.

We will construct a parameter space of the above data,
but we avoid the case of $m\,=\,2$ and $\gamma_1\,=\,\gamma_2\,=\,1$
and postpone its proof later. This is because this case needs an extra argument.

Excluding the case where $m\,=\,2$ and $\gamma_1\,=\,\gamma_2\,=\,1$,
we first construct the parameter space of the above data
with the further restricted conditions that
\begin{itemize}
\item $\overline{E}^{(i)}_s\,\not\cong\, \overline{E}^{(j)}_s$ for $i\,\neq\, j$, and

\item all the extensions in \eqref{fundamental piece of extensions of stable bundles} and
\eqref{fundamental piece of extensions of semistable sheaves}
do not split.
\end{itemize}
Set
\[
 N\ =\
 \left\{ (\overline{E}^{(i)}_s) \in \prod_{s=1}^m \prod_{i=1}^{\gamma_s} 
 {\mathcal N}^e (r^{(i)}_s,d^{(i)}_s)^{\circ}
 \,\middle|\, 
 \bigotimes_{s=1}^m\bigotimes_{i=1}^{\gamma_s} 
 \det(\overline{E}^{(i)}_s) \cong L \right\},
\]
where ${\mathcal N}^e (r^{(i)}_s, d^{(i)}_s)^{\circ}$
is the moduli space of stable vector bundles on $X$ of rank $r^{(i)}_s$
and of degree $d^{(i)}_s$. Since $\dim {\mathcal N}^e (r^{(i)}_s, d^{(i)}_s)^{\circ}
\,=\,(r^{(i)}_s)^2(g-1)+1$, we have
$\dim N\,=\,\sum_{s=1}^m \sum_{i=1}^{\gamma_s} \big( (r^{(i)}_s)^2(g-1)+1\big) -g$.
Take a quasi-finite covering $N'\,\longrightarrow\, N$
whose image consists of those points such that
$\overline{E}^{(i)}_s\,\not\cong\,\overline{E}^{(i')}_{s'}$ for $(i,\,s)\,\neq\,(i',\,s')$.
We may take a universal family of vector bundles
$\{\overline{\mathcal E}^{(i)}_s\}^{1\leq i\leq \gamma_s}_{1\leq s\leq m}$
on $X\times N'$ over $N'$ such that
$\bigotimes_{s=1}^m\bigotimes_{i=1}^{\gamma_s} 
\det \big( \overline{\mathcal E}^{(i)}_s \big)
\,\cong\, L\otimes{\mathcal L'}$ for some line bundle ${\mathcal L'}$ on $N'$.
After replacing $N'$ with a disjoint union of locally closed subsets,
we may further assume that
\begin{itemize}
\item the relative $\Ext$-sheaves
$\Ext^p_{X\times N'/N'} \big( \overline{\mathcal E}^{(i)}_s ,\, \overline{\mathcal E}^{(j)}_s \big)$
are locally free sheaves on $N'$ for $1\,\leq\, s\,\leq\, m$,
$p\,=\,0,\,1$ and any $j\,<\,i$, and

\item the canonical maps
$\Ext^p_{X\times N'/N'} \big( \overline{\mathcal E}^{(i)}_s ,\, \overline{\mathcal E}^{(j)}_s \big)\big|_z
\,\longrightarrow\, \Ext^p_{X\times\{z\}}\big( \overline{\mathcal E}^{(i)}_s|_{X\times\{z\}} , \,
\overline{\mathcal E}^{(j)}_s|_{X\times\{z\}} \big)$
are isomorphisms for all points $z\in N'$.
\end{itemize}
Set
\[
 P^{(2)}_s\ :=\
 \mathbb{P}_*\Ext^1_{X\times N'/N'} 
\big( \overline{\mathcal E}^{(2)}_s ,\, \overline{\mathcal E}^{(1)}_s \big)
\,=\, \Proj \mathrm{Sym}\left(\Ext^1_{X\times N'/N'} 
\big( \overline{\mathcal E}^{(2)}_s ,\, \overline{\mathcal E}^{(1)}_s \big)^{\vee}\right)
\]
for every $1\,\leq\, s\,\leq\, m$, where 
$\mathrm{Sym}\left(\Ext^1_{X\times N'/N'} 
\big( \overline{\mathcal E}^{(2)}_s ,\, \overline{\mathcal E}^{(1)}_s \big)^{\vee}\right)$
is the symmetric algebra of $\Ext^1_{X\times N'/N'} 
\big( \overline{\mathcal E}^{(2)}_s ,\, \overline{\mathcal E}^{(1)}_s \big)^{\vee}$
over ${\mathcal O}_{N'}$. Then there is a universal extension
\[
0 \,\longrightarrow\, \overline{\mathcal E}^{(1)}_s\, \longrightarrow\, {\mathcal E}^{(2)}_s
\,\longrightarrow\, \overline{\mathcal E}^{(2)}_s\otimes{\mathcal O}_{P^{(2)}_s}(1)
\, \longrightarrow\, 0
\]
on $X\times P^{(2)}_s$. Once $P^{(2)}_s,\ldots,P^{(i)}_s$
and ${\mathcal E}^{(2)}_s,\,\cdots,\,{\mathcal E}^{(i)}_s$
are defined, we set
\[
P^{(i+1)}_s\ =\ 
\mathbb{P}_*\Ext^1_{X\times_{N'} P^{(i)}_s/P^{(i)}_s} 
\big( \overline{\mathcal E}^{(i+1)}_s ,\, {\mathcal E}^{(i)}_s \big).
\]
There is a universal extension
\[
0 \,\longrightarrow\, {\mathcal E}^{(i)}_s\,\longrightarrow\, {\mathcal E}^{(i+1)}_s
\,\longrightarrow\, \overline{\mathcal E}^{(i+1)}_s\otimes{\mathcal O}_{P^{(i+1)}_s}(1)
\, \longrightarrow\, 0
\]
on $X\times P^{(i+1)}_s$. Set $P_s\,:=\,P^{(\gamma_s)}_s$ for $1\,\leq\, s\,\leq\, m$,
$P\,:=\, P_1\times_{N'}\times\cdots\times_{N'} P_m$
and $\overline{\mathcal E}_s\,:=\,{\mathcal E}^{(\gamma_s)}_s\otimes_{O_{P_s}}{\mathcal O}_P$.
After replacing $P$ with a disjoint union of locally closed subsets,
we may assume that
\begin{itemize}
\item the relative $\Ext$-sheaves
$\Ext^p_{X\times P/P} \big( \overline{\mathcal E}_s,\,\overline{\mathcal E}_{s'}\big)$
are all locally free sheaves on $P$ for $p\,=\,0,\,1$ and $s'\,<\,s$, and

\item the canonical homomorphisms
$\Ext^p_{X\times P/P} \big( \overline{\mathcal E}_s,\,\overline{\mathcal E}_{s'}\big)\big|_z
\,\longrightarrow \,\Ext^p_{X\times P/P} \big( \overline{\mathcal E}_s\big|_{X\times\{z\}},\,
\overline{\mathcal E}_{s'}\big|_{X\times\{z\}}\big)$
are isomorphisms for all points $z\,\in\, P$.
\end{itemize}
Set
\[
 Q_2\ :=\
 \mathbb{P}_*\Ext^1_{X\times P/P} \big( \overline{\mathcal E}_2,\,\overline{\mathcal E}_1\big)
 \ =\ 
 \Proj\mathrm{Sym} \left( 
 \Ext^1_{X\times P/P} \big( \overline{\mathcal E}_2,\,\overline{\mathcal E}_1\big)^{\vee} \right).
\]
Then there is a universal extension
$0\,\longrightarrow\, \overline{\mathcal E}_1 \,\longrightarrow \,
 {\mathcal E}_2 \,\longrightarrow\, \overline{\mathcal E}_2 \otimes{\mathcal O}_{Q_2}(1)
\, \longrightarrow\, 0
$
on $X\times Q_2$. Once $Q_2,\,\cdots,\, Q_s$
and ${\mathcal E}_2,\,\cdots,\, {\mathcal E}_s$ are defined, set
\[
Q_{s+1}\ :=\
\mathbb{P}_* \Ext^1_{X\times Q_s/Q_s} \big( \overline{\mathcal E}_{s+1},\, {\mathcal E}_s\big).
\]
Then there are universal extensions
$0\,\longrightarrow\, {\mathcal E}_s \,\longrightarrow\,{\mathcal E}_{s+1} 
\, \longrightarrow\, \overline{\mathcal E}_{s+1} \otimes{\mathcal O}_{Q_{s+1}}(1)
\, \longrightarrow\, 0$
for $1\,\leq\, s\,\leq\, m-1$.
Set
\[
 Q\, :=\, \left\{ z\,\in\, Q_m\,\,\middle|\,\,\,
 \text{$ {\mathcal E}_m|_{X\times \{z\}}$ is $n_0$-regular} \right\},
 \qquad
 {\mathcal E}\,:=\,{\mathcal E}_m|_{X\times Q}.
\]
Let $Y_Q$ be the flag bundle over $Q$ whose fiber over any $q\,\in\, Q$
is the parameter space of the filtrations
\[
{\mathcal E}|_{x_i\times q}\,=\,l^{(i)}_0\,\supset\, l^{(i)}_1\,\supset\,\cdots\,\supset\,
l^{(i)}_{r-1}\,\supset \,l^{(i)}_r\,=\,0 \qquad (1\,\leq\, i\,\leq\, n).
\]
Then there is a universal family of filtrations
$\ell$ so that $({\mathcal E},\,\ell)$
becomes a flat family of quasi-parabolic bundles
on $X\times Y_Q$ over $Y_Q$.
Let $Z_Q$ be the reduced closed subscheme of $Y_Q$
consisting of the points $y$ such that
$\dim\End\left( ({\mathcal E} ,\,\ell) |_{X\times y}\right)\,\geq\, 2$.

We want to prove that the dimension of $Z_Q$ is at most $(r^2-1)(g-1)+nr(r-1)/2-2$.
Recall that $\dim N'\,=\,-g+\sum_{s=1}^m \sum_{i=1}^{\gamma_s} \big( (r^{(i)}_s)^2(g-1)+1\big)$.
Since there are exact sequences
\[
 \Ext^1(\overline{E}^{(i)}_s ,\, E^{(j-1)}_s)
\,\longrightarrow\, \Ext^1(\overline{E}^{(i)}_s ,\, E^{(j)}_s)
 \longrightarrow \Ext^1(\overline{E}^{(i)}_s , \overline{E}^{(j)}_s)
\]
for $1\,\leq\, j\,<\,i$, the dimension of
$\mathbb{P}_* (\Ext^1(\overline{E}^{(i)}_s ,\, E^{(i-1)}_s))$ is at most
$
-1+\sum_{j<i} \dim \Ext^1(\overline{E}^{(i)}_s,\,\overline{E}^{(j)}_s)
$.
Furthermore, the Riemann--Roch theorem implies that
$\dim \Ext^1(\overline{E}^{(i)}_s,\,\overline{E}^{(j)}_s)\,=\,
r^{(i)}_s r^{(j)}_s(g-1)$,
because $\overline{E}^{(i)}_s$ and $\overline{E}^{(j)}_s$ are stable vector bundles
of the same slope and $E^{(i)}_s\,\not\cong\, E^{(j)}_s$.
Therefore, the dimension of the fibers of $P^{(i+1)}_s\,=\,
\mathbb{P}_*\Ext^1_{X\times_{N'} P^{(i)}_s/P^{(i)}_s} 
\big( \overline{\mathcal E}^{(i+1)}_s ,\, {\mathcal E}^{(i)}_s \big)$
over $P^{(i)}_s$ is at most $-1+\sum_{j=1}^{i-1} r^{(i)}_s r^{(j)}_s (g-1)$,
which implies that the dimension of the fibers of
$P_s\,=\,P^{(\gamma_s)}_s$ over $N'$ is at most
\begin{equation}\label{partial sum of fiber dimension}
 \sum_{i=2}^{\gamma_s} \Big( -1
 +\sum_{j=1}^{i-1} r^{(i)}_s r^{(j)}_s (g-1) \Big)
\, =\, 1-\gamma_s+ \sum_{1\leq j<i\leq \gamma_s} r^{(i)}_s r^{(j)}_s(g-1).
\end{equation}
Since the extensions in
\eqref{fundamental piece of extensions of semistable sheaves}
do not split,
we can see --- by an argument similar to the above --- that
the dimension of the fibers of $Q$ over $P_1\times_{N'}\times\cdots\times_{N'} P_m$
is at most
\begin{equation}\label{second sum of fiber dimension}
\sum_{s=2}^m \Big( -1+\sum_{t=1}^{s-1} 
\dim\Ext^1(\overline{E}_t, \overline{E}_s) \Big)
\,=\,
1-m+\sum_{1\leq s<t\leq m} \dim \Ext^1(\overline{E}_t, \overline{E}_s).
\end{equation}
By Lemma \ref{key inequality coming from Riemann--Roch},
we have the inequality
$$\dim \Ext^1_X(\overline{E}_t,\,\overline{E}_s)\,\leq\, 
\max\{ r_sr_t(2g-3),\, r_s r_t g-1\}\,\leq\, 2r_s r_t(g-1)-1.$$
Using the equality
$r_s\,=\,r^{(1)}_s+\cdots+r^{(\gamma_s)}_s$ we get the following:
\begin{align}\label{main inequality in higher genus case}
\begin{split}
 \dim Q
 \:&\:\leq\:
 -g+\sum_{s=1}^m \sum_{i=1}^{\gamma_s} \big( (r^{(i)}_s)^2(g-1)+1\big)
 +\sum_{s=1}^m \Big(1-\gamma_s
 +\sum_{1\leq j<i \leq \gamma_s} r^{(i)}_sr^{(j)}_s(g-1) \Big)
 \\
 &\qquad
 +1-m+\sum_{1\leq s<t\leq m} (2r_s r_t (g-1)-1)
 \\
 &\leq\,\,
 -g+m+\sum_{s=1}^m (r^{(1)}_s+\cdots+r^{(\gamma_s)}_s)^2(g-1)
 -\sum_{s=1}^m (\gamma_s-1)
 \\
 &\qquad +1-m+\sum_{1\leq s<t\leq m} 2r_s r_t (g-1) -\frac{m(m-1)}{2}
 \\
 &=\,\,
 (r^2-1)(g-1)-\frac{m(m-1)}{2}
 -\sum_{s=1}^m (\gamma_s-1).
\end{split}
\end{align}
Taking into account the condition $m\,\geq\, 2$, we have $\dim Q\,\leq\, (r^2-1)(g-1)-2$,
because we avoid the case where $m\,=\,2$ and $\gamma_1\,=\,\gamma_2\,=\,1$.
Since the dimension of the fibers of $Y_Q$ over $Q$
is $nr(r-1)/2$, and $Z_Q$ is contained in $Y_Q$,
we have $\dim Z_Q\,\leq\, \dim Q+nr(r-1)/2\,\leq\, (r^2-1)(g-1)+nr(r-1)/2-2$.

Consider the case where one of the extensions
\eqref{fundamental piece of extensions of stable bundles} and
\eqref{fundamental piece of extensions of semistable sheaves} splits,
while excluding the case of $m\,=\,2$ and $\gamma_1\,=\,\gamma_2\,=\,1$ again.
In this case, we replace
$P^{(i+1)}_s\,=\,\mathbb{P}_*\Ext^1_{X\times P^{(i)}_s/P^{(i)}_s}
(\overline{\mathcal E}^{(i+1)}_s,\, {\mathcal E}^{(i)}_s )$
with $P^{(i+1)}_s\,=\,P^{(i)}_s$
or replace
$Q_{s+1}\,=\,\mathbb{P}_*\Ext^1_{X\times Q_s/Q_s}(\overline{\mathcal E}_{s+1},\,{\mathcal E}_s)$
with $Q_{s+1}\,=\,Q_s$ (depending on which extension splits).
So the replacement of the estimation of \eqref{second sum of fiber dimension}
does not affect the calculation in \eqref{main inequality in higher genus case}.
Thus the inequality
$\dim Q\,\leq\, (r^2-1)(g-1)-2$ still holds, and we get that
$\dim Z_Q\,\leq\, (r^2-1)(g-1)+nr(r-1)/2-2$.

Next consider the case where $\overline{E}^{(i)}_s\,\cong\, \overline{E}^{(j)}_s$
for some $i\,\neq\, j$.
In the calculation of \eqref{partial sum of fiber dimension},
we should replace
$\dim\Ext^1_X(\overline{E}^{(i)}_s ,\, \overline{E}^{(j)}_s)\,=\,r^{(i)}_s r^{(j)}_s(g-1)$
with $\dim\Ext^1_X(\overline{E}^{(i)}_s , \,\overline{E}^{(j)}_s)
\,=\,r^{(i)}_s r^{(j)}_s(g-1)+1$ in the term related to the above pair $(i,\, j)$.
However, we replace the condition
$\overline{E}^{(i)}_s\,\not\cong\, \overline{E}^{(j)}_s$
with the condition $\overline{E}^{(i)}_s\,\cong\, \overline{E}^{(j)}_s$
in the definition of $N'$.
So the calculation of \eqref{main inequality in higher genus case}
is still valid and we get the inequality $\dim Z_Q\,\leq\, (r^2-1)(g-1)+nr(r-1)-2$.

Consider now the remaining case
where $m\,=\,2$ and $\gamma_1\,=\,\gamma_2\,=\,1$. In this case,
$Q$ is a parameter space of the extensions
\[
 0\,\longrightarrow\, E_1 \,\longrightarrow\, E\, \longrightarrow\, E_2\, \longrightarrow\, 0,
\]
where $E_1,\, E_2$ are stable vector bundles such that $\mu(E_1)\,=\,\mu_1\,>\,\mu_2\,=\,\mu(E_2)$.
In the calculation of \eqref{main inequality in higher genus case},
we have $\dim Q\,\leq\, (r^2-1)(g-1)-1$ in this case.
So we have $\dim Z_Q\,\leq\, (r^2-1)(g-1)+nr(r-1)/2-1$.
Note that an automorphism $\textbf{g}$ of $E$ makes the diagram
\[
 \begin{CD}
 0 @>>> E_1 @>>> E @>>> E_2 @>>> 0 \\
 & & @V \textbf{g}_1 VV @V \textbf{g} VV @V \textbf{g}_2 VV \\
 0 @>>> E_1 @>>> E @>>> E_2 @>>> 0
 \end{CD}
\]
commutative and we have
$\textbf{g}_1\,=\, c_1\mathrm{id}_{E_1}$ and $\textbf{g}_2\,=\,c_2\mathrm{id}_{E_2}$
for some $c_1,\, c_2\,\in\, k^{\times}$.

Consider the case where $\Hom(E_2,\,E_1)\,=\,0$ for generic members
$(E_1,\, E_2)$ of $N'$. In that case, the dimension of the locus $\Hom(E_2,\,E_1)\,\neq\, 0$
in $Z_Q$ is at most $(r^2-1)(g-1)+nr(r-1)-2$.
For a general member $(E_1,\, E_2)$ of $N'$,
the automorphisms $\textbf{g}$ of $E$ is given by
$(c_1,\, c_2)\,\in\, k^{\times}\times k^{\times}$
satisfying the conditions $\textbf{g}_1\,=\,c_1\mathrm{id}_{E_1}$ and
$\textbf{g}_2\,=\, c_2\mathrm{id}_{E_2}$.
Let $v\,=\, v_1e_1+\cdots+v_re_r$
be a generator of $l^{(1)}_{r-1}$ with respect to a chosen basis
$e_1,\,\cdots,\, e_r$ of $E|_{x_1}$
such that $e_1,\,\cdots,\, e_{r_1}$ generates $E_1|_{x_1}$.
Applying the automorphisms of $E$ of the above form,
we can normalize $l^{(1)}_{r-1}$ so that
a generator $v\,=\,v_1e_1+\cdots+v_re_r$ of $l^{(1)}_{r-1}$ 
satisfies $v_i\,=\,v_{r_1+j}$ or $v_iv_{r_1+j}\,=\,0$ for some $i,\,j$ with
$1\,\leq\, i\,\leq\, r_1$ and $1\,\leq \,j\,\leq\, r_2$.
The reduced subscheme of $Z_Q$ defined by this condition
is of dimension at most $(r^2-1)(g-1)+nr(r-1)-2$.

Consider the case where
$\Hom(E_2,\,E_1)\,\neq\, 0$ for generic members $(E_1,\, E_2)$ of $N'$.
Then there are automorphisms of $E$ of the form
$c\cdot \mathrm{id}_E+h$ with $0\,\neq\, h\,\in\,\Hom(E_2,\, E_1)$.
After replacing $(E_1,\,E_2)$ with
$(E_1\otimes {\mathcal L}^{\otimes r_2},\, E_2\otimes{\mathcal L}^{\otimes -r_1})$,
for a generic member ${\mathcal L}\,\in\, \Pic^0_X$,
we may assume that $h|_{x_1}\,\neq\, 0$, because the locus of $N'$ satisfying $\Hom(E_2,\,
E_1(-x_1))\,=\,\Hom(E_2,\,E_1)$ is of dimension less than $\dim N'$.
After applying the automorphisms of $E$,
we may normalize a generator $v\,=\,v_1e_1+\cdots+v_re_r$
of $l^{(1)}_{r-1}$ such that $v_i\,=\,0$ for some $1\,\leq\, i\,\leq\, r_1$
or $l^{(1)}_{r-1}\,\subset \,E_1|_{x_1}$.
The locus of $Z_Q$ defined by this condition is of dimension
at most $(r^2-1)(g-1)+nr(r-1)-2$.

The disjoint union of all of the $Z_Q$'s in the above arguments
and the flat family of quasi-parabolic bundles 
given by $({\mathcal E},\,\ell)$'s
satisfy the assertion of the proposition.
\end{proof}

\begin{Prop}\label{dimension of non-simple parabolic bundles when genus=1}
Let $X$ be an elliptic curve over $k$, and let $L$ be a line bundle of degree $d$ on $X$.
Assume that one of the following holds:
\begin{itemize}
\item $n\,\geq\, 3$ and $r\,\geq\, 2$, or
\item $n\,=\,2$ and $r\,\geq\, 3$.
\end{itemize}
Then there exists a scheme $Z$ of finite type over $k$
and a flat family $(\widetilde{E},\,\widetilde{\boldsymbol{l}})$ of quasi-parabolic bundles
on $X\times Z$ over $Z$ such that
\begin{itemize}
\item
$\dim Z \, \leq \, r(r-1)n/2-2$,
\item
$\dim \End((\widetilde{E},\,\widetilde{\boldsymbol{l}})|_{X\times \{z\}})\,\geq\, 2$ for any $z\,\in\, Z$
\end{itemize}
and each member of the complement
$\left| {\mathcal N}^{n_0\text{\rm -reg}}_{\mathrm{par}}(L) \right|
\setminus \left| {\mathcal N}^{n_0\text{\rm -reg}}_{\mathrm{par}}(L)^{\circ} \right|$
is isomorphic to $(\widetilde{E},\,\widetilde{\boldsymbol{l}})|_{X\times\{z\}}$
for some point $z\,\in\, Z$.
\end{Prop}

\begin{proof}
As in the proof of Proposition \ref{dimension of non-simple parabolic bundles in higher genus case},
we may assume that $r$ and $d$ are coprime.
Take a member 
$(E,\,\boldsymbol{l})\,\in\, \left| {\mathcal N}^{n_0\text{\rm -reg}}_{\mathrm{par}}(L) \right|
\setminus \left| {\mathcal N}^{n_0\text{\rm -reg}}_{\mathrm{par}}(L)^{\circ} \right|$.
Since $\dim\End(E,\,\boldsymbol{l})\,\geq\, 2$, it follows that $\dim\End(E)\,\geq\, 2$.
As $r$ and $d$ are coprime, the vector bundle $E$ is not semistable. Let
\[
 0\,\subset\, E_1\,\subset\, E_2\,\subset\,\cdots\,\subset\, E_m\,=\,E
\]
be the Harder--Narasimhan filtration of $E$; note that $m\,\geq\, 2$ because $E$ is not semistable.
Setting $\overline{E}_1\,=\, E_1$ and $\overline{E}_s\,:=\, E_s/E_{s-1}$ for
$2\,\leq\, s \,\leq\, m$, the slopes $\mu_t\,:=\,\mu(\overline{E}_t)$ satisfy the inequalities
\begin{equation}\label{slopes of HN on elliptic curve}
\mu_1\,>\,\mu_2\,>\, \cdots\,>\, \mu_m,
\end{equation}
and we get extensions
\begin{equation}\label{extension of semistable sheaves on elliptic curve}
0\,\longrightarrow\, E_s \,\longrightarrow\, E_{s+1} \,\longrightarrow \,\overline{E}_{s+1}
\,\longrightarrow \,0.
\end{equation}
Note that we have $\Ext^1(\overline{E}_t,\, \overline{E}_s)\,\cong\,
\Hom(\overline{E}_s ,\, \overline{E}_t)^{\vee}\,=\,0$
for $s\,<\,t$, because $\mu(\overline{E}_t)\,<\,\mu(\overline{E}_s)$
and $\overline{E}_s,\, \overline{E}_t$ are semistable.
Hence it follows that $\Ext^1(\overline{E}_{s+1},\,E_s)\,=\, 0$.
So the extension \eqref{extension of semistable sheaves on elliptic curve}
must split, and we have a decomposition
\[
 E\,\cong\, \bigoplus_{s=1}^m \overline{E}_s.
\]
Let $\{\overline{F}^{(i)}_s\}_{i=1,\ldots,\gamma_s}$ be the set of
stable bundles arising in the direct summands of $\mathrm{gr}(\overline{E}_s)$.
Fix an index $i\, \in\, \{1,\cdots,\, \gamma_s\}$.
Let $G_i\,\subset \,\overline{E}_s$ be a maximal subbundle
satisfying the condition $\Hom(G_i,\, \overline{F}^{(i)}_s)\,=\,0$.
Then we have $\Hom(\overline{F}^{(j)}_s,\,\overline{E}_s/G_i)\,=\,0$ for any $j\,\neq\, i$,
because otherwise the pullback of 
$\overline{F}^{(j)}_s\,\subset\, \overline{E}_s/G_i$ by the surjection
$\overline{E}_s\,\longrightarrow\, \overline{E}_s/G_i$
contradicts the maximality of $G_i$.

Taking account that $\overline{E}_s$ is semistable,
we can see that $\overline{E}_s/G_i$ is semistable of slope $\mu(\overline{E}_s)$
and $\mathrm{gr}\big(\overline{E}_s/G_i\big)\,\cong\, \big(\overline{F}^{(i)}_s\big)^{\oplus u}$
for some positive integer $u$.
So we have
$\Ext^1(\overline{E}_s/G_i,\,G_i)\,\cong\,
\Hom(G,\,\overline{E}_s/G_i)^{\vee}\,=\,0$
and the extension $0\,\longrightarrow\, G_i \,\longrightarrow\, \overline{E}_s
\,\longrightarrow \,\overline{E}_s/G_i \,\longrightarrow\, 0$
must split. Repeating the same argument to $G_i$, we finally get a decomposition
\[
 \overline{E}_s\ \cong\ \bigoplus_{i=1}^{\gamma_s} F^{(i)}_s,
\]
where $F^{(i)}_s$ is a semistable bundle satisfying the condition
$\mathrm{gr}\big(F^{(i)}_s\big) \,\cong\, \big(\overline{F}^{(i)}_s\big)^{\oplus u}$
for a positive integer $u$.
Note that $\mu(F^{(i)}_s)\,=\,\mu_s$ for any $i$ and
these satisfy the inequalities in \eqref{slopes of HN on elliptic curve}.
We may further assume that
$\overline{F}^{(i)}_s\,\not\cong\, \overline{F}^{(j)}_s$ for $i\,\neq\, j$.
Note that we have
\begin{equation}\label{determinant bundle condition for g=1}
 \bigotimes_{s=1}^m \bigotimes_{i=1}^{\gamma_s} \det F^{(i)}_s\ \cong\ L.
\end{equation}
The moduli space of stable bundles parameterizing $\overline{F}^{(i)}_s$
is isomorphic to $\Pic^0_X\,\cong\, X$ for all $i,\,s$.
Since we have $\dim\Ext^1(\overline{F}^{(i)}_s,\,\overline{F}^{(i)}_s)\,=\,1$
for a stable vector bundle $\overline{F}^{(i)}_s$,
a successive non-split extension of $\overline{F}^{(i)}_s$
is unique up to an isomorphism.
So, once $\overline{F}^{(i)}_s$ is given,
then the extensions $F^{(i)}_s$ of $\overline{F}^{(i)}_s$'s
are parameterized by a finite set.
Taking into account the relation \eqref{determinant bundle condition for g=1},
the underlying vector bundles $E$ of $(E,\,\boldsymbol{l})$
can be parameterized by a scheme $W$ of finite type over $\Spec k$
whose dimension is $-1+\sum_{s=1}^m \gamma_s$.

Let
\[ 
 Y \,\longrightarrow\, W
\]
be the flag bundle
parameterizing the quasi-parabolic structures on the vector bundles $E$
corresponding to the points of $W$.
There is a universal family of quasi-parabolic bundles
$(\widetilde{E},\,\widetilde{\boldsymbol{l}})$ on $X\times Y$ over $Y$.
Since each fiber of $Y$ over $W$
is of dimension $nr(r-1)/2$, we have 
\[
\dim Y
\,=\,\dim W+ \frac {nr(r-1)} {2}
\:=\:\sum_{s=1}^m\gamma_s + \frac {nr(r-1)} {2} -1.
\]
Write $(E,\,\boldsymbol{l})\,:=\,(\widetilde{E},\,\widetilde{\boldsymbol{l}})|_{X\times y}$ 
for each point $y\,\in\, Y$.

\noindent
\textbf{Case A.}\,
Consider the case where the number of components in the decomposition
$E\, =\, \bigoplus_{s,i}F^{(i)}_s$ is at least three.
Choose a basis
$e^{(i)}_{x,s,1},\,\cdots,\, e^{(i)}_{x,s,r_{i,s}}$
of $F^{(i)}_s|_x$ at each point $x\,\in \,D$ for $1\,\leq\, s\,\leq\, m$ and $1\,\leq\, i\,\leq\,\gamma_s$.
Let
\[
\sum_{s=1}^m\sum_{i=1}^{\gamma_s}\sum_{p=1}^{r_{i,s}} v^{(i)}_{x_2,s,p} e^{(i)}_{x_2,s,p}
\]
be a generator of $l^{(2)}_{r-1}$, and let
\[
\sum_{s=1}^m \sum_{i=1}^{\gamma_s}\sum_{p=1}^{r_{i,s}} w^{(i)}_{x_2,s,p} e^{(i)}_{x_2,s,p}
\]
be a representative of a generator of $l^{(2)}_{r-2}/l^{(2)}_{r-1}$. The group $\Aut E$
of automorphisms of $E$ consists of the invertible elements of the ring of endomorphisms of $E$:
\[
 \End E\, =\,
 \bigg(\bigoplus_{s,i} \End (F^{(i)}_s)\bigg)
 \bigoplus \bigg(\bigoplus_{(s,i)\neq (t,j)}
 \Hom(F^{(i)}_s,\,F^{(j)}_t)\bigg).
\]
By the assumption, we can choose $F^{(i)}_s,\,F^{(i')}_{s'},\,F^{(j)}_t$ and $F^{(j')}_{t'}$
whose indices satisfy
$s'\,<\,s$, $t'\,<\,t$ and
$((s,\,i),\,(s',\,i'))\,\neq\,((t,\,j),\,(t',\,j'))$.
So $\Aut E$ contains
the three types of automorphisms:
\[
\prod_{s,i} k^*\mathrm{id}_{F^{(i)}_s},\, \quad
\mathrm{id}_E+\Hom(F^{(i)}_s,\,F^{(i')}_{s'}), \, \quad 
\mathrm{id}_E+\Hom(F^{(j)}_t,\,F^{(j')}_{t'}).
\]
Note that the restriction maps
\begin{align}
 & \Hom\big(F^{(i)}_s,\,F^{(i')}_{s'}\big)\,\longrightarrow \,
 \Hom\left(F^{(i)}_s\big|_{x_2},\, F^{(i')}_{s'}\big|_{x_2}\right),
\label{equation: evaluation of first homomorphim}
 \\
 & \Hom\big(F^{(j)}_t,\,F^{(j')}_{t'}\big)\,\longrightarrow\,
 \Hom\left(F^{(j)}_t\big|_{x_2},\,F^{(j')}_{t'}\big|_{x_2}\right)
 \label{equation: evaluation of second homomorphism}
\end{align}
are not zero for generic choices of
$F^{(i)}_s,\,F^{(i')}_{s'},\,F^{(j)}_t$ and $F^{(j')}_{t'}$.

If $F^{(i)}_s\,\neq\, F^{(j)}_t$ and $v^{(i)}_{x_2,s,p}\,\neq\, 0$ for some $p$,
then we may normalize a representative of a generator of $l^{(2)}_{r-2}/l^{(2)}_{r-1}$
such that $w^{(i)}_{x_2,s,p}\,=\,0$. Applying the action of
$\mathrm{id}_E+\Hom(F^{(i)}_s,\,F^{(i')}_{s'})$ and
$\mathrm{id}_E+\Hom(F^{(j)}_t,\,F^{(j')}_{t'})$,
we may ensure that
$v^{(i')}_{x_2,s',p}\,=\,0$ for some $p$
and $w^{(j')}_{x_2,t',q} w^{(j)}_{x_2,t,q'}\,=\,0$ for some $q,\, q'$.
The Zariski closed subset $Y'$ defined by this condition
is of dimension $\dim Y-2\,=\,\sum_{s=1}^m\gamma_s+nr(r-1)/2-3$.

Assume that $F^{(i)}_s\,\neq\, F^{(j)}_t$ and $v^{(i)}_{x_2,s,p}\,=\,0$ for all $p$.
If in addition the condition $v^{(i)}_{x_2,s',p'}\,=\,0$ holds for all $p'$,
then such a locus is of dimension at most
$\dim Y-2\,=\,\sum_{s=1}^m\gamma_s+nr(r-1)/2-3$.
So assume that $v^{(i)}_{x_2,s',p'}\,\neq\, 0$ for some $p'$.
Then we can normalize a representative of a generator of
$l^{(2)}_{r-2}/l^{(2)}_{r-1}$ such that $w^{(i)}_{x_2,s',p'}\,=\,0$. Applying the action of
$\mathrm{id}_E+\Hom(F^{(j)}_t,\,F^{(j')}_{t'})$,
we may have $w^{(j')}_{x_2,t',q} w^{(j)}_{x_2,t,q'}\,=\,0$
for some $q,\,q'$. The Zariski closed subset of $Y'$ defined by the condition
$v^{(i)}_{x_2,s,p}\,=\,0$ for all $p$ and $w^{(j')}_{x_2,t',q} w^{(j)}_{x_2,t,q'}\,=\,0$
for some $q,\, q'$ is of dimension at most
$\dim Y-2\,=\,\sum_{s=1}^m\gamma_s+nr(r-1)/2-3$.

Assume that $F^{(i)}_s\,=\,F^{(j)}_t$.
Then we have $F^{(i')}_{s'}\,\neq\, F^{(j')}_{t'}$ by the choices of
$(s,\, i),\, (s',\, i'),\, (t,\, j),\, (t',\, j')$.
If $v^{(i)}_{x_2,s,p}\,\neq\, 0$ for some $p$,
then applying automorphisms in 
$\mathrm{id}_E+\Hom(F^{(i)}_s,\,F^{(i')}_{s'})$
and 
$\mathrm{id}_E+\Hom(F^{(j)}_t,\,F^{(j')}_{t'})$,
we may ensure that $v^{(i')}_{x_2,s',p'}\,=\,v^{(i)}_{x_2,t',q'}\,=\,0$
for some $p',\,q'$.
The Zariski closed subset of $Y$
defined by this condition
is of dimension at most
$\dim Y-2\,=\,\sum_{s=1}^m\gamma_s+nr(r-1)/2-3$.
Assume that $v^{(i)}_{x_2,s,p}\,=\,0$ for all $p$, while $F^{(i)}_s\,=\,F^{(j)}_t$
is still assumed.
If in addition we have $v^{(i')}_{x_2,s',p'}\,=\,0$ for all $p'$,
then such a locus in $Y$ is of dimension at most
$\dim Y-2\,=\,\sum_{s=1}^m\gamma_s+nr(r-1)/2-3$.
So assume that $v^{(i')}_{x_2,s',p'}\,\neq\, 0$ for some $p'$.
Then we may normalize a representative of a generator of
$l^{(2)}_{r-2}/l^{(2)}_{r-1}$ so that the condition $w^{(i')}_{x_2,s',p'}\,=\,0$ holds.
Applying an automorphism in
$\mathrm{id}_E+\Hom(F^{(j)}_t,\,F^{(j')}_{t'})$,
we may have
$w^{(j')}_{x_2,t',q'} w^{(j)}_{x_2,t,q} \,=\,0$
for some $q,\,q'$. The locus of $Y$ defined by 
$v^{(i')}_{x_2,s',p'}\,=\, 0$ for all $p'$
and $w^{(j')}_{x_2,t',q'} w^{(j)}_{x_2,t,q} \,=\, 0$ for some $q,\, q'$
is of dimension at most
$\dim Y-2\,=\,\sum_{s=1}^m\gamma_s+nr(r-1)/2-3$.

Therefore, in all cases we can get a disjoint union $Y'$ of locally closed subsets
of $Y$ and a flat family of quasi-parabolic bundles $(\widetilde{E},\,\widetilde{\boldsymbol{l}})$
on $X\times Y'$ over $Y'$ such that
$\dim Y'\,\leq \,\sum_{s=1}^m\gamma_s+nr(r-1)/2-3$
and every member of 
$\left| {\mathcal N}^{n_0\text{-reg}}_{\mathrm{par}}(L) \right|
\setminus \left| {\mathcal N}^{n_0\text{-reg}}_{\mathrm{par}}(L)^{\circ} \right|$
can be transformed by the actions of
$\mathrm{id}_E+\Hom(F^{(i)}_s,\,F^{(i')}_{s'})$ and
$\mathrm{id}_E+\Hom(F^{(j)}_t,\,F^{(j')}_{t'})$ to a quasi-parabolic bundle
$(\widetilde{E},\,\widetilde{\boldsymbol{l}})|_{X\times y}$
for some $y\,\in\, Y'$.

Using the action of the group $\prod_{s,i} k^{\times} \mathrm{id}_{F^{(i)}_s}$
on a generator $v^{(1)}_{x_1,1,1}e^{(1)}_{x_1,1,1}+\cdots 
+v^{(m)}_{x_1,\gamma_m,r_{m,\gamma_m}}e^{(m)}_{x_1,\gamma_m,r_{m,\gamma_m}}$
of $l^{(1)}_{r-1}$ we may have
$(v^{(s)}_{x_1,i,p}-1)v^{(s)}_{x_1,i,p}\,=\,0$
for $1\,\leq\, s\,\leq\, m$, $1\,\leq\, i \,\leq\, \gamma_s$ and any $p$.
The Zariski closed subset $Z$ of $Y'$ defined by this condition satisfies the following:
$\dim Z\, =\, \dim Y'-(-1+\sum_{s=1}^m\gamma_s)
\leq nr(r-1)/2-2$.

\noindent
\textbf{Case B.}\,
Consider the case where $E\,=\,F_1\oplus F_2$
with $\mu(F_1)\,>\,\mu(F_2)$, $r_i\,=\,\rank F_i$
and each $F_i$ is a successive extension of one stable vector bundle.
In this case, we have $m\,=\,2$ and $\gamma_1\,=\,\gamma_2\,=\,1$. So we have $\dim W\,=\, 1$
and $\dim Y \,=\, 1+nr(r-1)/2$. Since $\mu(F_1)\,>\,\mu(F_2)$
and $F_1,\, F_2$ are semistable, it follows that $\Hom(F_2,\,F_1)\,=\,0$.
So we have $\dim\Hom(F_2,\,F_1)\,=\,\deg(F_2^{\vee}\otimes F_1)\,>\,0$
by the Riemann--Roch theorem, and
\begin{equation}\label{equation: three cases of Hom(F_2,F_1(-x))}
 \dim \Hom(F_2,\,F_1(-x))
 \,=\,
 \begin{cases}
\deg (F_2^{\vee}\otimes F_1(-x)) & \text{if $\mu(F_2)\,<\, \mu(F_1(-x))$}
\\
 \deg(F_2^{\vee}\otimes F_1(-x))+\dim\Hom(F_1(-x),F_2)
 & \text{if $\mu(F_2)\,=\,\mu(F_1(-x))$}
 \\
 0 & \text{if $\mu(F_2)\,>\,\mu(F_1(-x))$}
 \end{cases}
\end{equation}
for a point $x$ of $X$. In the case where $\mu(F_2)\,=\,\mu(F_1(-x))$,
we have either $\dim\Hom(F_2,\,F_1(-x_1))\,=\,0$ or
$\dim\Hom(F_2,\,F_1(-x_2))\,=\,0$ because $x_1\,\neq\, x_2$.
So, in all cases of \eqref{equation: three cases of Hom(F_2,F_1(-x))},
at least one of the maps
\[
 \Hom(F_2,\,F_1)\,\longrightarrow\, \Hom(F_2|_{x_1},\,F_1|_{x_1}),
 \quad\quad
 \Hom(F_2,\,F_1)\,\longrightarrow\, \Hom(F_2|_{x_2},\,F_1|_{x_2})
\]
is not zero. Choose a basis $e_{x_i,1},\,\cdots,\,e_{x_i,r_1}$ of $F_1|_{x_i}$
and a basis $e'_{x_i,1},\,\cdots,\, e'_{x_i,r_2}$ of $F_2|_{x_i}$. Take a generator
$v_{x_i,1}e_{x_i,1}+\cdots+v_{x_i,r_1}e_{x_i,r_1}+v'_{x_i,1}e'_{x_i,1}+\cdots+v'_{x_i,r_2}e'_{x_i,r_2}$
of $l^{(i)}_{r-1}$.
Applying the action of $1_E+\Hom(F_2,\, F_1)$,
we may have $v_{x_2,q}\,=\,0$ for some $q$, or $v'_{x_2,q'}\,=\,0$ for all $q'$.
Moreover, applying the action of
$k^{\times}\mathrm{id}_{F_1}\times k^{\times}\mathrm{id}_{F_2}$, we may have
$v_{x_1,p}\,=\,v'_{x_1,p'}$ for some $p,\,p'$, or $v_{x_1,p}v'_{x_1,p'}\,=\,0$
for some $p,\, p'$.
Let $Y'$ be a disjoint union of subvarieties of $Y$
where the following two conditions hod: $v'_{x_2,q'}v_{x_2,q}\,=\,0$ for some $q,\, q'$
and $(v_{x_1,p}-v'_{x_1,p'})v_{x_1,p}v'_{x_1,p'}\,=\,0$
for some $p,\, p'$. Then we have
\[
\dim Y'\ \leq\ r(r-1)n/2-1
\]
and every member of the complement
$\left| {\mathcal N}^{n_0\text{\rm -reg}}_{\mathrm{par}}(L) \right|
\setminus \left| {\mathcal N}^{n_0\text{\rm -reg}}_{\mathrm{par}}(L)^{\circ} \right|$
can be transformed by the actions of
$k^{\times}\cdot 1_{F_1}\times k^{\times}\cdot 1_{F_2}$
and $1_E+\Hom(F_2,\, F_1)$ to a quasi-parabolic bundle
$(\widetilde{E},\,\widetilde{\boldsymbol{l}})|_{X\times y}$
for some $y\,\in\, Y'$. Let $Y''$ be the Zariski closed subset of $Y'$
defined by
\[
 Y''\ :=\
 \left\{ y\,\in\, Y' \,\, \middle|\,\,\,
 \dim\End\big( (\widetilde{E},\widetilde{\boldsymbol{l}})|_{X\times y} \big)\,\geq\, 2\right\}.
\]
For each point $y\,\in\, Y''$,
write $(E,\,\boldsymbol{l})\,:=\, (\widetilde{E},\,\widetilde{\boldsymbol{l}})|_{X\times y}$.
Set
\[
 H\ :=\ \left\{ g\,\in\,\Aut E \,\,\middle|\,\,\,
 \text{$g|_{x_i}(l^{(i)}_{r-1})\,=\,l^{(i)}_{r-1}$ for $i\,=\,1,\,2$} \right\}.
\]
Then $H$ contains non-scalar automorphisms.

\noindent
\textbf{Case B-I.}\,
Consider the case where
$H\ \not\subset\ k^{\times}\mathrm{id}_E+\Hom(F_2,F_1)$.
Take $\textbf{g}\,\in\, H\setminus\big( k^{\times}\mathrm{id}_E+\Hom(F_2,F_1)\big)$.
Then we can write
\[
\textbf{g}\ =\ 
 \begin{pmatrix}
 c_1\mathrm{Id}_{F_1} & b \\
 0 & c_2\mathrm{Id}_{F_2}
\end{pmatrix},
\]
where $c_1,\, c_2\,\in\, k^{\times}$, $b\,\in\, \Hom(F_2,\,F_1)$
and $c_1\,\neq\, c_2$.

\noindent
\textbf{Case B-I-(a).}\,
Consider the case where $n\,\geq\, 3$.
Since $\textbf{g}|_{x_3}$ has distinct eigenvalues $c_1,\, c_2$,
the condition that $\textbf{g}|_{x_3}$ preserves $l^{(3)}_{r-1}$
implies that $\dim Y''\,\leq\, \dim Y-1\leq nr(r-1)/2-2$.

\noindent
\textbf{Case B-I-(b).}\,
Consider the case where $r\,\geq\, 3$.
In this case, we have either
$r_1\,=\,\rank F_1\,\geq\, 2$ or $r_2\,=\,\rank F_2\,\geq\, 2$.

\noindent
\textbf{Case B-I-(b)-(i).}\,
Consider the case of $r_2\,\geq\, 2$.
If $l^{(i)}_{r-1}\,\subset\, F_1|_{x_i}$ for $i\,=\,1$ or $i\,=\,2$,
then $Y'$ can be replaced by the locus satisfying this condition
and we get that $\dim Y'\,\leq\, nr(r-1)/2-2$.
So we may assume that $l^{(i)}_{r-1}\,\not\subset\, F_1|_{x_i}$
for $i\,=\,1$ or $i\,=\,2$.
Then we have
$\textbf{g}|_{l^{(i)}_{r-1}}\,=\,c_2\mathrm{id}_{l^{(i)}_{r-1}}$
and $\textbf{g}$ induces a linear map
$\overline{\textbf{g}}\,\colon\, E|_{x_i}/l^{(i)}_{r-1} \,\longrightarrow\, E|_{x_i}/l^{(i)}_{r-1}$.
Since the eigenvalues of $\overline{\textbf{g}}$ are $c_1,\, c_2$ and $l^{(i)}_{r-2}/l^{(i)}_{r-1}$
is preserved by $\overline{\textbf{g}}$,
it follows that $\dim Y''\,\leq\, \dim Y'-1\,\leq\, nr(r-1)/2-2$.

\noindent
\textbf{Case B-I-(b)-(ii).}\,
Consider the case where $r_1\,\geq\, 2$. If $l^{(i)}_{r-1}\,\subset\, F_2|_{x_i}$,
then such a locus in $Y'$ is of dimension at most
$nr(r-1)/2-2$. So we may assume that $l^{(i)}_{r-1}\,\not\subset\, F_2|_{x_i}$.
Since the induced map
$\overline{\textbf{g}}\,\colon\, E|_{x_i}/l^{(i)}_{r-1}\, \longrightarrow\, E|_{x_i}/l^{(i)}_{r-1}$
has distinct eigenvalues $c_1,\, c_2$
and $l^{(i)}_{r-2}/l^{(i)}_{r-1}$ is preserved by $\overline{\textbf{g}}$,
it follows that $\dim Y''\,\leq\, \dim Y'-1\,\leq\, nr(r-1)/2-2$.

\noindent
\textbf{Case B-II.}\,
Consider the case where $H$ is contained in 
$k^{\times}\mathrm{id}_E+\Hom(F_2,\,F_1)$.

\noindent
\textbf{Case B-II-(a).}\,
Assume that $n\,\geq\, 3$, in addition to 
$H\,\subset\, k^{\times} \mathrm{id}_E+\Hom(F_2,\,F_1)$.
We may assume that the composition of maps
\[
 \left\{ a\,\in\,\Hom(F_2,\,F_1) \,\big\vert\,\, \mathrm{id}_E+a\,\in\, H \right\}
 \,\hookrightarrow\,
 \Hom(F_2,\,F_1)\,\longrightarrow\,
 \Hom(F_2|_{x_3},\,F_1|_{x_3})
\]
is injective, because the non-injective locus in $Y'$ is of dimension
at most $\dim Y'-1\leq nr(r-1)/2-2$.
Choose a basis $e_{x_i,1},\, \cdots,\, e_{x_i,r_1}$ of $F_1|_{x_i}$
and a basis $e'_{x_i,1},\, \cdots,\, e'_{x_i,r_2}$
of $F_2|_{x_i}$ for $i\,=\,1,\,2$. Take a generator
$v_{x_i,1}e_{x_i,1}+\cdots+v_{x_i,r_1}e_{x_i,r_1}
+v'_{x_i,1}e'_{x_i,1}+\cdots+v'_{x_i,r_2}e'_{x_i,r_2}$
of $l^{(i)}_{r-1}$.
Applying an automorphism $\mathrm{id}_E+a\,\in\, H$ 
with $a\,\in\,\Hom(F_2,\, F_1)$ satisfying the condition $a|_{x_3}\,\neq\, 0$,
we can normalize $l^{(3)}_{r-1}$ so that
the condition $v_{x_3,p}\,=\, 0$ holds for some $p$, or
the condition $v'_{x_3,p'}\,=\, 0$ holds for all $p'$.
The Zariski closed subset of $Y''$
defined by this condition is of dimension at most
$\dim Y'-1 \leq nr(r-1)/2-2$.

\noindent
\textbf{Case B-II-(b).}\,
Consider the case where $r\,\geq\, 3$
while $H\,\subset\, k^{\times}\mathrm{id}_E+\Hom(F_2,F_1)$
is again assumed.
We may assume the injectivity of the homomorphism
\[
 \Hom(F_2,\,F_1)\,\longrightarrow\,
 \Hom(F_2|_{x_2},\,F_1|_{x_2}),
\]
because it holds for a generic point of $Y'$.
Take a basis $f_1,\,f_2,\, \cdots,\, f_r$ of $E|_{x_2}$
such that $f_1$ is a generator of $l^{(2)}_{r-1}$.
If there is an element $1_E+a\,\in\, H$ such that
$a\,\in\,\Hom(F_2,\,F_1)\setminus\{0\}$
and $\mathrm{Im}(a|_{x_2})\,\not\subset\, l^{(2)}_{r-1}$,
then, after applying such an automorphism,
we can normalize a representative $a_2f_2+\cdots+a_rf_r$
of a generator of $l^{(2)}_{r-2}/l^{(2)}_{r-1}$ so that
the condition $a_p\,=\, 0$ holds for some $p\,\geq\, 2$.
Such a locus in $Y$ is of dimension at most $r(r-1)n/2-2$.
If the condition $\mathrm{Im}(a|_{x_2})\,\subset\, l^{(2)}_{r-1}$ holds
for all $a\,\in\,\Hom(F_2,F_1)$ satisfying $\mathrm{id}_E+a\,\in\, H$,
then we have $l^{(2)}_{r-1}\,=\, \mathrm{Im}(a|_{x_2})$
for such an $a$ with $a\,\neq\, 0$.
So we may replace $Y'$ with a Zariski closed subset
whose dimension is at most
$r(r-1)n/2-(r-1)\,\leq\, r(r-1)n/2-2$,
because $r\,\geq\, 3$.

Therefore, in all cases the disjoint union $Z$ of all the locally closed subsets of $Y''$
in the above argument
and the pullback of flat families
$(\widetilde{E},\,\widetilde{\boldsymbol{l}})|_{X\times Z}$
satisfy the assertion of the proposition.
\end{proof}

\begin{Prop}\label{proposition: codimension genus=0 and n>3}
Assume that $X\,=\,\mathbb{P}^1_k$, $L$ is a line bundle on $\mathbb{P}^1_k$
and one of the following two holds:
\begin{itemize}
 \item[I.]\,\, $n\,\geq\, 5$ and $r\,\geq\, 2$, 
 \item[II.]\,\, $n\,=\,4$ and $r\,\geq\, 3$.
\end{itemize}
Then there exists a scheme $Z$ of finite type over $\Spec k$
and a flat family $(\widetilde{E},\,\widetilde{\boldsymbol{l}})$ of quasi-parabolic bundles
on $\mathbb{P}^1\times Z$ over $Z$ such that
\begin{itemize}
\item
$\dim Z \,\leq\, -r^2+r(r-1)n/2-1$, and
\item
$\dim \End\left( (\widetilde{E},\,\widetilde{\boldsymbol{l}})|_{\mathbb{P}^1\times\{z\}}\right)
\,\geq\, 2$ for any $z\,\in \,Z$,
\end{itemize}
and each member of the complement
$\left| {\mathcal N}^{n_0\text{\rm -reg}}_{\mathrm{par}}(L) \right|
\setminus \left| {\mathcal N}^{n_0\text{\rm -reg}}_{\mathrm{par}}(L)^{\circ} \right|$
is isomorphic to $(\widetilde{E},\,\widetilde{\boldsymbol{l}})|_{\mathbb{P}^1\times\{z\}}$
for some point $z\,\in\, Z$.
\end{Prop}

\begin{proof}
Take a quasi-parabolic bundle $(E,\,\boldsymbol{l})$.
Write
\[
 E\,=\,{\mathcal O}_{\mathbb{P}^1}(a_1)^{\oplus r_1}\oplus\cdots\oplus
 {\mathcal O}_{\mathbb{P}^1}(a_m)^{\oplus r_m}
\]
with $a_1\,<\,\cdots\,<\,a_m$. If 
$l^{(i)}_{r-1}\,\not\subset\,
{\mathcal O}_{\mathbb{P}^1}(a_2)^{\oplus r_2}\big|_{x_i}\oplus\cdots\oplus
{\mathcal O}_{\mathbb{P}^1}(a_m)^{\oplus r_m}\big|_{x_i}$
for some $i$, set
\[
\begin{aligned}
 E'\,:=\,&\ \ker\left( E\longrightarrow E|_{x_i}/l^{(i)}_{r-1} \right) \otimes{\mathcal O}_{\mathbb{P}^1}(x_i) \\
 \cong\, &\ 
 {\mathcal O}_{\mathbb{P}^1}(a_1)^{\oplus r_1-1}\oplus{\mathcal O}_{\mathbb{P}^1}(a_1+1)
 \oplus {\mathcal O}_{\mathbb{P}^1}(a_2)^{\oplus r_2}\oplus\cdots
 \oplus{\mathcal O}_{\mathbb{P}^1}(a_m)^{\oplus r_m}.
\end{aligned}
\]
Repeating such process of elementary transformations
and a twist by a line bundle,
we may replace $(E,\,\boldsymbol{l})$ with a quasi-parabolic bundle
which satisfies one of the following two conditions:
\begin{itemize}
\item[(A)]
$E\,\cong\,{\mathcal O}_{\mathbb{P}^1}^{\oplus r}$,

\item[(B)]
$E\,=\,{\mathcal O}_{\mathbb{P}^1}(a_1)^{\oplus r_1}\oplus\cdots\oplus
{\mathcal O}_{\mathbb{P}^1}(a_m)^{\oplus r_m}$ and
$l^{(i)}_{r-1}\subset
{\mathcal O}_{\mathbb{P}^1}(a_2)^{\oplus r_2}\big|_{x_i}\oplus\cdots\oplus
{\mathcal O}_{\mathbb{P}^1}(a_m)^{\oplus r_m}\big|_{x_i}$
for any $i$,
where $a_1\,<\,a_2\,<\,\cdots\,<\,a_m$.
\end{itemize}

\noindent
Case (A).\,\, Consider the case where
$E\,\cong\,{\mathcal O}_{\mathbb{P}^1}^{\oplus r}$.

We will construct a parameter space for non-simple quasi-parabolic bundles $(E,\,\boldsymbol{l})$
satisfying $E\,\cong\,{\mathcal O}_{\mathbb{P}^1}^{\oplus r}$.
Let $e_1,\,\cdots,\,e_r$ be the basis of $E$ obtained by pulling back the canonical basis of
${\mathcal O}_{\mathbb{P}^1}^{\oplus r}$ via the isomorphism
$E\,\xrightarrow{\,\sim\,}\, {\mathcal O}_{\mathbb{P}^1}^{\oplus r}$.
We may assume that $l^{(1)}_*$ is given by
$l^{(1)}_j\,=\,\langle e_1,\,\cdots,\,e_{r-j} \rangle$ for $j\,=\,0,\,\cdots,\,r-1$,
after applying an automorphism of $E$.
Applying automorphisms of $E$ fixing $l^{(1)}_*$,
we may further assume that $l^{(2)}_*$ is given by
$l^{(2)}_j\,=\,\langle e_{\sigma(1)},\,\cdots,\, e_{\sigma(r-j)} \rangle$
for $j\,=\,0,1,\,\cdots,\,r-1$,
where $\sigma$ is a permutation of $\{1,\,\cdots,\,r\}$.
Let $w_1e_1+\cdots+w_re_r$ be a generator of $l^{(3)}_{r-1}$.
Applying a diagonal automorphism of $E$, which automatically fixes $l^{(1)}_*$ and $l^{(2)}_*$,
we may assume that either $w_i\,=\,1$ holds or $w_i\,=\, 0$ holds for any $i$.
Then the group of automorphisms of $E$
fixing $l^{(1)}_*$, $l^{(2)}_*$ and $l^{(3)}_{r-1}$
becomes
\[
 B''\,=\,
 \left\{
 (a_{ij}) \,\in\, {\rm GL}_r(k) \, \middle| \:
 \begin{array}{l}
 \text{$a_{ij}\,=\,0$ and $a_{\sigma(i)\sigma(j)}\,=\,0$ for $i\,>\,j$
 and there is $c\,\in\,k^{\times}$} \\
 \text{satisfying $a_{ii}w_i+\sum_{j\neq i} a_{ij}w_j\,=\,c\,w_i$ for any $i$}
 \end{array}
 \right\}.
\]
Since $\dim\End(E,\,\boldsymbol{l})\,\geq\, 2$ by the assumption, it follows that either
there is some $(i,\,j)$ with $i\,<\,j$ and $\sigma(i)\,<\,\sigma(j)$,
or there is some $i$ satisfying the condition $w_i\,=\,0$.

\noindent
Case (A)-I.
First assume that $n\,\geq\, 5$ and $r\,\geq\, 2$.

\noindent
(A)-I-(i)
Consider the case where $w_{i_1}\,=\,0$ for some $i_1$.
Then there are automorphisms $(a_{ij})$ in $B''$ such that
$a_{i_1i_1}\,=\,c\,\in\,k^{\times}$, $a_{ii}\,=\,1$ for $i\neq i_1$
and $a_{ij}\,=\,0$ for all $i\,\neq\, j$.
Applying these automorphisms to
a generator $v\,=\,v_1e_1+\cdots+v_re_r$
of $l^{(4)}_{r-1}$ normalize $v$ so that either
\begin{itemize}
\item[$\bullet$]
$v_{i_1}\,=\,0$, or

\item[$\bullet$]
$v_{i_1}\,\neq\, 0$ and $v_{i'}\,=\,0$ for any $i'\,\neq\, i_1$, or

\item[$\bullet$]
$v_{i_1}=v_{i'}\,\neq\, 0$ for some $i'\,\neq\, i_1$
\end{itemize}
holds.
So there is a parameter space of $l^{(4)}_{r-1}$
whose dimension is at most $r-1-1\,=\,r-2$.

\noindent
(A)-I-(ii).
Consider the case where $w_i\,=\,1$ for every $i$.
Then there are some $i_1\,<\,i_2$ with $\sigma(i_1)\,<\,\sigma(i_2)$,
because $\dim B''\,\geq \,2$. So there are automorphisms
$(a_{ij})$ in $B''$ of the form
$a_{\sigma(i_1)\sigma(i_1)}\,=\,c\in k^{\times}\setminus\{1\}$,
$a_{\sigma(i_1)\sigma(i_2)}\,=\,1-c$, $a_{ii}\,=\,1$ for $i\,\neq\, \sigma(i_1)$
and $a_{ij}\,=\,0$ if $i\,\neq\, j$ and $(i,\,j)\,\neq\,(\sigma(i_1),\,\sigma(i_2))$.
Applying these automorphisms to a generator $v\,=\,v_1e_1+\cdots+v_re_r$
of $l^{(4)}_{r-1}$ normalize $v$ so that one of the following holds: 
\begin{itemize}
\item $v_{\sigma(i_2)}\,=\,0$, or

\item $v_{\sigma(i_1)}\,=\,v_{\sigma(i_2)}\,\neq\, 0$, or

\item $v_{\sigma(i_1)}\,=\,0$, $v_{\sigma(i_2)}\,\neq\, 0$. 
\end{itemize}
So there is a parameter space of $l^{(4)}_{r-1}$
whose dimension is at most $r-1-1\,=\,r-2$.

In both cases (A)-I-(i) and (A)-I-(ii), consider the group of automorphisms
\[
 B'''\,:=\,
 \left\{ g\in B'' \,\, \middle|\, \:
 \text{$g$ fixes $l^{(1)}_*,\,l^{(2)}_*,\,l^{(3)}_{r-1}$
 and $l^{(4)}_{r-1}$} \right\}.
\]
Since $(E,\,l)$ is not simple,
there is an automorphism $\textbf{g}$ in $B'''$ other than a scalar endomorphism.
Then the parameter space of $l^{(5)}_{r-1}$ preserved by $\textbf{g}$
is of dimension at most $r-1-1$. Thus there is a parameter space of
$l^{(1)}_*,\,\ldots,\,l^{(n)}_*$
whose dimension is at most
\[
 \sum_{j=1}^{r-2} j
 +2\bigg( r-2+\sum_{j=1}^{r-2}j \,\bigg)
 +(n-5)\sum_{j=1}^{r-1} j
 \,=\,-r^2+1+\frac{1}{2}r(r-1)n-2.
\]

Case (A)-II.\,
Assume that $n\,=\,4$ and $r\,\geq\, 3$.

(A)-II-(i).\,
Assume that $w_{i_1}\,=\,0$ for some $i_1$.
Then there are automorphisms
$(a_{ij})$ in $B''$ of the form
$a_{ii}\,=\,1\,\in\,k^{\times}$ for $i\,\neq\, i_1$,
$a_{i_1i_1}\,=\,c\in\,k^{\times}$
and $a_{ij}\,=\,0$ for all $i\,\neq\, j$.
For a representative $v\,=\,v_1e_1+\cdots+v_re_r\,\in\, l^{(3)}_{r-2}$ of a generator 
of $l^{(3)}_{r-2}/l^{(3)}_{r-1}$,
we may assume, after adding an element of $l^{(3)}_{r-1}$, that
$v_{i_2}\,=\,0$ for some $i_2\,\neq\, i_1$. Applying an automorphism in $B''$ of the above form,
normalize $v$ so that one of the following holds:
\begin{itemize}
\item $v_{i_1}\,=\,0$, or

\item $v_{i_1}\,\neq\, 0$ and $v_{i'}\,=\,0$ for any $i'\,\neq\, i_1$, or

\item
$v_{i_1}\,=\,v_{i_3}\,\neq\, 0$ for some $i_3\,\neq\, i_1,\,i_2$.
\end{itemize}
So there is a parameter space of $l^{(3)}_{r-2}$
whose dimension is at most $r-2-1\,=\,r-3$.

(A)-II-(ii).\,
Assume that $w_i\,=\,1$ for any $i$.
Then there are some
$i_1\,<\,i_2$ with $\sigma(i_1)\,<\,\sigma(i_2)$ because $B''\,\neq\,k^{\times}\mathrm{id}$.
Then there are automorphisms $(a_{ij})$ in $B''$ of the form
$a_{\sigma(i_1) \sigma(i_1)}\,=\,c\,\in\,k^{\times}\setminus\{1\}$,
$a_{\sigma(i_1) \sigma(i_2)}\,=\,1-c$,
$a_{ii}\,=\,1$ for $i\,\neq\, \sigma(i_1)$
and $a_{ij}\,=\,0$ if $i\,\neq\, j$ and $(i,\,j)\,\neq\,(\sigma(i_1),\,\sigma(i_2))$.
For a representative $v\,=\,v_1e_1+\cdots+v_re_r\in l^{(3)}_{r-2}$
of a generator of $l^{(3)}_{r-2}/l^{(3)}_{r-1}$,
we may assume,
after adding an element of $l^{(3)}_{r-1}$, that $v_{i'}\,=\,0$
for some $i' \,\neq\, \sigma(i_1),\,\sigma(i_2)$.
Applying an automorphism in $B''$, normalize $v$ so that
\begin{itemize}
\item either $v_{\sigma(i_1)}\,=\,v_{\sigma(i_2)}$, or
\item
$v_{\sigma(i_1)}v_{\sigma(i_2)}\,=\,0$
\end{itemize}
holds. So there is a parameter space of $l^{(3)}_{r-2}$
whose dimension is at most $r-2-1\,=\,r-3$.

In both cases (A)-II-(i) and (A)-II-(ii), consider the group of automorphisms
\[
 B'''\,:=\,
 \left\{ \textbf{g}\,\in\, B'' \,\, \middle| \,\:
 \text{$\textbf{g}$ fixes $l^{(1)}_*,\,l^{(2)}_*,\,l^{(3)}_{r-1}$ and $l^{(3)}_{r-2}$} \right\}.
\]
Since $(E,\,l)$ is not simple,
there is an automorphism $\textbf{g}$ in $B'''$ other than a scalar automorphism.
The parameter space of $l^{(4)}_{r-1}$ preserved by $\textbf{g}$
is of dimension at most $r-1-1\,=\,r-2$.
Thus there is a parameter space of $l^{(1)}_*,\,\ldots,\,l^{(n)}_*$ whose dimension is at most
\begin{align*}
 (r-3)+\sum_{j=1}^{r-3} j
 \:+(r-2)+\sum_{j=1}^{r-2} j
 \: + \,\frac{1}{2}r(r-1) (n-4)
 &=\,
 -r^2-1+\frac{1}{2}r(r-1)n.
\end{align*}

\noindent
Case (B).\,
Consider the case where
$E\,\cong\, {\mathcal O}_{\mathbb{P}^1}(a_1)^{\oplus r_1}\oplus\cdots\oplus
{\mathcal O}_{\mathbb{P}^1}(a_m)^{\oplus r_m}$
with $a_1\,<\,a_2\,<\,\cdots\,<\,a_m$ and
$l^{(i)}_{r-1}\,\subset\, 
\big( {\mathcal O}_{\mathbb{P}^1}(a_2)^{\oplus r_2}\oplus\cdots\oplus
{\mathcal O}_{\mathbb{P}^1}(a_m)^{\oplus r_m}\big)\big|_{x_i}$
for any $i$.

We choose a basis
$e^{(i)}_{1,1},\,\cdots,\,e^{(i)}_{1,r_1},\,e^{(i)}_{2,1},\,\cdots,\,e^{(i)}_{2,r_2},\,\cdots,
\,e^{(i)}_{m,1},\,\cdots,\,e^{(i)}_{m,r_m}$ of
$E\big\vert_{x_i}$ corresponding to the given decomposition
$E\big\vert_{x_i}\,=\, {\mathcal O}_{\mathbb{P}^1}(a_1)\big|_{x_i}^{\oplus r_1}\oplus\cdots\oplus
{\mathcal O}_{\mathbb{P}^1}(a_m)\big|_{x_i}^{\oplus r_m}$.
For $1\, \leq\, p\, \leq\, m$, let
\[
 \pi^{(i)}_p\,\colon\, \ E|_{x_i} \ = \
 {\mathcal O}_{\mathbb{P}^1}(a_1)^{\oplus r_1}|_{x_i}\oplus\cdots\oplus
{\mathcal O}_{\mathbb{P}^1}(a_m)^{\oplus r_m}|_{x_i}
\ \longrightarrow \
{\mathcal O}_{\mathbb{P}^1}(a_p)^{\oplus r_p}|_{x_i}
\]
be the projection to the $p$-th direct summand. So any element
$v\,\in\, E\big\vert_{x_i}$
can be uniquely written as follows:
\[
 v\ =\ v_1+\cdots+v_m \quad\quad
 \big( \text{$v_p\,\in\, {\mathcal O}_{\mathbb{P}^1}(a_p)^{\oplus r_p}\big|_{x_i}$
 \, for \, $1\,\leq\, p\,\leq\, m$ } \big).
\]

We want to choose suitable generators
$v^{(i)}_{p^{(i)}(1),\,j^{(i)}(1)},\,\cdots,\,v^{(i)}_{p^{(i)}(s),\,j^{(i)}(s)}$
of $l^{(i)}_{r-s}$. First, define a number $p^{(i)}(1)$ with
$1\,\leq\, p^{(i)}(1)\,\leq\, m$ by setting
\[
 p^{(i)}(1)\,:=\,\min \left\{
 p\,\in\,\{1,\,\ldots,\,m\}\, \, \middle|\, \:
 \begin{array}{l}
 \text{the $p$-th component $v_p$ does not vanish} \\
 \text{for a generator
 $v\,=\,v_1+\cdots+v_m$ of $l^{(i)}_{r-1}$} 
 \end{array}
 \right\}
\]
for each $1\,\leq\, i\,\leq\, n$.
So we can choose an element
$v\,=\,v_{p^{(i)}(1)}+v_{p^{(i)}(1)+1}+\cdots+v_m\in l^{(i)}_{r-1}$
with $v_{p^{(i)}(1)}\,\neq\, 0$.
Put $j^{(i)}(1)\,=\,1$ and set
$v^{(i)}_{p^{(i)}(1),j^{(i)}(1)}\,:=\,v$.
Consider the projection
\[
 \pi^{(i)}_1\times\cdots\times\pi^{(i)}_p \, \colon
 \ E|_{x_i} \, = \,
 {\mathcal O}_{\mathbb{P}^1}(a_1)^{\oplus r_1}|_{x_i}\oplus\cdots\oplus
{\mathcal O}_{\mathbb{P}^1}(a_m)^{\oplus r_m}|_{x_i}
\ \longrightarrow \
{\mathcal O}_{\mathbb{P}^1}(a_1)^{\oplus r_1}|_{x_i}\oplus\cdots\oplus
{\mathcal O}_{\mathbb{P}^1}(a_p)^{\oplus r_p}|_{x_i}
\]
for $1\,\leq\, p\,\leq\, m$.
For $2\,\leq\, s\,\leq\, r-1$,
define $p^{(i)}(s)$, $j^{(i)}(s)$ and $v^{(i)}_{p^{(i)}(s),j^{(i)}(s)}$
inductively on $s$.
For each integer $s$ with $2\,\leq\, s\,\leq\, r-1$, define
$p^{(i)}(s)$ by the condition
\[
\begin{cases}
 (\pi^{(i)}_1\times\cdots\times\pi^{(i)}_{p^{(i)}(s)}) (l^{(i)}_{r-s+1})
 \,\,\subsetneq\,\,
 (\pi^{(i)}_1\times\cdots\times\pi^{(i)}_{p^{(i)}(s)}) (l^{(i)}_{r-s})
 & \text{and}
 \\
 (\pi^{(i)}_1\times\cdots\times\pi^{(i)}_p) (l^{(i)}_{r-s+1})
 \,\, =\,\,
 (\pi^{(i)}_1\times\cdots\times\pi^{(i)}_p) (l^{(i)}_{r-s})
 & \text{for $p\,<\,p^{(i)}(s)$}.
\end{cases}
\]
Set
\[
 j^{(i)}(s)\,:=\,
 1+\#\left\{ s'\in\{1,\ldots,s-1\} \, \,\middle|\, \: p^{(i)}(s')\,=\,p^{(i)}(s) \right\}.
\]
Then we can take an element $v^{(i)}_{p^{(i)}(s), j^{(i)}(s)}$ of $l^{(i)}_{r-s}$
such that
\[
(\pi^{(i)}_1\times\cdots\times\pi^{(i)}_{p^{(i)}(s)}) (v^{(i)}_{p^{(i)}(s),j^{(i)}(s)})
\,\,\notin\,\,
 (\pi^{(i)}_1\times\cdots\times\pi^{(i)}_{p^{(i)}(s)}) (l^{(i)}_{r-s+1}).
\]
By the construction it follows that
$l^{(i)}_{r-s}$ is generated by
$v^{(i)}_{p^{(i)}(1),j^{(i)}(1)},\,\cdots,\,v^{(i)}_{p^{(i)}(s),j^{(i)}(s)}$.

Applying an automorphism of $E$ given by an element of
\[
 B\,:=\,
 \left\{
 g\,=\,(a^{pq})_{1\leq p,q\leq m}\, \, \middle|\, \:
 \begin{array}{l}
 \text{$a^{pq}\,=\,(a^{pq}_{jj'})\,\in\,\Hom({\mathcal O}_{\mathbb{P}^1}(a_q)^{\oplus r_q},
 {\mathcal O}_{\mathbb{P}^1}(a_p)^{\oplus r_p})$
 for $p\,\geq\, q$,}
 \\
 \text{$a^{pp}\,=\,(a^{pp}_{jj'})\,\in\,\Aut({\mathcal O}_{\mathbb{P}^1}(a_p)^{\oplus r_p})$
 for $1\,\leq\, p\,\leq\, m$} \\
 \text{and $a^{pq}\,=\,0$ for $p\,<\,q$}
 \end{array}
 \right\},
\]
we may assume that
$v^{(1)}_{p,j}\,=\,e^{(1)}_{p,j}$ for $1\,\leq\, p\,\leq\, m$ and $1\,\leq\, j\,\leq\, r_p$.
Note that the group of automorphisms of $E$ fixing $l^{(1)}_*$ is
\[
 B'\,=\,\left\{
 g\,=\,
 (a^{pq})\,\in\, B
 \,\, \middle| \,\:
 \begin{array}{l}
 \text{$a^{p^{(1)}(s), p^{(1)}(s')} _{j^{(1)}(s),j^{(1)}(s')} \Big|_{x_1}\,=\,0$ for $s\,>\,s'$,
 and for each $1\,\leq\, p\,\leq\, m$,} \\
 \text{$\big(a^{pp}_{jj'}\big|_{x_1}\big)\,\in\,\Aut({\mathcal O}_{\mathbb{P}^1}(a_p)^{\oplus r_p}\big|_{x_1})$
 is an upper triangular matrix}
 \end{array}
 \right\}.
\]
If $p\,>\,q$, we can always take an element
$g\,=\,(a^{pq}_{jj'})$ of $B'$ such that
$a^{pq}_{jj'}\big|_{x_2}\,\neq\, 0$.
So, after applying an automorphism in $B'$ to $l^{(2)}_*$, it may be assumed that
the condition $v^{(2)}_{p,j}\,=\,e^{(2)}_{p,\sigma_p(j)}$
holds for $1\,\leq\, j\,\leq\, r_p$,
where $\sigma_p$ is a permutation of $\{1,\,\cdots,\,r_p\}$.
The generator $v^{(3)}_{p^{(3)}(1),p^{(3)}(1)}$ of $l^{(3)}_{r-1}$ can be written as
\[
v^{(3)}_{p^{(3)}(1),j^{(3)}(1)}\,\,=\,\,w_{1,1}e^{(3)}_{1,1}+\cdots+w_{m,r_m}e^{(3)}_{m,r_m}.
\]
Consider the diagonal automorphisms $\textbf{g}\,=\,(a^{pq}_{jj'})$ of $E$
given by $a^{pq}_{jj'}\,=\,0$ for $(p,\,j)\,\neq\,(q,\,j')$
and $a^{pp}_{jj}\,\in\,k^{\times}$ for any $(p,\,j)$.
After applying such automorphisms, normalize $v^{(3)}_{1,1}$ such that
either $w_{p,j}\,=\,1$ holds or $w_{p,j}\,=\,0$ holds for any $p,\,j$.
Note that the conditions $p^{(1)}(1)\,\geq \,2$, $p^{(2)}(1)\,\geq\, 2$
and $w_{1,1}\,=\,0$ hold, because of the assumption
that $l^{(i)}_{r-1}\,\subset\, 
\big( {\mathcal O}_{\mathbb{P}^1}(a_2)^{\oplus r_2}\oplus\cdots\oplus
{\mathcal O}_{\mathbb{P}^1}(a_m)^{\oplus r_m}\big)\big|_{x_i}$
for $i\,=\,1,\,2,\,3$.

\noindent
Case (B)-I.\,
If $n\,\geq \,5$,
then we can give a parameter space of $l^{(i)}_{r-1}$ 
whose dimension is at most $r-2$ for each $4\,\leq\, i \,\leq\, n$,
because $l^{(i)}_{r-1}\,\subset \,
\big( {\mathcal O}_{\mathbb{P}^1}(a_2)^{\oplus r_2}\oplus\cdots\oplus
{\mathcal O}_{\mathbb{P}^1}(a_m)^{\oplus r_m}\big)\big|_{x_i}$.
So there is a parameter space of $(E,\,l)$
whose dimension is at most
\begin{align*}
 \sum_{j=1}^{r-2} j
 +
 (n-3)\bigg((r-2)+\sum_{j=1}^{r-2}j \bigg)
 &\,=\,
 \frac{1}{2}r(r-1)(n-2)-(r-1)-(n-3)
 \\
 &=\,
 -r^2+1+\frac{1}{2}r(r-1)n-(n-3)
 \leq -r^2+1+\frac{1}{2}r(r-1)n-2.
\end{align*}

\noindent
Case (B)-II.\,
Assume that $n\,=\,4$ and $r\,\geq\, 3$.
Recall that $v^{(3)}_{p^{(3)}(1),p^{(3)}(1)}$ is a generator of $l^{(3)}_{r-1}$ and 
we can write $v^{(3)}_{p^{(3)}(1),p^{(3)}(1)}\,=\,\sum_{p,j} w_{p,j}e^{(3)}_{p,j}$ with $w_{1,1}\,=\, 0$. 
Take a representative $u\,=\,\sum_{p,j} u_{p,j} e^{(3)}_{p,j}\,\in\, l^{(3)}_{r-2}$
of a generator of $l^{(3)}_{r-2}/l^{(3)}_{r-1}$ with the normalized condition
$u_{p^{(3)}(1),j^{(3)}(1)}\,=\,0$.
Consider the diagonal automorphisms $\textbf{g}\,=\,(a^{pq}_{jj'})$ of $E$
determined by
$a^{pq}_{jj'}\,=\,0$ for $(p,\,j)\,\neq\, (q,\,j')$,
$a^{11}_{1,1}\,=\,c\,\in\,k^{\times}$
and $a^{p,p}_{j,j}\,=\,1\,\in\,k^{\times}$ for $(p,\,j)\,\neq\, (1,\,1)$.
Then such automorphisms
preserve $l^{(1)}_*$, $l^{(2)}_*$ and $l^{(3)}_{r-1}$.
Choose an index $(q,\,j')\,\neq\, (p^{(3)}(1),\,j^{(3)}(1)),\,(1,\,1)$.
Applying automorphisms of the above form to $u$, we may assume that one of the followings holds:
\begin{itemize}
\item
$u_{1,1}\,=\,0$, or 
\item
$u_{q,j'}\,=\,0$, or
\item
$u_{1,1}\,=\,u_{q,j'}\,\neq\, 0$.
\end{itemize}
So we can give a parameter space of such $l^{(3)}_{r-2}$
of dimension at most $r-3$.
Furthermore, the parameter space of
$l^{(4)}_{r-1}$ is at most $r-2$, because of the condition
$l^{(4)}_{r-1}\,\subset\,
\big( {\mathcal O}_{\mathbb{P}^1}(a_2)^{\oplus r_2}\oplus\cdots\oplus
{\mathcal O}_{\mathbb{P}^1}(a_m)^{\oplus r_m}\big)\big|_{x_4}$.
So we can give a parameter space of $(E,\,l)$ whose dimension is at most
\begin{align*}
 &r-3+\sum_{j=1}^{r-3} j 
 +(r-2)+\sum_{j=1}^{r-2} j 
\,=\,r^2-2r-1,
\end{align*}
which is equal to $-r^2+1+r(r-1)n/2-2$
as $n\,=\,4$.
\end{proof}

\begin{Prop}\label{proposition: codimension genus zero n=3 case}
Assume that $X\,=\,\mathbb{P}^1_k$, $n\,=\,3$, $r\,\geq\, 4$ and $L$ is a line bundle
on $\mathbb{P}^1_k$. 
Then there exists a scheme $Z$ of finite type over $\Spec k$,
and a flat family $(\widetilde{E},\,\widetilde{\boldsymbol{l}})$
of quasi-parabolic bundles on $\mathbb{P}^1\times Z$ over $Z$, such that
\begin{itemize}
\item
$\dim Z \,\leq\, (r^2-3r+2)/2-2$,
\item
$\dim\End\big( (\widetilde{E},\,\widetilde{\boldsymbol{l}}) |_{\mathbb{P}^1\times\{z\}} \big)
\,\geq\, 2$ for any $z\,\in\, Z$,
\end{itemize}
and each member of
$\left| {\mathcal N}^{n_0\text{\rm -reg}}_{\mathrm{par}}(L) \right|
\setminus \left| {\mathcal N}^{n_0\text{\rm -reg}}_{\mathrm{par}}(L)^{\circ} \right|$
is isomorphic to $(\widetilde{E},\,\widetilde{\boldsymbol{l}})|_{\mathbb{P}^1\times\{z\}}$
for some point $z\,\in\, Z$.
\end{Prop}

\begin{proof}
First we fix a universal constant
\begin{equation}\label{l0}
\lambda_0\,\in\, k\setminus\{0,\,1\}.
\end{equation}
As in the proof of Proposition \ref{proposition: codimension genus=0 and n>3},
we may assume that the quasi-parabolic bundles $(E,\,\boldsymbol{l})$ satisfy one of the following
conditions:
\begin{itemize}
\item[(A)]
$E\,\cong\,{\mathcal O}_{\mathbb{P}^1}^{\oplus r}$,
or
\item[(B)]
$E\,=\,{\mathcal O}_{\mathbb{P}^1}(a_1)^{\oplus r_1}\oplus\cdots\oplus
{\mathcal O}_{\mathbb{P}^1}(a_m)^{\oplus r_m}$ and
$l^{(i)}_{r-1}\subset
{\mathcal O}_{\mathbb{P}^1}(a_2)^{\oplus r_2}\big|_{x_i}\oplus\cdots\oplus
{\mathcal O}_{\mathbb{P}^1}(a_m)^{\oplus r_m}\big|_{x_i}$
for any $i$,
where $a_1\,<\,a_2\,<\,\cdots\,<\,a_m$.
\end{itemize}

\noindent
Case (A).
First consider the case where
$E\,\cong\,{\mathcal O}_{\mathbb{P}^1}^{\oplus r}$.

As in the proof of Proposition \ref{proposition: codimension genus=0 and n>3},
we may assume that
$l^{(1)}_*$ is determined by the standard basis
$e_1,\,\cdots,\,e_r$, and
$l^{(2)}_*$ is determined by the basis
$e_{\sigma(1)},\,\cdots,\,e_{\sigma(r)}$ for a permutation $\sigma$
of $\{1,\,\cdots,\,r\}$ while $l^{(3)}_{r-1}$ is generated by
$w\,=\,w_1e_1+\ldots+w_re_r$
with $w_i\,=\,1$ or $w_i\,=\,0$ for each $i$. Consider the following three cases:
\begin{itemize}
\item[(a)]
$w_{i_1}\,=\,w_{i_2}\,=\,0$
for some $i_1\,\neq\, i_2$,

\item[(b)]
$w_{i_1}\,=\,0$ for some $i_1$ and $w_i\,=\,1$ for any $i\neq i_1$, and

\item[(c)]
$w_i\,=\,1$ for any $i$.
\end{itemize}

\noindent
Case (A)-(a).
Assume that $w_{i_1}\,=\,w_{i_2}\,=\,0$ for $i_1\,\neq\, i_2$.
Fix indices $i_3,\,i_4$ such that $w_{i_3}\,=\,1$
and that $i_4\,\neq\, i_1,\,i_2,\,i_3$.
Consider the automorphisms 
$(a_{ij})$ of $E$ satisfying
$a_{i_1i_1}\,=\,c_{i_1}\,\in\,k^{\times}$,\,
$a_{i_2i_2}\,=\,c_{i_2}\,\in\,k^{\times}$,\,
$a_{ii}\,=\,1$ for $i\,\neq\, i_1,\,i_2$
and $a_{ij}\,=\,0$ for $i\,\neq\, j$.
Then such automorphisms preserve $l^{(1)}_*$, $l^{(2)}_*$ and $l^{(3)}_{r-1}$.
Normalize a representative $v\,=\,v_1e_1+\ldots+v_re_r\in l^{(3)}_{r-2}$ of a generator of
$l^{(3)}_{r-2}/l^{(3)}_{r-1}$ such that $v_{i_3}\,=\,0$,
after adding an element of $l^{(3)}_{r-1}$.
Applying the above type of automorphisms to $v$, we can assume that one of the following statements holds:
\begin{itemize}
\item
$v_{i_1}\,=\,v_{i_2}\,=\,0$,

\item
$v_{i_1}\,=\,v_{i_4}\,=\,0$,
\item
$v_{i_2}\,=\,v_{i_4}=0$,

\item
$v_{i_1}\,=\,0$ and $v_{i_2}\,=\,v_{i_4}\,\neq\, 0$,

\item
$v_{i_2}\,=\,0$ and $v_{i_1}\,=\,v_{i_4}\,\neq\, 0$, 

\item
$v_{i_1}\,=\,v_{i_2}\,\neq\, 0$ and $v_{i_4}\,=\,0$,

\item
$v_{i_1}\,=\,v_{i_2}\,=\,v_{i_4}\,\neq\, 0$.
\end{itemize}
So there is a parameter space of $l^{(3)}_{r-2}$
whose dimension is at most $r-2-2\,=\,r-4$.
Adding the data of $l^{(3)}_{r-3},\,\cdots,\,l^{(3)}_1$,
we can get a parameter space of $(E,\,\boldsymbol{l})$
whose dimension is at most
\[
 (r-4) + \sum_{j=1}^{r-3} j
 \ = \
 \frac {r^2-3r+2} {2} - 2.
\]

\noindent
Case (A)-(b).
Assume that $w_{i_1}\,=\,0$ for some $i_1$ and $w_i\,=\,1$ for any $i\,\neq\, i_1$.
Fix an index $i_2$ other than $i_1$.
For a representative $v\,=\,v_1e_1+\cdots+v_re_r\,\in\, l^{(3)}_{r-2}$
of a generator of $l^{(3)}_{r-2}/l^{(3)}_{r-1}$,
we may assume, after adding an element of $l^{(3)}_{r-1}$, that $v_{i_2}\,=\,0$.
Then we have one of the following three cases:
\begin{enumerate}
\item[(i)]
$v_{i_1}\,=\,0$,
\item[(ii)]
$v_{i_1}\,\neq\, 0$ and $v_i\,=\,0$ for any $i\,\neq\, i_1,\,i_2$,
\item[(iii)]
$v_{i_1}\,\neq\, 0$ and
$v_{i_3}\,\neq\, 0$ for some $i_3$ with $i_3\,\neq\, i_1$ and $i_3\,\neq\, i_2$.
\end{enumerate}

\noindent
(A)-(b)-(i).
Consider the case where $v_{i_1}\,=\,0$.
Then we can give a parameter space of
$l^{(3)}_{r-2}$ whose dimension is at most $r-2-1\,=\,r-3$.
Consider the automorphisms $\textbf{g}\,=\, (a_{ij})$ of $E$
given by $a_{i_1i_1}\,=\,c\,\in\,k^{\times}$,
$a_{ii}\,=\,1$ for $i\,\neq\, i_1$
and $a_{ij}\,=\,0$ for $i\,\neq\, j$. 
Then such automorphisms preserve
$l^{(1)}_*$,\, $l^{(2)}_*$,\, $l^{(3)}_{r-1}$ and also $l^{(3)}_{r-2}$.
Since $v\,\neq\, 0$, we may choose an index $i_3$ such that
$v_{i_3}\,\neq\, 0$ and $i_3\,\neq\, i_1,\,i_2$.
For a representative $u\,=\,u_1e_1+\cdots+u_re_r\in l^{(3)}_{r-3}$ of a generator of
$l^{(3)}_{r-3}/l^{(3)}_{r-2}$,
we may assume,
after adding an element in $l^{(3)}_{r-2}$, that $u_{i_2}\,=\,u_{i_3}\,=\,0$.
After applying the above type of automorphisms to $u$,
we may assume that one of the following statements holds:
\begin{itemize}
\item $u_{i_1}\,=\,0$,
\item
$u_{i_1}\,\neq\, 0$ and $u_i\,=\,0$ for $i\,\neq\, i_1,\,i_2,\,i_3$,
\item 
$u_{i_1}\,=\,u_{i_4}\,\neq\, 0$ for some $i_4\,\neq\, i_1,\,i_2,\,i_3$.
\end{itemize}
In all these cases, there is a parameter space of $l^{(3)}_{r-3}$ whose dimension
is at most $r-3-1\,=\,r-4$.
Adding the data of $l^{(3)}_{r-4},\,\cdots,\,l^{(3)}_1$
we can get a parameter space of $(E,\,\boldsymbol{l})$
whose dimension is at most
\[
 r-3+r-4 + \sum_{j=1}^{r-4} j
 \,=\,
 \frac {r^2-3r+2} {2} - 2.
\]

\noindent
(A)-(b)-(ii).
Consider the case where
$v_{i_1}\,\neq\, 0$ and $v_i\,=\,0$ for any $i\,\neq\, i_1,\,i_2$.
Then $l^{(3)}_{r-2}$ is uniquely determined.
So the dimension of the parameter space of such $(E,\,\boldsymbol{l})$ is at most
\[
 (r-3)+\sum_{j=1}^{r-4} j
 \ =\
 \frac {r^2-3r+2} {2} - r+2
 \, \leq \,
 \frac {r^2-3r+2} {2}-2 .
\] 

\noindent
(A)-(b)-(iii).
Consider the case where $v_{i_1}\,\neq\, 0$ and
$v_{i_3}\,\neq\, 0$ for some $i_3$ with $i_3\,\neq\, i_1$ and $i_3\,\neq\, i_2$.
Recall again that we normalize a representative $v\,\in\, l^{(3)}_{r-2}$
of a generator of $l^{(3)}_{r-2}/l^{(3)}_{r-1}$ such that $v_{i_2}\,=\,0$.
Suppose that the condition $\sigma(i)\,>\,\sigma(j)$ holds for any $i\,<\,j$.
Then the automorphisms of $E$ preserving $l^{(1)}_*$ and $l^{(2)}_*$
are only diagonal automorphisms
$\textbf{g}\,=\,(a_{ij})$, which satisfy the condition $a_{ij}\,=\,0$ for $i\,\neq\, j$.
If $\textbf{g}\,=\,(a_{ij})$ preserves
$l^{(3)}_{r-1}\,=\,\langle w \rangle$ and $l^{(3)}_{r-2}\,=\,\langle w,\,v\rangle$
in addition, then we have $a_{ii}\,=\,a_{jj}$ for $i,\,j\,\neq\, i_1$ and $a_{i_1i_1}\,=\,a_{i_3i_3}$.
So $\textbf{g}$ must be a constant scalar multiplication,
which is in contradiction to the assumption that 
$\dim \Aut(E,\boldsymbol{l})\,\geq \,2$.
Thus we have the following:
\begin{itemize}
\item
There are $i_0\,<\,j_0$ satisfying $\sigma(i_0)\,<\,\sigma(j_0)$.
\end{itemize}
So consider the following cases:
\begin{itemize}
\item[($\alpha$)]
$\sigma(j_0)\,=\,i_1$ and $\sigma(i_0)\,=\,i_2$,

\item[($\beta$)]
$\sigma(j_0)\,=\,i_1$ and $\sigma(i_0)\,\neq\, i_2$,

\item[($\gamma$)]
$\sigma(j_0)\,\neq\, i_1$ and $\sigma(i_0)\,=\,i_2$,

\item[($\delta$)]
$\sigma(j_0)\,\neq\, i_1$ and $\sigma(i_0)\,\neq\, i_1,\,i_2$,

\item[($\epsilon$)]
$\sigma(i_0)\,=\,i_1$ and
$\{ j\,\in\,\{1,\,\cdots,\,r\} \, \, \big\vert\, \: j\,>\,i,\; \sigma(j)\,>\,\sigma(i) \}\,=\,\emptyset$
 for any $i\,\neq\, i_0$.
\end{itemize}
More precisely, in the remaining case other than
($\alpha$), ($\beta$), ($\gamma$) and ($\delta$), we have $\sigma(i_0)\,=\,i_1$.
If there are $i'\,\neq\, i_0$ and $j'\,>\,i'$ satisfying $\sigma(j')\,>\,\sigma(i')$,
then we replace $(i_0,\,j_0)$ with $(i',\,j')$
and reduce to the case ($\alpha$), ($\beta$), ($\gamma$) or ($\delta$).
Otherwise, we may assume ($\epsilon$).

\noindent
(A)-(b)-(iii)-($\alpha$).
Assume that $\sigma(i_0)\,=\,i_2$ and $\sigma(j_0)\,=\,i_1$.
Consider the automorphisms $\textbf{g}\,=\,(a_{ij})$ of $E$
given by $a_{i_1i_1}\,=\,c\,\in\,k^{\times}$,
$a_{\sigma(i_0)\sigma(j_0)}\,=\,a_{i_2i_1}\,=\,a\,\in\,k$,
$a_{ii}\,=\,1$ for $i\,\neq\, i_1\,=\,\sigma(j_0)$
and $a_{ij}\,=\,0$ for $i\,\neq\, j$ satisfying $(i,\,j)\,\neq\,(i_2,\,i_1)$.
Then such automorphisms preserve
$l^{(1)}_*$, $l^{(2)}_*$ and $l^{(3)}_{r-1}$.
The coefficient of $e_{i_2}$ in
\[
 gv\,\,=\,
 v_1e_1+\ldots+cv_{\sigma(j_0)}e_{i_1}+\ldots+(v_{\sigma(i_0)}+av_{\sigma(j_0)})e_{i_2}
 +\ldots+v_{i_3}e_{i_3}+\ldots+v_re_r
\]
is $v_{\sigma(i_0)}+av_{\sigma(j_0)}\,=\,av_{\sigma(j_0)}$
because of $v_{\sigma(i_0)}\,=\,v_{i_2}\,=\,0$, and hence
the normalized representative of a generator of
$l^{(3)}_{r-2}/l^{(3)}_{r-1}$ becomes
\begin{align*}
 & \
 gv-av_{\sigma(j_0)}w
 \\
 &=\,
 (v_1-av_{\sigma(j_0)})e_1+\cdots+cv_{\sigma(j_0)}e_{i_1}+\cdots+0e_{i_2}
 +(v_{i_3}-av_{\sigma(j_0)})e_{i_3}
 +\cdots+(v_r-av_{\sigma(j_0)})e_r.
\end{align*}
If we choose an index $i_4$ other than $i_1,\,i_2,\,i_3$,
we may assume that one of the following two holds:
\begin{itemize}
\item
$v_{i_1}\,=\,v_{i_3}\,=\,v_{i_4}\,\neq\, 0$,

\item $v_{i_1}\,=\,v_{i_3}\,\neq\, 0$ and $v_{i_4}\,=\,0$.
\end{itemize}
So we can give a parameter space for such $l^{(3)}_{r-2}$
whose dimension is at most $r-4$.

\noindent
(A)-(b)-(iii)-($\beta$).
Assume that
$\sigma(j_0)\,=\,i_1$ and $\sigma(i_0)\,\neq \,i_2$.
Consider the automorphisms $\textbf{g}\,=\,(a_{ij})$ of $E$ given by
$a_{i_1i_1}\,=\,c\,\in\,k^{\times}$,
$a_{\sigma(i_0)\sigma(j_0)}\,=\,a\,\in\,k$,
$a_{ii}\,=\,1$ for $i\,\neq\, i_1\,=\,\sigma(j_0)$
and $a_{ij}\,=\,0$ for $i\,\neq\, j$ satisfying $(i,\,j)\,\neq\,(\sigma(i_0),\,\sigma(j_0))$.
Then such automorphisms preserve
$l^{(1)}_*$,\, $l^{(2)}_*$ and $l^{(3)}_{r-1}\,=\,\langle w\rangle$.
Since
\[
\textbf{g}v\,=\,
 v_1e_1+\cdots+(v_{\sigma(i_0)}+av_{\sigma(j_0)})e_{\sigma(i_0)}
 +\cdots+cv_{\sigma(j_0)}e_{i_1}+\cdots+0e_{i_2}+\cdots+v_re_r,
\]
we may assume that one of the following holds:
\begin{itemize}
\item
$\sigma(i_0)\,\neq\, i_3$ and $v_{\sigma(i_0)}\,=\,v_{\sigma(j_0)}\,=\,v_{i_3}\,\neq\, 0$,

\item
$\sigma(i_0)\,=\,i_3$ and
$v_{\sigma(i_0)}\,=\,v_{\sigma(j_0)}\,=\,v_{i_4}\,\neq\, 0$
for some $i_4$ other than $i_1,\,i_3,\,i_2$,

\item
$\sigma(i_0)\,=\,i_3$, $v_{\sigma(i_0)}\,=\,v_{\sigma(j_0)}$ and
$v_i\,=\,0$ for any $i$ other than $i_1(=\sigma(j_0)),\, i_3$.
\end{itemize}
So we can give a parameter space of such $l^{(3)}_{r-2}$
whose dimension is at most $r-4$.

\noindent
(A)-(b)-(iii)-($\gamma$).
Assume that $\sigma(j_0)\,\neq\, i_1$ and $\sigma(i_0)\,=\,i_2$.
In this case, consider the automorphisms
$\textbf{g}\,=\,(a_{ij})$ of the form $a_{i_1i_1}\,=\,c\,\in\,k^{\times}$,
$a_{i_2i_2}\,=\,a\,\in\,k^{\times}\setminus\{1\}$,
$a_{i_2\sigma(j_0)}\,=\,1-a$,
$a_{ii}\,=\,1$ for $i\,\neq\, i_1,\,i_2$
and $a_{ij}\,=\,0$ for any $i\,\neq \,j$ satisfying $(i,\,j)\,\neq\, (i_2,\,\sigma(j_0))$.
Then such a $\textbf{g}$ preserves $l^{(1)}_*$, $l^{(2)}_*$ and $l^{(3)}_{r-1}$.
Since the $e_{i_2}$-coefficient of
\[
\textbf{g}v\,=\,
 v_1+\cdots+cv_{i_1}e_{i_1}+\cdots+(av_{i_2}+(1-a)v_{\sigma(j_0)})e_{i_2}
 +\cdots+v_{i_3}e_{i_3}+\cdots+v_re_r
\]
is $av_{i_2}+(1-a)v_{\sigma(j_0)}\,=\,(1-a)v_{\sigma(j_0)}$,
we should replace $\textbf{g}v$ with its normalization
\begin{align*}
 &
 \ \textbf{g}v-(1-a)v_{\sigma(j_0)}w
 \\
 &=\,
 (v_1-(1-a)v_{\sigma(j_0)})e_1+\cdots+cv_{i_1} e_{i_1}
 +\cdots+0v_{i_2}+\cdots+av_{\sigma(j_0)} e_{\sigma(j_0)}
 +\cdots+(v_r-(1-a)v_{\sigma(j_0)})e_r.
\end{align*}
Fix an index $i_4$ other than $\sigma(j_0)$,\, $i_1$ and $i_2$.
After applying an automorphism of the above form,
we may assume that one of the following holds:
\begin{itemize}
\item
$\sigma(j_0)\,\neq\, i_3$ and $v_{i_1}\,=\,v_{i_3}\,=\,v_{\sigma(j_0)}\,\neq\, 0$, 

\item
$\sigma(j_0)\,\neq\, i_3$ and $v_{i_1}\,=\,v_{i_3}\,=\,\lambda_0 v_{\sigma(j_0)}\,\neq\, 0$
(see \eqref{l0} for $\lambda_0$),

\item $\sigma(j_0)\,\neq\, i_3$, $v_{i_1}\,=\,v_{i_3}\,\neq\, 0$
and $v_{\sigma(j_0)}\,=\,0$, 

\item $\sigma(j_0)\,=\,i_3$ and $v_{i_1}\,=\,v_{\sigma(j_0)}\,=\,v_{i_4}\,\neq\, 0$,

\item $\sigma(j_0)\,=\,i_3$,\, $v_{i_1}\,=\,v_{\sigma(j_0)}$
and $v_{i_4}\,=\,0$.
\end{itemize}
So we can give a parameter space of such $l^{(3)}_{r-2}$
whose dimension is at most $r-4$.

\noindent
(A)-(b)-(iii)-($\delta$).
Assume that $\sigma(j_0)\,\neq\, i_1$ and $\sigma(i_0)\,\neq\, i_1,\,i_2$.
Consider the automorphisms $\textbf{g}\,=\,(a_{ij})$ of $E$ given by
$a_{i_1i_1}\,=\,c\,\in\,k^{\times}$,
$a_{\sigma(i_0)\sigma(i_0)}\,=\,a\,\in\,k^{\times}$,
$a_{\sigma(i_0),\, \sigma(j_0)}\,=\,1-a$,\,
$a_{ii}\,=\,1$ for $i\,\neq\, i_1,\, \sigma(i_0)$
and $a_{ij}\,=\,0$ for $i\,\neq\, j$ satisfying $(i,\,j)\,\neq \,(\sigma(i_0),\,\sigma(j_0))$.
Then such automorphisms preserve $l^{(1)}_*$, \, $l^{(2)}_*$ and $l^{(3)}_{r-1}$,
and we have
\[
\textbf{g}v\,=\, 
 v_1e_1+\cdots+cv_{i_1}e_{i_1}+\cdots+(av_{\sigma(i_0)}+(1-a)v_{\sigma(j_0)})e_{\sigma(i_0)}
 +\cdots+0e_{i_2} +\cdots+v_re_r.
\]
In the case where $v_{\sigma(i_0)}\,=\,v_{\sigma(j_0)}$, we can normalize $v$ so that
the condition $v_{i_1}\,=\,v_{i_3}$ holds.
In the case where $v_{\sigma(i_0)}\,\neq\, v_{\sigma(j_0)}$,
we can normalize $v$ so that one of the following holds:
\begin{itemize}
\item
$\sigma(i_0)\,\neq \,i_3$ and 
$v_{i_1}\,=\,v_{\sigma(i_0)}\,=\,v_{i_3}\,\neq\, 0$, 

\item $\sigma(i_0)\,=\,i_3$ and 
$v_{i_1}\,=\,v_{\sigma(i_0)}\,=\,v_{i_4}\neq 0$
for some $i_4$ other than $i_1,\,i_2,\,i_3$,

\item
$\sigma(i_0)\,=\,i_3$ and
$v_{i_1}\,=\,v_{\sigma(i_0)}\,=\,\lambda_0 v_{i_4}\,\neq\, 0$
for some $i_4$ other than $i_1,\,i_2,\,i_3$ (see \eqref{l0} for $\lambda_0$),

\item $\sigma(i_0)\,=\,i_3$,\, $v_{i_1}\,=\,v_{\sigma(i_0)}\,\neq\, 0$ and $v_{i_4}\,=\,0$ for some $i_4$ other than $i_1,\,i_2,\,i_3$.
\end{itemize}
So we can give a parameter space of such $l^{(3)}_{r-2}$ whose dimension is at most $r-4$.

\noindent
(A)-(b)-(iii)-($\epsilon$).
Assume that $\sigma(i_0)\,=\,i_1$ and that $\{ j\,>\,i \, \,| \: \sigma(j)\,>\,\sigma(i)\}\,=\,\emptyset$
for all $i\,\neq\, i_0$. Then the group of automorphisms of $E$ preserving
$l^{(1)}_*$,\, $l^{(2)}_*$ and $l^{(3)}_{r-1}\,=\,\langle w \rangle$ becomes
\[
 B''\,=\,
 \left\{ g\,=\,
 \begin{pmatrix}
 a_{11} & \cdots & a_{1r} \\
 0 & \ddots & \vdots \\
 0 & 0 & a_{rr}
 \end{pmatrix}
 \, \middle| \:
 \begin{array}{l}
 \text{$a_{ij}\,=\,0$ for $i\,\neq \,j$ satisfying $i\,\neq\, i_1$}\\
 \text{$a_{i_1i_1}\,\in\,k^{\times}$, $a_{ii}\,=\,c\,\in\,k^{\times}$ for $i\,\neq\, i_1$ and} \\
 \displaystyle \sum_{j>i_0,\, \sigma(j)>\sigma(i_0)=i_1} a_{i_1\sigma(j)} \,=\,0 
 \end{array}
 \right\}.
\]

Suppose that for any index $j_1$ satisfying $j_1\,\neq\, j_0$ and 
$i_0\,<\,j_1$, we have $\sigma(i_0)\,>\,\sigma(j_1)$.
Then any automorphism $\textbf{g}$ in $B''$ becomes diagonal. 
In other words, $\textbf{g}\,=\,(a_{ij})$ satisfies the following conditions: $a_{ij}\,=\,0$ for $i\,\neq \,j$
and there is a $c\,\in\,k^{\times}$ such that $a_{ii}\,=\,c$ for $i\,\neq\, i_1$. 
If $\textbf{g}$ further preserves $l^{(3)}_{r-2}$,
then we have $a_{i_1i_1}\,=\,a_{i_3i_3}\,=\,c$,
because $v_{i_1}\,\neq\, 0$,\, $v_{i_3}\,\neq\, 0$ and $v_{i_2}\,=\,0$.
Thus $\textbf{g}$ must be a constant scalar multiplication, which is a contradiction because
$(E,\,l)$ is not simple.

So there is an index $j_1$ with $j_1\,\neq\, j_0$ satisfying
the conditions $i_0\,<\,j_1$ and $\sigma(i_0)\,<\,\sigma(j_1)$.
Consider the automorphisms $\textbf{g}\,=\,(a_{ij})$ of the form
$a_{i_1i_1}\,=\,c'\,\in\,k^{\times}$,\, $a_{i_1\sigma(j_0)}\,=\,a\,=\,-a_{i_1\sigma(j_1)}\,\in\,k$,\,
$a_{ii}\,=\,1$ for $i\,\neq\, i_1$ and $a_{ij}\,=\,0$
for any $i\,\neq\, j$ satisfying $(i,\,j)\,\neq\,(i_1,\, \sigma(j_0)),\, (i_1,\,\sigma(j_1))$.
Recall that the representative $v\,=\,\sum_{i=1}^r v_ie_i$
of a generator of $l^{(3)}_{r-2}/l^{(3)}_{r-1}$ is normalized so that $v_{i_2}\,=\,0$. 
We further normalize a representative
$u\,=\,\sum_{i=1}^r u_i e_i\in l^{(3)}_{r-3}$
of a generator of $l^{(3)}_{r-3}/l^{(3)}_{r-2}$ so that $u_{i_2}\,=\,u_{i_3}\,=\,0$. We may assume that
$\{ \sigma(j_0),\,\sigma(j_1)\}\,\neq\,\{i_2,\,i_3\}$,
because otherwise we can replace $i_2$ or $i_3$
by another index $i_4$ other than $i_1,\,i_2,\,i_3$ according to
whether $v_{i_4}\,=\,0$ or $v_{i_4}\,\neq\, 0$.
So assume that $\sigma(j_0)\,\neq\, i_2,\,i_3$.
Applying an automorphism $\textbf{g}$ of the above form to $v$ and $u$, we have
\begin{align*}
 \textbf{g}v&\,=\,
 v_1e_1+\cdots+(c'v_{i_1}+av_{\sigma(j_0)}-av_{\sigma(j_1)})e_{i_1}+\cdots
 \\
 &\quad +0e_{i_2}+\cdots+v_{\sigma(j_0)}e_{\sigma(j_0)}
 +\cdots+v_{\sigma(j_1)}e_{\sigma(j_1)}+\cdots+v_re_r,
 \\
 \textbf{g}u&\,=\,
 u_1e_1+\cdots+(c'u_{i_1}+au_{\sigma(j_0)}-au_{\sigma(j_1)})e_{i_1}+\cdots
 \\
 &\quad
 +0e_{i_2}+\cdots+0e_{i_3}+\cdots+u_{\sigma(j_0)}e_{\sigma(j_0)}+\cdots
 +u_{\sigma(j_1)}e_{\sigma(j_1)}+\cdots+u_re_r.
\end{align*}
So we may assume that one of the following holds:
\begin{itemize}
\item
$v_{\sigma(j_0)}\,=\,v_{\sigma(j_1)}$ and $u_{\sigma(j_0)}\,=\,u_{\sigma(j_1)}$,

\item $v_{\sigma(j_0)}\,=\,v_{\sigma(j_1)}$, $u_{\sigma(j_0)}\,\neq\, u_{\sigma(j_1)}$
and $u_{i_1}\,=\,u_{\sigma(j_0)}-u_{\sigma(j_1)}$, 

\item $u_{\sigma(j_0)}\,=\,u_{\sigma(j_1)}$, $v_{\sigma(j_0)}\,\neq\, v_{\sigma(j_1)}$
and $v_{i_1}\,=\,v_{\sigma(j_0)}-v_{\sigma(j_1)}$, 

\item $v_{\sigma(j_0)}-v_{\sigma(j_1)}\,\neq\, 0$, $u_{\sigma(j_0)}-u_{\sigma(j_1)}\,\neq\, 0$,
$v_{i_1}(u_{\sigma(j_0)}-u_{\sigma(j_1)})-u_{i_1}(v_{\sigma(j_0)}-v_{\sigma(j_1)})\,=\,0$
and $v_{i_1}\,=\,v_{\sigma(j_0)}-v_{\sigma(j_1)}$,

\item $v_{\sigma(j_0)}-v_{\sigma(j_1)}\,\neq\, 0$,
$u_{\sigma(j_0)}-u_{\sigma(j_1)}\,\neq\, 0$,
$v_{i_1}(u_{\sigma(j_0)}-u_{\sigma(j_1)})-u_{i_1}(v_{\sigma(j_0)}-v_{\sigma(j_1)})\,\neq\, 0$,
$v_{i_1}\,=\,v_{\sigma(j_0)}-v_{\sigma(j_1)}$ and
$u_{i_1}\,=\,\lambda_0(u_{\sigma(j_0)}-u_{\sigma(j_1)})$.
\end{itemize}
In each of the above cases, we can give a parameter space of $l^{(3)}_{r-2}$ and $l^{(3)}_{r-3}$
whose dimension is at most $r-3+r-4\,=\,2r-7$.

In all cases of (A)-(b)-(iii), we can give a parameter space of $(E,\,\boldsymbol{l})$ whose dimension is at most
\[
 (r-2)+(r-3)-2+(r-4)+(r-5)+\cdots+1
 \,=\, \frac {r^2-3r+2} {2} - 2.
\]

\noindent
Case (A)-(c). Consider the case where $w_\ell\,=\,1$ for any $\ell$.

\noindent
Case A-(c)-(i).\, Assume further that there are
$i_1\,<\,j_1$ and $i_2\,<\,j_2$ satisfying the conditions
$\sigma(i_1)\,\neq\, \sigma(i_2)$,\, $\sigma(i_1)\,<\,\sigma(j_1)$ and $\sigma(i_2)\,<\,\sigma(j_2)$.
Let $B''$ be the group of automorphisms of $E$ preserving $l^{(1)}_*,\,l^{(2)}_*$ and $l^{(3)}_{r-1}$.
Then $B''$ contains two types of automorphisms $(a_{ij})$,\, $(b_{ij})$ such that
\begin{itemize}
\item
$a_{\sigma(i_1)\,\sigma(i_1)}\,=\,c\,\in\,k^{\times}$,\,
$a_{\sigma(i_1)\,\sigma(j_1)}\,=\,1-c$,\, $a_{ii}\,=\,1$ for $i\,\neq \,\sigma(i_1)$ and $a_{ij}\,=\,0$ for $i\,\neq\, j$
satisfying the condition $(i,\,j)\,\neq\, (\sigma(i_1),\,\sigma(j_1))$,

\item
$b_{\sigma(i_2)\,\sigma(i_2)}\,=\,c'\,\in\,k^{\times}$,\, $b_{\sigma(i_2)\,\sigma(j_2)}\,=\,1-c'$,\,
$b_{ii}\,=\,1$ for $i\,\neq\, \sigma(i_2)$ and $b_{ij}\,=\,0$ for $i\,\neq\, j$
satisfying the condition $(i,\,j)\,\neq\,(\sigma(i_2),\,\sigma(j_2))$.
\end{itemize}
For a representative $v\,=\,v_1e_1+\cdots+v_re_r$ of a generator of
$l^{(3)}_{r-2}/l^{(3)}_{r-1}$,
we may assume --- after adding an element of $l^{(3)}_{r-1}$ --- that
$v_{i'_2}\,=\,0$ for some $i'_2$ such that
$i'_2\,\neq \,\sigma(i_1),\,\sigma(j_1),\,\sigma(i_2)$.
We may further assume that $i'_2\,\neq\, \sigma(j_2)$
if $\sigma(i_1),\,\sigma(j_1),\,\sigma(i_2),\,\sigma(j_2)$ are not distinct.
Applying automorphisms of the above type,
we may assume that one of the following holds:
\begin{itemize}
\item
$j_1\,=\,j_2$, $v_{\sigma(j_1)}\,=\,0$ and 
$(v_{\sigma(i_1)}-v_{\sigma(i_2)})v_{\sigma(i_1)}v_{\sigma(i_2)}\,=\,0$,

\item
$j_1\,=\,j_2$, $v_{\sigma(j_1)}\,\neq\, 0$
and 
$(v_{\sigma(j_1)}-v_{\sigma(i_1)})(v_{\sigma(j_1)}-\lambda_0v_{\sigma(i_1)})
\,=\, (v_{\sigma(j_1)}-v_{\sigma(i_2)})(v_{\sigma(j_1)}-\lambda_0v_{\sigma(i_2)})\,=\,0$,

\item $j_1\,\neq \,j_2$ and
$v_{\sigma(j_1)}(v_{\sigma(i_1)}-v_{\sigma(j_1)})\,=\,v_{\sigma(j_2)}(v_{\sigma(i_2)}-v_{\sigma(j_2)})\,=\,0$,

\item $j_1\,\neq\, j_2$, $0\,\neq\, v_{\sigma(j_1)}\,=\,\lambda_0v_{\sigma(i_1)}$
and $(v_{\sigma(i_2)}-v_{\sigma(j_2)})v_{\sigma(j_2)}\,=\, 0$,

\item $j_1\,\neq\, j_2$, $(v_{\sigma(i_1)}-v_{\sigma(j_1)})v_{\sigma(j_1)}\,=\,0$ and 
$0\,\neq\, v_{\sigma(j_2)}\,=\,\lambda_0v_{\sigma(i_2)}$,

\item $j_1\,\neq\, j_2$,
$v_{\sigma(i_1)}\,=\,\lambda_0v_{\sigma(j_1)}\,\neq\, 0$
and $v_{\sigma(i_2)}\,=\,\lambda_0v_{\sigma(j_2)}\,\neq\, 0$.
\end{itemize}
So we can give a parameter space for $l^{(3)}_{r-2}$
whose dimension is at most $r-2-2\,=\,r-4$.
Adding the data of $l^{(3)}_{r-3},\,\cdots,\,l^{(3)}_1$,
we can get a parameter space for $(E,\,\boldsymbol{l})$ whose dimension is at most
\[
(r-4)+\sum_{j=1}^{r-3} j \ =\ \frac {r^2-3r+2} {2} - 2.
\]

\noindent
Case A-(c)-(ii).\,
Consider the rest case of A-(c).
So there is at most one $i_0$ such that there is $j\,>\,i_0$
for which $\sigma(i_0)\,<\,\sigma(j)$.
Recall that we assumed that $w_i\,=\,1$ for any $i$.
Then the automorphism group $B''$ of $E$
preserving $l^{(1)}_*$,\, $l^{(2)}_*$ and $l^{(3)}_{r-1}$
becomes
\[
 B''=
 \left\{\textbf{g}\,=\, (a_{ij}) \, \middle| \:
 \begin{array}{l}
 \bullet \text{ there is a $c\in\,k^{\times}$ such that
 $a_{ii}=c$ for $i\neq \sigma(i_0)$,} \\
 \bullet \text{ for any $i\neq \sigma(i_0)$, $a_{ij}=0$ for $i\neq j$ and} \\
 \bullet \ \displaystyle
 a_{\sigma(i_0)\sigma(i_0)}+\sum_{\genfrac{}{}{0pt}{}{j>i_0}{\sigma(j)>\sigma(i_0)}}
 a_{\sigma(i_0)\sigma(j)}=c
 \end{array}
\right\}.
\]
Since there are non-scalar automorphisms in $B''$, there is some $j_0\,>\,i_0$ for which
$\sigma(j_0)\,>\,\sigma(i_0)$. 
Choosing $i'_2$ other than $\sigma(i_0)$ and $\sigma(j_0)$,
we can normalize a representative
$v\,=\,v_1e_1+\cdots+v_re_r$ of a generator of
$l^{(3)}_{r-2}/l^{(3)}_{r-1}$ so that $v_{i'_2}\,=\,0$.
Consider the automorphisms
$\textbf{g}\,=\,(a_{ij})$ of $E$ given by
$a_{ii}\,=\,1$ for $i\,\neq \,\sigma(i_0)$,\,
$a_{\sigma(i_0)\sigma(i_0)}\,=\,c\,\in\,k^{\times}$,
$a_{\sigma(i_0)\sigma(j_0)}\,=\,1-c$
and $a_{ij}\,=\,0$ for any $i\neq j$ such that $(i,\,j)\,\neq\,(\sigma(i_0),\,\sigma(j_0))$.
Then such automorphisms preserve $l^{(1)}_*$,\, $l^{(2)}_*$ and $l^{(3)}_{r-1}$.
Choose $i'_3$ other than $\sigma(i_0),\,\sigma(j_0),\, i'_2$.
Since $v$ is sent to
\[
 v_1e_1+\cdots+(cv_{\sigma(i_0)}+(1-c)v_{\sigma(j_0)})e_{\sigma(i_0)}
 +\cdots+v_{\sigma(j_0)}e_{\sigma(j_0)}+\cdots+v_re_r
\]
by the automorphism $\textbf{g}$, we can assume that one of the following holds:
\begin{itemize}
\item[($\alpha$)]
$v_{\sigma(i_0)}\,=\,v_{\sigma(j_0)}$,

\item[($\beta$)]
$v_{\sigma(i_0)}\,\neq\, v_{\sigma(j_0)}$ and $v_{\sigma(j_0)}\,=\,v_{i'_3}$,

\item[($\gamma$)]
$v_{\sigma(i_0)}\,\neq\, v_{\sigma(j_0)}$, \,$v_{\sigma(j_0)}\,\neq\, v_{i'_3}$ and
$v_{\sigma(i_0)}\,=\,v_{i'_3}$.
\end{itemize}
If in addition we have $v_{i'_3}\,=\,0$, then we can give a parameter space for
such $l^{(3)}_{r-2}$ whose dimension is at most
$r-4$. So we assume that $v_{i'_3}\,\neq\, 0$.

\noindent
A-(c)-(ii)-($\alpha$).\, Assume that the condition
$v_{\sigma(i_0)}\,=\,v_{\sigma(j_0)}$ holds.
Recall that we are assuming that $i'_2\,\neq \,\sigma(i_0),\,\sigma(j_0)$
and $i'_3\,\neq\, \sigma(i_0),\,\sigma(j_0),\,i'_2$.
Furthermore, we are normalizing $v$ so that $v_{i'_2}\,=\,0$. Consider
the automorphisms $\textbf{g}\,=\,(a_{ij})$ given by $a_{ii}\,=\,1$ for $i\,\neq\,\sigma(i_0)$,
$a_{\sigma(i_0)\sigma(i_0)}\,=\,c\,\in\, k^{\times}$,
$a_{\sigma(i_0)\sigma(j_0)}\,=\,1-c$ and $a_{ij}\,=\,0$ for any $i\,\neq \,j$
such that $(i,\,j)\,\neq\,(\sigma(i_0),\,\sigma(j_0))$.
Then such automorphisms $g$ preserve not only $l^{(1)}_*$,\, $l^{(2)}_*$ and $l^{(3)}_{r-1}$
but also $v$. Consider a normalized representative
$u\,=\,u_1e_1+\cdots+u_re_r$ of a generator of
$l^{(3)}_{r-3}/l^{(3)}_{r-2}$ such that
$u_{i'_2}\,=\,u_{i'_3}\,=\,0$. Then $u$ is sent to
\[
 u_1e_1+\cdots+(cu_{\sigma(i_0)}+(1-c)u_{\sigma(j_0)})e_{\sigma(i_0)}
 +\cdots+u_{\sigma(j_0)}e_{\sigma(j_0)}+0u_{i'_2}+0u_{i'_3}+\cdots+u_re_r
\]
by the above automorphism $\textbf{g}$.
Replacing $u$ by some $\textbf{g}u$, we may assume that one of the following holds
\begin{itemize}
\item
$u_{\sigma(i_0)}\,=\,u_{\sigma(j_0)}$,

\item $u_{\sigma(i_0)}\,\neq\, u_{\sigma(j_0)}$ and $u_{\sigma(j_0)}\,=\,0$, 

\item $u_{\sigma(i_0)}\,\neq\, u_{\sigma(j_0)}$,\, $u_{\sigma(j_0)}\,\neq\, 0$
and $u_{\sigma(i_0)}\,=\,0$.
\end{itemize}
So we can give a parameter space for $(E,\,\boldsymbol{l})$ whose dimension is at most
\[
 (r-3)+(r-4)+\sum_{j=1}^{r-4} j\ =\ \frac {r^2-3r+2} {2} - 2.
\]

\noindent
A-(c)-(ii)-($\beta$).\, Assume that 
$v_{\sigma(i_0)}\,\neq\, v_{\sigma(j_0)}$ and
$v_{\sigma(j_0)}\,=\,v_{i'_3}$.
Recall that we are assuming that $v_{i'_3}\,\neq\, 0$.
After applying an automorphism in $B''$,
we may assume that $v_{\sigma(i_0)}\,=\,\lambda_0 v_{i'_3}$.
So we can give a parameter space for such $l^{(3)}_{r-2}$
of dimension at most $r-4$.
Then we can give a parameter space for $(E,\,\boldsymbol{l})$ whose dimension is at most
\[
 (r-4)+\sum_{j=1}^{r-3} j \, =\,\frac {r^2-3r+2} {2} - 2.
\]

\noindent
A-(c)-(ii)-($\gamma$).\, Assume that
$v_{\sigma(i_0)}\,\neq\, v_{\sigma(j_0)}$,\, $v_{\sigma(j_0)}\,\neq\, v_{i'_3}$ and
$v_{\sigma(i_0)}\,=\,v_{i'_3}$. Note that there are non-scalar automorphisms
$\textbf{g}\,=\,(a_{ij})\,\in\, B''$ preserving $l^{(3)}_{r-2}$.
Recall that there is a $c\,\in\,k^{\times}$ such that $a_{ii}\,=\,c$ for $i\,\neq\, \sigma(i_0)$.
Since $\textbf{g}v\,\in\,\langle v,\,w\rangle$, and the coefficient of $e_{i'_2}$ in
\begin{align*}
 \textbf{g}v 
 &\,=\,
 cv_1e_1+\cdots+
 \bigg( a_{\sigma(i_0)\sigma(i_0)}v_{\sigma(i_0)}+
 \sum_{\genfrac{}{}{0pt}{}{j>i_0}{\sigma(j)>\sigma(i_0)}} a_{\sigma(i_0)\sigma(j)}v_{\sigma(j)} 
 \bigg) e_{\sigma(i_0)}
 \\
 &\quad\quad
 +\cdots+cv_{\sigma(j_0)}e_{\sigma(j_0)}+\cdots+cv_{i'_3}e_{i'_3}+\cdots+cv_re_r
\end{align*}
is zero, we must have $\textbf{g}v\,=\,cv$. Comparing the coefficients of $e_{\sigma(i_0)}$, we have
\[
 a_{\sigma(i_0)\sigma(i_0)}v_{\sigma(i_0)}+
 \sum_{\genfrac{}{}{0pt}{}{j>i_0}{\sigma(j)>\sigma(i_0)}} a_{\sigma(i_0)\sigma(j)}v_{\sigma(j)} 
 \,=\,c \, v_{\sigma(i_0)}.
\] 
Combining with the equality
$\displaystyle
a_{\sigma(i_0)\sigma(i_0)}+\sum_{\genfrac{}{}{0pt}{}{j>i_0}{\sigma(j)>\sigma(i_0)}}
a_{\sigma(i_0)\sigma(j)}\,=\,c$,
it follows that
\[
\sum_{\genfrac{}{}{0pt}{}{j>i_0}{\sigma(j)>\sigma(i_0)}}
a_{\sigma(i_0)\sigma(j)}(v_{\sigma(j)}-v_{\sigma(i_0)})\,=\,0.
\]
So there is $j_1\,\neq\, j_0$ for which $j_1\,>\,i_0$ and $\sigma(j_1)\,>\,\sigma(i_0)$.

If $v$ satisfies the condition $v_{\sigma(j_0)}\,=\,v_{\sigma(j_1)}$,
then, taking into account the condition ($\gamma$),
we can give a parameter space for such $l^{(3)}_{r-2}$ of dimension 
at most $r-2-2\,=\,r-4$.

So we assume that $v_{\sigma(j_0)}\,\neq\, v_{\sigma(j_1)}$.
For $a\,\in\,k^{\times}$,
we can construct an automorphism $\textbf{g}\,=\,(a'_{ij})\,\in\, B''$ satisfying the
following conditions:
\begin{itemize}
\item
$a'_{ii}\,=\,1$ for $i\,\neq\,\sigma(i_0)$,

\item
$a'_{ij}\,=\,0$ for any $i\,\neq\, j$ for which 
$(i,\,j)\,\neq\,(\sigma(i_0),\,\sigma(j_0)),\,\, (\sigma(i_0),\,\sigma(j_1))$,

\item $a'_{\sigma(i_0)\sigma(i_0)}\,=\,a$,\, $a'_{\sigma(i_0)\sigma(j_0)}\,=\,b\,\in\,k$,\,
$a'_{\sigma(i_0)\sigma(j_1)}\,=\,b'\,\in\,k$,

\item $a+b+b'\,=\,1$ and $av_{\sigma(i_0)}+bv_{\sigma(j_0)}+b'v_{\sigma(j_1)}\,=\,v_{\sigma(i_0)}$.
\end{itemize}
Indeed, if $a\,\in\,k^{\times}$ is given,
then $b'$ is determined by the equality
\[
 (a-1)(v_{\sigma(i_0)}-v_{\sigma(j_0)})\ =\ b'(v_{\sigma(j_0)}-v_{\sigma(j_1)})
\]
and $b$ is determined by the condition $b\,=\,1-a-b'$.
Recall that we normalized $v_{i'_2}\,=\,0$
and we are assuming that $i'_3\,\neq\, i'_2,\,\sigma(i_0),\,\sigma(j_0)$.
Consider a representative $u\,=\,u_1e_1+\cdots+u_re_r\in l^{(3)}_{r-3}$
of a generator of $l^{(3)}_{r-3}/l^{(3)}_{r-2}$
satisfying the normalized condition $u_{i'_2}\,=\,u_{i'_3}\,=\,0$.
Then the $e_{i'_2}$-coefficient and the $e_{i'_3}$-coefficient of
\[
 \textbf{g}u\,=\,
 u_1e_1+\cdots+(au_{\sigma(i_0)}+bu_{\sigma(j_0)}+b'u_{\sigma(j_1)})e_{\sigma(i_0)}
 +\cdots+u_{\sigma(j_0)}+\cdots+0e_{i'_2}+0e_{i'_3}+\cdots+u_re_r
\]
vanish, and the $e_{\sigma(i_0)}$-coefficient of $\textbf{g}u$ is
\begin{align*}
 au_{\sigma(i_0)}+bu_{\sigma(j_0)}+b'u_{\sigma(j_1)}
 &\,=\,
 au_{\sigma(i_0)}+(1-a-b')u_{\sigma(j_0)}+b'u_{\sigma(j_1)}
 \\
 &=\,
 a(u_{\sigma(i_0)}-u_{\sigma(j_0)})
 +u_{\sigma(j_0)}
 -(a-1)\frac{v_{\sigma(i_0)}-v_{\sigma(j_0)}} {v_{\sigma(j_0)}-v_{\sigma(j_1)}} 
 (u_{\sigma(j_0)}-u_{\sigma(j_1)}).
\end{align*}
If $\displaystyle
u_{\sigma(i_0)}-u_{\sigma(j_0)}\,\neq\,
\frac{v_{\sigma(i_0)}-v_{\sigma(j_0)}} {v_{\sigma(j_0)}-v_{\sigma(j_1)}}
(u_{\sigma(j_0)}-u_{\sigma(j_1)})$,
then we can normalize $u$ so that $u_{\sigma(i_0)}\,=\,u_{\sigma(j_0)}$.
So we can give a parameter space for such $l^{(3)}_{r-3}$ whose dimension is at most
$r-4$. If the equality $\displaystyle
u_{\sigma(i_0)}-u_{\sigma(j_0)}\,=\,
\frac{v_{\sigma(i_0)}-v_{\sigma(j_0)}} {v_{\sigma(j_0)}-v_{\sigma(j_1)}}
(u_{\sigma(j_0)}-u_{\sigma(j_1)})$ holds,
then we can give a parameter space for such $l^{(3)}_{r-3}$ whose dimension is at most $r-4$.
Therefore, in all cases we can give a parameter space of $(E,\,\boldsymbol{l})$ whose dimension is at most
\[
 (r-3)+(r-4)+\sum_{j=1}^{r-4} j \ =\ \frac {r^2-3r+2} {2} - 2.
\]

\noindent
Case (B). \ Consider the case where
$E\,=\,{\mathcal O}_{\mathbb{P}^1}(a_1)^{\oplus r_1}\oplus\cdots\oplus
{\mathcal O}_{\mathbb{P}^1}(a_m)^{\oplus r_m}$ 
with $a_1\,<\,a_2\,<\,\cdots\,<\,a_m$ and 
$l^{(i)}_{r-1}\,\subset\,
{\mathcal O}_{\mathbb{P}^1}(a_2)^{\oplus r_2}\big|_{x_i}\oplus\cdots\oplus
{\mathcal O}_{\mathbb{P}^1}(a_m)^{\oplus r_m}\big|_{x_i}$
for $1\,\leq\, i\,\leq \,n$.

As in the proof of Proposition \ref{proposition: codimension genus=0 and n>3},
we choose a basis $e^{(i)}_{j,1},\,\cdots,\,e^{(i)}_{j,r_j}$ of
${\mathcal O}_{\mathbb{P}^1}(a_j)^{\oplus r_j}\big|_{t_i}$
for each $i,\,j$ and we choose suitable generators
$v^{(i)}_{p^{(i)}(1),j^{(i)}(1)},\,\cdots,\, v^{(i)}_{p^{(i)}(s),j^{(i)}(s)}$
of $l^{(i)}_{r-s}$. We may further assume that
$l^{(i)}_{r-s}$ is generated by 
$e^{(i)}_{p^{(i)}(1),j^{(i)}(1)},\,\cdots,\,e^{(i)}_{p^{(i)}(s),j^{(i)}(s)}$
for $i\,=\,1,\, 2$.
Since diagonal automorphisms $\textbf{g}\,=\,(a^{p,q}_{j,j'})$ of $E$ given by
$a^{pp}_{jj}\,\in\,k^{\times}$ and $a^{pq}_{jj'}\,=\,0$ for $(p,\,j)\,\neq\,(q,\,j')$
preserve $l^{(1)}_*$ and $l^{(2)}_*$, we can normalize the generator
\[
 v^{(3)}_{p^{(3)}(1),j^{(3)}(1)}
\ =\ w_{1,1}e^{(3)}_{1,1}+\ldots +w_{m,r_m}e^{(3)}_{m,r_m}
\]
of $l^{(3)}_{r-1}$ so that either $w_{p,j}\,=\,1$ or $w_{p,j}\,=\,0$ for any $p,\,j$.
Note that $w_{1,j}\,=\,0$ for $1\,\leq\, j\,\leq \,r_1$
by the assumption of Case (B). There are the following two possible cases:
\begin{itemize}
\item[(i)] $r_1\,\geq\, 2$,

\item[(ii)] $r_1\,=\,1$.
\end{itemize}

\noindent
Case (B)-(i).\, Assume that the condition $r_1\,\geq\, 2$ holds.
After adding an element of $l^{(3)}_{r-1}$, we can assume that
a representative $v\,=\,v_{1,1}e^{(3)}_{1,1}+\cdots+v_{m,r_m}e^{(3)}_{m,r_m}$
of a generator of $l^{(3)}_{r-2}/l^{(3)}_{r-1}$
satisfies the condition $v_{p^{(3)}(1),j^{(3)}(1)}\,=\,0$.
Consider the automorphisms $\textbf{g}\,=\,(a^{p,q}_{j,j'})$ of $E$ given by
$a^{1,1}_{j,j}\,=\,c_j\,\in\,k^{\times}$ for $1\,\leq\, j\,\leq\, r_1$,\,
$a^{p,p}_{j,j}\,=\,c'\,\in\,k^{\times}$ for $p\,\geq\, 2$ and $1\,\leq\, j\,\leq\, r_p$,\,
$a^{p,q}_{j,j'}\,=\,0$ for $(p,\,j)\,\neq\,(q,\,j')$.
Then such automorphisms preserve $l^{(1)}_*$,\, $l^{(2)}_*$ and $l^{(3)}_{r-1}$.
Since
\[
\textbf{g}v\,=\,
 c_1v_{1,1}e^{(3)}_{1,1}+\cdots+c_{r_1}v_{1,r_1}e^{(3)}_{1,r_1}
 +c'v_{2,1}e^{(3)}_{2,1}+\cdots+0e^{(3)}_{p^{(3)}(1),j^{(3)}(1)}+\cdots+c'v_{m,r_m}e^{(3)}_{m,r_m},
\]
we can assume that either $v_{1,j}\,=\,1$ or $v_{1,j}\,=\,0$ holds for any $p,\,j$. 
If $r_1\,>\,2$, then the parameter space for such generators of
$l^{(3)}_{r-2}/l^{(3)}_{r-1}$ is of dimension at most $r-4$.
If $r_1\,=\,2$, we may further assume that for some $(p',\,j')\,\neq\,(p^{(3)}(1),\,j^{(3)}(1))$
with $p'\,\geq\, 2$ the following holds: either $v_{p',j'}\,=\,1$ or $v_{p',j'}\,=\,0$.
So we can give a parameter space for $l^{(3)}_{r-2}$ whose dimension is at most $r-4$.
Adding the data $l^{(3)}_{r-3},\,\cdots,\,l^{(3)}_1$,
we can give a parameter space for $(E,\,\boldsymbol{l})$ whose dimension is at most
\[
 (r-4)+(r-3)+(r-4)+(r-5)+\cdots+1\ =\ \frac {r^2-3r+2} {2} - 2.
\]

\noindent
Case (B)-(ii). \, Assume that $r_1\,=\,1$.
We again take a representative $v\,=\,v_{1,1}e^{(3)}_{1,1}+\cdots+v_{m,r_m}e^{(3)}_{m,r_m}$
of a generator of $l^{(3)}_{r-2}/l^{(3)}_{r-1}$
so that $v_{p^{(3)}(1),j^{(3)}(1)}\,=\,0$. We may assume that one of the following holds:
\begin{itemize}
\item[($\alpha$)] $v_{1,1}\,=\,0$,
\item[($\beta$)] $v_{1,1}\,\neq\, 0$ and $(p^{(1)}(1),\,j^{(1)}(1))\,\neq\,(p^{(2)}(1),j^{(2)}(1))$,
\item[($\gamma$)] $v_{1,1}\,\neq\, 0$ and $(p^{(1)}(1),\,j^{(1)}(1))\,=\,(p^{(2)}(1),\,j^{(2)}(1))$.
\end{itemize}

\noindent
(B)-(ii)-($\alpha$)\, Assume that $v_{1,1}\,=\,0$ holds.
After adding an element of $l^{(3)}_{r-1}$,
we can normalize $v\,=\,v_{1,1}e^{(3)}_{1,1}+\cdots+v_{m,r_m}e^{(3)}_{m,r_m}$
so that $v_{p^{(3)}(1),j^{(3)}(1)}\,=\,0$. Consider the automorphisms
$\textbf{g}\,=\,(a^{p,q}_{j,j'})$ of $E$ given by $a^{1,1}_{11}\,=\,c_1\,\in\,k^{\times}$,\,
$a^{p,p}_{j,j}\,=\,c_2\,\in\,k^{\times}$ for $p\,\geq\, 2$
and $a^{p,q}_{j,j'}\,=\,0$ for $(p,\,j)\,\neq\,(q,\,j')$.
Then such automorphisms preserve not only $l^{(1)}_*$,\, $l^{(2)}_*$ and $l^{(3)}_{r-1}$
but also $l^{(3)}_{r-2}$. Choose $(q_1,\,j_1)$ such that $q_1\,\geq\, 2$,\, $v_{q_1,j_1}\,\neq\, 0$
and $(q_1,\,j_1)\,\neq\,(p^{(3)}(1),\,j^{(3)}(1))$.
We can normalize a representative $u\,=\,u_{1,1}e^{(3)}_{1,1}+\ldots+u_{m,r_m}e^{(3)}_{m,r_m}$
of a generator of $l^{(3)}_{r-3}/l^{(3)}_{r-2}$ by adding an element of $l^{(3)}_{r-2}$
such that $u_{p^{(3)}(1),j^{(3)}(1)}\,=\,u_{q_1,j_1}\,=\,0$.
Take an index $(q_2,\,j_2)$ other than $(1,\,1)$,\, $(p^{(3)}(1),\,j^{(3)}(1))$ and $(q_1,\,j_1)$.
Since
\begin{align*}
\textbf{g}u
 &\,=\,
 c_1u_{1,1}e^{(3)}_{1,1}+c_2u_{2,1}e^{(3)}_{2,1}+\ldots +0e^{(3)}_{p^{(3)}(1),j^{(3)}(1)} +\ldots
 \\
 &\quad\quad
+0e^{(3)}_{q_1,j_1}+\cdots+c_2u_{q_2,j_2}e^{(3)}_{q_2,j_2}+\ldots +c_2u_{m,r_m}e^{(3)}_{m,r_m},
\end{align*}
we may assume that one of the following holds:
\begin{itemize}
\item $u_{1,1}\,=\,u_{q_2,j_2}\,\neq\, 0$,

\item $u_{1,1}\,=\,0$,

\item $u_{q_2,j_2}\,=\,0$.
\end{itemize}
So we can give a parameter space for $l^{(3)}_{r-2},\,l^{(3)}_{r-3}$
whose dimension is at most $(r-3)+(r-4)$.
Adding the data $l^{(3)}_{r-4},\,\cdots,\,l^{(3)}_1$,
we can give a parameter space for $(E,\,\boldsymbol{l})$ whose dimension is at most
\[
 (r-3)+(r-4) +\sum_{j=1}^{r-4} j
\, =\,\frac {r^2-3r+2} {2} - 2.
\]

\noindent
(B)-(ii)-($\beta$).\, Assume that the conditions $v_{1,1}\,\neq\, 0$ and
$(p^{(1)}(1),\,j^{(1)}(1))\,\neq\,(p^{(2)}(1),\,j^{(2)}(1))$ hold.
After replacing the indices $i\,=\,1$ and $2$ if necessary, we may assume
that $(p^{(3)}(1),\,j^{(3)}(1))\,\neq\, (p^{(1)}(1),\,j^{(1)}(1))$.
Consider the automorphisms $\textbf{g}\,=\,(a^{p,q}_{j,j'})$ of $E$
given by
\begin{enumerate}
\item[]
$a^{1,1}_{1,1}\,=\,c_1\,\in\,k^{\times}$
and $a^{p,p}_{j,j}\,=\,c_2\,\in\,k^{\times}$ for $p\,\geq \,2$ and $1\,\leq\, j\,\leq\, r_p$,
\item[]
$a^{p^{(1)}(1),1}_{j^{(1)}(1),1}\,\in\,
\Hom({\mathcal O}_{\mathbb{P}^1}(a_1),\,{\mathcal O}_{\mathbb{P}^1}(a_{p^{(1)}(1)}))$
satisfying $a^{p^{(1)}(1),1}_{j^{(1)}(1),1}\big|_{x_2}\,=\,0$ 
and
\item[]
$a^{p,q}_{j,j'}\,=\,0$ for any $(p,\,j)\,\neq\,(q,\,j')$ satisfying 
$((p,\,j), \, (q, \, j') )\,\neq\,((p^{(1)}(1),\,j^{(1)}(1)), \, (1,\, 1))$.
\end{enumerate}
Such automorphisms preserve $l^{(1)}_*$,\, $l^{(2)}_*$ and 
$l^{(3)}_{r-1}\,=\,\big\langle v^{(3)}_{p^{(3)}(1),j^{(3)}(1)} \big\rangle$.
We can normalize the representative
$v\,=\,v_{1,1}e^{(3)}_{1,1}+\ldots +v_{m,r_m}e^{(3)}_{m,r_m}$ of
a generator of $l^{(3)}_{r-2}/l^{(3)}_{r-1}$ after adding an element of $l^{(3)}_{r-1}$
such that $v_{p^{(3)}(1),j^{(3)}(1)}\,=\,0$.
Choose $(q_1,\,j_1)$ such that $q_1\,\geq\, 2$ and 
$(q_1,\,j_1)\,\neq\, (p^{(1)}(1),\,j^{(1)}(1)),\, (p^{(3)}(1),\,j^{(3)}(1))$.
For $\textbf{g}\,=\,(a^{p,q}_{j,j'})\,\in\, B'$, we have
\begin{align*}
\textbf{g}v
 &\,=\,
 c_1v_{1,1}e^{(3)}_{1,1}+\ldots+
 \left( a^{p^{(1)}(1),1}_{j^{(1)}(1),1}v_{1,1} +c_2v_{p^{(1)}(1),j^{(1)}(1)} \right)
 e^{(3)}_{p^{(1)}(1),j^{(1)}(1)} +\ldots \\
 & \quad\quad
 +\cdots+ 0e^{(3)}_{p^{(3)}(1),j^{(3)}(1)}
 +\ldots+c_2v_{q_1,j_1}e^{(3)}_{q_1,j_1}+\ldots+c_2v_{m,r_m}e^{(3)}_{m,r_m}.
\end{align*}
So we can normalize $v$ so that one of the following statements holds:
\begin{itemize}
\item
$v_{1,1}\,=\,v_{p^{(1)}(1),j^{(1)}(1)}\,=\,v_{q_1,j_1}\,\neq\, 0$,
\item
$v_{1,1}\,=\,v_{p^{(1)}(1),j^{(1)}(1)}\,\neq\, 0$ and $v_{q_1,j_1}\,=\,0$.
\end{itemize}
Thus we can give a parameter space for $l^{(3)}_{r-2}$
whose dimension is at most $r-4$.
Adding the data $l^{(3)}_{r-2},\,\cdots,\,l^{(3)}_1$,
we can give a parameter space for $(E,\,\boldsymbol{l})$ whose dimension is at most
\[
 (r-4)+(r-3)+(r-4)+\cdots+1
\ =\ \frac {r^2-3r+2} {2} - 2.
\]

\noindent
(B)-(ii)-($\gamma$).\,
Assume that the following conditions hold: $v_{1,1}\,\neq\, 0$ and
$(p^{(1)}(1),\,j^{(1)}(1))\,=\,(p^{(2)}(1),j^{(2)}(1))$.
By the definition, $v^{(3)}_{p^{(3)}(1),j^{(3)}(1)}
\,=\,w_{2,1}e^{(3)}_{2,1}+\ldots+w_{m,r_m}e^{(3)}_{m,r_m}$
is a fixed generator of $l^{(3)}_{r-1}$,
and for any $p,\,j$, we have either $w_{p,j}\,=\,1$ or $w_{p,j}\,=\,0$.
We can choose $(q_1,\,j_1)$ such that $w_{q_1,\,j_1}\,=\,1$.
Take $(q_2,\,j_2)$ such that $q_2\,\geq\, 2$ and $(q_2,\,j_2)\,\neq\,(q_1,\,j_1),\, (p^{(1)}(1),\,j^{(1)}(1))$.
After replacing $(q_1,\,j_1)$ and $(q_2,\,j_2)$ if necessary, we can assume
that one of the following statements holds:
\begin{itemize}
\item[($\gamma$-1)]
$w_{q_2,j_2}\,=\,0$,

\item[($\gamma$-2)]
$w_{p^{(1)}(1),j^{(1)}(1)}\,=\,0$,

\item[($\gamma$-3)]
$w_{q_2,j_2}\,=\,w_{p^{(1)}(1),j^{(1)}(1)}\,=\,1$
and $w_{p,j}\,=\,0$ for any $(p,\,j)\,\neq\,(q_2,\, j_2),\, (p^{(1)}(1),\,j^{(1)}(1))$,

\item[($\gamma$-4)]
$w_{q_2,\,j_2}\,=\,w_{p^{(1)}(1),j^{(1)}(1)}\,=\,1$,\, $(q_1,\,j_1)\,\neq\,(q_2,\,j_2),\, (p^{(1)}(1),\,j^{(1)}(1))$
and $q_2\,>\,p^{(1)}(1)$,

\item[($\gamma$-5)]
$w_{q_2,j_2}\,=\,w_{p^{(1)}(1),j^{(1)}(1)}\,=\,1$,\, $(q_1,\,j_1)\,\neq\,(q_2,\, j_2),\, (p^{(1)}(1),\,j^{(1)}(1))$
and $q_2\,\leq \,p^{(1)}(1)$.
\end{itemize}

\noindent
(B)-(ii)-($\gamma$-1).\, Assume that the condition $w_{q_2,j_2}\,=\,0$ holds.
Consider the diagonal automorphisms
$\textbf{g}\,=\,(a^{p,q}_{j,j'})$ of $E$ given by
$a^{1,1}_{1,1}\,=\,c_1\,\in\,k^{\times}$,\,
$a^{q_2,q_2}_{j_2,j_2}\,=\,c_2\,\in\,k^{\times}$,\,
$a^{p,p}_{j,j}\,=\,c_3\,\in\,k^{\times}$
for $(p,\,j)\,\neq\, (1,\, 1),\, (q_2,\, j_2)$ and
$a^{p,q}_{j,j'}\,=\,0$ for $(p,\,j)\,\neq\,(q,\,j')$.
Then such automorphisms preserve
$l^{(1)}_*$,\, $l^{(2)}_*$ and $l^{(3)}_{r-1}$.
Consider a representative
$v\,=\,v_{1,1}e^{(3)}_{1,1}+\ldots+v_{m,r_m}e^{(3)}_{m,r_m}$
of a generator of $l^{(3)}_{r-2}/l^{(3)}_{r-1}$
with the normalizing condition $v_{q_1,j_1}\,=\,0$.
Applying the above type of automorphisms to $v$, we have
\[
\textbf{g}v\,=\,c_1v_{1,1}e^{(3)}_{1,1}+c_3v_{2,1}e^{(3)}_{2,1}+\cdots+0e^{(3)}_{q_1,j_1}
 +\cdots+c_2v_{q_2,j_2}e^{(3)}_{q_2,j_2}
 +\cdots+c_3v_{m,r_m}e^{(3)}_{m,r_m}.
\]
So we can normalize $v$ so that one of the following holds:
\begin{itemize}
\item
$v_{1,1}\,=\,v_{q_2,j_2}\,=\,v_{p^{(1)}(1),j^{(1)}(1)}\,\neq\, 0$,

\item $v_{1,1}\,=\,v_{q_2,j_2}\,\neq\, 0$ and $v_{p^{(1)}(1),j^{(1)}(1)}\,=\,0$,

\item $v_{1,1}\,=\,v_{p^{(1)}(1),j^{(1)}(1)}\,\neq\, 0$ and $v_{q_2,j_2}\,=\,0$,

\item $v_{q_2,j_2}\,=\,v_{p^{(1)}(1),j^{(1)}(1)}\,=\,0$.
\end{itemize}
So we can give a parameter space for $l^{(3)}_{r-2}$ whose dimension is at most $r-4$.

\noindent
(B)-(ii)-($\gamma$-2).\,
Assume that the condition $w_{p^{(1)}(1),j^{(1)}(1)}\,=\,0$ holds.
In this case, we have $(q_1,\,j_1)\,\neq\,(p^{(1)}(1),\,j^{(1)}(1))$,
because $w_{q_1,j_1}\,=\,1\,\neq\, 0$. Consider the automorphisms
$\textbf{g}\,=\,(a^{p,q}_{j,j'})$ of $E$ given by
$a^{1,1}_{1,1}\,=\,c_1\,\in\,k^{\times}$,\,
$a^{p^{(1)}(1),p^{(1)}(1)}_{j^{(1)}(1),j^{(1)}(1)}\,=\,c_2\,\in\,k^{\times}$,\,
$a^{p,p}_{j,j}\,=\,c_3\,\in\,k^{\times}$ for $(p,\,j)\,\neq\, (1,\, 1),\, (p^{(1)}(1),\,j^{(1)}(1))$
and $a^{p,q}_{j,j'}\,=\,0$ for $(p,\,j)\,\neq\,(q,\,j')$.
Then such automorphisms $\textbf{g}$ preserve
$l^{(1)}_*$,\, $l^{(2)}_*$ and $l^{(3)}_{r-1}$.
Normalize the representative
$v\,=\,v_{1,1}e^{(3)}_{1,1}+\ldots+v_{m,r_m}e^{(3)}_{m,r_m}$
of a generator of $l^{(3)}_{r-2}/l^{(3)}_{r-1}$
so that $v_{q_1,j_1}\,=\,0$.
Applying the above type of automorphism $\textbf{g}$ to $v$, we have
\[
\textbf{g}v\,=\,
 c_1v_{1,1}e^{(3)}_{1,1}+\cdots+c_2v_{p^{(1)}(1),j^{(1)}(1)}+\cdots+0e^{(3)}_{q_1,j_1}
 +\cdots+c_3v_{q_2,j_2}e^{(3)}_{q_2,j_2}+\cdots+c_3v_{m,r_m}e^{(3)}_{m,r_m}.
\]
So we can assume that one of the following holds: 
\begin{itemize}
\item $v_{1,1}\,=\,v_{p^{(1)}(1),j^{(1)}(1)}\,=\,v_{q_2,j_2}\,\neq\, 0$,

\item $v_{1,1}\,=\,v_{p^{(1)}(1),j^{(1)}(1)}\,\neq\, 0$ and $v_{q_2,j_2}\,=\,0$,

\item $v_{1,1}\,=\,v_{q_2,j_2}\,\neq\, 0$ and $v_{p^{(1)}(1),j^{(1)}(1)}\,=\,0$, 

\item $v_{p^{(1)}(1),j^{(1)}(1)}\,=\,v_{q_2,j_2}\,=\,0$.
\end{itemize}
So we can give a parameter space for $l^{(3)}_{r-2}$ whose dimension is at most $r-4$.

\noindent
(B)-(ii)-($\gamma$-3).\, Assume that $w_{q_2,j_2}\,=\,w_{p^{(1)}(1),j^{(1)}(1)}\,=\,1$ 
and $w_{p,j}\,=\,0$ for any $(p,\,j)$
other than $(q_2,\,j_2),\,(p^{(1)}(1),\,j^{(1)}(1))$.
In this case, we have $(q_1,\,j_1)\,=\,(p^{(1)}(1),\,j^{(1)}(1))$
because $w_{q_1,j_1}\,=\,1\,\neq\, 0$.
Consider the diagonal automorphisms
$\textbf{g}\,=\,(a^{p,q}_{j,j'})$ of $E$ given by
$a^{1,1}_{1,1}\,=\,c_1\,\in\,k^{\times}$,\,
$a^{q_2,q_2}_{k_2,j_2}\,=\,a^{p^{(1)}(1),p^{(1)}(1)}_{j^{(1)}(1),j^{(1)}(1)}\,=\,c_2\,\in\,k^{\times}$,\,
$a^{p,p}_{j,j}\,=\,c_3\,\in\,k^{\times}$
for $(p,\,j)\,\neq\,(1,\, 1),\, (q_2,\, j_2),\, (p^{(1)}(1),\, j^{(1)}(1))$
and $a^{p,q}_{j,j'}\,=\,0$ for $(p,\,j)\,\neq\, (q,\,j')$.
Such an automorphism $\textbf{g}$ preserves
$l^{(1)}_*$,\, $l^{(2)}_*$ and $l^{(3)}_{r-1}$.
We normalize again the representative $v\,=\,v_{1,1}e^{(3)}_{1,1}+\ldots+v_{m,r_m}e^{(3)}_{m,r_m}$
of a generator of $l^{(3)}_{r-2}/l^{(3)}_{r-1}$ such that $v_{q_1,j_1}\,=\,0$.
Further, fix an index $(q_3,\,j_3)$ other than $(1,\,1),\, (q_1,\, j_1),\, (q_2,\, j_2)$.
Applying the above type of automorphisms to $v$, we have
\[
\textbf{g}v\,=\,
 c_1v_{1,1}e^{(3)}_{1,1}+\cdots+0e^{(3)}_{q_1,j_1}+\cdots+c_2v_{q_2,j_2}e^{(3)}_{q_2,j_2}
 +\cdots+c_3v_{q_3,j_3}e^{(3)}_{q_3,j_3}+\cdots+c_3v_{m,r_m}e^{(3)}_{m,r_m}.
\]
So we may assume that one of the following holds:
\begin{itemize}
\item
$v_{1,1}\,=\,v_{q_2,j_2}\,=\,v_{q_3,j_3}\,\neq \,0$,

\item $v_{1,1}\,=\,v_{q_2,j_2}\,\neq\, 0$ and $v_{q_3,j_3}\,=\,0$,

\item $v_{1,1}\,=\,v_{q_3,j_3}\,\neq\, 0$ and $v_{q_2,j_2}\,=\,0$,

\item $v_{q_2,j_2}\,=\,v_{q_3,j_3}\,=\,0$.
\end{itemize}
Then we can give a parameter space for $l^{(3)}_{r-2}$ whose dimension is at most $r-4$.

\noindent
(B)-(ii)-($\gamma$-4).\, Assume that the following three conditions hold: $w_{q_2,j_2}\,=\,
w_{p^{(1)}(1),j^{(1)}(1)}\,=\,1$,\, $(q_1,\,j_1)\,\neq\,(q_2,\,k_2),\, (p^{(1)}(1),\,j^{(1)}(1))$
and $q_2\,>\,p^{(1)}(1)$. In this case, we have $a_1\,\leq\, a_{q_2}-2$
and we can take sections $\alpha$ of
$\Hom({\mathcal O}_{\mathbb{P}^1}(a_1),\,{\mathcal O}_{\mathbb{P}^1}(a_{q_2}))$
such that $\alpha|_{x_1}\,=\,\alpha|_{x_2}\,=\,0$ but
$\alpha|_{x_3}$ is arbitrary.
Recall that $v\,=\,v_{1,1}e^{(3)}_{1,1}+\cdots+v_{m,r_m}e^{(3)}_{m,r_m}$
gives a representative of a generator of $l^{(3)}_{r-2}/l^{(3)}_{r-1}$ with $v_{1,1}\,\neq\, 0$.
We impose the normalizing condition $v_{q_1,j_1}\,=\,0$
after adding an element of $l^{(3)}_{r-1}$ to $v$.
Consider the automorphisms $\textbf{g}\,=\,(a^{p,q}_{j,j'})$ of $E$ given by
\begin{itemize}
\item
$a^{1,1}_{1,1}\,=\,c_1\,\in\,k^{\times}$, \
$a^{p,p}_{j,j}\,=\,c_2\,\in\,k^{\times}$ for $(p,\,j)\,\neq\,(1,\,1)$,

\item
$a^{q_2,1}_{k_2,1}\,=\,\alpha\,\in\, \Hom({\mathcal O}_{\mathbb{P}^1}(a_1),\, {\mathcal O}_{\mathbb{P}^1}(a_{q_2}))$
satisfying $\alpha|_{x_1}\,=\,0$,\, $\alpha|_{x_2}\,=\,0$ and

\item $a^{p,q}_{j,j'}\,=\,0$ for any $(p,j,\,q,k)$ such that
$(p,\,j)\,\neq\,(q,\,j')$ and $(p,\,j,\,q,\, j)\,\neq\, (q_2,\,k_2,\,1,\,1)$.
\end{itemize}
Then such automorphisms preserve $l^{(1)}_*$,\, $l^{(2)}_*$ and $l^{(3)}_{r-1}$.
Applying such an automorphism,
the representative $v$ of a generator of $l^{(3)}_{r-2}/l^{(3)}_{r-1}$
is sent to
\begin{align*}
\textbf{g}v
 &\,=\,
 c_1v_{1,1}e^{(3)}_{1,1}+\cdots
 + c_2v_{p^{(1)}(1),j^{(1)}(1)}e^{(3)}_{p^{(1)}(1),j^{(1)}(1)}
 +\cdots \\
 &\quad\quad
 +0e^{(3)}_{q_1,j_1}+\cdots
 + \left( \alpha|_{x_3} v_{1,1} +c_2 v_{q_2,j_2} \right) e^{(3)}_{q_2,j_2}
 +\cdots+c_2v_{m,r_m}.
\end{align*}
So we may assume that one of following two holds:
\begin{itemize}
\item
$v_{1,1}\,=\,v_{p^{(1)}(1),j^{(1)}(1)}\,=\,v_{q_2,j_2}\,\neq\, 0$,
\item
$v_{1,1}\,=\,v_{q_2,j_2}\,\neq\, 0$ and $v_{p^{(1)}(1),j^{(1)}(1)}\,=\,0$.
\end{itemize}
So we can give a parameter space for $l^{(3)}_{r-2}$ whose dimension is at most $r-4$.

\noindent
(B)-(ii)-($\gamma$-5).\, Assume that the following three conditions hold:
$w_{q_2,j_2}\,=\,w_{p^{(1)}(1),j^{(1)}(1)}\,=\,1$,\,
$(q_1,\,j_1)\,\neq\,(q_2,\,j_2),\,(p^{(1)}(1),\,j^{(1)}(1))$
and $q_2\,\leq\, p^{(1)}(1)$.
Consider the automorphisms $\textbf{g}\,=\, (a^{p,q}_{j,j'})$ of $E$ given by
\begin{itemize}
\item $a^{1,1}_{1,1}\,=\,c_1\,\in\,k^{\times}$, \
$a^{p^{(1)}(1),p^{(1)}(1)}_{j^{(1)}(1),j^{(1)}(1)}\,=\,c_2\,\in\,k^{\times}$, \
$a^{p,p}_{j,j}\,=\,c_3\,\in\,k^{\times}$ for $(p,\,j)\,\neq\,(1,\,1),\,(p^{(1)}(1),\,j^{(1)}(1))$,

\item $a^{p^{(1)}(1),q_2}_{j^{(1)}(1),j_2}\,=\,b\,\in\, 
\Hom({\mathcal O}_{\mathbb{P}^1}(a_{q_2}),\,
{\mathcal O}_{\mathbb{P}^1}(a_{p^{(1)}(1)}))$
such that $c_2w_{p^{(1)}(1),j^{(1)}(1)}+b|_{x_3} w_{q_2,j_2}\,=\,c_3w_{p^{(1)}(1),j^{(1)}(1)}$
and

\item $a^{p,q}_{j,j'}\,=\,0$ for any $(p,\,j,\,q,\,j)$ such that
$(p,\,j)\,\neq\, (q,\,j')$ and $(p,\,j,\,q,\,j')\,\neq\, (p^{(1)}(1),\,j^{(1)}(1),\,q_2,\,j_2)$.
\end{itemize}
Note that we can always choose $b\,\in\, H^0({\mathcal O}_{\mathbb{P}^1}(a_{p^{(1)}(1)}-a_{q_2}))$
satisfying the condition that $b|_{x_3} w_{q_2,j_2}\,=\,(c_3-c_2)w_{p^{(1)}(1),j^{(1)}(1)}$
for any given $c_2,\,c_3\,\in\,k^{\times}$.
Such automorphisms preserve $l^{(1)}_*$,\, $l^{(2)}_*$ and $l^{(3)}_{r-1}$.
Applying such an automorphism,
the representative $v\,\in\, l^{(3)}_{r-2}$ of a generator of
$l^{(3)}_{r-2}/l^{(3)}_{r-1}$ is sent to
\begin{align*}
\textbf{g}v
 &\,=\,
 c_1v_{1,1}e^{(3)}_{1,1}+c_3v_{2,1}e^{(3)}_{2,1}+\cdots
 +\left( c_2v_{p^{(1)}(1),j^{(1)}(1)}+b|_{x_3} v_{q_2,j_2} \right) e^{(3)}_{p^{(1)}(1),j^{(1)}(1)}
 +\cdots 
 \\
 &\quad\quad
 +0e^{(3)}_{q_1,j_1}+\cdots+c_3 v_{q_2,j_2} e^{(3)}_{q_2,j_2}
 +\cdots+c_3v_{m,r_m}.
\end{align*}
So we can assume that one of the following holds:
\begin{itemize}
\item $v_{1,1}\,=\,v_{p^{(1)}(1),j^{(1)}(1)}\,=\,v_{q_2,j_2}\,\neq\, 0$,

\item $v_{1,1}\,=\,v_{p^{(1)}(1),j^{(1)}(1)}\,\neq\, 0$ and $v_{q_2,j_2}\,=\,0$,

\item $v_{1,1}\,=\,v_{q_2,j_2}\,\neq\, 0$ and $v_{p^{(1)}(1),j^{(1)}(1)}\,=\,0$,

\item $v_{p^{(1)}(1),j^{(1)}(1)}\,=\,v_{q_2,j_2}\,=\,0$.
\end{itemize}
Then we can give a parameter space for $l^{(3)}_{r-2}$ whose dimension is at most $r-4$.

In all cases of (B)-(ii)-($\gamma$),
by adding the data $l^{(3)}_{r-3},\,\cdots,\,l^{(3)}_1$ to the parameter space of $l^{(3)}_{r-2}$,
we can give a parameter space for $(E,\,\boldsymbol{l})$ whose dimension is at most
\[
 (r-4)+(r-3)+\sum_{j=1}^{r-4} j
 \,=\,\frac {r^2-3r+2} {2} - 2.
\]
This completes the proof.
\end{proof}

Define the open subset 
${\mathcal M}^{n_0\text{\rm -reg}}_{\mathrm{PC}}(\boldsymbol{\nu},\nabla_L)^{\circ}$
of $\mathcal{M}_{\mathrm{PC}}^{n_0\text{-\rm reg}}(\boldsymbol{\nu},\nabla_L)$
\begin{equation}\label{simple underlying parabolic bundle PC moduli}
{\mathcal M}^{n_0\text{\rm -reg}}_{\mathrm{PC}}(\boldsymbol{\nu},\nabla_L)^{\circ}
\,:=\,
\left\{ (E,\,\nabla,\,\boldsymbol{l}) \,\in\, \mathcal{M}_{\mathrm{PC}}^{n_0\text{-\rm reg}}
(\boldsymbol{\nu},\nabla_L) \,\,\middle|\,\,\,
\dim \big( \End(E,\,\boldsymbol{l}) \big) \,=\, 1 \right\}
\end{equation}
which consists of $\boldsymbol{\nu}$-parabolic connections
$(E,\,\nabla,\,\boldsymbol{l})$ with the determinant isomorphic to $(L,\nabla_L)$
such that the underlying quasi-parabolic bundle
$(E,\,\boldsymbol{l})$ is simple.

\begin{Prop}\label{proposition: codimension in de Rham moduli}
Let $X$ be a smooth projective curve of genus $g$
over an algebraically closed field $k$, and let $L$ be a line bundle on $X$.
Let $r$ and $n$ be positive integers such that $r$ is not divisible by the characteristic of $k$
and one of the following holds:
\begin{itemize}
\item $n\,\geq\, 1$ and $r\,\geq\, 2$ are arbitrary if $g\,\geq\, 2$,

\item $n\,\geq\, 2$, $r\geq 2$ and $n+r\,\geq\, 5$ if $g\,=\,1$,

\item $n\,\geq \,3$, $r\geq 2$ and $n+r\,\geq \,7$ if $g\,=\,0$.
\end{itemize}
Then the following holds:
\[
 \mathrm{codim}_{\mathcal{M}_{\mathrm{PC}}^{n_0\text{-\rm reg}}(\boldsymbol{\nu},\nabla_L)}
 \left( \mathcal{M}_{\mathrm{PC}}^{n_0\text{-\rm reg}}(\boldsymbol{\nu},\nabla_L)
 \setminus \mathcal{M}_{\mathrm{PC}}^{n_0\text{-\rm reg}}(\boldsymbol{\nu},\nabla_L)^{\circ} \right)
 \,\,\geq\,\, 2.
\]
\end{Prop}

\begin{proof}
By Proposition \ref{dimension of non-simple parabolic bundles in higher genus case},
Proposition \ref{dimension of non-simple parabolic bundles when genus=1},
Proposition \ref{proposition: codimension genus=0 and n>3} and
Proposition\ref{proposition: codimension genus zero n=3 case}
there is a scheme $Z$ of finite type over $k$
and a flat family $(\widetilde{E},\,\widetilde{\boldsymbol{l}})$
of quasi-parabolic bundles on $X\times Z$ over $Z$
such that
\begin{enumerate}
\item{} $\dim\End((\widetilde{E},\,\widetilde{\boldsymbol{l}})|_{X\times z})\,\geq\, 2$
for any point $z\,\in\, Z$,

\item{} $\dim Z\leq (r^2-1)(g-1)+nr(r-1)/2-2$, and

\item each quasi-parabolic bundle in
$\left| {\mathcal N}^{n_0\text{-reg}}_{\mathrm{PC}}(L) \right|
\setminus \left| {\mathcal N}^{n_0\text{-reg}}_{\mathrm{PC}}(L)^{\circ} \right|$
is isomorphic to $(\widetilde{E},\,\widetilde{\boldsymbol{l}})|_{X\times\{z\}}$
for some point $z\,\in\, Z$.
\end{enumerate}
We may assume that there is an isomorphism
$\varphi\,\colon\, \det(\widetilde{E})\,\xrightarrow{\,\,\sim\,\,\,}\,
L\otimes{\mathcal L}$ for some line bundle
${\mathcal L}$ on $Z$.

Let
\begin{equation}\label{aes}
 0 \,\longrightarrow\, {\mathcal End}(\widetilde{E})\,\longrightarrow\,
 \mathrm{At}(\widetilde{E})\, \xrightarrow{\,\,\mathrm{symb}_1\,\,\,}\,
 T_{X\times Z/Z}\, \longrightarrow\, 0
\end{equation}
be the relative Atiyah exact sequence, where $\mathrm{At}(\widetilde{E})$ is
the Atiyah bundle for $\widetilde E$. Setting $\mathrm{At}_D(\widetilde{E})$ to be the pullback of
$T_{X\times Z/Z}(-D\times Z)$ by the surjection
$\mathrm{At}(\widetilde{E})\,\longrightarrow\, T_{X\times Z/Z}$ in \eqref{aes},
we get a short exact sequence
\[
0 \,\longrightarrow\, {\mathcal End}(\widetilde{E})\, \longrightarrow
\, \mathrm{At}_D(\widetilde{E}) \,\xrightarrow{\,\,\mathrm{symb}_1\,\,\,}\, 
 T_{X\times Z/Z}(-D\times Z)\, \longrightarrow\, 0.
\]
By \cite[Theorem 7.7.6]{Groth1},
there exists a coherent sheaf
${\mathcal H}$ on $Z$
and a functorial isomorphism
\[
 \pi_{S*}\left( \mathrm{At}_D(\widetilde{E})
 \otimes\Omega^1_{X\times Z/Z}(D_Z)
 \otimes_{{\mathcal O}_Z} {\mathcal Q}\right)
\ \cong\ {\mathcal Hom}_{{\mathcal O}_S}
 \big({\mathcal H}\otimes_{{\mathcal O}_Z}{\mathcal O}_S
 {\mathcal Q} \big)
\]
for any morphism $S\,\longrightarrow\, Z$
and any coherent sheaf ${\mathcal Q}$ on $S$.
Set
${\mathcal V}\,:=\,\Spec \left( \mathrm{Sym}^*( {\mathcal H} )\right)$.
Then there is a universal section
$\widetilde{\Psi}\,\colon\, T_{X\times{\mathcal V}/{\mathcal V}}(-D\times{\mathcal V})
\,\longrightarrow \,\mathrm{At}_D(\widetilde{E})$.
Note that the composition of maps $\mathrm{symb}\circ\widetilde{\Psi}$
defines a global section of ${\mathcal O}_{X\times{\mathcal V}}$,
which is a section of ${\mathcal O}_{\mathcal V}$.
Let ${\mathcal V}'$ be the closed subscheme of ${\mathcal V}$ defined by
the condition $\mathrm{symb}_1\circ\widetilde{\Psi}\,=\,1$.
Then the restriction $\widetilde{\Psi}|_{{\mathcal V}'}$
defines a universal relative connection
\[
\widetilde{\nabla}\,\,\colon\,\,
\widetilde{E}_{\mathcal V'}\,\,\longrightarrow\,\,
\widetilde{E}_{\mathcal V'}\otimes\Omega_{X\times{\mathcal V'}/{\mathcal V'}}
(D_{\mathcal V'}).
\]
Let $B$ be the maximal closed subscheme of ${\mathcal V'}$
such that
$\big( \res_{x_i\times{\mathcal V'}}(\widetilde{\nabla}) -\nu^{(i)}_j\mathrm{id}\big)
(\widetilde{l}^{(i)}_j)_{\mathcal V'} \,\subset\, (\widetilde{l}^{(i)}_{j+1})_{\mathcal V'}$
for any $i,\,j$ and $(\varphi\otimes\mathrm{id})\circ\widetilde{\nabla}\circ\varphi^{-1}
\,=\,\nabla_L\otimes\mathrm{id}_{\mathcal L}$.
Set
\begin{align*}
 \widetilde{\mathcal D}^{\mathrm{par}}_{\mathfrak{sl},0}
 &\,:=\,
 \left\{ u\,\in\, {\mathcal End}(\widetilde{E}_B)\,\, \middle| \:
 \text{$\Tr(u)\,=\,0$ and
 $u|_{x_i\times B}(\widetilde{l}^{(i)}_j)_B\,\subset\, (\widetilde{l}^{(i)}_j)_B$
 for any $i,\,j$}
 \right\},
 \\
 \widetilde{\mathcal D}^{\mathrm{par}}_{\mathfrak{sl},1}
 &\,:=\,
 \left\{ u\,\in\, {\mathcal End}(\widetilde{E}_B)\otimes K_X(D)\ \middle|\, \:
 \text{$\Tr(u)\,=\,0$ and $\mathrm{res}_{x_i\times B}(u)(\widetilde{l}^{(i)}_j)_B
 \,\subset\, (\widetilde{l}^{(i)}_{j+1})_B$ for any $i,\,j$}
 \right\},
 \\
 \nabla_{\widetilde{\mathcal D}^{\mathrm{par}}_{\mathfrak{sl},\bullet}}
 &\colon \:
\widetilde{\mathcal D}^{\mathrm{par}}_{\mathfrak{sl},0}\, \longrightarrow\,
\widetilde{\mathcal D}^{\mathrm{par}}_{\mathfrak{sl},1},\ \ \,
u\, \longmapsto\, \widetilde{\nabla}\circ u-(u\otimes\mathrm{id})\circ\widetilde{\nabla}.
\end{align*}
There is a canonically induced morphism
\[
 B \ \longrightarrow\ Z
\]
whose fiber over a point $z$ 
is an affine space isomorphic to
$H^0\big(X,\, \widetilde{\mathcal D}^{\mathrm{par}}_{\mathfrak{sl},1}
|_{X\times\{z\}}\big)$. Set
\[
 B^{\circ}\,:=\,
 \left\{ x\,\in\, B\, \, \middle|\, \:
 \text{$(\widetilde{E},\,\widetilde{\nabla},\,\widetilde{\boldsymbol{l}})\big\vert_{X\times x}$\, is simple}
 \right\}.
\]
Then there is a canonically induced morphism
\[
 q\,\, \colon \,\, B^{\circ} \,\,\longrightarrow\,\,
 \mathcal{M}_{\mathrm{PC}}^{n_0\text{-\rm reg}}(\boldsymbol{\nu},\nabla_L).
\]
By the construction, the complement
$\mathcal{M}_{\mathrm{PC}}^{n_0\text{-\rm reg}}(\boldsymbol{\nu},\nabla_L)
\setminus \mathcal{M}_{\mathrm{PC}}^{n_0\text{-\rm reg}}(\boldsymbol{\nu},\nabla_L)^{\circ} $
coincides with the image $q(B^{\circ})$.
So it suffices to show that for every irreducible component
$B'$ of $B^{\circ}$, the closure
$\overline{q(B')}$
has dimension at most $2(r^2-1)(g-1)+r(r-1)n-2$.

For each point $b\,\in \,B'$, consider the group
$\mathrm{Aut}((\widetilde{E},\,\widetilde{\boldsymbol{l}},\,\det \widetilde{E})|_{X\times \{b\}})$
of automorphisms of $\widetilde{E}|_{X\times\{b\}}$ preserving 
$\widetilde{\boldsymbol{l}}_{D\times\{b\}}$
and $\det \widetilde{E}|_{X\times\{b\}}$.
Then the tangent space of 
$\Aut((\widetilde{E},\,\widetilde{\boldsymbol{l}},\,\det \widetilde{E})|_{X\times \{b\}})$ is isomorphic to 
$H^0(X,\, \widetilde{\mathcal D}^{\mathrm{par}}_{\mathfrak{sl},0}
|_{X\times\{b\}})$.
For a point $b$ of $B'$, there is the orbit map
\[
\mathrm{Aut}((\widetilde{E},\,\widetilde{\boldsymbol{l}},\,\det \widetilde{E})|_{X\times\{b\}}) 
\,\, \longrightarrow\,\, B',\,\,\,
\; g \,\longmapsto \,g\cdot b,
\]
whose differential
$H^0\big(X,\, \widetilde{\mathcal D}^{\mathrm{par}}_{\mathfrak{sl},0}
|_{X\times\{b\}}\big)\, 
\xrightarrow{\,\, \nabla_{\widetilde{\mathcal D}^{\mathrm{par}}_{\mathfrak{sl},\bullet}}\,\,\,}\,
H^0\big(X,\,\widetilde{\mathcal D}^{\mathrm{par}}_{\mathfrak{sl},1}|_{X\times\{b\}}\big)
$
is injective because 
$(\widetilde{E},\,\widetilde{\nabla},\,\widetilde{\boldsymbol{l}})|_{X\times\{b\}}$
is simple. Since the fiber $q^{-1}(x)$ over a point $x$ of
$\mathcal{M}_{\mathrm{PC}}^{n_0\text{-\rm reg}}(\boldsymbol{\nu},\nabla_L)$
contains an orbit for the action of
$\mathrm{Aut}((\widetilde{E},\,\widetilde{\boldsymbol{l}},\,\det \widetilde{E})|_{X\times\{b\}})$,
we have 
\[
\dim q^{-1}(x)\,\,\geq\,\,
\dim H^0(X,\, \widetilde{\mathcal D}^{\mathrm{par}}_{\mathfrak{sl},0}
|_{X\times\{b\}}).
\]
Note that we have 
$\big(\widetilde{\mathcal D}^{\mathrm{par}}_{\mathfrak{sl},0}\big)^{\vee}
\otimes K_X \,\cong\, \widetilde{\mathcal D}^{\mathrm{par}}_{\mathfrak{sl},1}$,
and
\[
\dim H^0\big(X,\widetilde{\mathcal D}^{\mathrm{par}}_{\mathfrak{sl},1}
|_{X\times \{b\}}\big)
-\dim H^0\big(X,\widetilde{\mathcal D}^{\mathrm{par}}_{\mathfrak{sl},0}
|_{X\times\{b\}}\big)\,=\,(r^2-1)(g-1)+nr(r-1)/2
\]
by the Riemann--Roch theorem.
If we choose $x$ to be a generic point of $q(B')$, then we have
\begin{align*}
 \dim \overline{q(B')}
 &\,=\,
 \dim B'-\dim q^{-1}(x)
 \\
 &\leq\,
 \dim B'
 -\dim H^0\big(X,\, \widetilde{\mathcal D}^{\mathrm{par}}_{\mathfrak{sl},0}
 |_{X\times \{b\}}\big)
 \\
 &\leq\,
 \dim Z
 +\dim H^0\big(\widetilde{\mathcal D}^{\mathrm{par}}_{\mathfrak{sl},1}
 |_{X\times\{b\}}\big)
 -\dim H^0\big(\widetilde{\mathcal D}^{\mathrm{par}}_{\mathfrak{sl},0}
 |_{X\times\{b\}}\big)
 \\
 &=\,\dim Z +(r^2-1)(g-1)+r(r-1)n/2
 \\
 &\leq\,
 2(r^2-1)(g-1)+nr(r-1)-2.
\end{align*}
Since $q(B^{\circ})=\mathcal{M}_{\mathrm{PC}}^{n_0\text{-\rm reg}}(\boldsymbol{\nu},\nabla_L)
\setminus \mathcal{M}_{\mathrm{PC}}^{n_0\text{-\rm reg}}(\boldsymbol{\nu},\nabla_L)^{\circ} $
is a union of $q(B')$'s, the proof is completed.
\end{proof}

Define the open subset
$\mathcal{M}^{n_0\text{-\rm reg}}_{\mathrm{Higgs}}(\boldsymbol{\mu},\,\Phi_L)^{\circ}$ of
$\mathcal{M}^{n_0\text{-\rm reg}}_{\mathrm{Higgs}}(\boldsymbol{\mu},\,\Phi_L)$ by
\begin{equation}\label{underlying parabolic bundle is simple Higgs moduli}
\mathcal{M}^{n_0\text{-\rm reg}}_{\mathrm{Higgs}}(\boldsymbol{\mu},\,\Phi_L)^{\circ}
:= \left\{ (E,\,\Phi,\,\boldsymbol{l}) \in
\mathcal{M}^{n_0\text{-\rm reg}}_{\mathrm{Higgs}}(\boldsymbol{\mu},\,\Phi_L) \,\middle|\, \dim\big( \End(E,\,\boldsymbol{l}) \big)=1 \right\}
\end{equation}
which consists of $\boldsymbol{\mu}$-parabolic Higgs bundles
$(E,\,\Phi,\,\boldsymbol{l})$ with the determinant isomorphic to $(L,\Phi_L)$
such that the underlying quasi-parabolic bundle $(E,\,\boldsymbol{l})$ is simple.

The proof of the following proposition uses an argument similar to one in the proofs of 
Proposition \ref{proposition: codimension in de Rham moduli}.

\begin{Prop}\label{Higgs moduli codimension for underlying parabolic bundle}
Let $X$ be a smooth projective curve of genus $g$ over an algebraically closed field $k$,
and let $L$ be a line bundle on $X$ with a homomorphism
$\Phi_L \,\colon\, L\,\longrightarrow\, L\otimes K_X(D)$.
Take positive integers $r,n$ and a tuple $\boldsymbol{\mu}=(\mu^{(i)}_j)^{1\leq i\leq n}_{0\leq j\leq r-1}\in\,k^{nr}$
such that $\res_{x_i}(\Phi_L)\,=\,\sum_{j=0}^{r-1}\mu^{(i)}_j$
for any $i$.
Assume that $r$ is not divisible by the characteristic of $k$ and one of the following holds:
\begin{itemize}
\item $n\,\geq\, 1$ and $r\,\geq\, 2$ are arbitrary if $g\,\geq\, 2$,

\item $n\,\geq\, 2$, $r\geq 2$ and $n+r\,\geq\, 5$ if $g\,=\,1$,

\item $n\,\geq \,3$, $r\geq 2$ and $n+r\,\geq \,7$ if $g\,=\,0$.
\end{itemize}
Then $\mathrm{codim}_{\mathcal{M}^{n_0\text{-\rm reg}}_{\mathrm{Higgs}}(\boldsymbol{\mu},\,\Phi_L)}
\left( \mathcal{M}^{n_0\text{-\rm reg}}_{\mathrm{Higgs}}(\boldsymbol{\mu},\,\Phi_L)
\setminus \mathcal{M}^{n_0\text{-\rm reg}}_{\mathrm{Higgs}}(\boldsymbol{\mu},\,\Phi_L)^{\circ}\right)
\,\,\geq\,\, 2$.
\end{Prop}

\begin{proof}
By Proposition \ref{dimension of non-simple parabolic bundles in higher genus case},
\ref{dimension of non-simple parabolic bundles when genus=1},
\ref{proposition: codimension genus=0 and n>3},
\ref{proposition: codimension genus zero n=3 case},
there is a scheme $Z$ of finite type over $\Spec k$
and a flat family $(\widetilde{E},\,\widetilde{\boldsymbol{l}})$
of quasi-parabolic bundles on $X\times Z$ over $Z$
such that
\begin{enumerate}
\item $\dim Z\,\leq\, (r^2-1)(g-1)+nr(r-1)/2-2$,

\item $\dim\End((\widetilde{E},\,\widetilde{\boldsymbol{l}})|_{X\times z})\,\geq\, 2$
for all $z\,\in\, Z$, and

\item{} each quasi-parabolic bundle in the complement
$\left| {\mathcal N}^{n_0\text{-reg}}_{\mathrm{par}}(L) \right|
\setminus \left| {\mathcal N}^{n_0\text{-reg}}_{\mathrm{par}}(L)^{\circ} \right|$
is isomorphic to $(\widetilde{E},\,\widetilde{\boldsymbol{l}})|_{X\times\{z\}}$
for some $z\,\in\, Z$.
\end{enumerate}
Define
\begin{align*}
 \widetilde{\mathcal D}^{\mathrm{par}}_{\mathfrak{sl},0}
 &\,:=\,
 \left\{ u\,\in\, {\mathcal End}(\widetilde{E}) \, \middle| \:
 \text{$u|_{x_i\times Z_{\alpha}}(\widetilde{l}^{(i)}_j)\,\subset\, \widetilde{l}^{(i)}_j$
 for any $i,\,j$}\right\}
 \\
 \widetilde{\mathcal D}^{\mathrm{par}}_{\mathfrak{sl},1}
 &\,:=\,
 \left\{ u\,\in\, {\mathcal End}(\widetilde{E})\otimes K_X(D) \, \middle| \:
 \text{$\mathrm{res}_{x_i\times Z_{\alpha}}(u)
 (\widetilde{l}^{(i)}_j)\,\subset\, \widetilde{l}^{(i)}_{j+1}$ for any $i,\,j$}
 \right\} .
\end{align*}
By \cite[Theorem 7.7.6]{Groth1}, there is a coherent sheaf ${\mathcal H}$ on
$Z$ together with a functorial isomorphism
\[
 \Hom({\mathcal H}\otimes_{{\mathcal O}_Z}{\mathcal O}_S,\,Q)\,\cong\,
 H^0(X\times S,\,\widetilde{\mathcal D}^{\mathrm{par}}_{\mathfrak{sl},1}
 \otimes_{{\mathcal O}_Z}Q)
\]
for any Noetherian scheme $S$ over $Z$ and any coherent sheaf $Q$ on $S$.
For $\mathbb{V}({\mathcal H})\,:=\,
\Spec \big( \mathrm{Sym}^*({\mathcal H})\big)$,
there is a universal family of Higgs fields
$\widetilde{\Phi}\,\in\, 
H^0(X\times S,\,\widetilde{\mathcal D}^{\mathrm{par}}_{\mathfrak{sl},0}
\otimes K_X(D)\otimes_{{\mathcal O}_Z}
{\mathcal O}_{\mathbb{V}({\mathcal H})})$
on $(\widetilde{E},\,\widetilde{\boldsymbol{l}})\otimes {\mathcal O}_{\mathbb{V}({\mathcal H})}$.
We may assume that
$\det(\widetilde{E})\,\cong\, L\otimes{\mathcal P}$
for some line bundle ${\mathcal P}$ on $Z$.
Let $B$ be the maximal locally closed subscheme of
$\mathbb{V}({\mathcal H})$
such that the composition of the homomorphisms
\[
 L\otimes{\mathcal P}_B
 \,\xrightarrow{\,\,\sim\,}\,\det(\widetilde{E})_B
 \,\xrightarrow{\,\,\Tr\widetilde{\Phi}\,}\,
 \det(\widetilde{E})\otimes K_X(D)_B\,
 \xrightarrow{\,\,\sim\,}\, L\otimes{\mathcal P}_B
\]
coincides with $\nabla_L\otimes{\mathcal P}_B$ and
$(\res_{x_i\times Z} (\widetilde{\Phi}) - \mu^{(i)}_\ell ) (\widetilde{l}^{(i)}_\ell)
\,\subset\, \widetilde{l}^{(i)}_{\ell +1}$ for any $i,\,\ell$ and also
$(\widetilde{E},\,\widetilde{\boldsymbol{l}},\,\widetilde{\Phi})|_{X\times b}$ is simple for any $b\,\in\, B$.
Then the family $(\widetilde{E},\,\widetilde{\boldsymbol{l}},\,\widetilde{\Phi})_{B}$
defines a morphism
\begin{equation}\label{equation: morphism to special Dolbeault moduli}
 B \ \longrightarrow\
 \mathcal{M}^{n_0\text{-\rm reg}}_{\mathrm{Higgs}}(\boldsymbol{\mu},\,\Phi_L)
\end{equation}
whose image coincides with
$\mathcal{M}^{n_0\text{-\rm reg}}_{\mathrm{Higgs}}(\boldsymbol{\mu},\,\Phi_L)
\setminus \mathcal{M}^{n_0\text{-\rm reg}}_{\mathrm{Higgs}}(\boldsymbol{\mu},\,\Phi_L)^{\circ}$.
Note that the fibers of the morphism in \eqref{equation: morphism to special Dolbeault moduli}
contain orbits of the action by the automorphism group of
$(\widetilde{E},\,\widetilde{\boldsymbol{l}},\,\det(\widetilde{E}))_z$
whose dimension is that of 
$H^0(X,\,\widetilde{\mathcal D}^{\mathrm{par}}_{\mathfrak{sl},0}|_{X\times z})$.
So we have
\begin{align*}
 \dim\im\left(B\longrightarrow 
 \mathcal{M}^{n_0\text{-\rm reg}}_{\mathrm{Higgs}}(\boldsymbol{\mu},\,\Phi_L)\right)
 &\leq\,
 \dim B-
 \dim H^0\big(X,\,\widetilde{\mathcal D}^{\mathrm{par}}_{\mathfrak{sl},0}|_{X\times z}\big)
 \\
 &\leq \, \dim Z+
 \dim H^0 \big( X,\,\widetilde{\mathcal D}^{\mathrm{par}}_{\mathfrak{sl},1}|_{X\times z}\big)
 -
 \dim H^0\big(X,\,\widetilde{\mathcal D}^{\mathrm{par}}_{\mathfrak{sl},0}|_{X\times z}\big)
 \\
 &=\, \dim Z+ (r^2-1)(g-1)+nr(r-1)/2
 \\
 &\leq \, \dim \mathcal{M}^{n_0\text{-\rm reg}}_{\mathrm{Higgs}}(\boldsymbol{\mu},\,\Phi_L)-2.
\end{align*}
Since 
$\mathcal{M}^{n_0\text{-\rm reg}}_{\mathrm{Higgs}}(\boldsymbol{\mu},\,\Phi_L)
\setminus \mathcal{M}^{n_0\text{-\rm reg}}_{\mathrm{Higgs}}(\boldsymbol{\mu},\,\Phi_L)^{\circ}$
coincides with the image of the morphism in \eqref{equation: morphism to special Dolbeault moduli},
the proof is complete.
\end{proof}

As a corollary of the above theorem,
we can also get a result, \cite[Theorem 4.2, (c)]{B-Y}, by Boden and Yokogawa.

\begin{Cor}\label{Cor: irreducibility of moduli of parabolic Higgs bundles}
Under the same assumption as in
Proposition \ref{proposition: codimension in de Rham moduli} and
Proposition \ref{Higgs moduli codimension for underlying parabolic bundle},
the moduli spaces
$\mathcal{M}_{\mathrm{PC}}^{n_0\text{-\rm reg}}(\boldsymbol{\nu},\nabla_L)$
and
$\mathcal{M}^{n_0\text{-\rm reg}}_{\mathrm{Higgs}}(\boldsymbol{\mu},\,\Phi_L)$
are irreducible.
\end{Cor}

\begin{proof}
We only prove the irreducibility for
$\mathcal{M}^{n_0\text{-\rm reg}}_{\mathrm{Higgs}}(\boldsymbol{\mu},\,\Phi_L)$
as the proof is same for
$\mathcal{M}_{\mathrm{PC}}^{n_0\text{-\rm reg}}(\boldsymbol{\nu},\nabla_L)$.
The open subspace
$\mathcal{M}^{n_0\text{-\rm reg}}_{\mathrm{Higgs}}(\boldsymbol{\mu},\,\Phi_L)^{\circ}$
is isomorphic to an affine space bundle over the moduli space
${\mathcal N}_{\rm par}^{n_0\text{-reg}}(L)$ of $n_0$-regular simple quasi-parabolic bundles
with the determinant $L$.
Since ${\mathcal N}_{\rm par}^{n_0\text{-reg}}(L)$ is irreducible, it follows that
$\mathcal{M}^{n_0\text{-\rm reg}}_{\mathrm{Higgs}}(\boldsymbol{\mu},\,\Phi_L)^{\circ}$
is also irreducible.
Recall that the moduli space
$\mathcal{M}^{n_0\text{-\rm reg}}_{\mathrm{Higgs}}(\boldsymbol{\mu},\,\Phi_L)$
is smooth of equi-dimension by
Proposition \ref{Prop: smoothness and dimension of moduli of Hiiggs bundles}.
So $\mathcal{M}^{n_0\text{-\rm reg}}_{\mathrm{Higgs}}(\boldsymbol{\mu},\,\Phi_L)$
is connected and thus irreducible, because
$\dim \left( \mathcal{M}^{n_0\text{-\rm reg}}_{\mathrm{Higgs}}(\boldsymbol{\mu},\,\Phi_L)
\setminus\mathcal{M}^{n_0\text{-\rm reg}}_{\mathrm{Higgs}}(\boldsymbol{\mu},\,\Phi_L)^{\circ}\right)
\,<\,\dim \mathcal{M}^{n_0\text{-\rm reg}}_{\mathrm{Higgs}}(\boldsymbol{\mu},\,\Phi_L)$
by Proposition \ref{Higgs moduli codimension for underlying parabolic bundle}.
\end{proof}

\begin{Rem}
The proof of Corollary \ref{Cor: irreducibility of moduli of parabolic Higgs bundles}
is in fact valid under a weaker assumption than that of 
Theorem \ref{Higgs moduli codimension for underlying parabolic bundle}.
Indeed, it is valid under the same assumption as that of \cite[Theorem 2.2]{Inaba-1}.
\end{Rem}

\subsection{The moduli space is not affine}\label{subsection: transcendence degree}

We use the notation of Section \ref{subsection: moduli of parabolic connections and Higgs bundles}.
In this subsection, $k$ is assumed to be an algebraically closed field
of arbitrary characteristic unless otherwise noted.

Let $X$ be a smooth projective curve over $k$ of genus $g$.
Fix a line bundle $L$ of degree $d$ on $X$ equipped with a logarithmic connection
$\nabla_L\,\colon\, L \,\longrightarrow\, L\otimes K_X(D)$, and also fix a string of local exponents
$\boldsymbol{\nu}\,=\, (\nu^{(i)}_j)\,\in\, k^{nr}$
such that $\res_{x_i}(\nabla_L)\,=\,\sum_{j=0}^{r-1}\nu^{(i)}_j$ for any $i$.
We assume the following:
\begin{equation}\label{assumption of local exponents}
\sum_{i=1}^n\sum_{\ell =1}^s \nu^{(i)}_{j^{(i)}_\ell}
\,\notin\,\mathrm{Im}\left(\mathbb{Z}\rightarrow k\right)
\quad \text{for any choice of $s$ elements
$\{ j^{(i)}_1,\,\cdots,\, j^{(i)}_s \}$ in $\{1,\,\cdots,\,r\}$.}
\end{equation}
Under the assumption in \eqref{assumption of local exponents},
any $\boldsymbol{\nu}$-parabolic connection is irreducible, and hence it is
$\boldsymbol{\alpha}$-stable for any parabolic weight $\boldsymbol{\alpha}$.
So we have
$\mathcal{M}_{\mathrm{PC}}^{\boldsymbol{\alpha}}(\boldsymbol{\nu},\,\nabla_L)
\,=\,\mathcal{M}_{\mathrm{PC}}(\boldsymbol{\nu},\,\nabla_L)$.
In this subsection we will show that
the moduli space $\mathcal{M}_{\mathrm{PC}}^{\boldsymbol{\alpha}}(\boldsymbol{\nu},\,\nabla_L)$
is not affine. This will be done by comparing the transcendence degree of
the ring of global algebraic functions on the moduli space
$\mathcal{M}_{\mathrm{PC}}^{\boldsymbol{\alpha}}(\boldsymbol{\nu},\nabla_L)$
of parabolic connections with the transcendence degree of
the ring of global algebraic functions on the moduli space of parabolic Higgs bundles.

Consider the moduli space
\[
 {\mathcal M}^{\boldsymbol{\alpha}}_{\rm Higgs}(d)\, =\,
 \left\{
 (E,\,\Phi,\,\boldsymbol{l})\,\, \middle|\,
 \begin{array}{l}
 \text{$(E,\,\boldsymbol{l})$ is a quasi-parabolic bundle of rank $r$ and degree $d$,} \\
 \text{$\Phi\colon E\to E\otimes K_X(D)$ is an ${\mathcal O}_X$-homomorphism such that} \\
 \text{$\res_{x_i}(\Phi)(l^{(i)}_j)\,\subset\, l^{(i)}_j$ for any $i,\,j$,
 and $(E,\,\Phi,\,\boldsymbol{l})$ is $\boldsymbol{\alpha}$-stable}
 \end{array}
 \right\}
\]
of $\boldsymbol{\alpha}$-stable parabolic Higgs bundles.
Setting
\[
 \Lambda_{\rm Higgs}
 =
 \Bigg\{
 \boldsymbol{\mu}=(\mu^{(i)}_j)^{1\leq i\leq n}_{0\leq j\leq r-1} 
 \, \in \, k^{nr}
 \ \Bigg| \ \sum_{i=1}^n\sum_{j=0}^{r-1}\mu^{(i)}_j=0 \Bigg\},
\]
we have a canonical morphism
\[
{\mathcal M}^{\boldsymbol{\alpha}}_{\rm Higgs}(d)
\,\, \longrightarrow \,\, \Lambda_{\rm Higgs}
\]
whose fiber over any $\boldsymbol{\mu}\,\in\, \Lambda_{\rm Higgs}$
is the moduli space
${\mathcal M}^{\boldsymbol{\alpha}}_{\rm Higgs}(\boldsymbol{\mu})$
of $\boldsymbol{\alpha}$-stable $\boldsymbol{\mu}$-parabolic Higgs bundles.
For a parabolic Higgs bundle
$(E,\,\Phi,\,\boldsymbol{l})\,\in\, {\mathcal M}^{\boldsymbol{\alpha}}_{\rm Higgs}(d)$,
consider the homomorphism
\[
 T\,\mathrm{Id}_E-\Phi \,\,\colon\, \,
 E \otimes k[T] \,\longrightarrow\,
 E\otimes \mathrm{Sym}^*(K_X(D))\otimes k[T],
\]
where $T$ is an indeterminate. We can write
\[
 \det(T\,\mathrm{Id}_E-\Phi)
\, =\,T^r+s_1T^{r-1}+\cdots+s_{r-1}T+s_r
\]
with $s_j\,\in\, H^0(X,\,K_X^{\otimes j}(jD))$.
Note that $s_1\,=\,-\Tr(\Phi)$. Set
\[
 W\,\,:=\,\,\bigoplus_{j=1}^r H^0(X,\,K_X^{\otimes j}(jD)).
\]
Using the above constructed $(s_1,\,\cdots,\,s_r)$, we get a morphism
\begin{equation}\label{equation: Hitchin map}
 H \,\,\colon\,\, {\mathcal M}^{\boldsymbol{\alpha}}_{\rm Higgs}(d)
 \,\, \longrightarrow \,\, W,
\end{equation}
which is called the Hitchin map.
A remarkable property of the Hitchin map is that it is proper, which was proved by Hitchin, Simpson and Nitsure.
We use the parabolic version of it which was proved by Yokogawa.

\begin{Thm}[{\cite{Hi}, \cite{Sim2}, \cite{Nit}, \cite{Yokogawa}}]\label{Theorem: Hitchin map is proper}
Under the assumption that $\boldsymbol{\alpha}$-semistability implies $\boldsymbol{\alpha}$-stability,
the Hitchin map
$H\,\,\colon\, \,{\mathcal M}^{\boldsymbol{\alpha}}_{\rm Higgs}(d)\,\, \longrightarrow\,\, W$
in \eqref{equation: Hitchin map} is a proper morphism.
\end{Thm}

Set
\[
 A_{\rm Higgs}\,\,:=\,\,
 \left\{ \boldsymbol{a}\,=\,(a^{(i)}_j)^{1\leq i\leq n}_{1\leq j\leq r} 
 \, \in \, k^{nr} \ \middle| \ 
 \sum_{i=1}^n a^{(i)}_1=0
 \right\}.
\]
Using the correspondence $(s_\ell)_{1\leq \ell\leq r}\, \longmapsto\,
(\res_{x_i}(s_\ell))^{1\leq i\leq n}_{1\leq\ell\leq r}$,
we define a morphism
\[
 W \, \longrightarrow \, A_{\rm Higgs}
\]
which is a linear surjection under any of the following conditions:
\begin{enumerate}
\item[(i)]
$n\,\geq\, 1$ when $g\,\geq\, 2$,

\item[(ii)]
$n\,\geq\, 2$ when $g\,=\,1$,

\item[(iii)] $n\,\geq\, 3$ when $g\,=\,0$. 
\end{enumerate}
There is also a morphism
\[
 \Lambda_{\rm Higgs} \, \longrightarrow \, A_{\rm Higgs}
\]
that associates the coefficients of 
$\prod_{j=0}^{r-1}(t-\mu^{(i)}_j)$.
Then the Hitchin map induces a morphism
\begin{equation}\label{base change of Hitchin map}
{\mathcal M}^{\boldsymbol{\alpha}}_{\rm Higgs}(d)
\ \longrightarrow \
W\times_{ A_{\rm Higgs} } \Lambda_{\rm Higgs},
\end{equation}
which is proper by Theorem \ref{Theorem: Hitchin map is proper}.

Fix a line bundle $L$ on $X$ of degree $d$, and consider the closed subvariety
\[
{\mathcal M}^{\boldsymbol{\alpha}}_{\rm Higgs}(L)
\, :=\,
\left\{ (E,\,\Phi,\,\boldsymbol{l}) \, \in \, {\mathcal M}^{\boldsymbol{\alpha}}_{\rm Higgs}(d)
\,\ \middle| \,\ \det(E)\,\cong\, L \right\}
\]
of ${\mathcal M}^{\boldsymbol{\alpha}}_{\rm Higgs}(d)$.
Then the restriction of the map in \eqref{base change of Hitchin map}
\begin{equation}\label{base change of Hitchin map with fixed determinant}
{\mathcal M}^{\boldsymbol{\alpha}}_{\rm Higgs}(L)
\ \longrightarrow \
W\times_{ A_{\rm Higgs} } \Lambda_{\rm Higgs} 
\end{equation}
is also a proper morphism.

Generic fibers of the Hitchin map were investigated by Logares and Martens 
in \cite[Proposition 2.2]{Log-Mar}. The following result is likely to be well-known
to the experts. We give a proof of it using the arguments given by
Alfaya and G\'omez in \cite[Lemma 3.2]{Alf-Gom}.

\begin{Cor}\label{Cor: base change of Hitchin map is surjective} 
Assume that $\boldsymbol{\alpha}$-semistability implies $\boldsymbol{\alpha}$-stability.
Also, assume that one of the following statements holds:
\begin{enumerate}
\item[(i)]
$n\,\geq\, 1$ if $g\,\geq\, 2$,
\item[(ii)]
$n\,\geq\, 2$ if $g\,=\,1$,
\item[(iii)]
$n\,\geq\, 3$ if $g\,=\,0$. 
\end{enumerate}
Then the morphism
${\mathcal M}^{\boldsymbol{\alpha}}_{\rm Higgs}(L)
\ \longrightarrow\ W\times_{ A_{\rm Higgs} } \Lambda_{\rm Higgs}$
in \eqref{base change of Hitchin map with fixed determinant} is surjective.
\end{Cor}

\begin{proof}
It suffices to prove that the morphism in \eqref{base change of Hitchin map with fixed determinant} is dominant,
because it is proper. Take any
$(s\,=\,(s_\ell),\,\boldsymbol{\mu})\,\in\, W\times_{ A_{\rm Higgs} } \Lambda_{\rm Higgs}$.
Consider the corresponding spectral curve $X_s\, \subset\, \mathbb{P}({\mathcal O}_X\oplus K_X(D))$ which is defined by the equation
\[
 y^r+s_1y^{r-1}+\cdots+s_{r-1}y+s_r\,=\,0,
\]
where $y$ is the section of ${\mathcal O}_{\mathbb{P}({\mathcal O}_X\oplus K_X(D))}(1)$ corresponding to the inclusion map
${\mathcal O}_X\,\hookrightarrow \, {\mathcal O}_X\oplus K_X(D)$.
Take a section $\tau\,\in\, H^0(X,\, K_X^{\otimes r}(rD))$
which has at most simple zeroes; since $K_X^{\otimes r}(rD)$ is very ample
by the assumption in the corollary, such a section exists.
Then the spectral curve $y^r-\tau\,=\,0$ has no singular points.

Since the smoothness is an open condition,
there is an open subset $U\,\subset\, W\times_{A_{\rm Higgs}}\Lambda_{\rm Higgs}$ such that
the spectral curve $X_s$ is smooth for every $s\,\in \,U$. Take a line bundle ${\mathcal L}$ on $X_s$
such that the locally free sheaf $E\,:=\,\pi_*({\mathcal L})$
has its determinant $\det(E)$ isomorphic to $L$,
where $\pi\,\colon\, X_s\, \longrightarrow\, X$ is the natural projection.
By the Beauville--Narasimhan--Ramanan correspondence \cite[Proposition 3.6]{B-N-R}, 
there is a Higgs field $\Phi\,\colon\, E \,\longrightarrow\, E\otimes K_X(D)$
induced by the action of $y$ on ${\mathcal L}$.
Shrinking $U$ if necessary, we may further assume that
$\mu^{(i)}_0,\,\cdots,\,\mu^{(i)}_{r-1}$ are mutually distinct for any fixed $i$.
Then we can associate a unique parabolic structure $l$ on $E$
compatible with $\Phi$.
Since $(E,\,\Phi,\,\boldsymbol{l})$ is irreducible by its construction,
it is evidently $\boldsymbol{\alpha}$-stable.
So we have $(E,\,\Phi,\,\boldsymbol{l})\,\in\, {\mathcal M}^{\boldsymbol{\alpha}}_{\rm Higgs}(L)$
which is sent to $(s,\,\boldsymbol{\mu})$ under the morphism in
\eqref{base change of Hitchin map with fixed determinant}.
Thus the morphism in \eqref{base change of Hitchin map with fixed determinant} 
is dominant because its image contains the dense open subset 
$U$ of $W\times_{ A_{\rm Higgs} } \Lambda_{\rm Higgs}$.
\end{proof}

As a consequence of Theorem \ref{Theorem: Hitchin map is proper}
and Corollary \ref{Cor: base change of Hitchin map is surjective},
we can determine the transcendence degree of the ring of global algebraic functions 
on the moduli space of parabolic Higgs bundles.

\begin{Cor}\label{Cor: global algebraic functions on Higgs moduli space}
Let $L$ be a line bundle on $X$ with a Higgs field
$\Phi_L\,\colon\, L \,\longrightarrow\, L\otimes K_X(D)$.
Take $\boldsymbol{\mu}\,=\,(\mu^{(i)}_j)\,\in \,\Lambda_{\rm Higgs}$
satisfying the condition $\res_{x_i}(\Phi_L)\,=\,\sum_{j=0}^{r-1} \mu^{(i)}_j$
for all $i$. Then, under the same assumption as in
Theorem \ref{Higgs moduli codimension for underlying parabolic bundle},
the transcendence degree of the ring of
global algebraic functions on the moduli space of
parabolic Higgs bundles is given by the following:
\[
\mathrm{tr.deg}_{\,k}\,
\Gamma \big( {\mathcal M}^{\boldsymbol{\alpha}}_{\rm Higgs}(\boldsymbol{\mu}, \Phi_L) ,\,
{\mathcal O}_{ {\mathcal M}^{\boldsymbol{\alpha}}_{\rm Higgs}(\boldsymbol{\mu}, \Phi_L) } \big)
\, =\, (r^2-1)(g-1)+\frac{1}{2}nr(r-1).
\]
\end{Cor}

\begin{proof}
The closed subvariety
\[
 Y\,:=\, \left\{
 (s\,=\,(s_\ell )_{1\leq\ell\leq r-1},\, \boldsymbol{\mu}) \,\in\, W\times_{A_{\rm Higgs}}
\{ \boldsymbol{\mu} \} \ \middle| \ s_1\,=\,-\Phi_L \right\} 
\]
of $W\times_{A_{\rm Higgs}} \Lambda_{\rm Higgs}$ is isomorphic to an affine space.
Its inverse image
${\mathcal M}^{\boldsymbol{\alpha}}_{\rm Higgs}(L) \times_{ W\times_{ A_{\rm Higgs} } \Lambda_{\rm Higgs} } Y$
for the morphism in \eqref{base change of Hitchin map with fixed determinant}
is nothing but the moduli space
${\mathcal M}^{\boldsymbol{\alpha}}_{\rm Higgs}(\boldsymbol{\mu}, \Phi_L)$
of $\boldsymbol{\alpha}$-stable $\boldsymbol{\mu}$-parabolic Higgs bundles
with determinant $(L,\,\Phi_L)$.
By Corollary \ref{Cor: base change of Hitchin map is surjective}, the base change
\[
 H'\,\,\colon\,\, {\mathcal M}^{\boldsymbol{\alpha}}_{\rm Higgs}(\boldsymbol{\mu}, \Phi_L)
\,\, \longrightarrow \,\, Y
\]
is also a proper and surjective morphism. So the ring homomorphism
${\mathcal O}_Y\,\longrightarrow\,
H'_*{\mathcal O}_{{\mathcal M}^{\boldsymbol{\alpha}}_{\rm Higgs}(\boldsymbol{\mu},\Phi_L)}$
is injective, and
$H'_*{\mathcal O}_{{\mathcal M}^{\boldsymbol{\alpha}}_{\rm Higgs}(\boldsymbol{\mu},\Phi_L)}$
is a finite algebra over ${\mathcal O}_Y$.
Therefore,
$\Gamma({\mathcal M}^{\boldsymbol{\alpha}}_{\rm Higgs}(\boldsymbol{\mu}, \Phi_L),\, 
{\mathcal O}_{{\mathcal M}^{\boldsymbol{\alpha}}_{\rm Higgs}(\boldsymbol{\mu},\Phi_L)})$
is a finite algebra over $\Gamma(Y,\,{\mathcal O}_Y)$
whose Krull dimension is
\[
 \dim Y
\, =\,-n(r-1)+\sum_{j=2}^r 
 \dim H^0\left(X,\, K_X^{\otimes j}(jD)\right)
\]
\[
=\, -n(r-1)+\sum_{j=2}^r \left( (2g-2)j+jn+(1-g) \right) 
\, =\,(r^2-1)(g-1)+\frac{nr(r-1)} {2}.
\]
Since 
$\Gamma({\mathcal M}^{\boldsymbol{\alpha}}_{\rm Higgs}(\boldsymbol{\mu}, \Phi_L), \,
{\mathcal O}_{{\mathcal M}^{\boldsymbol{\alpha}}_{\rm Higgs}(\boldsymbol{\mu},\Phi_L)})$
is a finitely generated algebra over $k$,
its transcendence degree over $k$ coincides with its Krull dimension.
\end{proof}

\begin{Prop}\label{Prop: Deligne--Hitchin family}
There is a projective flat morphism
\[
 \overline{\mathcal M'} \ \longrightarrow \ \mathbb{A}^1\ =\ \Spec k[t],
\]
and a $\mathbb{A}^1$-relative very ample divisor $Y\,\subset\, \overline{\mathcal M'}$,
such that the complement ${\mathcal M'}\,:=\,\overline{\mathcal M'}\setminus Y$ satisfies the following:
\[
 {\mathcal M'}_h
\ \cong\
 \begin{cases}
\mathcal{M}_{\mathrm{PC}}^{n_0\text{-\rm reg}}(\boldsymbol{\nu},\nabla_L)^{\circ} & 
(0\,\neq \,h \,\in\,\mathbb{A}^1) \\
\mathcal{M}^{n_0\text{-\rm reg}}_{\mathrm{Higgs}}(\boldsymbol{0},\,0)^{\circ} &(h\,=\,0).
\end{cases}
\]
\end{Prop}

\begin{proof}
Let $\mathcal{N}^{n_0\text{-reg}}_{\mathrm{par}}(L)$ be the moduli space of
simple $n_0$-regular quasi-parabolic bundles $(E,\,\boldsymbol{l})$
with $\det E\,\cong\, L$.
Let $(\widetilde{E},\,\widetilde{\boldsymbol{l}})$ be the universal family over
$X\times \mathcal{N}^{n_0\text{-reg}}_{\mathrm{par}}(L)$.
As in the proof of Proposition \ref{2020.2.24.14.52},
we can construct the relative Atiyah bundle
$\mathrm{At}_D(\widetilde{E})$, which fits in the exact sequence
\[
 0 \,\longrightarrow\, {\mathcal End}(\widetilde{E}) \otimes K_X(D)
 \,\longrightarrow \,\mathrm{At}_D(\widetilde{E}) \otimes K_X(D)
\, \longrightarrow\, {\mathcal O}_{X\times \mathcal{N}^{n_0\text{-reg}}_{\mathrm{par}}(L) }
\, \longrightarrow\, 0.
\]
Recall the construction of the homomorphism \eqref{eq:2021_12_17_12_9}
in the proof of Proposition \ref{2020.2.24.14.52},
which defines a surjection
\[
 \mathrm{At}_D(\widetilde{E})\otimes K_X(D)
 \ \longrightarrow \
 \big( \mathrm{At}_D(\widetilde{E})\otimes K_X(D) \big) \big/
 \big( \mathrm{At}(\widetilde{E})\otimes K_X \big)
 \ \xrightarrow{\,\,\,\sim\,\,} \
 {\mathcal End}(\widetilde{E})\big|_{D\times \mathcal{N}^{n_0\text{-reg}}_{\mathrm{par}}(L) }.
\]
Let $\mathrm{At}_D(\widetilde{E},\,\widetilde{\boldsymbol{l}})
\,\subset\, \mathrm{At}_D(\widetilde{E})$ be the pullback of the subsheaf
\[
\Big\{ a\,\in\,{\mathcal End}(\widetilde{E})\big|_{D\times \mathcal{N}^{n_0\text{-reg}}_{\mathrm{par}}(L) }
\, \Big| \: 
a|_{x_i\times\mathcal{N}^{n_0\text{-reg}}_{\mathrm{par}}(L) }
(\widetilde{l}^{(i)}_j) \,\subset\, \widetilde{l}^{(i)}_j \text{ for any $i,\,j$} \Big\}
\,\,\subset\, {\mathcal End}(\widetilde{E})\big|_{D\times \mathcal{N}^{n_0\text{-reg}}_{\mathrm{par}}(L) }
\]
by the above surjection.

Since $\det(\widetilde{E})\,\cong\, L\otimes{\mathcal P}$ for a line bundle
${\mathcal P}$ on $\mathcal{N}^{n_0\text{-reg}}_{\mathrm{par}}(L)$,
it follows that $\mathrm{At}_D(\det(\widetilde{E}))\,\cong\,
\mathrm{At}_D(L)\otimes{\mathcal O}_{X\times \mathcal{N}^{n_0\text{-reg}}_{\mathrm{par}}(L)}$.
There is an exact sequence
\[
 0 \,\longrightarrow\, {\mathcal O}_{X\times \mathcal{N}^{n_0\text{-reg}}_{\mathrm{par}}(L)}
\, \longrightarrow \,
 \mathrm{At}_D(L)\otimes{\mathcal O}_{X\times \mathcal{N}^{n_0\text{-reg}}_{\mathrm{par}}(L)}
\, \xrightarrow{\,\,\,\mathrm{symb}_1\,\,}
 T_X(-D)\otimes{\mathcal O}_{X\times \mathcal{N}^{n_0\text{-reg}}_{\mathrm{par}}(L)}
\, \longrightarrow\, 0,
\]
which admits a section
$T_X(-D)\otimes{\mathcal O}_{X\times \mathcal{N}^{n_0\text{-reg}}_{\mathrm{par}}(L)}
\,\longrightarrow\,
\mathrm{At}_D(L)\otimes{\mathcal O}_{X\times \mathcal{N}^{n_0\text{-reg}}_{\mathrm{par}}(L)}$
induced by $\nabla_L$.
So its image determines a subbundle of
$\mathrm{At}_D(L)\otimes{\mathcal O}_{X\times \mathcal{N}^{n_0\text{-reg}}_{\mathrm{par}}(L)}$.
Let $\mathrm{At}_D(\widetilde{E},\,\widetilde{\boldsymbol{l}},\,\nabla_L)$
be the pullback of this subbundle by the homomorphism
\begin{equation}\label{equation: connection induced on determinant bundle}
 \mathrm{At}_D(\widetilde{E},\,\widetilde{\boldsymbol{l}})
\, \longrightarrow\, \mathrm{At}_D(\det(\widetilde{E}))
\, \cong\, \mathrm{At}_D(L\otimes{\mathcal O}_{X\times \mathcal{N}^{n_0\text{-reg}}_{\mathrm{par}}(L)})
\end{equation}
defined by
$D\,\longmapsto\, D\wedge\mathrm{Id}\wedge\cdots\wedge\mathrm{Id}+\cdots+
\mathrm{Id}\wedge\cdots\wedge\mathrm{Id}\wedge D$.
If we set
\[
\widetilde{\mathcal D}^{\mathrm{par}}_{\mathfrak{sl},0}
 \,:=\,
 \left\{ a\,\in\, {\mathcal End}(\widetilde{E})\,
 \middle|\,\,
 \text{$\mathrm{Tr}(a)\,=\,0$ and 
 $a|_{x_i\times\mathcal{N}^{n_0\text{-reg}}_{\mathrm{par}}(L)}(\widetilde{l}^{(i)}_j)
 \,\subset\, \widetilde{l}^{(i)}_j$
 for any $i,\,j$} \right\},
\]
then the subbundle
$\mathrm{At}_D(\widetilde{E},\,\widetilde{\boldsymbol{l}},\,\nabla_L)
\,\subset\, \mathrm{At}_D(\widetilde{E},\,\widetilde{\boldsymbol{l}})$
fits in the exact sequence
\[
 0 \,\longrightarrow\, 
 {\mathcal End}_{\mathrm{par},\mathfrak{sl}}(\widetilde{E},\,\widetilde{\boldsymbol{l}})
 \,\longrightarrow\,
\mathrm{At}_D(\widetilde{E},\,\widetilde{\boldsymbol{l}},\,\nabla_L)
\, \xrightarrow{\,\,\mathrm{symb}_1\,\,}\,
 T_X(-D)\otimes {\mathcal O}_{X\times\mathcal{N}^{n_0\text{-reg}}_{\mathrm{par}}(L)}
\, \longrightarrow\, 0.
\]
If we set
\[
 \widetilde{\mathcal D}^{\mathrm{par}}_{\mathfrak{sl},1}
 \,:=\,
 \left\{ a\,\in \,{\mathcal End}(\widetilde{E})\otimes K_X(D)\,
 \middle|\,\,
 \text{$\mathrm{Tr}(a)\,=\,0$ and 
 $\mathrm{res}_{x_i\times\mathcal{N}^{n_0\text{-reg}}_{\mathrm{par}}(L)}(a) (\widetilde{l}^{(i)}_j)
 \,\subset\, \widetilde{l}^{(i)}_{j+1}$
 for any $i,\,j$} \right\},
\]
then, by Serre duality,
\[
H^1\big(X,\,\widetilde{\mathcal D}^{\mathrm{par}}_{\mathfrak{sl},0}
\big|_{X\times \{p\}}\otimes K_X(D)\big)^{\vee}
\,\cong\,
H^0\big(X,\, \widetilde{\mathcal D}^{\mathrm{par}}_{\mathfrak{sl},1}
\big|_{X\times \{p\}}\otimes T_X(-D)\big)
\,\subset\, H^0\big(X,\, \widetilde{\mathcal D}^{\mathrm{par}}_{\mathfrak{sl},0}
\big|_{X\times \{p\}}\big),
\]
which in fact becomes zero because the underlying quasi-parabolic bundle
$(\widetilde{E},\,\widetilde{\boldsymbol{l}})|_{X\times\{p\}}$
is simple. Let 
$$\pi \,\,\colon\,\, X\times {\mathcal N}^{n_0\text{-reg}}_{\rm par}(L)
\,\,\longrightarrow \,\, {\mathcal N}^{n_0\text{-reg}}_{\rm par}(L)$$
be the projection.
Then we have $R^1\pi_*
\left(\widetilde{\mathcal D}^{\mathrm{par}}_{\mathfrak{sl},0}
\otimes K_X(D) \right)\,=\,0$,
and get a short exact sequence
\begin{align*}
 0 \,\longrightarrow\, 
 \pi_* \left( \widetilde{\mathcal D}^{\mathrm{par}}_{\mathfrak{sl},0} \otimes K_X(D) \right)
 \,\longrightarrow\,
 & \pi_*\left(\mathrm{At}_D(\widetilde{E},\,\widetilde{\boldsymbol{l}},\,\nabla_L)\otimes K_X(D)\right)
 \,\xrightarrow{\,\mathrm{symb}_1\otimes\mathrm{id}\,\,}\,
 \pi_*\left( {\mathcal O}_{X\times\mathcal{N}^{n_0\text{-reg}}_{\mathrm{par}}(L)}\right)
 \,\longrightarrow \,0.
\end{align*}
Note that
$\pi_*\left( {\mathcal O}_{X\times\mathcal{N}^{n_0\text{-reg}}_{\mathrm{par}}(L)}\right)
\,\cong\,{\mathcal O}_{\mathcal{N}^{n_0\text{-reg}}_{\mathrm{par}}(L)}$.
Consider the homomorphism
\[
 \Psi_t\,\colon\,\,
 \left(\pi_*\left(\mathrm{At}_D(\widetilde{E},\,\widetilde{\boldsymbol{l}},\,\nabla_L)\otimes K_X(D)\right)
 \, \oplus \, {\mathcal O}_{\mathcal{N}^{n_0\text{-reg}}_{\mathrm{par}}(L)} \right)
 \otimes k[t]\,\,\longrightarrow\,\,
 {\mathcal O}_{\mathcal{N}^{n_0\text{-reg}}_{\mathrm{par}}(L)}\otimes k[t]
\]
on $\mathcal{N}^{n_0\text{-reg}}_{\mathrm{par}}(L)\times\Spec k[t]$
defined by
$$(u,\,f)\ \longmapsto\ (\mathrm{symb}_1\otimes\mathrm{id}_{K_X(D)})(u)-tf$$ for 
$u\,\in\,\pi_*\left(\mathrm{At}_D(\widetilde{E},\,\widetilde{\boldsymbol{l}},\,\nabla_L)\otimes K_X(D)\right)
\otimes k[t]$ and
$f\,\in\, {\mathcal O}_{\mathcal{N}^{n_0\text{-reg}}_{\mathrm{par}}(L)}\otimes k[t]$.
Then $\ker\Psi_t$ is a locally free sheaf on 
$\mathcal{N}^{n_0\text{-reg}}_{\mathrm{par}}(L)\times\Spec k[t]$, and we have
\[
 \ker\Psi_t\otimes k[t]/(t-h)\,\,\cong\,\,
 \begin{cases}
 \pi_*\left(\mathrm{At}_D(\widetilde{E},\,\widetilde{\boldsymbol{l}},\,\nabla_L)\otimes K_X(D)\right)
 & (h\neq 0) \\
 \pi_* \left(\widetilde{\mathcal D}^{\mathrm{par}}_{\mathfrak{sl},0}
 \otimes K_X(D) \right)
 \oplus {\mathcal O}_{\mathcal{N}^{n_0\text{-reg}}_{\mathrm{par}}(L)}
 & (h\,=\,0).
 \end{cases}
\]
Define the projective bundle
\[
 \mathbb{P}_*(\ker\Psi_t)\,:=\,
 \Proj \left( \mathrm{Sym} \left( (\ker\Psi_t)^{\vee} \right) \right)
\]
over $\mathcal{N}^{n_0\text{-reg}}_{\mathrm{par}}(L)\times\Spec k[t]$.
There is a tautological line-subbundle
\[
 {\mathcal O}_{\mathbb{P}_*(\ker\Psi_t)}(-1)
\,\, \hookrightarrow\,
 \ker\Psi_t\otimes{\mathcal O}_{\mathbb{P}_*(\ker\Psi_t)}.
\]
Consider the sections
\begin{align*}
 y \,\colon\,
 {\mathcal O}_{\mathbb{P}_*(\ker\Psi_t)}(-1)
 &\hookrightarrow
 \ker\Psi_t\otimes{\mathcal O}_{\mathbb{P}_*(\ker\Psi_t)}
 \\
 &\hookrightarrow
 \pi_*\left(\mathrm{At}_D(\widetilde{E},\,\widetilde{\boldsymbol{l}},\,\nabla_L)\otimes K_X(D)\right)
 \otimes {\mathcal O}_{\mathbb{P}_*(\ker\Psi_t)}
 \, \oplus \, {\mathcal O}_{\mathbb{P}_*(\ker\Psi_t)}
 \,\longrightarrow\,
 {\mathcal O}_{\mathbb{P}_*(\ker\Psi_t)}
 \\
 \widetilde{\nu}^{(i)}_j \,\colon\,
 {\mathcal O}_{\mathbb{P}_*(\ker\Psi_t)}(-1)
 &\hookrightarrow\,
 \ker\Psi_t\otimes{\mathcal O}_{\mathbb{P}_*(\ker\Psi_t)}
 \\
 &\hookrightarrow
 \pi_*\left(\mathrm{At}_D(\widetilde{E},\,\widetilde{\boldsymbol{l}},\,\nabla_L)\otimes K_X(D)\right)
 \otimes {\mathcal O}_{\mathbb{P}_*(\ker\Psi_t)}
 \, \oplus \, {\mathcal O}_{\mathbb{P}_*(\ker\Psi_t)}
 \\
 &\longrightarrow
 \pi_*\left(\mathrm{At}_D(\widetilde{E},\,\widetilde{\boldsymbol{l}},\,\nabla_L)\otimes K_X(D)\right)
 \otimes {\mathcal O}_{\mathbb{P}_*(\ker\Psi_t)} \\
 &\xrightarrow{\mathrm{res}_D} \
 \pi_*(\widetilde{\mathcal D}^{\mathrm{par}}_{\mathfrak{sl},0}
 |_{D\times {\mathbb{P}_*(\ker\Psi_t)}} )
 \ \longrightarrow \
 \pi_* ( {\mathcal End}(\widetilde{l}^{(i)}_j/\widetilde{l}^{(i)}_{j+1})
 \otimes{\mathcal O}_{\mathbb{P}_*(\ker\Psi_t)} )
 \,=\,{\mathcal O}_{\mathbb{P}_*(\ker\Psi_t)}.
\end{align*}
Let $I$ be the ideal sheaf of the graded algebra
$\mathrm{Sym} \left( (\ker\Psi_t)^{\vee} \right)$
over $\mathcal{N}^{n_0\text{-reg}}_{\mathrm{par}}(L)$,
which is generated by 
$\left\{ \widetilde{\nu}^{(i)}_j-\nu^{(i)}_j ty \, \middle| \, 1\leq i\leq n , 0\leq j\leq r-1 \right\}$.
Set
\[
 \overline{\mathcal M'}\,:=\,
 \Proj \big(\mathrm{Sym}\left( \ker\Psi_t ^{\vee}\right)/I \big)
 \,\subset\, \mathbb{P}_*(\ker\Psi_t).
\]
Then there is a canonical structure morphism
\[
\overline{\mathcal M'} \ \longrightarrow\ \Spec k[t].
\]
Let $Y\,\subset\,\overline{\mathcal M'}$ be the effective divisor
defined by the equation $y\,=\,0$. Setting
${\mathcal M'}\,:=\,\overline{\mathcal M'}\setminus Y$,
we see by the construction that
${\mathcal M'}_h\,\cong\, \mathcal{M}_{\mathrm{PC}}^{n_0\text{\rm -reg}}(\boldsymbol{\nu},\nabla_L)^{\circ}$
for $h\,\neq\, 0$ and
${\mathcal M'}_0\,\cong\, \mathcal{M}^{n_0\text{-\rm reg}}_{\mathrm{Higgs}}(\boldsymbol{0},\,0)^{\circ}$.
\end{proof}

\begin{Thm}\label{Thm: global algebraic functions on connection moduli}
Let $X$ be a smooth projective curve of genus $g$
over an algebraically closed field $k$ of arbitrary characteristic,
and let $D\,=\,\sum_{i=1}^n x_i$ be a reduced effective divisor on $X$.
Fix a line bundle $L$ over $X$ with a connection
$\nabla_L \,\colon\, L\,\longrightarrow\, L\otimes K_X(D)$.
Take positive integers $r$ and $d$ 
such that $r\,\geq\, 2$, $n\,\geq\, 1$
and $r$ is not divisible by the characteristic of $k$.
Assume that one of the following statements holds:
\begin{itemize}
\item
$n\,\geq\, 1$ and $r\,\geq\, 2$ are arbitrary if $g\,\geq \,2$,

\item $n\,\geq\, 2$ and $n+r\,\geq\, 5$ if $g\,=\,1$,

\item $n\,\geq\, 3$ and $n+r\,\geq\, 7$ if $g\,=\,0$.
\end{itemize}
Also, assume that the exponent
$\boldsymbol{\nu}=(\nu^{(i)}_j)^{1\leq i\leq n}_{0\leq j\leq r-1}$
satisfies the condition $\res_{x_i}(\nabla_L)\,=\,\sum_{j=0}^{r-1}\nu^{(i)}_j$ for any $i$, and furthermore,
$
 \sum_{i=1}^n\sum_{\ell=1}^s \nu^{(i)}_{j^{(i)}_\ell}
$
is not contained in the image of $\mathbb{Z}$ in $k$
for any integer $1\,\leq\, s\,<\, r$ and any choice
of $s$ elements $\{j^{(i)}_1,\,\cdots,\,j^{(i)}_s\}$
in $\{1,\,\ldots,\,r\}$, for each $1\,\leq\, i\,\leq\, n$.
Then the transcendence degree of the global algebraic functions on
the moduli space
$\mathcal{M}_{\mathrm{PC}}^{\boldsymbol{\alpha}}(\boldsymbol{\nu})$
of $\boldsymbol{\alpha}$-stable $\boldsymbol{\nu}$-parabolic connections
satisfies the inequality
\[
 \mathrm{tr.deg}_{\,k} \, \Gamma\big(
 \mathcal{M}_{\mathrm{PC}}^{\boldsymbol{\alpha}}(\boldsymbol{\nu}),\,
 {\mathcal O}_{ \mathcal{M}_{\mathrm{PC}}^{\boldsymbol{\alpha}}(\boldsymbol{\nu}) } \big)
\, \,\leq \,\, r^2(g-1)-g+1+\frac{nr(r-1)} {2}.
\]
\end{Thm}

\begin{proof}
Note that
$
 \Gamma \left(
 \mathcal{M}^{n_0\text{-\rm reg}}_{\mathrm{Higgs}}(\boldsymbol{0},\,0)^{\circ},\,
 {\mathcal O}_{ \mathcal{M}^{n_0\text{-\rm reg}}_{\mathrm{Higgs}}(\boldsymbol{0},\,0)^{\circ}}\right)
\, =\,
 \Gamma \left(
 \mathcal{M}^{n_0\text{-\rm reg}}_{\mathrm{Higgs}}(\boldsymbol{0},\,0) ,\,
 {\mathcal O}_{\mathcal{M}^{n_0\text{-\rm reg}}_{\mathrm{Higgs}}(\boldsymbol{0},\,0) } \right)
$
by Proposition \ref{Higgs moduli codimension for underlying parabolic bundle}.
Since we can extend the Hitchin map in \eqref{equation: Hitchin map}
to a morphism
$
 \mathcal{M}^{n_0\text{-\rm reg}}_{\mathrm{Higgs}}(\boldsymbol{0},\,0)
 \,\longrightarrow\, W$, we have the inclusion maps
\[
 \Gamma(W,\,{\mathcal O}_W)
 \,\subset\, \Gamma \left(
 \mathcal{M}^{n_0\text{-\rm reg}}_{\mathrm{Higgs}}(\boldsymbol{0},\,0) ,\,
 {\mathcal O}_{ \mathcal{M}^{n_0\text{-\rm reg}}_{\mathrm{Higgs}}(\boldsymbol{0},\,0) } \right)
\, \subset\,
 \Gamma \left(
 \mathcal{M}^{\boldsymbol{\alpha'}}_{\mathrm{Higgs}}(\mathbf{0},\,0),\,
 {\mathcal O}_{\mathcal{M}^{\boldsymbol{\alpha'}}_{\mathrm{Higgs}}(\mathbf{0},\,0) } \right),
\]
where we take $\boldsymbol{\alpha'}$ generic so 
that $\boldsymbol{\alpha'}$-semistability implies $\boldsymbol{\alpha'}$-stability.
Then using Corollary \ref{Cor: global algebraic functions on Higgs moduli space} it follows that
$ \Gamma \left(
 \mathcal{M}^{n_0\text{-\rm reg}}_{\mathrm{Higgs}}(\boldsymbol{0},\,0)^{\circ},\,
 {\mathcal O}_{ \mathcal{M}^{n_0\text{-\rm reg}}_{\mathrm{Higgs}}(\boldsymbol{0},\,0)^{\circ}}\right)
 \,=\,
\Gamma \left(
 \mathcal{M}^{n_0\text{-\rm reg}}_{\mathrm{Higgs}}(\boldsymbol{0},\,0) ,\,
 {\mathcal O}_{ \mathcal{M}^{n_0\text{-\rm reg}}_{\mathrm{Higgs}}(\boldsymbol{0},\,0) } \right)$
is a finitely generated $k$-algebra whose Krull dimension is
$r^2(g-1)-g+1+nr(r-1)/2$. 

We use the notation in the proof of Proposition \ref{Prop: Deligne--Hitchin family}.
Note that
$\mathcal{M}^{n_0\text{-\rm reg}}_{\mathrm{Higgs}}(\boldsymbol{0},\,0)^{\circ}$
is isomorphic to the cotangent bundle over
$\mathcal{N}^{n_0\text{-reg}}_{\mathrm{par}}(L)$.
So we have
$
 \mathcal{M}^{n_0\text{-\rm reg}}_{\mathrm{Higgs}}(\boldsymbol{0},\,0)^{\circ}
 \,\cong\,
 \mathrm{Spec} \Big( \mathrm{Sym}^* \Big( \pi_* \big( 
 \widetilde{\mathcal D}^{\mathrm{par}}_{\mathfrak{sl},1}\big)^{\vee}\Big) \Big)$,
which implies that
\begin{equation*} 
 \Gamma \left(
 \mathcal{M}^{n_0\text{-\rm reg}}_{\mathrm{Higgs}}(\boldsymbol{0},\,0)^{\circ} ,\,
 {\mathcal O}_{ \mathcal{M}^{n_0\text{-\rm reg}}_{\mathrm{Higgs}}(\boldsymbol{0},\,0)^{\circ} } \right)
 \,\cong\,
 \bigoplus_{m=0}^{\infty}
 H^0\left( \mathcal{N}^{n_0\text{-reg}}_{\mathrm{par}}(L) ,\,
 \mathrm{Sym}^m \left( \pi_* \big( 
\widetilde{\mathcal D}^{\mathrm{par}}_{\mathfrak{sl},1} \big) ^{\vee} \right) \right).
\end{equation*}
Note that there is a short exact sequence
\[
 0\,\longrightarrow\,
 \pi_*\big(\widetilde{\mathcal D}^{\mathrm{par}}_{\mathfrak{sl},1} \big)
 \,\longrightarrow\,
 \pi_* \big( 
\widetilde{\mathcal D}^{\mathrm{par}}_{\mathfrak{sl},0}
 \otimes K_X(D) \big)\,
 \stackrel{q}{\longrightarrow}\,
 \bigoplus_{i,j}
 \pi_*\big({\mathcal End}(\widetilde{l}^{(i)}_j/\widetilde{l}^{(i)}_{j+1})\big)
 \,\longrightarrow \,0.
\]
We can see that the above homomorphism $q$ determines
the equalities $\big( \widetilde{\nu}^{(i)}_j-\nu^{(i)}_jty \big)|_{t=0}$
on the fiber $\mathbb{P}(\ker\Psi_t^{\vee}\otimes\mathbb{C}[t]/(t))$
over $t\,=\,0$. Taking the dual of the above exact sequence,
\begin{align*}
 \left(\mathrm{Sym}\left( \ker\Psi_t ^{\vee}\right)/I\right) \otimes k[t]/(t) 
 &\,\cong\,
 \mathrm{Sym}\left(\pi_*\big(
\widetilde{\mathcal D}^{\mathrm{par}}_{\mathfrak{sl},1} \big)^{\vee}
 \oplus {\mathcal O}_{\mathcal{N}^{n_0\text{-reg}}_{\mathrm{par}}(L)} \right)
 \,\cong\,
 \bigoplus_{d=0}^{\infty} \bigoplus_{d_1+d_2=d} \mathrm{Sym}^{d_1}
 \left(\pi_*\big( \widetilde{\mathcal D}^{\mathrm{par}}_{\mathfrak{sl},1}
 \big)^{\vee}\right) T^{d_2},
\end{align*}
where $T$ is a variable corresponding to the second component of
$\ker\Psi_t\otimes k[t]/(t)
\,=\,\pi_*\big( \widetilde{\mathcal D}^{\mathrm{par}}_{\mathfrak{sl},1} \big)^{\vee}
\oplus {\mathcal O}_{{\mathcal N}^{n_0\text{-\rm reg}}_{\mathrm{par}}(L)}$.
So the ring of global sections of 
this sheaves of algebras
over $\mathcal{N}^{n_0\text{-reg}}_{\mathrm{par}}(L)$ becomes
a polynomial ring
\[
 \Gamma \Big(
 \mathcal{N}^{n_0\text{-reg}}_{\mathrm{par}}(L) ,\,
 \left(\mathrm{Sym}\left( \ker\Psi_t ^{\vee}\right)/I\right) \otimes k[t]/(t) \Big)
 \,\cong\,
 \Gamma \left(
 \mathcal{M}^{n_0\text{-\rm reg}}_{\mathrm{Higgs}}(\boldsymbol{0},\,0)^{\circ} ,\,
 {\mathcal O}_{ \mathcal{M}^{n_0\text{-\rm reg}}_{\mathrm{Higgs}}(\boldsymbol{0},\,0)^{\circ} } \right)
 [T]
\]
over $\Gamma \left(
 \mathcal{M}^{n_0\text{-\rm reg}}_{\mathrm{Higgs}}(\boldsymbol{0},\,0)^{\circ} ,\,
 {\mathcal O}_{ \mathcal{M}^{n_0\text{-\rm reg}}_{\mathrm{Higgs}}(\boldsymbol{0},\,0)^{\circ} } \right)$.
In particular, 
$ \dim ((\mathrm{Sym}^m (\ker\Phi_t^{\vee})/I_m)\otimes k[t]/(t))$
becomes a polynomial in $m$ of degree
$$\text{Krull-dim}\,\Gamma \left(
 \mathcal{M}^{n_0\text{-\rm reg}}_{\mathrm{Higgs}}(\boldsymbol{0},\,0)^{\circ} ,\,
 {\mathcal O}_{ \mathcal{M}^{n_0\text{-\rm reg}}_{\mathrm{Higgs}}(\boldsymbol{0},\,0)^{\circ} } \right)
 \ =\ r^2(g-1)-g+1+\frac{nr(r-1)} {2}.$$

Let $\left(\mathrm{Sym}\left( \ker\Psi_t^{\vee} \right)/I \right)_{(y)}$
be the subalgebra of the localized graded algebra
$\left(\mathrm{Sym}\left( \ker\Psi_t^{\vee} \right)/I \right)_y$
consisting of homogeneous elements of degree zero.
Then we have 
\[
 \mathcal{M}_{\mathrm{PC}}^{n_0\text{-reg}}(\boldsymbol{\nu},\nabla_L)^{\circ}\ \cong\
 \Spec \left( \left(\mathrm{Sym}\left( \ker\Psi_t^{\vee} \right)/I \right)_{(y)}
 \otimes k[t]/(t-h) \right)
\]
for $h\,\neq\, 0$.
By the assumption in \eqref{assumption of local exponents} on the choice of
the exponent $\boldsymbol{\nu}$,
and by Proposition \ref{proposition: codimension in de Rham moduli},
we have
\begin{align*}
 \Gamma( \mathcal{M}_{\mathrm{PC}}^{\boldsymbol{\alpha}}(\boldsymbol{\nu},\nabla_L),\,
{\mathcal O}_{\mathcal{M}_{\mathrm{PC}}^{\boldsymbol{\alpha}}(\boldsymbol{\nu},\nabla_L)} )
 &\,=\,
\Gamma( \mathcal{M}_{\mathrm{PC}}^{n_0\text{-reg}}(\boldsymbol{\nu},\nabla_L),\,
 {\mathcal O}_{\mathcal{M}_{\mathrm{PC}}^{n_0\text{-reg}}(\boldsymbol{\nu},\nabla_L)} )
 \\
 &=\,
 \Gamma( \mathcal{M}_{\mathrm{PC}}^{n_0\text{-reg}}(\boldsymbol{\nu},\nabla_L)^{\circ} ,\,
 {\mathcal O}_{\mathcal{M}_{\mathrm{PC}}^{n_0\text{-reg}}(\boldsymbol{\nu},\nabla_L)^{\circ} })
 \\
 &=\,
\Gamma \left(\mathcal{N}^{n_0\text{-reg}}_{\mathrm{par}}(L) , \,
\left(\mathrm{Sym}\left( \ker\Psi_t^{\vee} \right)/I \right)_{(y)}
\otimes k[t]/(t-h) \right).
\end{align*}
By Lemma \ref{lemma: upper semi-continuity of global sections}
which is proved later in Section \ref{section: appendix}, the function 
\begin{align*}
 h\; \longmapsto \; 
 &\dim H^0\left(\mathcal{N}^{n_0\text{-reg}}_{\mathrm{par}}(L) , \,
 \mathrm{Sym}^m(\ker\Psi_t^{\vee})/I_m\big|_{t=h}\right)
 \\
 &=\,
 \dim H^0\left(\mathcal{N}^{n_0\text{-reg}}_{\mathrm{par}}(L) , \,
 \mathrm{Sym}^m(\ker\Psi_t^{\vee})/I_m\otimes k[t]/(t-h)\right)
\end{align*}
is upper semi-continuous in $h$.
So we have
\begin{align*}
 &\dim H^0\left( \mathcal{N}^{n_0\text{-reg}}_{\mathrm{par}}(L) , \,
 \mathrm{Sym}^m(\ker\Psi_t^{\vee})/I_m\big|_{t=h}\right)
 \,\leq\,
 \dim H^0\left(\mathcal{N}^{n_0\text{-reg}}_{\mathrm{par}}(L), \,
 \mathrm{Sym}^m(\ker\Psi_t^{\vee})/I_m\big|_{t=0}\right)
\end{align*}
for $h\,\neq\, 0$. 

Let $d$ be the transcendence degree of
$\Gamma( \mathcal{M}_{\mathrm{PC}}^{\boldsymbol{\alpha}}(\boldsymbol{\nu},\nabla_L),\,
{\mathcal O}_{\mathcal{M}_{\mathrm{PC}}^{\boldsymbol{\alpha}}(\boldsymbol{\nu},\nabla_L)} )$
over $k$.
Then we have
\begin{align*}
 d&\,=\,
 \mathrm{tr.deg}_{\,k}\Gamma \left( \mathcal{N}^{n_0\text{-reg}}_{\mathrm{par}}(L), \,
 \left(\mathrm{Sym}\left( \ker\Psi_t \right)/I\right)_{(y)}
 \otimes k[t]/(t-h) \right)
 \\
 &=\,
 \mathrm{tr.deg}_{\,k}
 \Gamma \left( \mathcal{N}^{n_0\text{-reg}}_{\mathrm{par}}(L),\,
 \left(\mathrm{Sym}\left( \ker\Psi_t \right)/I\right)
 \otimes k[t]/(t-h) \right) -1.
\end{align*}
Take homogeneous elements $x_1,\,\cdots,\,x_d$ of
$\Gamma \left( \mathcal{N}^{n_0\text{-reg}}_{\mathrm{par}}(L), \,
 \left(\mathrm{Sym}\left( \ker\Psi_t \right)/I\right)
 \big|_{t=h} \right)$
such that $\{x_1,\,\cdots,\,x_d,\,y\}$ is 
a transcendence basis of
$\Gamma \left( \mathcal{N}^{n_0\text{-reg}}_{\mathrm{par}}(L), \,
 \left(\mathrm{Sym}\left( \ker\Psi_t \right)/I\right)
 \big|_{t=h} \right)$
over $k$. Let $S$ be the graded subalgebra of
$\Gamma \left( \mathcal{N}^{n_0\text{-reg}}_{\mathrm{par}}(L), \,
 \left(\mathrm{Sym}\left( \ker\Psi_t \right)/I\right)
 \big|_{t=h} \right)$
generated by $x_1,\,\cdots,x_d,\,y$.
Then
\begin{align*}
 \dim S_m
 &\,\leq\,
 \dim H^0\left(\mathcal{N}^{n_0\text{-reg}}_{\mathrm{par}}(L), \,
 \mathrm{Sym}^m(\ker\Psi_t^{\vee})/I_m\big|_{t=h} \right)
 \\
 &\leq\,
 \dim H^0\left(\mathcal{N}^{n_0\text{-reg}}_{\mathrm{par}}(L), \,
 \mathrm{Sym}^m(\ker\Psi_t^{\vee})/I_m\big|_{t=0} \right).
\end{align*}
Since $S_m$ is a polynomial in $m$ of degree
$d$ for $m\,\gg\, 0$, it follows that $d\,\leq\, r^2(g-1)-g+1+nr(r-1)/2$.
\end{proof}

\begin{Rem}
A statement similar to Theorem \ref{Thm: global algebraic functions on connection moduli} can be considered for connections without pole.
When $X$ is a curve over the field of complex numbers whose genus
is greater than $2$, then there are only constant global algebraic functions on the de Rham moduli space of connections without pole by 
\cite[Corollary 4.4]{BR}.
So the inequality similar to Theorem \ref{Thm: global algebraic functions on connection moduli}
becomes strict in that case.
On the other hand, if $X$ is defined over the base field of positive characteristic, it is proved in \cite[Theorem 1.1]{Groech}
that the Hitchin map for the de Rham moduli space connections without pole is \'etale locally equivalent to that on the Dolbeault moduli space.
So the ring of global algebraic functions on the de Rham moduli space
has the same transcendence degree as that of the ring of global algebraic
functions on the Dolbeault moduli space in that case.
The Hitchin map for the logarithmic de Rham moduli space
over the base field of positive characteristic is introduced in \cite{dC-H-Z}. 
\end{Rem}

The following is an immediate consequence of Theorem \ref{Thm: global algebraic functions on connection 
moduli}.

\begin{Cor}\label{Cor: connection moduli is not affine}
The moduli space
$\mathcal{M}_{\mathrm{PC}}^{\boldsymbol{\alpha}}(\boldsymbol{\nu},\nabla_L)$ 
of $\boldsymbol{\alpha}$-stable $\boldsymbol{\nu}$-parabolic connections is not affine.
\end{Cor}

From now on, consider the case of $k\,=\,\mathbb{C}$.

Since the fundamental group
$\pi_1(X\setminus D,\,*)$ is finitely presented, the space of representations
\[
 \Hom(\pi_1(X\setminus D,*),\, \mathrm{GL}(r,\mathbb{C}))
\]
can be realized as an affine variety.
Take generators $\alpha_1,\,\beta_1,\,\cdots,\,\alpha_g,\,\beta_g$
of the fundamental group $\pi_1(C,\, *)$,
and choose a loop $\gamma_i$ around each $x_i$ with respect to the base point $*$.
Then the fundamental group $\pi_1(X\setminus D,\,*)$ is generated by
$\alpha_1,\,\beta_1,\,\cdots,\,\alpha_g,\,\beta_g,\,\gamma_1,\,\cdots,\,\gamma_n$
with the single relation $[\alpha_1,\,\beta_1]\cdots[\alpha_g,\,\beta_g]\gamma_1\cdots\gamma_n\,=\,1$.
The space of representations of $\pi_1(X\setminus D,\,*)$ can be realized as the affine variety
\begin{align*}
 &\Hom(\pi_1(X\setminus D,\,*),\,\, \mathrm{GL}(r,\mathbb{C}))
 \\
 &=\,
 \left\{ (A_1,B_1,\ldots,A_g,B_g,C_1,\ldots,C_n)\,\in\, \mathrm{GL}(r,\mathbb{C})^{2g+n}
\,\, \middle|\,\, A_1^{-1}B_1^{-1}A_1B_1\cdots A_g^{-1}B_g^{-1}A_gB_gC_1\cdots C_n\,=\,I_r \right\}.
\end{align*}
Note that the connection $\nabla_L$ on the line bundle $L$ induces
a one dimensional representation $\rho_{\nabla_L}$ of $\pi_1(X\setminus D,\,*)$.
Define a tuple $(b^{(i)}_j)$ by $b^{(i)}_j\,:=\,e^{-2\pi\sqrt{-1}\nu^{(i)}_j}$,
and consider the closed subvariety
\[
 Y\,=\,
 \left\{ ((A_k,B_k),(C_i)) \,\in\, 
 \Hom(\pi_1(X\setminus D,*),\mathrm{GL}(r,\mathbb{C}))
\, \middle|\, 
 \begin{array}{l}
 \text{$\rho_{\nabla_L}(\alpha_k)=\det(A_k)$ and
 $\rho_{\nabla_L}(\beta_k)=\det(B_k)$}
 \\
 \text{ for $1\leq i\leq g$ and}
 \\
 \det(TI_r-C_i)=\prod_{j=0}^{r-1} (T-b^{(i)}_j) 
 \ \text{for $1\leq i \leq n$} 
 \end{array}\right\}
\]
of $\Hom(\pi_1(X\setminus D,\,*),\,\mathrm{GL}(r,\mathbb{C}))$.
There is a canonical action of $\mathrm{GL}(r,\mathbb{C})$
on $Y$ given by the adjoint action of $\mathrm{GL}(r,\mathbb{C})$ on
itself, and we can take the corresponding categorical quotient
\begin{equation}\label{character variety with fixed local monodromy}
 \mathrm{Ch}_{X\setminus D,(b^{(i)}_j)}\,\,:=\,\,
 Y/\!\!/\mathrm{GL}(r,\mathbb{C})
\,\,=\,\,\Spec \Gamma(Y,\, {\mathcal O}_Y)^{\mathrm{GL}(r,\mathbb{C})}.
\end{equation}
Under the genericity assumption in \eqref{assumption of local exponents} of the eigenvalues of the residues,
this quotient is in fact a geometric quotient, and we have a Riemann--Hilbert morphism
\[
 \mathrm{RH} \,\,\colon\,\,
 \mathcal{M}_{\mathrm{PC}}^{\boldsymbol{\alpha}}(\boldsymbol{\nu},\nabla_L)
 \,\, \longrightarrow\,\,
 \mathrm{Ch}_{X\setminus D,(b^{(i)}_j)}.
\]
By \cite{Inaba-1}, the above Riemann--Hilbert morphism $\mathrm{RH}$
is a proper and surjective holomorphic map, which is generically an isomorphism.
So $\mathcal{M}_{\mathrm{PC}}^{\boldsymbol{\alpha}}(\boldsymbol{\nu},\nabla_L)$
gives an analytic resolution of singularities of
$\mathrm{Ch}_{X\setminus D,(b^{(i)}_j)}$.
Since the character variety
$\mathrm{Ch}_{X\setminus D , (b^{(i)}_j)}$
is affine by its definition,
it is evident that
\begin{equation}\label{equation: global algebraic functions on character variety}
 \mathrm{tr.deg}_{\mathbb{C}} \Gamma(\mathrm{Ch}_{X\setminus D , (b^{(i)}_j)},\,
 {\mathcal O}_{\mathrm{Ch}_{X\setminus D , (b^{(i)}_j)}} )
 \,=\,
 \dim \mathrm{Ch}_{X\setminus D , (b^{(i)}_j)} \,=\,
 2(r^2-1)(g-1)+r(r-1)n.
\end{equation}

By Theorem \ref{Thm: global algebraic functions on connection moduli}
and \eqref{equation: global algebraic functions on character variety}
(or by Corollary \ref{Cor: connection moduli is not affine}),
we have the following:

\begin{Cor}\label{Cor: logarithmic Riemann-Hilbert is not algebraic}
The Riemann-Hilbert morphism
$\mathrm{RH}\,\colon\,
\mathcal{M}_{\mathrm{PC}}^{\boldsymbol{\alpha}}(\boldsymbol{\nu},\nabla_L)
\,\longrightarrow\, \mathrm{Ch}_{X\setminus D,(b^{(i)}_j)}$
is not an algebraic morphism.
\end{Cor}

\section{Appendix}\label{section: appendix}

Let $k$ be an algebraically closed field of arbitrary characteristic.
We will prove a lemma on the upper semi-continuity of the dimension
of global sections of vector bundles on an
algebraic space containing a projective variety over $k$.

Recall that an algebraic space ${\mathcal X}$ of finite type over $\Spec k$
is said to be locally separated over $\Spec k$
if there is a scheme $U$ of finite type over $\Spec k$
together with an \'etale surjective morphism
$U\,\longrightarrow\,{\mathcal X}$ such that
$U\times_{\mathcal X}U$ is a locally closed subscheme of
$U\times_{\Spec\mathbb{C}}U$.
A locally separated algebraic space ${\mathcal X}$ 
of finite type over $\Spec k$ is irreducible
if the underlying topological space $|{\mathcal X}|$ is irreducible.
In other words, any two non-empty open subspaces
$U_1,\,U_2\,\subset\,{\mathcal X}$ intersect: $U_1\cap U_2\,\neq\,\emptyset$.

\begin{Lem}\label{lemma: upper semi-continuity of global sections}
Let ${\mathcal X}$ be a locally separated smooth irreducible algebraic space,
which is of finite type over $\Spec k$.
Assume that $\overline{X}$ is an open subspace of 
${\mathcal X}$ such that $\overline{X}$ is isomorphic to a smooth projective variety over $k$.
Let $T$ be an affine variety, and let ${\mathcal F}$ be a locally free sheaf of finite rank on
${\mathcal X}\times T$.
For each point $t\in T$,\, denote by
$\Gamma({\mathcal X}\times\{t\} ,\, {\mathcal F}|_{{\mathcal X}\times\{t\}})$
the space of global sections of the restriction
${\mathcal F}|_{{\mathcal X}\times\{t\}}$.
Then the function
\[
T \ \longrightarrow\ \mathbb{Z}_{\geq 0},\ \ \ \, t \ \longmapsto \ 
 \dim \Gamma( {\mathcal X}\times\{t\} ,\, {\mathcal F}|_{{\mathcal X}\times\{t\}})
\]
is upper semi-continuous.
\end{Lem}

\begin{proof}
Since the upper semi-continuity is a local property on $T$,
we may replace $T$ with a neighborhood at any point of $T$.
Take a finite number of smooth affine varieties
$\{U_i\}_{i=1}^n$ and an \'etale surjective morphism
\[
 f\,\,\colon\,\, \overline{X} \sqcup \coprod_{i=1}^n U_i
 \ \longrightarrow \ {\mathcal X},
\]
whose restriction to $\overline{X}$ coincides with the given inclusion
map $f|_{\overline{X}}\,\,\colon\, \overline{X}\,\hookrightarrow\, {\mathcal X}$.
After shrinking $U_i$ and $T$, we may assume that
${\mathcal F}|_{U_i\times T}\,\cong\,{\mathcal O}_{U_i\times T}^{\oplus r}$ for every $i$.
Since ${\mathcal X}$ is irreducible, we have
$\overline{X}\cap(\bigcap_{i=1}^n f(U_i))\,\neq\,\emptyset$.
So there is a non-empty affine open subset
$V\,\subset\, \overline{X}\cap\bigcap_{i=1}^n f(U_i)$.
Take a non-empty smooth affine variety
$\widetilde{V}$ with \'etale morphisms 
$\widetilde{V}\,\longrightarrow\, V$, and $\widetilde{f}_i\,\colon\, \widetilde{V}\,\longrightarrow\, U_i$
for $1\,\leq\, i\,\leq\, n$, such that the diagram
\[
 \begin{CD}
\widetilde{V} @>\widetilde{f}_i >> U_i \\
 @VVV @VV f V \\
 V @>>> {\mathcal X}
 \end{CD}
\]
is commutative for every $1\,\leq\, i\,\leq\, n$.

Let $\widetilde{X}$ be the normalization of $\overline{X}$
in the field $K(\widetilde{V})$ of rational functions on $\widetilde{V}$.
Then $\widetilde{X}$ is a projective variety
with a canonical commutative diagram
\[
 \begin{CD}
 \widetilde{V} @>>> \widetilde{X} \\
 @VVV @VVV \\
 V @>>> \overline{X}.
 \end{CD}
\]
After shrinking $\widetilde{V}$ if necessary,
$\widetilde{V} \,\longrightarrow \,\widetilde{X}$ is an open immersion.
We can take a very ample divisor $D\,\subset\, \overline{X}$ such that
$\overline{X}\setminus V \,\subset\, D$.
Choose a very ample divisor $\widetilde{D}$ on $\widetilde{X}$
such that the inclusion $\widetilde{X}\setminus \widetilde{V}\,\subset\, \widetilde{D}$ holds set theoretically
and that $D\times_{\overline{X}}\widetilde{X}\,\subset\, \widetilde{D}$.

We can construct a projective variety $P_i$ with a very ample divisor
$D_i\,\subset\, P_i$ such that $P_i\setminus D_i$ is isomorphic to $U_i$.
We can also take a very ample divisor $D'_i\,\subset\, P_i$ such that
$P_i\setminus \widetilde{f}_i(\widetilde{V})\,\subset\, D'_i$ holds set theoretically
and that $D'_i\,=\,D_i+B_i$ holds
for a divisor $B_i$ without any common component with $D_i$.

For $i\,<\,j$, the fiber product $U_i
\times_{\mathcal X}U_j$ is a smooth quasi-affine scheme over $\Spec k$.
So we can construct a projective scheme
$P_{ij}$ over $\Spec k$,
which contains $U_i\times_{\mathcal X}U_j$ as a Zariski open subscheme.
Choose a very ample divisor $D_{ij}\,\subset\, P_{ij}$ such that
$P_{ij}\setminus(U_i\times_{\mathcal X}U_j)\,\subset\, D_{ij}$.

Since $\overline{X}$ is projective, and $D$ is very ample,
we can take a sufficiently large integer $l$ such that
$H^p(\overline{X}\times\{t\},\, {\mathcal F}|_{\overline{X}\times\{t\}}(lD))\,=\,0$
for all $p\,\geq\, 1$ and $t\,\in\, T$.
After shrinking $T$, the space of sections $\Gamma({\mathcal F}|_{\overline{X}\times T}(lD))$ 
is a free $\Gamma({\mathcal O}_T)$--module of finite rank
and the map $\Gamma({\mathcal F}|_{\overline{X}\times T}(lD))\otimes k(t)
\,\longrightarrow\,
\Gamma({\mathcal F}|_{\overline{X}\times \{t\}}(lD))$
is bijective for any $t\,\in\, T$,
where $k(t)$ is the residue field of ${\mathcal O}_{T,t}$.

Choose generators $s_1,\,\cdots,\,s_N$ of $\Gamma({\mathcal F}|_{\overline{X}\times T}(lD))$.
Consider the pullbacks of these sections by the morphism
$P_i\setminus D'_i\,\hookrightarrow\,
\widetilde{f}_i(\widetilde{V})\,\xrightarrow{\,\,(f|_{U_i})|_{\widetilde{f}_i(\widetilde{V})}\,\,\,}\,
V\,\hookrightarrow\, \overline{X}$
and denote them by
\[
 s_1|_{P_i\setminus D'_i},\ldots,s_N|_{P_i\setminus D'_i}
\, \,\in\,\, \Gamma({\mathcal F}|_{(U_i\setminus (U_i\cap D'_i))\times T})
\, \,\cong \,\,\Gamma({\mathcal O}_{(P_i\setminus D'_i)\times T}^{\oplus r}).
\]
There is a sufficiently large integer $l_i$ such that
each $s_1|_{P_i\setminus D'_i},\,\cdots,\,s_N|_{P_i\setminus D'_i}$
can be lifted to a section of
$\Gamma({\mathcal O}_{P_i\times T}(l_iD'_i))$.

After shrinking $T$, the space of sections $\Gamma({\mathcal O}_{P_i\times T}(l_iD_i))$
is a free $\Gamma({\mathcal O}_T)$--module of finite rank.
Fix a basis $t^{(i)}_1,\,\cdots,\,t^{(i)}_{N_i}$ of it.
Let $t^{(i)}_\ell\big|_{\widetilde{X}\setminus\widetilde{D}}$ be the pullback of 
$t^{(i)}_\ell$ by the composition of the maps
$$\widetilde{X}\setminus\widetilde{D}\,\hookrightarrow\, \widetilde{V}\,\longrightarrow\,
\widetilde{f}_i(\widetilde{V})\,\hookrightarrow\, U_i\,=\,P_i\setminus D_i\,\hookrightarrow\, P_i .$$
Then there is an integer $\widetilde{l}\,\geq\, l$ such that
all $t^{(i)}_1\big|_{\widetilde{X}\setminus\widetilde{D}},\, \cdots,\, t^{(i)}_{N_i}\big|_{\widetilde{X}\setminus\widetilde{D}}$ 
can be lifted to sections 
$\widetilde{t}^{(i)}_1,\,\cdots, \,\widetilde{t}^{(i)}_{N_i}$ of 
$\Gamma({\mathcal F}_{\widetilde{X}\times T}(\widetilde{l}\,\widetilde{D}))$
for $1\,\leq \,i\,\leq\, n$.

Consider the pullback $t^{(i)}_{\gamma}\big|_{P_{ij}\setminus D_{ij}}$ of $t^{(i)}_{\gamma}$ by the
composition of maps
\[
P_{ij}\setminus D_{ij}\,\hookrightarrow\, U_i\times_NU_j\,\longrightarrow\, U_i\,\hookrightarrow \,P_i.
\]
If we choose $l_{ij}$ sufficiently large, all
$t^{(i)}_1\big|_{P_{ij}\setminus D_{ij}},\,\cdots,\,t^{(i)}_{N_i}\big|_{P_{ij}\setminus D_{ij}}$
can be lifted to sections $t^{(i)}_{j,1},\,\cdots,\,t^{(i)}_{j,N_i}$ of
$\Gamma({\mathcal O}_{P_i\times T}(l_{ij}D_{ij})^{\oplus r})$.
We may also assume that all
$t^{(j)}_1\big|_{P_{ij}\setminus D_{ij}},\,\cdots,\,t^{(j)}_{N_j}\big|_{P_{ij}\setminus D_{ij}}$
can be lifted to sections $t^{(j)}_{i,1},\,\cdots,\,t^{(j)}_{i,N_j}$ of
$\Gamma({\mathcal O}_{P_i\times T}(l_{ij}D_{ij})^{\oplus r})$.

Take a resolution
\[
 {\mathcal L}_1 \,\xrightarrow{\,\, \partial_{{\mathcal L}_{\bullet}}\,\,\,}\, {\mathcal L}_0\,
\xrightarrow{\,\,\,\psi\,\,\,}\, {\mathcal F}|_{\overline{X}\times T}^{\vee}
\, \longrightarrow \,0,
\]
where ${\mathcal L}_i\,=\, {\mathcal O}_{\overline{X}\times T}(-m_i)^{\oplus R_i}$
for $i\,=\,1,\,2$ and $m_i\,\gg\, 1$. 
After shrinking $T$, both $\Gamma({\mathcal L}_i^{\vee})$ and
$\Gamma({\mathcal L}_i^{\vee}(lD))$ are free $\Gamma({\mathcal O}_T)$--modules
for $i\,=\,1,\,2$.
Let ${\mathcal F}|_{\widetilde{X}\times T}$ and ${\mathcal L}_i|_{\widetilde{X}\times T}$
respectively be the pullbacks of ${\mathcal F}|_{\overline{X}\times T}$ and
${\mathcal L}_i$ by the morphism
$\widetilde{X}\times T\,\longrightarrow \,\overline{X}\times T$.
Then there is the following commutative diagram with exact rows:
\[
 \begin{CD}
 0 @>>> \Gamma({\mathcal F}|_{\overline{X}\times T})
 @>>> \Gamma({\mathcal L}_0^{\vee})
 @>\Gamma(\partial_{{\mathcal L}_{\bullet}}^{\vee}) >> \Gamma({\mathcal L}_1^{\vee})
 \\
 & & @VVV @VVV @VVV \\
 0 @>>> \Gamma({\mathcal F}|_{\widetilde{X}\times T}(\widetilde{l}\,\widetilde{D}))
 @>>> \Gamma({\mathcal L}_0^{\vee}|_{\widetilde{X}\times T}(\widetilde{l}\,\widetilde{D}))
 @>>> \Gamma({\mathcal L}_1^{\vee}|_{\widetilde{X}\times T}(\widetilde{l}\,\widetilde{D})).
 \end{CD}
\]
Consider the homomorphism
\[
 \Phi\, \,\colon\,\,
 \Gamma({\mathcal L}_0^{\vee}) \oplus 
 \bigoplus_{i=1}^n \Gamma({\mathcal O}_{P_i\times T}(l_iD_i)^{\oplus r})
\,\, \longrightarrow\,\,
 \Gamma({\mathcal L}_1^{\vee})\oplus
 \Gamma(\widetilde{\mathcal L}_0^{\vee}(\widetilde{l}\widetilde{D}))^{\oplus n} \oplus
 \bigoplus_{i<j} \Gamma({\mathcal O}_{P_{ij}\times T}(l_{ij}D_{ij})^{\oplus r}),
\]
defined by
\[
\Big( \alpha,\, \Big( \sum_{\gamma=1}^{N_i} c^{(i)}_{\gamma} t^{(i)}_{\gamma}
\Big) \Big)\, \longmapsto\,
 \Big(
 \Gamma(\partial_{{\mathcal L}_{\bullet}}^{\vee}) (\alpha), \,
 \Big(\iota(\alpha)
 -\sum_{\gamma=1}^{N_i} c^{(i)}_{\gamma} \Gamma(\psi^{\vee}) (\widetilde{t}^{(i)}_{\gamma}) \Big)_i
 , \, \Big( \sum_{\gamma=1}^{N_i} c^{(i)}_{\gamma} t^{(i)}_{j,\gamma} 
 -\sum_{\gamma=1}^{N_j} c^{(j)}_{\gamma} t^{(j)}_{i,\gamma} \Big)_{i<j} \Big),
\]
where $\iota\,\colon\, \Gamma({\mathcal L}_0^{\vee})
\,\longrightarrow\, \Gamma({\mathcal L}_0^{\vee}(\widetilde{l}\,\widetilde{D}))$
is the canonical inclusion map and
$$\Gamma(\psi^{\vee})\,\colon\,
\Gamma({\mathcal F}|_{\widetilde{X}\times T}(\widetilde{l}\,\widetilde{D}))
\, \longrightarrow\,
\Gamma({\mathcal L}_0^{\vee}|_{\widetilde{X}\times T}(\widetilde{l}\,\widetilde{D}))$$
is the map induced by $\psi$.

\noindent
{\bf Claim.}
$\Gamma({{\mathcal X}\times\{t\},\,\, \mathcal F}|_{{\mathcal X}\times\{t\}})
\,\,=\,\, \ker (\Phi\otimes k(t))$
for any $t\,\in \,T$.

\noindent
{\it Proof of Claim.}\,
Take a section $s\,\in\, \Gamma({{\mathcal X}\times\{t\},\, \mathcal F}|_{{\mathcal X}\times\{t\}})$.
Its restriction $s|_{\overline{X}\times\{t\}}$
is a section of $\Gamma(\overline{X},\,{\mathcal F}|_{\overline{X}\times\{t\}})
\,\subset\, \Gamma(\overline{X},{\mathcal F}|_{\overline{X}\times\{t\}}(lD))$. From
the choice of $l_i$,
the pullback
$(f|_{U_i\cap\widetilde{f}_i(\widetilde{V})})^*(s|_{\overline{X}\times\{t\}})$
can be lifted to a section $\sigma_i$ of
$\Gamma({\mathcal O}_{P_i\times\{t\}}(l_iD'_i))$.
On the other hand, we have $(s|_{U_i})|_{\widetilde{f}_i(\widetilde{V})}
\,=\, (f|_{U_i\cap\widetilde{f}_i(\widetilde{V})})^*(s|_{\overline{X}\times\{t\}})$.
Since $s|_{U_i}$ does not have pole along $B_i$, it follows that $\sigma_i$ belongs to
$\Gamma({\mathcal O}_{P_i\times\{t\}}(l_iD_i))$.
So we get an element
$\big( \psi^{\vee}\big(s|_{\overline{X}\times\{t\}}\big) ,\, (\sigma_i)_i \big)$
of $\Big( \Gamma({\mathcal L}_0^{\vee}) \oplus 
\bigoplus_{i=1}^n \Gamma({\mathcal O}_{P_i\times T}(l_iD_i)^{\oplus r})
\Big)\otimes k(t)$.
By the construction, we have
$\Phi \big( \psi^{\vee}\big(s|_{\overline{X}\times\{t\}}\big) ,\, (\sigma_i)_i \big) \,=\,0$.
So we get the inclusion map
$\Gamma({{\mathcal X}\times\{t\}, \,\mathcal F}|_{{\mathcal X}\times\{t\}})
\,\subset\, \ker (\Phi\otimes k(t))$.

To prove the reverse direction, take a section
$(\alpha,\,(s_i))\,\in\, \ker (\Phi\otimes k(t))$.
Since $\Gamma(\partial_{{\mathcal L}_{\bullet}})(\alpha)\,=\,0$,
there is a section $s\,\in\, \Gamma({\mathcal F}|_{\overline{X}\times\{t\}})$
such that $\psi^{\vee}(s)\,=\,\alpha$. Considering the middle component of
$\Phi(\alpha,\, (s_i))\,=\,0$, we obtain the equality
$s|_{\overline{X}\times_{\mathcal X}U_i}\,=\,s_i|_{\overline{X}\times_{\mathcal X}U_i}$,
because the maps
$\Gamma({\mathcal F}|_{\widetilde{X}\times T}(\widetilde{l}\widetilde{D}))
\,\longrightarrow\, \Gamma({\mathcal L}_0^{\vee}|_{\widetilde{X}\times T}(\widetilde{l}\widetilde{D}))$
and
$\Gamma({\mathcal F}|_{(\overline{X}\times_{\mathcal X}U_i)\times\{t\}})
\,\longrightarrow\,
\Gamma({\mathcal F}|_{(\widetilde{X}\setminus\widetilde{D})\times\{t\}})$
are injective. So $(s,\,(s_i))$ is in the kernel of
\[
\Gamma({\mathcal F}|_{(\overline{X}\sqcup\coprod_{i=1}^n U_i)\times\{t\}})
\ \longrightarrow\ \Gamma({\mathcal F}|_{(\overline{X}\sqcup\coprod_{i=1}^n U_i)\times_{\mathcal X}
 (\overline{X}\sqcup\coprod_{i=1}^n U_i)\times\{t\}}),
\]
which is in fact
$\Gamma({{\mathcal X}\times\{t\},\, \mathcal F}\big|_{{\mathcal X}\times\{t\}})$.
So we also have the inclusion
$\ker (\Phi\otimes k(t)) \,\subset\,
\Gamma({{\mathcal X}\times\{t\}, \mathcal F}\big|_{{\mathcal X}\times\{t\}})$.

This proves the claim.

Since the Claim holds, it suffices to show that
$$\left\{ t \,\in \,T \,\, \middle|\, \: \dim\ker(\Phi\otimes k(t))\,\geq \,d\right\}$$
is Zariski closed for any $d\,\in\,\mathbb{Z}_{\geq 0}$.
Note that
\[
 \dim \ker (\Phi\otimes k(t))
\, =\,\rank_{\Gamma({\mathcal O}_T)} 
 \Big( \Gamma({\mathcal L}_0^{\vee}) \oplus 
 \bigoplus_{i=1}^n \Gamma({\mathcal O}_{P_i\times T}(l_iD_i)^{\oplus r}) \Big)
 -\rank (\Phi\otimes k(t)).
\]
Since the subset of $T$ given by locus of all points satisfying the condition
$$\rank(\Phi\otimes k(t))\,\leq\,
\rank_{\Gamma({\mathcal O}_T)} 
\Big( \Gamma({\mathcal L}_0^{\vee}) \oplus 
\bigoplus_{i=1}^n \Gamma({\mathcal O}_{P_i\times T}(l_iD_i)^{\oplus r}) \Big) -d$$
is Zariski closed, the proof of the lemma is complete.
\end{proof}

\section*{Acknowledgements}

We are very grateful to the referee for detailed comments to improve the exposition.
The work was initiated during a visit of IB to Kobe University. He is grateful to Kobe University for hospitality
and great working conditions.
During the conference `The 13th MSJ-SI' ``Differential Geometry and Integrable Systems",
all of the authors fortunately had opportunities to meet together which pushed the work to its completion.
They are grateful to Professor Masashi Yasumoto for an excellent organization of the conference.

\end{document}